\documentclass[reqno]{amsart}

\usepackage[mathscr]{eucal}
\usepackage[all]{xy}
\usepackage{mathrsfs}
\usepackage{xypic}
\usepackage{amsfonts}
\usepackage{amsmath}
\usepackage{amsthm}
\usepackage{amssymb}
\usepackage{latexsym}

\theoremstyle{plain}
\newtheorem{thm}[equation]{Theorem}
\newtheorem{prop}[equation]{Proposition}
\newtheorem{lem}[equation]{Lemma}
\newtheorem{cor}[equation]{Corollary}

\theoremstyle{definition}
\newtheorem{defn}[equation]{Definition}
\newtheorem{prob}[equation]{Problem}
\newtheorem{conjecture}[equation]{Conjecture}
\newtheorem{question}[equation]{Question}

\theoremstyle{remark}
\newtheorem{examp}[equation]{Example}
\newtheorem{examps}[equation]{Examples}
\newtheorem{rem}[equation]{Remark}
\newtheorem{claim}[equation]{Claim}

\numberwithin{equation}{subsection}
\newcommand{\pf}{\noindent{\em Proof.}\;}

\newcommand{\Lmod}[1]{#1\text{-}\mathsf{mod}}
\newcommand{\Rmod}[1]{\mathsf{mod}\text{-}#1}
\newcommand{\hdot}{{\:\raisebox{3pt}{\text{\circle*{1.5}}}}}

\newcommand{\idot}{{\:\raisebox{1pt}{\text{\circle*{1.5}}}}}

\DeclareMathOperator{\Tor}{\mathrm{Tor}}
\DeclareMathOperator{\Ext}{\mathrm{Ext}}
\DeclareMathOperator{\Rad}{\mathrm{Rad}}

\DeclareMathOperator{\rk}{\mathrm{rk}}
\DeclareMathOperator{\Fun}{\mathrm{Fun}}

\newcommand{\pare}[1]{\bigl( #1\bigr)}

\def\map{\longrightarrow}
\newcommand{\dis}{\displaystyle}

\newcommand{\beq}{\begin{equation}\label}
\newcommand{\eeq}{\end{equation}}
\DeclareMathOperator{\Spec}{\mathrm{Spec}}
\DeclareMathOperator{\pr}{pr}
\newcommand{\iso}{{\;\;\stackrel{_\sim}{\longrightarrow}\;\;}}

\def\ccirc{{{}_{^{\,^\circ}}}}
\DeclareMathOperator{\Proj}{\mathrm{Proj}}
\DeclareMathOperator{\GL}{\mathrm{GL}}
\DeclareMathOperator{\mres}{mod}

\newcommand{\cO}{\mathcal{O}}
\newcommand{\cC}{{\mathscr C}}

\newcommand{\g}[1]{\mathfrak{#1}}
\newcommand{\scr}[1]{\mathscr{#1}}
\newcommand{\cal}[1]{\mathcal{#1}}

\newcommand{\Perf}{{D_{\mathrm{perf}}}}
\newcommand{\Dcoh}{{D^b_{\mathrm{coh}}}}
\newcommand{\Coh}{\mathrm{Coh}} 

\newcommand{\en}{\enspace}
\DeclareMathOperator{\End}{\mathrm{End}}
\newcommand{\op}{\mathrm{op}}
\DeclareMathOperator{\Hom}{\mathrm{Hom}}
\newcommand{\inj}{\hookrightarrow}
\DeclareMathOperator{\Der}{\mathrm{Der}}

\newcommand{\erem}{\hfill$\lozenge$\end{rem}}

\newcommand{\dder}{{\mathbb{D}\mathbf{er}}}
\renewcommand{\o}{\otimes}

\newcommand{\Gm}{{\mathbb{G}}_{\mathbf{m}}}
\newcommand{\ac}{\text{\it{act}}}
\newcommand{\ev}{\mathrm{ev}}
\newcommand{\odd}{\mathrm{odd}}
\newcommand{\id}{\mathrm{id}}
\newcommand{\Id}{\mathrm{Id}}
\newcommand{\e}{\varepsilon} 
\newcommand{\bimod}[1]{#1\text{-}\mathsf{bimod}} 
\newcommand{\Mat}{\mathrm{Mat}}
\newcommand{\Tr}{\mathrm{Tr}}
\newcommand{\tr}{\mathrm{tr}}
\newcommand{\diag}{\mathrm{diag}}
\newcommand{\Ab}{\mathsf{Ab}}
\newcommand{\Set}{\mathsf{Sets}}
\newcommand{\OO}{{\mathbb{O}}}
\newcommand{\im}{\mathrm{Im}}
\newcommand{\f}{\varphi}

\newcommand{\eer}{\hfill$\lozenge$\end{rem}}
\newcommand{\eex}{\hfill$\lozenge$\end{examp}}
\newcommand{\spad}{\heartsuit}
\newcommand{\eps}{\varepsilon}
\DeclareMathOperator{\Inn}{\mathrm{Inn}}
\DeclareMathOperator{\Aut}{\mathrm{Aut}}

\newcommand{\st}{\,\mid\,}
\DeclareMathOperator{\Sym}{\mathrm{Sym}}
\DeclareMathOperator{\clif}{\mathrm{Cliff}}
\DeclareMathOperator{\ad}{\mathrm{ad}}
\newcommand{\sgn}{\tt{sign}}
\renewcommand{\L}{{\mathscr{L}}}
\def\bplus{{\mbox{$\bigoplus$}}}
\def\btimes{{\mbox{$\bigotimes$}}}
\DeclareMathOperator{\Ker}{\mathrm{Ker}}
\DeclareMathOperator{\Coker}{\mathrm{Coker}}
\newcommand{\tooo}{{\;{-\!\!\!-\!\!\!-\!\!\!-\!\!\!\longrightarrow}\;}}

\def\calA{{\mathcal{A}}}
\def\calE{{\mathcal{E}}}
\def\calT{{\mathscr{T}}}
\def\calD{{\mathscr{D}}}
\def\scrc{{\mathscr{C}}}
\def\calL{{\mathscr{L}}}
\def\tinv{\mbox{$\frac{1}{t}$}}
\def\sZ{{\mathsf{Z}}}

\newcommand{\sminus}{\smallsetminus}
\newcommand{\mto}{\longmapsto}

\newcommand{\inv}{^{-1}}
\newcommand{\vi}{${\sf {(i)}}\;$}
\newcommand{\vii}{${\sf {(ii)}}\;$}
\newcommand{\viii}{${\sf {(iii)}}\;$}
\newcommand{\iv}{${\sf {(iv)}}\;$}
\newcommand{\vv}{${\sf {(v)}}\;$}

\newcommand{\sset}{\subset}

\newcommand{\into}{{}^{\,}\hookrightarrow^{\,}}
\newcommand{\too}{\,\longrightarrow\,}
\newcommand{\onto}{{}^{\,}\twoheadrightarrow^{\,}}

\newcommand{\oper}{\operatorname}
\newcommand{\C}{\Bbbk}
\newcommand{\R}{\mathbb R}
\newcommand{\N}{\mathbb N}
\newcommand{\Z}{\mathbb Z}
\newcommand{\Q}{\mathbb Q}
\newcommand{\HH}{{H\!H}}
\newcommand{\DR}{{\oper{DR}}}

\newcommand{\comO}{\Omega_{\mathrm{com}}}
\newcommand{\ncO}{\Omega_{\mathrm{nc}}}
\newcommand{\ncP}{{\mathcal P}_{\mathrm{nc}}}
\newcommand{\Vect}{\mathsf{Vect}}
\newcommand{\Th}{\Theta}
\newcommand{\ncT}{\Theta_{\mathrm{nc}}}
\newcommand{\oA}{\overline{A}}
\newcommand{\Lie}{\mathrm{Lie}}
\newcommand{\alt}{\mathrm{alt}}
\newcommand{\pder}[2]{\frac{\partial #1}{\partial #2}}
\DeclareMathOperator{\gr}{\mathrm{gr}}
\DeclareMathOperator{\pic}{\mathrm{Pic}}
\def\ip<#1,#2>{\left\langle#1,#2\right\rangle}
\def\sp<#1>{\left\langle#1\right\rangle}

\newcommand{\Triv}{\mathrm{Triv}}
\newcommand{\Harm}{\mathrm{Harm}}
\newcommand{\ch}{\mathrm{ch}}
\newcommand{\cch}{{\mathsf{c}}}
\newcommand{\ncW}{W_{\mathrm{nc}}}

\newcommand{\cs}{\mathrm{cs}}

\newcommand{\Rep}{\mathrm{Rep}}
\newcommand{\beps}{\C[\eps]/\langle\!\langle\eps^2\rangle\!\rangle}
\newcommand{\llb}{\langle\!\langle}
\newcommand{\rrb}{\rangle\!\rangle}

\newcommand{\lotimes}{\,{\stackrel{L}\otimes}}

\newcommand{\Eu}{\mathsf{eu}}

\def\ip<#1,#2>{\left\langle#1,#2\right\rangle}

\newcommand{\basic}{\mathrm{basic}}

\def\pbox{\parbox[t]{140mm}}
\def\npb{\noindent$\bullet\quad$\parbox[t]{120mm}}
\def\pb{$\bullet\quad$\parbox[t]{120mm}}
\def\hp{\hphantom{x}}

\newcommand{\ee}{^{\operatorname{e}}}

\newcommand{\rep}{{\mathsf{rep}}}
\newcommand{\aand}{\quad\text{and}\quad}

\newcommand{\der}{{\mathtt{Der}^{\,}}}

\newcommand{\derp}{\mathtt{Der}_{_{\!{\mathscr{P}}}}}

\newcommand{\tdash}{\mbox{\tiny{-}}}

\newcommand{\la}{\lambda}
\newcommand{\om}{\omega}
\newcommand{\Om}{\Omega}
\newcommand{\wh}{\widehat}

\newcommand{\PP}{{\mathcal{P}}}

\newcommand{\bul}{^\hdot\!}
\newcommand{\bulp}{^\hdot_{_{\!\PP\!}}}
\newcommand{\sss}{{\mathsf{ss}}}

\def\oo{{\mathcal O}}

\def\U{{\mathcal U}}
\def\up{{\mathcal U}^{\PP\!}\!}
\def\upa{{{\mathcal U}^{\PP\!}\!A}}
\def\pp{_{_{\!\PP}}}

\def\B{{\sf B}}

\def\T{{\mathsf T}}
\def\Tc{{\mathsf{\check{T}}}}
\def\tc{\check{T}}

\def\se{{\mathsf{S}}}
\def\sr{{\mathsf{R}}}
\def\sd{{\mathsf{D}}}
\def\su{{\mathsf{U}}}
\def\k{{\Bbbk}}
\def\s{{\mathbb{S}}}

\def\bpa{{{\mathsf{B}}_{\bullet}^\PP\!A}}
\def\ccirc{{{}_{^{^\circ}}}}

\def\hp{\hphantom{x}}

\begin{document}
\setlength{\parindent}{6mm}
\setlength{\parskip}{3pt plus 5pt minus 0pt}


\centerline{\Large {\bf Lectures on Noncommutative Geometry}}

\vskip 5mm
\centerline{\large {\sc Victor Ginzburg}}
\vskip 3mm
\centerline{University of Chicago, Department of Mathematics}
\centerline{{\tt ginzburg@math.uchicago.edu}}
\vskip 5mm


\hphantom\qquad{\textbf{Abstract.}\quad
These Lectures are based on a course
on noncommutative geometry given by the author in 2003. 
The lectures contain some standard material, such as
Poisson and Gerstenhaber algebras, deformations,
 Hochschild cohomology, Serre functors,
etc. We also discuss many less known as well as some new results
such as noncommutative
Chern-Weil theory,  noncommutative symplectic geometry,
noncommutative differential forms and double-tangent bundles.

\bigskip

\centerline{\sf Table of Contents}
\vskip -1mm

$\hspace{30mm}$ {\footnotesize \parbox[t]{115mm}{
\hp${}_{}$\!
\hp\!1.{ $\;\,$} {\tt Introduction} \newline
\hp2.{ $\;\,$} {\tt Morita equivalence}\newline
\hp3.{ $\;\,$} {\tt Derivations and Atiyah algebras}\newline
\hp4.{ $\;\,$} {\tt The Bar complex}\newline
\hp5.{ $\;\,$} {\tt Hochschild homology and cohomology}\newline
\hp6.{ $\;\,$} {\tt Poisson brackets and Gerstenhaber algebras}\newline
\hp7.{ $\;\,$} {\tt Deformation quantization}\newline
\hp8.{ $\;\,$} {\tt K\"ahler differentials}\newline
\hp9.{ $\;\,$} {\tt The Hochschild-Kostant-Rosenberg Theorem}\newline
10.{ $\;\,$}  {\tt Noncommutative differential forms}\newline
12. {$\;\,$}  {\tt Noncommutative Calculus}\newline
13. {$\;\,$}  {\tt The Representation functor}\newline
14. {$\;\,$}  {\tt Double-derivations and the double tangent bundle}\newline
15. {$\;\,$}  {\tt Noncommutative symplectic geometry}\newline
16. {$\;\,$}  {\tt Kirillov-Kostant bracket}\newline
17. {$\;\,$}  {\tt Review of (commutative) Chern-Weil theory}\newline
18. {$\;\,$}  {\tt Noncommutative Chern-Weil theory}\newline
19. {$\;\,$}  {\tt Chern character in $K$-theory}\newline
20. {$\;\,$}  {\tt Formally smooth algebras}\newline
21. {$\;\,$}  {\tt Serre functors and duality}\newline
22. {$\;\,$}  {\tt Geometry over an operad}\newline
}}
}

\section{Introduction} This is an expanded version of a course
given by the author at the University of Chicago in Winter 2003.
A preliminary draft of lecture notes was prepared by Daniel Hoyt
at the time of the course. The present version is about twice
as large as  the original notes and  it contains a lot of additional
material.
However, the reader should be warned at the outset that the 
version at hand is still
a very rough draft which is by no means complete.
 I decided to make the text public `as is'
since there is a real danger that a more
perfect and more complete version will not appear
in a foreseeable future.

In these lectures we will not attempt to present a systematic treatment of noncommutative
   geometry since we don't think such a theory presently exists. 
   Instead, we will try to convey an (almost random) list of beautiful 
   concrete examples and  general guiding principles which seem certain
   to be part of any future theory, even though we don't know what that
   theory is going to be. Along the way, 
we will try to formulate many open questions and problems.

To avoid misunderstanding, we should caution the reader 
 that the name `noncommutative
   geometry' is quite ambiguous; different people 
attach to it different meanings. 
Typically, by noncommutative
(affine) algebraic   geometry one understands studying
noncommutative algebras from the point of view
of their similarity to coordinate rings
of affine algebraic varieties.
More generally,
(a not necessarily affine) noncommutative geometry
studies  (some interesting) abelian, resp. triangulated,
categories which share  some properties of the
abelian category of coherent sheaves on
a (not necessarily affine) scheme, resp. 
the corresponding derived category.

It is important  to make 
a distinction between what may be called noncommutative
   geometry `in the small', and  noncommutative
   geometry `in the large'. The former is
a {\em generalization} of the conventional
`commutative' algebraic geometry to
the noncommutative world. The objects that one studies 
here should be thought of as noncommutative
deformations, sometimes referred to as {\em quantizations},
 of their commutative counterparts. 
A typical example of this approach is the way of
thinking about the universal enveloping algebra
of a finite dimensional Lie algebra $\g g$
as a deformation of the symmetric algebra $S(\g g)$,
which is isomorphic to the polynomial algebra.

As opposed to the noncommutative
   geometry `in the small',
 noncommutative
   geometry `in the large' is  {\em not}
a generalization of commutative theory.
The world of noncommutative
   geometry `in the large' does not contain
commutative world as a special case,
but is only similar,  parallel, to it.
The concepts and results that one develops here,
do  {\em not}
specialize to their commutative analogues.
Consider for instance
the notion of {\em smoothness}  that exists
both in commutative algebraic geometry and
in noncommutative algebraic geometry  `in the large'.
A commutative algebra $A$ may be smooth
in the sense of  commutative algebraic geometry,
and at the same time be non-smooth
from the point of view of noncommutative geometry  `in the large'.
An explantation of this phenomenon comes from operad theory,
see e.g. \cite{MSS}, \cite{GiK}, \cite{Ka1}.
Each of the mathematical worlds that we study is governed
by an appropriate operad. Commutative geometry is
 governed
by the operad of commutative (associative) algebras,
while noncommutative geometry  `in the large' is governed
by the operad of associative not necessarily commutative
algebras. In this sense, it would be more appropriate
to speak of `associative geometry' instead of
what we call  noncommutative geometry  `in the large'.
There are other geometries arising from operads
of Lie algebras, Poisson algebras, etc.

Many interesting and important topics of noncommutative 
geometry are completely left out in these notes. For example, I
have not discussed noncommutative {\em projective}
geometry at all. The interested reader is referred to
\cite{St} and \cite{SvB} for an excellent reviews. Some additional
references are given in the bibliography at the end
of the lectures.

\subsection*{Basic notation.}
Throughout this text we fix $\C$, an algebraically closed field
of characteristic zero, which may be assumed without loss of generality
to be the field of complex numbers. By an 
{\em algebra} we always mean an associative, not
necessarily commutative, unital
$\k$-algebra.  If $A$ is an algebra, we denote by
$\Lmod{A}$, $\Rmod{A}$ and $\bimod{A}$ the categories of left
$A$-modules, right $A$-modules and $A$-bimodules, respectively.
We write $\otimes=\otimes_\C$.

\subsection*{Acknowledgments.}{\small{
I would like to thank D. Boyarchenko for providing  proofs of
several results claimed in the course without proof
 and to D. Hoyt for making preliminary notes of the course.

I am very much indebted to Maxim Kontsevich for 
generously sharing with me his unpublished ideas.
I have  benefited from  many useful discussions with
Bill Crawley-Boevey, Pavel Etingof,
Mikhail Kapranov,  Toby Stafford, Boris Tsygan, and
Michel Van den Bergh.

This work was partially supported by the NSF.}}

\section{Morita Equivalence}
\subsection{Categories and functors.}
We remind the reader some basic concepts involving
 categories and functors. 

Recall that  a category where
$\Hom$-spaces are  abelian groups  and such that
the notion of 
direct sum
(also called {\em coproduct}) of a family of objects is defined
 is called  an additive category. 

A functor $F\colon\cC_1\to\cC_2$ is
\emph{fully faithful} if $F$ is an isomorphism on every set of
morphisms, and that $F$ 
is \emph{essentially surjective} if for every object $X\in\cC_2$, 
there is some $Y\in\cC_1$ such that $X$ and $F(Y)$ are isomorphic.

The most commonly used way to establish an equivalence of categories is
provided by the following

\begin{lem}\label{faithful}
Let $\cC_1$ and $\cC_2$ be two abelian categories, and let
$F\colon\cC_1\to\cC_2$ be an exact, fully faithful, essentially
surjective functor.  Then $F$ is an equivalence of categories.
\qed
\end{lem}

 Let $\Set$ be the category of sets.
For any category $\cC$, functors from $\cC$ to $\Set$ form
a category $\Fun(\cC,\Set)$.
Now,  any object $X\in\cC$,
gives rise to the functor
$\Hom_\cC(X, -)\colon \cC\to\Set.$
It is straightforward to check that 
the assignment $X\mto \Hom_\cC(X, -)$ extends to a natural
contravariant functor $\cC\map \Fun(\cC,\Set)$.

\begin{lem}[Yoneda lemma] \label{yoneda} The functor
 $\cC\map \Fun(\cC,\Set)$ induces isomorphisms
on $\Hom$'s, in other words, it is a {\em fully faithful}
functor.
\end{lem}

Write $\Ab$ for the  category of abelian groups.
\begin{defn} An object $P$ of  an abelian category $\cC$ is
said to be
\emph{projective} if the functor $\Hom_\cC(P,-)\colon\cC\to\Ab$ is
exact.
\end{defn}

In other words, $P$ is projective if given a short exact sequence
$$
0\to M'\to M\to M''\to0
$$
in $\cC$, we have that
$$
0\to\Hom_\cC(M',P)\to\Hom_\cC(M,P)\to\Hom_\cC(M'',P)\to0
$$ is exact in $\Ab$.  An object $G$ of $\cC$ is called a \emph{generator} if
$\Hom_\cC(G,A)$ is nonzero for every nonzero object $A$ of $\cC$.

\begin{defn}\label{compact} Let  $\cC$ be an additive category
 with arbitrary direct sums (also referred to as {\em coproduct}).
An object $X$ of $\cC$ is called \emph{compact} if,
for an arbitrary set of objects of $\cC$ and a  morphism
$
f\colon X\to\bigoplus_{\alpha\in I}M_\alpha,
$
 there exists some finite 
set $F\subset I$ such that $\im f$ is a subobject of $\bigoplus_{\alpha\in F}M_\alpha$.
\end{defn}

An easy consequence of the definition of compactness is the following

\begin{lem}\label{L:DirSumComp}
An object $X$ in an abelian category $\cC$ (with arbitrary direct sums)
is compact if and 
only if the functor $\Hom_\cC(X,-)$ commutes with arbitrary direct sums, that is
$$
\Hom_\cC\left(X\,,\,\oplus_{\alpha\in\Lambda}\,Y_\alpha\right)=
\oplus_{\alpha\in\Lambda}\,\Hom_\cC(X,Y_\alpha).\qquad\Box
$$
\end{lem}

\begin{lem} Let $A$ be a ring and $M$ an $A$-module.

\vi If $M$ is a finitely generated $A$-module, then $M$ is a
compact object of $\Lmod{A}$.

\vii If $M$ is projective and is a compact object of $\Lmod{A}$,
then $M$ is finitely generated.
\end{lem}
\begin{proof} \vi is obvious. To prove (ii), assume $M$ is projective, and
choose any surjection $p:A^{\oplus I}\twoheadrightarrow M$, where
$I$ is a possibly infinite set. There exists a section
$s:M\hookrightarrow A^{\oplus I}$. If $M$ is compact, the image of
$s$ must lie in a submodule $A^{\oplus J}\subseteq A^{\oplus I}$
for some finite subset $J\subseteq I$. Then $p\big\vert_{A^{\oplus
J}}$ is still surjective, which shows that $M$ is finitely
generated.
\end{proof}

The following result provides a very useful criterion for
an  abelian category to be equivalent to the category
of left modules over a ring.

\begin{prop}\label{P:AbCatEqCond}
Let $\cC$ be an abelian category with arbitrary 
direct sums.  Let $P\in\cC$ be a compact projective 
generator and set $B=(\End_\cC P)^\op$.  Then  the functor
$\Hom_\cC(P,-)$
yields an equivalence of categories between $\cC$ and $\Lmod B$.
\end{prop}

\begin{proof}[Proof of Proposition \ref{P:AbCatEqCond}]
We will show that $F(X)=\Hom_\cC(P,X)$ is fully faithful and then
apply Lemma \ref{faithful}.  

We wish to
show that there is an identification between $\Hom_\cC(X,M)$ and
$\Hom_{\Lmod B}(F(X),F(M))$ for all $M\in\cC$.  Since $P$ is a
generator, we deduce that $\Hom_\cC(P,M)\ne0$.  Define
$$
\f\colon\bigoplus_{f\in\Hom_\cC(P,M)}P\to M
$$
by $\f(p_f):=\sum_{f\in\Hom_\cC(P,M)}f(p_f)$ (this sum is finite).  Let $L=\im\f$.  Then $L$ is a submodule of $M$, and if it is not all of $M$ then $M/L$ is nonzero, hence there is some nonzero map $P\to M/L$.  But this map lifts to a map to $M$, and the image of this lift must include points not in $L$, which is a contradiction.  So, every $M\in\Lmod A$ can be written as a quotient of $P^{\oplus T}$ for some cardinal $T$.  Let $K$ denote the kernel, and take $P^{\oplus S}$ surjecting onto $K$.  Then composing this with the inclusion of $K$ into $P^{\oplus T}$ yields the exact sequence
$$
P^{\oplus S}\to P^{\oplus T}\to M\to0.
$$
Since $P$ is projective, $\Hom_\cC(P,-)$ is exact.  Hence,
$$
\Hom_\cC(P,P^{\oplus S})\to\Hom_\cC(P,P^{\oplus T})\to\Hom_\cC(P,M)\to0
$$
is exact.  Since $P$ is finitely generated (i.e., compact), it commutes
with arbitrary direct sums (see Lemma \ref{L:DirSumComp}.  So, using 
exactness of the above sequence it suffices to check that $\Hom_\cC(X,P)=\Hom(F(X),F(P))$.  Since $F(P)=\End_AP=B^\op$, this is automatic.
\end{proof}

We need two more definitions.

\begin{defn}
Let $\cC_1$ and $\cC_2$ be two categories, and let $F,G\colon\cC_1\to\cC_2$ be two functors.  A \emph{morphism} $\phi\colon F\to G$ is a natural transformation, i.e., a collection of morphisms
$$
\phi_X\colon F(X)\to G(X)
$$
for each $X\in\cC_1$ such that for any morphism $f\colon X\to Y$ ($Y\in\cC_1$), the following diagram commutes:
$$
\xymatrix{F(X)\ar[r]^{\phi_X}\ar[d]_{F(f)}&G(X)\ar[d]^{G(f)}\\
F(Y)\ar[r]^{\phi_Y}&G(Y)}.
$$
In particular, if $\cC$ is an abelian category, we have the identity morphism $\id_\cC\colon\cC\to\cC$.  We define the \emph{center} $\sZ(\cC)=\End(\id_\cC)$.
\end{defn}

\begin{examp}
It is a worthwhile exercise to 
check that $\sZ(\Coh X)\simeq{\mathcal{O}}(X)$ for any algebraic variety
(where 
${\mathcal{O}}(X)$ is the  ring of global regular functions on $X$).
\eex 

\begin{lem}\label{Z(A)}
Let $A$ be a associative algebra.  Then $\sZ(\Lmod A)=\sZ_A$.
\end{lem}

\begin{proof}
Choose an element $z\in \sZ_A$.  Define an endomorphism of $\id_\cC$ by setting $\phi_X$ to be the action of $z$ on the module $X$.  Since $z$ is central, this is a module homomorphism and it commutes with all module maps, that is, the diagram
$$
\xymatrix{X\ar[r]^{\phi_X}\ar[d]_f&X\ar[d]^f\\
Y\ar[r]^{\phi_Y}&Y}.
$$
Conversely, suppose an endomorphism $\phi\colon\id_\cC\to\id_\cC$ is given.  Set $z=\phi_A(1_A)$.  We need to check that this is indeed central.  Choose any $a\in A$, and define the left $A$-module map $f\colon A\to A$ by $f(x)=xa$.  Then since $\phi$ is a morphism, we know that $f\circ\phi_A=\phi_A\circ f$.  So,
\begin{align*}
za&=\phi_A(1_A)a=(f\circ\phi_A)(1_A)\\
&=(\phi_A\circ f)(1_A)=\phi_A(a)=a\phi_A(1_A)=az.
\end{align*}
So, $z\in \sZ_A$.  It is clear that this association is an algebra homomorphism.
\end{proof} 

\subsection{Algebras and spaces.}
One of the cornerstones of geometry is the equivalence of categories of
spaces and categories of algebras.  For example, the
Gelfand theorem asserts an (anti)-equivalence 
between the category of locally compact Hausdorff spaces 
with proper maps and the category of commutative $C^*$-algebras.  
Similarly, in algebra, one has an (anti)-equivalence 
between the category of affine algebraic varieties
and the category of finitely generated commutative algebras
without nilpotent elements.
A first step to such an equivalence in 
algebraic geometry is by
associating to each scheme $X$ its structure sheaf $\cO_X$.  However,
this is unsatisfactory since $\cO_X$ explicitly refers to the space $X$,
since it is a sheaf on~$X.$

One way to resolve this difficulty is
to forget about  scheme $X$ altogether, and to consider
instead the abelian category $\Coh(X)$ of coherent
sheaves over $X$.  
The space $X$ can be reconstructed in a natural way from this category.

A first step in considering such abelian categories is to look at the
category  $\Lmod A$ of 
left  modules over an associative algebra $A$.  One would hope that
the category $\Lmod A$ should uniquely determine $A$ up to isomorphism.
  It turns out that
this is the case for {\em commutative}  associative algebras,
but {\em not} for arbitrary  associative algebras.

\begin{defn}
Let $A$ and $B$ be  associative, not necessarily commutative, algebras.
  Then we say that $A$ and $B$ are
\emph{Morita equivalent}, if
 there is an
 equivalence of categories between $\Lmod A$ and $\Lmod B$.
\end{defn}

Morita equivalence of two commutative algebras is particularly simple.

\begin{prop}\label{P:CommMorita}
Commutative algebras $A$ and $B$ are Morita equivalent if and only if 
they are isomorphic.
\end{prop}

\begin{proof}[Proof of Proposition \ref{P:CommMorita}]
Clearly isomorphic algebras are Morita equivalent.  
Now suppose that $A$ and $B$ are  Morita equivalent
associative algebras.  Then $\Lmod A\sim\Lmod B$, so certainly we have
that
 $\sZ(\Lmod A)\simeq \sZ(\Lmod B)$.  
  If $A$ and $B$ are both commutative, then
by Lemma \ref{Z(A)} we have $A=\sZ_A$ and $B=\sZ_B$.
\end{proof}

\begin{examp} 
Given an algebra $A$ and an
integer $n\geq 1$, let $\Mat_n(A)$, be 
the algebra of $n\times n$-matrices with entries in $A$.
A typical example of Morita equivalence involving noncommutative
algebras is provided by the following
easy result.

\begin{lem}\label{morita_mat} For any algebra $A$ and any
integer $n\geq 1$, the algebras  $A$ and  $\Mat_n(A)$ are Morita equivalent.
\end{lem}

Note that even if $A$ is a commutative algebra, the algebra
 $\Mat_{n}(A)$ is not commutative for any $n>1$.

\begin{proof} 
Let $e_{11}\in \Mat_n(A)$ be the matrix with entry $1$
at the $1\times 1$ spot, and zeros elsewhere.
It is clear that the algebra $e_{11}\cdot \Mat_n(A)\cdot e_{11}$ is isomorphic to
$A$.
On the other hand, it is clear that multiplying the
matrices with arbitrary entry $a\in A$  at the $1\times 1$ spot, and
zeros elsewhere
by other matrices from $\Mat_n(\C)\sset \Mat_n(A)$, one can
obtain, taking linear combinations, every $A$-valued 
$n\times n$-matrix. Thus, we have shown that
$\Mat_n(A)=\Mat_n(A)\cdot e_{11}\cdot \Mat_n(A)$ and
that $A\simeq e_{11}\cdot \Mat_n(A)\cdot e_{11}$.
At this point, the result follows
from the general criterion of Corollary \ref{AeA} below.
\end{proof}
\end{examp}

Thus, in general, the algebra $A$ can not be recovered from
the corresponding abelian category $\Lmod{A}$. Therefore,
in order for a concept in noncommutative geometry to have an intrinsic
meaning, that concept must be Morita invariant.
In particular, the question of which properties of an algebra
 are Morita invariant  becomes very important.

\subsection{Morita theorem.} The main result about Morita equivalent
algebras is provided by the following\footnote{The exposition below
follows the notes prepared by M. Boyarchenko.}

\begin{thm}\label{t:morita}
Let $A$ and $B$ be two rings, and $F:\Lmod{A}\to\Lmod{B}$ an
additive right exact functor. Then there exists a $(B,A)$-bimodule
$Q$, unique up to isomorphism, such that $F$ is isomorphic to the
functor
\[
\Lmod{A}\too\Lmod{B},\quad M\longmapsto Q\otimes_A M.
\]
\end{thm}
\begin{proof}
The uniqueness of $Q$ (if it exists) is clear,
since we have $Q=Q\o_AA=F(A)$. To prove existence,
let $Q=F(A)$; by assumption, this is a left $B$-module. Moreover,
for every $a\in A$, the operator $\rho_a$ of right multiplication
by $a$ is an endomorphism of $A$ as a left $A$-module, whence we
obtain a ring homomorphism $A^{op}\to\End_B(Q)$, $a\mapsto
F(\rho_a)$. This homomorphism makes $Q$ into a $(B,A)$-bimodule.

Now, for every $M\in\Lmod{A}$, we define, functorially, a
$B$-module homomorphism $Q\otimes_A M\to F(M)$. Let us first
define a $\Z$-bilinear map $\phi_M:Q\times M\to F(M)$. An element
$m\in M$ gives rise to an $A$-module homomorphism $\rho_m:A\to M$,
$a\mapsto a\cdot m$. We define $\phi_M(q,m)=F(\rho_m)(q)$. Now
since $F(\rho_m)$ is a $B$-module homomorphism, the map $\phi_M$
commutes with left multiplication by elements of $B$. Also, if
$a\in A$, then
\[\phi_M(qa,m)=F(\rho_m)\pare{F(\rho_a)(q)}=F(\rho_m\circ\rho_a)(q)=
F(\rho_{am})(q)=\phi_M(a,qm),\]
whence $\phi_M$ descends to a left $B$-module homomorphism
$\psi_M:Q\otimes_A M \to F(M)$.

It is obvious that $\psi_M$ is functorial with respect to $M$.
Moreover, by construction, $\psi_M$ is an isomorphism whenever $M$
is free. In general, we use the fact that both $F$ and the functor
$Q\otimes_A -$ are exact. Given any left $A$-module $M$, choose an
exact sequence $F_1\to F_0\to M\to 0$, where $F_0$ and $F_1$ are
free $A$-modules, and apply both functors to this sequence. Using
the morphism of functors, we get a commutative diagram, and the
Five Lemma finishes the proof.
\end{proof}

\begin{cor}\label{c:left}
Two rings, $A$ and $B$, are Morita equivalent if and only if there
exist an $(A,B)$-bimodule $P$ and a $(B,A)$-bimodule $Q$ such that
$P\otimes_B Q\cong A$ as $A$-bimodules and $Q\otimes_A P\cong B$
as $B$-bimodules. Under this assumption, we have
$$\End_{\Lmod{A}}(P)=B^{op},\quad\op{and}\quad 
\End_{\Lmod{B}}(Q)=A^{op}.
$$

Moreover, $P$ is projective as an $A$-module and $Q$ is projective
as a $B$-module.
\end{cor}

\begin{proof}
Equivalences of categories are exact functors and preserve
projective objects.
\end{proof}
\begin{cor}\label{c:bimod}
If $A$ and $B$ are Morita equivalent rings, then the categories
$\Rmod{A}$ and $\Rmod{B}$ are also equivalent. Moreover, there is
a natural equivalence of categories $\bimod{A}\to\bimod{B}$ which
takes $A$ to $B$ (with their natural bimodule structures).
\end{cor}
\pf
Let $P$ and $Q$ be as in the previous corollary. For the first
statement, use the functors $-\otimes_A P$ and $-\otimes_B Q$. For
the second statement, use the functor
\[
Q\otimes_A - \otimes_A P : \bimod{A}\too\bimod{B}.
\qquad\Box\]
\medskip

Let $A$ be a ring and $e=e^2\in A$ an idempotent.
Clearly, $eAe$ is a subring in $A$ and $Ae$ is naturally 
an $(A,eAe)$-bimodule and $eA$ is an $(eAe,A)$-bimodule. 
Note that $e$ is the unit of the ring $A$. We
see that the 
inclusion map $eAe\hookrightarrow A$ is usually not a ring
homomorphism, and $A$ is not an $eAe$-module in a natural way.

For any left $A$-module $M$, the space $eM$ has a natural
 left $eAe$-module structure.
\begin{cor}\label{AeA}
Let $A$ be a ring and $e\in A$ an idempotent. 
The functor $\Lmod{A}\map\Lmod{B},\,
M\mto eM$ is a  Morita equivalence if and only if $AeA=A$.
In this case, an inverse  equivalence is given by
the functor $N\mto Ae\o_{eAe}N.$
\end{cor}
\begin{rem}
Observe also that 
there exist rings $R$ such that $R\times R\cong R$ (for example,
take $R$ to be the product of an infinite number of copies of some
fixed ring). Now for such a ring, let $A=R\times R$ and
$e=(1,0)\in A$. Then $eAe=R\cong A$, and in particular, $eAe$ is
Morita equivalent to $A$, but we have $AeA=R\times\{0\}\neq A$.
\end{rem}
\begin{proof}
Assume that $AeA=A$. We claim
that $eA\otimes_A Ae\cong eAe$ as $eAe$-bimodules (in fact, this
holds regardless of the assumption on $e$). Note that $eA$ is a
direct summand of $A$ as a right $A$-module: $A=eA\oplus (1-e)A$.
Similarly, $A=Ae\oplus A(1-e)$. Now the multiplication map
$A\otimes_A A\to A$ is obviously an $A$-bimodule isomorphism, and
it clearly restricts to an isomorphism $eA\otimes_A Ae
\stackrel{\sim}\map
eAe$ of $eAe$-bimodules.

Now we also claim that $Ae\otimes_{eAe}eA\cong A$ as
$A$-bimodules. We have the natural $A$-bimodule map
$Ae\otimes_{eAe}eA\stackrel{m}\map A$ given by multiplication. We can write
down an explicit inverse for this map. Namely, since $AeA=A$,
there exist elements $a_j,b_j\in A$ such that $1=\sum_j a_j e
b_j$. We define a map of abelian groups $f:A\to Ae\otimes_{eAe}eA$
by $f(a)=\sum_j aa_j e\otimes eb_j$. It is obvious that $m\circ
f=\id_A$. We must check that $f\circ m$ is also the identity. We
have
\begin{align*}
(f\circ m)&\left( \sum_k c_ke\otimes ed_k\right) = f\left( \sum_k
c_k e d_k\right) = \sum_{j,k} (c_k ed_k a_j e)\otimes eb_j \\
&= \sum_{j,k} c_k e\otimes(ed_k a_j e e b_j) = \sum_k c_k
e\otimes \left[ ed_k\cdot\sum_j a_j e b_j\right] = \sum_k c_k
e\otimes e d_k.
\end{align*}

We leave the rest of  the proof to the reader.
\end{proof}

\begin{rem}
Let us write $\Lmod{A}_f$ for the category of finitely generated
$A$-modules. In the situation of Theorem \ref{t:morita}, it is
clear that the functor $F$ takes $\Lmod{A}_f$ into $\Lmod{B}_f$ if
and only if $Q$ is finitely generated as a $B$-module. We claim
that this is always the case whenever $F$ is an equivalence; in
particular, if $A$ and $B$ are Morita equivalent rings, then the
categories $\Lmod{A}_f$ and $\Lmod{B}_f$ are equivalent. The proof
is based on the notion of a compact object.
In particular, if $F:\Lmod{A}\to\Lmod{B}$ is an equivalence of
abelian categories, then $F(A)$ must be a compact object of
$\Lmod{B}$, but it is also projective, whence finitely generated.
\end{rem}

\section{Derivations and Atiyah algebras}\subsection{}
We recall the definitions of derivations and super-derivations.

Let $A$ be an algebra, and let $M$ be an $A$-bimodule.  A $\C$-linear map $\delta\colon A\to M$ is called a \emph{derivation} if it satisfies the Leibniz rule, that is, if
$$
\delta(a_1a_2)=a_1\delta(a_2)+\delta(a_1)a_2\quad\text{for all}\quad a_1,a_2\in A.
$$

 We let $\Der(A,M)$ denote the $\C$-vector space of all derivations from
$A$ to $M$.
For any  $m\in m$, the
map $\ad m:\ M\to M\,,\, a\mapsto ma- am$ gives a derivation
of $M$. The derivations of the form $\ad m\,,\,m\in M,$ are called
{\em inner} derivations. 
We write  $\Inn(A,M)\sset \Der(A,M)$ for the space
of inner derivations, and 
$\sZ(M)=\{m\in M\mid am=ma,\;\forall a\in A\},$
for the `{\em center}' of $A$-bimodule $M$.
Thus, we have the following exact sequence
\beq{inn_der}
0\map \sZ(M)\map M\stackrel{\ad}\map
\Der(A,M)\map \Der(A,M)/\Inn(A,M)\map 0.
\end{equation}

\begin{rem} If the algebra $A$ is commutative,
then any left $A$-module $M$ may be viewed as
an $A$-bimodule, with right action being
given by $m\cdot a:= a\cdot m,\,\forall a\in A, m\in M$
(this formula only gives a right $A$-module structure
if $A$ is commutative).
Bimodules of that type are called {\em symmetric}.
Thus, given a left $A$-module  viewed
as a symmetric $A$-bimodule,
one may consider derivations $A\to M$.
\end{rem}

In the special case that $M=A$, we abbreviate $\Der(A,A)$ to
$\Der(A)$, resp. $\Inn(A,A)$ to $\Inn(A)$.
  It is an easy calculation that $\Der(A)$ is a Lie algebra
under the commutator, and that inner derivations
 form a Lie ideal $\Inn(A)$ in  $\Der(A)$.
It is also straightforward to check that the assignment
$a\mto \ad a$ is a Lie algebra map
(with respect to the commutator bracket on $A$),
whose kernel is the center, $\sZ_A$, of the algebra $A$.

The inter-relationships between $A, \sZ_A, $ and $\Der(A)$
are summirized in the following result which says that
$\Der(A)$ is a {\em Lie algebroid} on $\Spec\sZ_A$, see
Sect. \ref{algebroid_sec}.
\begin{prop}\label{der_algebroid}\vi The Lie algebra $\Der(A)$ acts on $A$, and this action
preserves the center $\sZ_A\sset A$.

\vii For any $z\in \sZ_A$ and $\theta\in \Der(A)$,
the map $z\theta: a\mapsto z\cdot\theta(a)$ is again a derivation of
$A$. The assignment $\theta\mapsto z\theta$ makes $\Der(A)$ a
$\sZ_A$-module.

\viii For any $\theta,\delta\in \Der(A)$ and $z\in \sZ_A$,
one has
$$ [z\theta,\delta]=z[\theta,\delta]-\delta(z)\theta.\qquad\Box$$
\end{prop}

\begin{examp} Let 
$A=\C[X]$, be the coordinate ring of a smooth
affine variety $X$.
 For each $x\in X$, let $\g m_x=\{f\in A\st f(x)=0\}$
be the corresponding  maximal ideal in $A$.
Then $\C_x=A/\g m_x$ is a 1-dimensional $A$-module (which we will also
view as an $A$-bimodule).

Given $x\in X$, let
$T_xX$ be the tangent space at $x$.
For any tangent vector $\xi\in T_xX$,
differentiating the function
$f$ with respect to $\xi$ gives
a map $A\ni f\mto (\xi{f})(x)=df(\xi)|_x\in\C$.
This map is a derivation $\partial_\xi: A\to \C_x.$
It is an elementary result of commutative algebra that
this way one gets an isomorphism
$$
T_xX\iso\Der(A,\C),\quad\xi\mto\partial_\xi.
$$

On $X$, we have the tangent bundle
$TX\to X$. We write $\calT_X$ for the tangent sheaf,
the sheaf of algebraic sections of the tangent bundle.
This sheaf is locally free  since $X$ is smooth, 
the sections of  $\calT_X$  are nothing but  algebraic 
vector fields on $X$. Write $\calT(X):=\Gamma(X,\calT_X)$ for
the vector space of (globally defined) algebraic vector fields on $X$.
Commutator of vector fields
makes  $\calT(X)$ a Lie algebra.

For any algebraic vector field $\xi$ on $X$,
the map $\C[X]\ni f\mapsto \xi f$ gives a derivation of $A=\C[X]$.
This way one obtains a canonical Lie algebra isomorphism
$$
\calT(X)\iso \Der(A), \quad\text{where}\quad A=\C[X].
$$

Note that the product of a function and a vector field is 
again a well-defined vector field, in accordance with
part (ii) of Proposition \ref{der_algebroid}.
\eex

\begin{rem}\label{intuitive_der} Let $A$ be an associative algebra.
 A derivation $\delta: A\to A$ may be
thought of,
heuristically, 
as a generator of an `infinitesimal' one-parameter group
$\eps\mapsto\exp(\eps\cdot\delta)=\id_A+\eps\cdot\delta
+\frac{1}{2} \eps^2\cdot\delta\ccirc\delta+\ldots,$
of automorphisms of $A$.
To formalize this, introduce  the ring $\beps,$
called the ring of `dual numbers',
and form the tensor product algebra $\beps\otimes A.$
 To any linear map $f: A\to A$
we associate the map
$$F:\ A \to \beps\otimes A,\; a\mto F(a):= a+ \eps\cdot
f(a).
$$
Here, we think of $\eps$ as an `infinitesimally small' parameter,
so, up to higher powers of $\eps$,
one has $F=\Id+\eps\cdot f \sim \exp(\eps\cdot f).$
Thus,
 the map $F$ may be thought of as a family of
maps $A\to A$  `infinitesimally close' to the identity.
Now the equation $F(a\cdot a')=F(a)\cdot F(a'),\,\forall
a,a'\in A$,
saying that our family is a family of algebra homomorphisms,
when expressed in terms of $f$, reads (recall that $\eps^2=0$):
\beq{inf_aut}
a\cdot a'+ \eps\cdot f(a\cdot a')=
a\cdot a'+ \eps\cdot \bigl(a\cdot f(a')+f(a)\cdot a')\bigr),\;\forall
a,a'\in A.
\end{equation}
Equating the coefficients in front of $\eps$,
we see that \eqref{inf_aut} reduces to the condition for
$f$  to be a derivation, as promised.

Let $\Aut(A)$ denote the group of all (unit preserving) automorphisms
of the algebra $A$.
Thus, if we think of  $\Aut(A)$ as some sort of Lie group,
then the `Lie algebra' of that group is given by 
$$\Lie \Aut(A)= \Der(A).
$$
Moreover, one can argue that the Lie bracket on the left-hand side
of this formula corresponds to the commutator bracket on the space of
derivations on the right.
\eer

\subsection{Square-zero construction.} We would like to extend the
intuitive point of view explained in Remark \ref{intuitive_der}
to a more general case where $\theta\in\Der(A,M)$
for an arbitrary $A$-bimodule $M$. This can be achieved
by the following general construction.

Suppose we are considering some class of algebraic structure,
 be it associative algebras, Lie algebras, Poisson algebras, etc.  
We also wish to discuss modules over these algebras.  There is a natural way of defining what the ``correct'' module structure is for a given type of algebra called the \emph{square zero} construction.  Suppose that $A$ is some sort of algebra over $\C$ and that $M$ is a $\C$-vector space.  We wish to give $M$ the type of module structure appropriate to the structure of $A$.  Consider the vector space $A\oplus M$.  Then giving a correct ``bimodule'' structure on $M$ is equivalent to giving an algebra structure on $A\oplus M$ such that
\begin{enumerate}
\item[(i)]{the projection $A\oplus M\to A$ is an algebra map, and}
\item[(ii)]{$M^2=0$ and $M$ is an ideal.}
\end{enumerate}
It is easy to see this principle at work.  
If $A$ is an associative algebra and $M$ is an $A$-bimodule, then
 $A\oplus M$ is an associative algebra under the multiplication
 $(a\oplus m)(a'\oplus m'):=aa'\oplus(am'+ma')$.  
Indeed, $M^2=0$, $M$ is an ideal, and the projection $A\oplus M\to A$ is
 an algebra map. 
 This is a rather trivial case, however.

\begin{lem}
A linear map $\theta\colon A\to M$ is a derivation if and only if the
map 
$A{\sharp}M\to A{\sharp}M$ given by $(a,m)\mapsto(a,m+\theta(a))$
 is an algebra automorphism.
\end{lem}

Thus, for an arbitrary bimodule $M$ we may think
of derivations $\theta\in \Der(A,M)$ as
`infinitesimal automorphisms' of the algebra
$A{\sharp}M$.

\subsection{Super-derivations.}
Now suppose that $A$ is $\Z$-graded, that is, there is a direct sum
decomposition (as $\C$-vector spaces) $A=\bigoplus_{i\in\Z}A_i$ such
that $A_iA_j\subset A_{i+j}$ for all $i,j\in\Z$.  
We put
$$
A_\ev=\bigoplus\nolimits_{i\in\Z}\,A_{2i}\quad\text{and}\quad
A_\odd=\bigoplus\nolimits_{i\in\Z}\,A_{2i+1}.
$$
A linear map $f: A\to A$ is said to be {\em even}, resp., {\em odd},
if  $f(A_\pm)\subset A_\pm,$  resp., $f(A_\pm)\subset A_\mp$,
where the plus sign stands for `even' and the minus sign stands for `odd'.

\begin{defn} An odd $\C$-linear map $f\colon A\to A$ is 
called either a \emph{super-derivation} or an \emph{odd derivation},
if it satisfies the graded Leibniz rule, that is,
$$
f(a_1a_2)=f(a_1)a_2+(-1)^{\deg a_1}a_1f(a_2),
$$
where $a_1$ is a homogeneous element of degree $k=\deg a_1$, that is,
 $a_1\in A_k$.  
\end{defn}

It is an easy calculation to check that the $\Z/(2)$-graded vector space of all
even/odd derivations 
forms a Lie super-algebra under the super-commutator:
$$[f,g]:= f\ccirc g \mp\, g\ccirc f,$$
where the plus sign is taken if both $f$ and $g$ are odd,
and the minus sign in all other cases.
In particular,  for an odd derivation $f: A\to A$,
we see that
$$f\circ f=\mbox{$\frac{1}{2}$}[f,f]$$
is an {\em even} derivation.

Below,  we will frequently use the following result,
which follows from  the Leibniz formula by an obvious induction.
\begin{lem}\label{triv}
Let $f,g: A\to A$ be to  derivations (both even, or odd)
of an associative algebra $A$. Let $S\sset A$ be a set of
algebra generators for $A$.
Then, we have
$$ f(s)=g(s),\enspace\forall s\in S\quad\Longrightarrow\quad f=g.\qquad\Box$$
\end{lem}

The most commonly used application of the Lemma is the following
\begin{cor}\label{dif_cor}
Let $d: A \to A$ be an odd derivation of a graded algebra $A$.
If  $S\sset A$ is a set of algebra generators for $A$, and
 $d^2(s)=0$, for any $s\in S$, then $d^2=0$.\qed
\end{cor}

\subsection{The tensor algebra of a bimodule.}\label{tensor_algebra}
Let $A$ be an associative algebra. 
Given two  $A$-bimodules $M$ and $N$ one has a well-defined $A$-bimodule
$M\otimes_AN$. In particular, for  an $A$-bimodule $M$ we put 
$T^n_AM:= M\otimes_AM\otimes_A\ldots\otimes_AM$ ($n$ times),
in particular, $T^1_AM=M$, and we put formally $T^0_AM:=A$.
The direct sum
$$ T^\hdot_AM:= \oplus_{i\geq 0}\;T^i_AM
$$
acquires an obvious  graded  associative  (unital) algebra structure,
called the
tensor algebra of an $A$-bimodule $M$.

Given a $\C$-vector space $V$, we will often write $V^{\otimes n}$ instead of 
$T^n_\C V$.

The tensor algebra construction may be usefully interpreted
as an adjoint functor. Specifically,
fix an  associative unital algebra $A$, and consider the category
$A\textsf{-algebras},$ whose objects are pairs $(B,f)$, where
$B$ is  an  associative unital algebra and $f: A\to B$ is
an algebra morphism such that $f(1)=1$.
Morphisms in $A\textsf{-algebras}$ are defined as algebra maps
$\varphi: B\to B'$ making the following natural diagram commute
$$
\xymatrix{
&A\ar[dl]_{f}\ar[dr]^{f'}&\\
B\ar[rr]^{\varphi}&&B'
}
$$
Note that an algebra morphism $A\to B$ makes $B$ an $A$-bimodule.
Thus, we get a functor $A\textsf{-algebras}\too A\textsf{-bimodules}$.
It is straighforward to check that
the assignment $M\mto T^\hdot_AM$ gives a right adjoint to that functor,
i.e., 
there is a canonical adjunction isomorphism
$$\Hom_{_{A\textsf{-bimodules}}}(M,B)\iso\Hom_{_{A\textsf{-algebras}}}(T^\hdot_AM,B),
$$
for any $A$-bimodule $M$ and an $A$-algebra $B$.

Now, fix an  algebra morphism $A\to B$, and view
$B$ as an $A$-bimodule. Let $M$ be another $A$-bimodule.
Any  morphism $\delta:M\to B,$
 of $A$-bimodules induces, by the adjunction isomorphism above,
an algebra homomorphism $T^\hdot_AM\to B$. 
This way, the algebra $B$ may be regarded as a $T^\hdot_AM$-bimodule.

We leave to the reader to prove the following
\begin{lem}\label{der_extend} The  $A$-bimodule map  $\delta:M\to B$
 can be uniquely extended to a  derivation, resp., super-derivation,
$$T(\delta):\ T_A^\hdot M \to B\quad\text{such that}\quad
T(\delta)\big|_{T^0_AM}=f,\quad\text{and}\quad
T(\delta)\big|_{T^1_AM}=\delta.\qquad\Box
$$
\end{lem}

\subsection{Picard group of a category.} Given a category $\scrc$ we
let $\pic(\scrc)$ be the group of all autoequivalences of  $\scrc$,
the group structure being given by composition.

For example, let $A$ be an associative algebra. We need the following
\begin{defn} A finitely-generated $A$-bimodule
$Q$ is said to be {\em invertible} if there exists a finitely-generated $A$-bimodule
$P$ such that $Q\otimes_A P\cong A\cong P\otimes_A Q$.
Then, $P$ is called an {\em inverse} of $Q$.
\end{defn}
It is clear that (the isomorphism classes of) invertible
 $A$-bimodules form a group under the tensor
product operation $P, P'\mapsto P\otimes_A P'$. The unit element
of this group is the isomorphism class of  $A$,
viewed as an  $A$-bimodule.

Now,  let $\scrc:= \Lmod{A}$ be the
category of left $A$-modules.
We know by Morita theory, see Theorem \ref{t:morita},
 that any equivalence $\Lmod A\to \Lmod A$
has the form $M\rightsquigarrow Q\otimes_A M$, for an invertible
$A$-bimodule
$Q$. Thus, we obtain
\begin{cor}\label{inv_bimod} The group $\pic(\Lmod A)$ is canonically
isomorphic to the group of (isomorphism classes of) invertible
 $A$-bimodules.\qed
\end{cor}

\begin{rem} In the traditional commutative algebraic geometry,
given an algebraic variety $X$, one writes
$\pic(X)$ for the abelian group formed by
the (isomorphism classes of) line bundles on $X$.
Now assume $X$ is  irreducible and  affine, and
put  $A=\C[X]$. Then any line bundle $\calL$ on $X$ gives
a rank 1 projective {\em left} $A$-module $L:=\Gamma(X,\calL)$.
This way, the group $\pic(X)$ gets identified
with the group formed by  rank 1 projective {\em left} $A$-modules
equipped with the tensor product structure. That gives,
in view of Corollary \ref{inv_bimod},
a certain justification for the name `Picard group' of a category.
 \eer

Recall   the group $\Aut(A)$  of automorphisms of the algebra $A$.
Given an algebra automorphism $g\in \Aut(A)$ and
an  $A$-bimodule $P$, we define
 a new $A$-bimodule $P^g$ by "twisting" the natural left action
on $P$ via $g$, i.e., by letting
$a\otimes a'\in A\ee$ act on $p\in P$ by the formula
$p\mto g(a)\cdot p\cdot a'$. It is clear that, given two automorphisms
$f,g: A\to A$, there is a canonical isomorphism $P^{fg}\cong
(P^g)^f.$ 

Let $\Aut_{\sf{in}}(A)\sset \Aut(A)$ denote the (normal)
subgroup formed by {\em inner} automorphisms of $A$,
i.e., automorphisms of the form $a\mapsto u\cdot a\cdot u\inv$,
where $u\in  A^\times$ is an invertible element.
We have the following result.
\begin{prop} The assignment $g\mapsto A^g$ yields
a group homomorphism $\Aut(A)$
$\to\pic(\Lmod A)$.
This  homomorphism descends to a well-defined
and injective  homomorphism $\Aut(A)/\Aut_{\sf{in}}(A)\into \pic(\Lmod A)$.
\end{prop}

\begin{proof} Let $u\in A^\times$ be an invertible element
and $g=g_u: a\mapsto u\cdot a\cdot u\inv$, the corresponding
inner automorphism of $A$.
It is straightforward to verify that the
map $\varphi: x\mapsto u\cdot x$ yields an isomorphism
of $A$-bimodules $\varphi: A\iso A^{g_u}$. Conversely,
given $g\in \Aut(A)$ and an  $A$-bimodule isomorphism
$\varphi: A\iso A^g$, put $u:=\varphi(1)$.
Then we have $\varphi(a)=\varphi(1\cdot a)=\varphi(1)\cdot a=u\cdot a$.
It follows that $g=g_u$, moreover, $u$ is invertible
since $u\inv=g\inv(1)$.
This completes the proof.
\end{proof}

It is instructive to think of the group $\pic(\Lmod A)$
as some sort of Lie group. The `Lie algebra' of that group
should therefore be formed by $A$-bimodules that
are `infinitesimally close' to $A$, the unit of
$\pic(\Lmod A)$. It may be argued that any right
$A$-module that is `infinitesimally close'  to a rank one free
 right
$A$-module is itself isomorphic to  a rank one free
 right
$A$-module. Thus, any object in
$\Lie \pic(\Lmod A)$ may be viewed as being
the 
right $A$-module $A$ on which the standard 
left multiplication-action is  `infinitesimally deformed'.
Denote this  `deformed' left action of  an element $a\in A$ 
by $b\mapsto a*b$. The deformed left action must commute with
the standard right action, hence, for any $a,b,c\in A$, we must have
$(a*b)\cdot c=a*(b\cdot c)$. This forces the deformed
action to be of the form $a*b=F(a)\cdot b$, where
$F: A \to A$ is a certain linear map that should be
`infinitesimally close' to the identity $\Id: A\to A$.
We express the latter condition by
introducing the ring $\beps$ of dual numbers
and writing the map $F$ in the form
$F(a)= a+ \eps\cdot \psi(a)$, as we have already done earlier.
Now the condition that the assignment 
$A\times A \to \beps\otimes A\,,\,a,b\mto a*b=F(a)\cdot b$ gives a (left) action
 becomes  the equation $F(a\cdot a')=F(a)\cdot F(a'),\,\forall
a,a'\in A$. The last equation is equivalent to the condition
that $\psi$ is a derivation, see \eqref{inf_aut}.

Furthermore, arguing similarly,
one finds that any {\em inner} derivation $\ad a\in \Inn(A)$ gives rise
to a bimodule of the form $A^g$ where $g:  b\mapsto (1+\eps\cdot
a)\cdot b\cdot(1+\eps\cdot a)\inv$ is an 
`infinitesimal inner automorphism'.
Thus, we conclude that the `Lie algebra' of the group $\pic(\Lmod A)$
is given by the formula
$$\Lie\pic(\Lmod A)=\Der(A)/\Inn(A).$$

\subsection{Atiyah algebra of a vector bundle.}\label{atiyah_alg}
 Let $X$ be an affine variety,
and $\calE$ an algebraic vector bundle on $X$,
locally trivial in the Zariski topology.
Giving $\calE$ is equivalent to giving
$E:=\Gamma(X,\calE)$, the vector space of global
sections of $\calE$, which is a finitely generated {\em projective}
 $\C[X]$-module.

Write  $\End(\calE)$ for  the (noncommutative) associative algebra
of $\C[X]$-linear endomorphisms of the vector bundle $\calE$.
The algebra  $\End(\calE)$ acts naturally on
$E:=\Gamma(X,\calE)$.
The center of the algebra $\End(\calE)$ is the subalgebra
$\C[X]=\C[X]\cdot\Id_\calE\into \End(\calE)$, formed by
{\em scalar endomorphisms}.

The  following result provides a basic example of Morita equivalence
\begin{thm} The algebras $\C[X]$ and  $\End(\calE)$ are Morita
equivalent.
Specifically, the following functors
\begin{align*}
\xymatrix{
{\Lmod {\C[X]}}\ar@/^/[rrr]|{S}&&& {\Lmod
{(\End(\calE))}}\ar@/^/[lll]|{T}
}
\\
S: M\mto E\otimes_{\C[X]} M,\aand
T:\ F \mto \Hom_{_{\End(\calE)}}(E,F)
\end{align*}
provide mutually inverse equivalences.
\end{thm}

\begin{defn}\label{atiyah_def} The Lie algebra
$\calA(\calE):=\Der(\End(\calE)),$
the derivations of the associative algebra $\End(\calE),$ is 
called the Atiyah algebra of $\calE$.
\end{defn}

To obtain a more explicit description of the Atiyah algebra,
assume that the affine variety $X$ is smooth.
Let $\calD(\calE)$ be
the (associative) algebra of
algebraic differential operators acting
on sections of $\calE$.
This algebra has a natural increasing filtration
$0=\calD_{-1}(\calE)\sset
\calD_0(\calE)\sset\calD_1(\calE)\sset\ldots,$ by the
order of differential operator.
In particular, we have $\calD_0(\calE)=\End(\calE)$.

Write $\gr_\idot\calD(\calE) =\bigoplus_{i\geq 0}\,
\calD_i(\calE)/\calD_{i-1}(\calE)$ for 
 the  associated graded algebra.
Assigning the {\em principal symbol} to a  differential operator
gives rise to a canonical graded algebra isomorphism
$$\sigma:\ \gr\calD(\calE)\iso \End(\calE)\otimes_{_{\C[X]}}
\Sym^\hdot\calT(X).
$$
\begin{rem} Note that the algebra $\gr\calD(\calE)$
is {\em not} commutative unless $\calE$ has rank one,
i.e.,  unless $\calE$ is a line bundle.
\eer

The top row of the  diagram below 
is a natural short exact sequence
involving the principal symbol map on the
 space of first order differential operators.
$$
\xymatrix{
0\ar[r]&\End(\calE)\ar[r]&\calD_1(\calE)
\ar[r]^<>(0.5){\sigma}& \End(\calE)\otimes_{_{\C[X]}}\calT(X)\ar[r]&
0\\
0\ar[r]&\End(\calE)\ar[r]\ar@{=}[u]^{\id}
&
{\;\calD_1^\spad(\calE)\;}\,\ar@{_{(}->}[u]\ar[r]&
{\left(
{}^{\quad\text{  \footnotesize{Scalar}}}_
{\text{\footnotesize{endomorphisms}}}
\!\right)}\bigotimes_{\C[X]}\calT(X)\ar@{_{(}->}[u]^{\imath\otimes\id}\ar[r]
&
0}$$
The space
 $\calD_1^\spad(\calE)$ in the bottom
row is formed by first order differential operators with
{\em scalar} principal symbol.

The Theorem below shows that the space  $\calD_1^\spad(\calE)$
is closely related to the Atiyah algebra
$\calA(\calE)$.

\begin{thm}\label{der_atiyah} \vi The space  $\calD_1^\spad(\calE)$
is a Lie subalgebra in the associative algebra $\calD(\calE)$
(with respect to the commutator bracket).

\vii The bottom row in the diagram above is an extension of
Lie algebras; in particular, 
$\calD_0(\calE)=\End(\calE)\sset \calD_1^\spad(\calE)$
is a Lie ideal.

\viii The adjoint action of an element $u\in \calD_1^\spad(\calE)$
in the ideal $\End(\calE)$ gives a derivation, $\ad u,$ of
$\End(\calE)$ viewed as an  associative algebra.

\iv The assignment
$u\mapsto \ad u$ 
is a Lie algebra homomorphism whose
kernel is the space $\C[X]\sset \End(\calE)=\calD_0(\calE),$
of scalar endomorphisms. In particular, this space  $\C[X]$
is a Lie ideal in $\calD_1^\spad(\calE)$.

\vv The adjoint action described in {\sf{(iii)}} gives rise to the following
Lie algebra exact sequence
$$
0\map \C[X]\too 
\calD_1^\spad(\calE)\;\stackrel{u\mapsto \ad u}\too\;
\Der\End(\calE)
\map 0.
$$
Thus, the map $u\mapsto \ad u$ induces an isomorphism
$\calD_1^\spad(\calE)/\C[X]\iso \calA(\calE).$
\end{thm}
\begin{proof} Observe first that the Theorem may be seen as a
unification of two special cases:

\npb{For $\calE=\oo_X,$ the trivial rank one bundle, 
we have $ \End(\calE)=\C[X]$, and the
derivations of the latter algebra are  given
by vector fields on $X$, by definition.}

\npb{If $X=\{pt\}$ is a single point, then we have $ \End(\calE)=\End(E)
\simeq\Mat_n\C,$ is a matrix algebra, and any derivation
of the matrix algebra is well-known to be inner.}

Thus, our Theorem says, essentially, that in the general case
of an arbitrary vector bundle on a variety $X$,
the Lie algebra of  derivations of $ \End(\calE)$
is an extention of the Lie algebra of vector fields
by inner derivations.

To prove the Theorem, we use Morita invariance
of  Hochschild cohomology of an algebra, cf. 
\S5 below. In particular, applying this for
the 1-st  Hochschild cohomology of  two Morita equivalent
algebras, $A$ and $B$, we get $\Der(A)/\Inn(A)\simeq \Der(B)/\Inn(B).$

Now, take $A=\C[X]$ and $B=\End\calE$.
Since $A$ is commutative, we have $\Inn(A)=0$,
hence $\Der(A)/\Inn(A)=\Der(A)=\calT(X)$.
On the other hand, inner derivations of the algebra
$B$ form the Lie algebra $B/\sZ_B\simeq \End(\calE)/\C[X]$.
Thus, the isomorphism of the previous paragraph yields
a short exact sequence
$$ 
\frac{\End(\calE)}{\C[X]}=\Inn(B)\into \Der(B)
\onto \frac{\Der(B)}{\Inn(B)}=\frac{\Der(A)}{\Inn(A)}=\calT(X).
$$
We leave to the reader to verify that this  short exact sequence
coincides with the bottom row of the
diagram preceeding Theorem \ref{der_atiyah}.
\end{proof}

\section{The Bar Complex}
\subsection{Free product of algebras}\label{free_product} Given two associative algebras,
$A$ and $B$, let $A*B$ be the $\C$-vector space whose basis is formed by
 words in elements of $A$ and $B$, with additional
relations that adjacent elements from the same algebra are multiplied
together. If, in addition, the algebras have units $1_A\in A$ and
$1_B\in B$,
then we impose the relation that multiplication  by either unit
acts as identity.

More generally, given two algebras $A,B$, and algebra
imbeddings $\imath: C\into A,$ and $\jmath: C\into B,$
such that $\imath(1_C)=1_A,$ and $\jmath(1_C)=1_B$,
one defines $A\,{*_{_C}}B,$
the free product of  $A$ and $B$ over $C$, as the
following
unital associative algebra
$$ A\,{*_{_C}}B:=\frac{T_\C(A\oplus B)}{\left\langle\!\!\left\langle
\begin{array}{c}
a\otimes a' = a\cdot a'\,,\,
b\otimes b'=
b\cdot b'\\
\imath(c)=\jmath(c)\,,\,
1_A=1=1_B
\end{array}
\right\rangle\!\!\right\rangle}_{a,a'\in A,\,
b,b'\in B,\,c\in C}
$$
where $\langle\!\langle\ldots\rangle\!\rangle$
denotes the two-sided ideal generated by the indicated relations.

The operation of free product of associative algebras 
plays a role somewhat analogous to the role of tensor product
for {\em commutative}  associative algebras.

\subsection{}
Throughout, we let $A$ be an associative $\C$-algebra (with $1$ as
usual).  
We wish to associate to $A$ a sequence of homology groups which will
play the noncommutative role of (co-)homology of a space.  In order to
do this, we wish to construct a particular 
resolution of $A$ by free $A$-bimodules.  
Before beginning, 
we remark that an $A$-bimodule is the same thing as a left $A\otimes_\C
A^\op$-module; if $M$ is an $A$-bimodule, then 
we define the action of $A\o A^\op$ on $M$ by
$
(a_1\otimes a_2^\op)m:=a_1ma_2.
$

We consider the following complex of $A$-bimodules
$$
\xymatrix{\cdots\ar[r]^<>(.5){b}&A^{\otimes  4}\ar[r]^<>(.5){b}&A^{\otimes
3}\ar[r]^<>(.5){b}
&A^{\otimes 2}=A\o A\ar[r]^<>(.5){m}&A\ar[r]&0},
$$
where $m\colon A\o A\to A$ is the multiplication on $A$ and $b\colon A^{\otimes  (n+1)}\to A^{\otimes   n}$ is given by
\begin{align*}
b(a_0\otimes a_1\otimes\cdots\otimes a_n)&:=a_0a_1\otimes a_2\otimes\cdots\otimes a_n+\\
&\qquad+\sum_{j=1}^{n-1}(-1)^ja_0\otimes\cdots\otimes(a_ja_{j+1})\otimes\cdots\otimes a_n.
\end{align*}

It is a tedious (but simple) calculation that $b^2=0$. 

\begin{defn} We set $A\ee:=A\o A^\op$, and for any $i=0,1,2,\ldots,$
put $\B_iA:= A\otimes A^{\otimes i}\otimes A$,
a free $A\ee$-module generated by the $\C$-vector space
$A^{\otimes i}$.  This way, the complex above can be writen
as the following \emph{bar complex}
$$\B_\idot A\colon\;\,
\big[\xymatrix{\cdots\ar[r]^<>(.5){b}&\B_2A\ar[r]^<>(.5){b}&\B_1A
\ar[r]^<>(.5){b}
&\B_0A=A\o A}\big]\,\onto A.
$$
\end{defn}

  We claim that this sequence is exact,
i.e., that the bar complex provides a free
$A$-bimodule resolution of $A$, viewed as an $A$-bimodule.
 To show this, we will construct a chain homotopy $h\colon A^{\otimes   i}\to A^{\otimes  (i+1)}$ such that $b\circ h+h\circ b=\id$.  Usual homological algebra then implies exactness.  We will first construct $h$ ``by hand,'' then give alternate descriptions of the bar complex which will make this definition more natural (and the proof of the homotopy easier).  In particular, we define
$$
h(a_1\otimes\cdots\otimes a_i)=1_A\otimes a_1\otimes\cdots\otimes a_i.
$$

We check that $h$ indeed satisfies $b\circ h+h\circ b=\id$ 
on the degree $2$ piece.  That is, we will show that $h(b(a_1\otimes a_2))+b(h(a_1\otimes a_2))=a_1\otimes a_2$.  Now, by definition $b(a_1\otimes a_2)=a_1a_2$, and $h(a_1a_2)=1_A\otimes a_1a_2$.  For the second term, $h(a_1\otimes a_2)=1_A\otimes a_1\otimes a_2$, and by the definition of $b$ we have
$$
b(1_A\otimes a_1\otimes a_2)=1_Aa_1\otimes a_2-1_A\otimes a_1a_2=a_1\otimes a_2-1_A\otimes a_1a_2.
$$
So,
$$
(h\circ b+b\circ h)(a_1\otimes a_2)=1_A\otimes a_1a_2+a_1\otimes a_2-1_A\otimes a_1a_2=a_1\otimes a_2,
$$
as desired.

\subsection{Second construction of the bar complex (after Drinfeld).}\label{drin_bar}
Let $A*\C[\e]$ be the free product of the algebra $A$ and the
polynomial algebra $\C[\e]$ in one variable $\e$.
An element of this free product can be written in the form
$a_1\e^{n_1}a_2\e^{n_2}\cdots a_k$ for elements $a_1,\ldots,a_k\in A$
and non-negative integers $n_1,\ldots,n_{k-1}$ (if the last factor is
from $\C[\e]$, simply ``pad'' with $1_A$, and similarly for the first
factor).  However, we can rewrite $\e^{n_j}=1_A\e1_A\cdots1_A\e1_A$,
where we have $n_j$ factors of $\e$.  
So, we can always write any element of $A*\C[\e]$ in the form $a_1\e a_2\e\cdots \e a_k$ for some $a_1,\ldots,a_k\in A$.  So, the $\e$ plays no role other than a separator, and we shall replace it by a bar, that is, we define $a_1\mid a_2\mid\cdots\mid a_k:=a_1\e a_2\cdots\e a_k$.  This notation is the genesis of the name ``bar complex.''

Now, we put a grading on $A*\C[\e]$ by declaring $\deg a=0$ for all $a\in A$ and $\deg\e=-1$.  We will now make $A*\C[\e]$ a differential graded algebra (DGA from now on) by defining a super-differential $d\colon A*\C[\e]\to A*\C[\e]$.  Recall that a super-differential is simply a super-derivation $d$ satisfying $d^2=0$.  We define $d$ on generators, namely, we set $da=0$ for all $a\in A$ and $d\e=1_A$, and we extend it to $A*\C[\e]$ uniquely by requiring that it obey the graded Leibniz rule.  If we now identify $A\e A\e\cdots\e A$ ($n$ factors of $A$) with $A^{\otimes   n}$ in the obvious fashion, we obtain an identification of $\B_\idot A$ with $A*\C[\e]$, and the super-differential $d$ on $A*\C[\e]$ becomes the bar differential.

Observe next that $d^2(\e)=d(1)=0$, and also $d^2(a)=0$, for
any $a\in A$, by definition.
Hence, Corollary \ref{dif_cor} yields $d^2=0$. We conclude that
$d$ is a differential on $A*\C[\e]$.

In this context, the proof of exactness becomes trivial.  Since $A*\C[\e]$ is a DGA, we can calculate its cohomology in the usual fashion.  The claim that $A*\C[\e]$ is exact (i.e., that $\B_\idot A$ is exact) is equivalent to claiming that $H^\hdot(A*\C[\e])=0$.  Now, by definition we have that $d\e=1_A$, hence $1_A$ is a coboundary and therefore zero in cohomology.  But since $A*\C[\e]$ is a DGA, $H^\hdot(A*\C[\e])$ is an algebra (it's even graded, but that is not important for our purposes).  In particular, the cohomology class of $1_A$ acts as a multiplicative unit for $H^\hdot(A*\C[\e])$, and since $[1_A]=0$, we see that $H^\hdot(A*\C[\e])=0$.  We can of course mimic the previous proof and construct a (co-)chain homotopy.  In the notation of $A*\C[\e]$, the homotopy $h$ is defined by $h(u)=\e u$ for all $u\in A*\C[\e]$.  Then $d(\e u)=d\e u-\e\,du=1_Au-\e\,du$.  So, $(d\circ h+h\circ d)(u)=u-\e\,du+\e\,du=u$, as required.

\subsection{Third construction of the bar complex.}
The following construction only applies (as presented) in the case where $A$ is finite dimensional.  Let $A^*=\Hom_\C(A,\C)$ be the dual vector space of $A$.  Then we form its non-unital tensor algebra
$$
T^+(A^*):=A^*\oplus(A^*\o A^*)\oplus\cdots\oplus(A^*)^{\otimes   k}\oplus\cdots.
$$
This is the free associative algebra on $A^*$ with no unit.

\begin{prop}
Giving an associative algebra structure on $A$ is equivalent to 
giving a map $d\colon T^+(A^*)\to T^+(A^*)$ such that

\vi $\;d$ is a super-derivation of degree $1$, i.e., $d((A^*)^{\otimes
k})
\subset (A^*)^{\otimes  (k+1)}$; and

\vii $\;d^2=0$.
\end{prop}

\begin{proof}
First, suppose we have any linear map $d\colon A^*\to A^*\o 
A^*$.  Then we can always extend this map to a super-derivation on
$T^+(A^*)$ by applying the super-Leibniz rule.  Now, if we are given a
multiplication map $m\colon A\o A\to A$, then by taking 
transposes we obtain a map $m^{_{\!\top}}\colon
A^*\to (A\o A)^*$.  Since $A$ is finite dimensional,
$(A\o A)^*$ and $A^*\o A^*$ are canonically isomorphic.
So, we can regard the transpose of multiplication as a map
$m^{_{\!\top}}\colon A^*\to A^*\o A^*$.  We let
$d=m^{_{\!\top}}$ and extend this to a super-derivation as described
above.  Of course, if we are already given a super-derivation of
$T^+(A^*)$, then we can transpose its restriction to 
$A^*\o A^*$ to obtain a 
multiplication on $A$ (since we can identify $A^{**}$ with $A$ and 
$(A^*\o A^*)^*$ with $A\otimes A$).

We will now show that the associative law for $m$ is equivalent to
$d^2=0$ where $m$ and $d$ are related as in the previous paragraph.
Since $T^+(A^*)$ is generated by elements of $A^*$, we need only show
that $d^2\colon A^*\to A^*\o A^*\o A^*$ is the zero 
map--the super-Leibniz rule will take care of the rest.  So, suppose we
are given some linear functional $\lambda\in A^*$.  Then $d\lambda\in
A^*\o A^*$, which we have 
canonically identified with $(A\o A)^*$.  In particular,
$d\lambda(a\otimes b)=m^{_{\!\top}}\lambda(a\otimes
b)=\lambda(m(a,b))=\lambda(ab)$ (where we will write $m(a,b)=ab$ for
simplicity's sake).  Further, since $d\lambda\in A^*\o A^*$, we
can find $\mu_1,\ldots,\mu_n,\nu_1,\ldots,\nu_n\in A^*$ so that
$d\lambda=
\sum_{i=1}^n(\mu_i\otimes\nu_i)$.  Then we have that
$$
d\lambda(a\otimes b)=\sum_{i=1}^n(\mu_i\otimes\nu_i)(a\otimes b)=\sum_{i=1}^n\mu_i(a)\nu_i(b).
$$

Now, we wish to consider $d(d\lambda)$.  Since this lies in $A^*\o A^*\o A^*$, we can view this is a linear functional on $A^{\otimes  3}$.  Then, using the fact that $d$ is a super-derivation, we obtain
\begin{align*}
d(d\lambda)(a\otimes b\otimes c)&=d\left[\sum_{i=1}^n(\mu_i\otimes\nu_i)\right](a\otimes b\otimes c)\\
&=\sum_{i=1}^n(d\mu_i\otimes\nu_i-\mu_i\otimes d\nu_i)(a\otimes b\otimes c)\\
&=\sum_{i=1}^nd\mu_i(a\otimes b)\nu_i(c)-\sum_{i=1}^n\mu_i(a)d\nu_i(b\otimes c)\\
&=\sum_{i=1}^n\mu_i(ab)\nu_i(c)-\sum_{i=1}^n\mu_i(a)\nu_i(bc)\\
&=\lambda((ab)c)-\lambda(a(bc))=\lambda((ab)c-a(bc)).
\end{align*}
So, if $d^2=0$, we see that $\lambda((ab)c-a(bc))=0$ for every
functional $\lambda$, hence $(ab)c=a(bc)$ (i.e., $m$ is associative).
Conversely, if $m$ is associative, 
we see that $d^2=0$.
\end{proof}

This is not quite the bar complex, since for one instance the
differential goes in the wrong direction.  But since $d\colon
T^i(A^*)\to T^{i+1}(A^*)$, we can consider the transpose of $d$,
$d^\top\colon(T^{i+1}(A^*))^*\to(T^i(A^*))^*$.  Since $A$ was assumed
finite dimensional, we can identify $T^i(A^*)^*$ and $T^i(A)$, so
finally $d^\top\colon T^{i+1}(A)\to 
T^i(A)$ is the desired complex.  One advantage of this approach is that
similar constructions can be made for different algebraic structures.

The reader is invited to 
consider how such a construction can be performed for, say Lie algebras.

\subsection{Reduced Bar complex.}
It turns out that the Bar complex $\B_\idot(A)$
contains  a large acyclic subcomplex.
Specifically, for each $n>1$,
in $\B_n{A}= A\otimes A^{\otimes n}\otimes A$ consider
the following $A$-subbimodule 
$$
\Triv_n{A}:= \sum_{i=1}^n A\otimes \bigl(A^{\otimes (i-1)}\otimes\C\otimes
A^{\otimes (n-i)}\bigr)\otimes A\;\sset \;A\otimes A^{\otimes n}\otimes A.
$$
It is easy to see from the formula for the bar-differential
that $b(\Triv_n{A})\sset\Triv_{n-1}{A}$.
Thus, $\Triv_\idot{A}$ is a subcomplex in $\B_\idot{A}$.

Define the {\em  reduced bar complex}
to be the quotient
$\overline{\B}_\idot{A}:=\B_\idot{A}/\Triv_\idot{A}$.
By definition,  the $n^{\text{th}}$ term 
of the  reduced bar complex is
$$\overline{\B}_n{A}=\B_n{A}/\Triv_n{A}=A\otimes\bar A^{\otimes
n}\otimes A,
$$
where $\bar A=A/\C$ as a vector space.  
The differential is the one induced by $d$ on the bar complex.  

The reduced bar  complex has the following interpretation in terms
of the free product construction $\B_\idot{A}=A*\C[\e]$,
see \S\ref{drin_bar}. 
Observe that since
$\deg\e=-1$, the Leibniz formula for an odd derivation yields
$$d(\e^2)=(d\e)\cdot\e-\e\cdot(d\e)=1\cdot\e-\e\cdot 1=0.$$
It follows readily that the two-sided ideal
$\langle\!\langle\e^2\rangle\!\rangle\sset A*\C[\e]$,
generated by the element $\e^2$ is $d$-stable,
i.e.,
we have $d\bigl(\langle\!\langle\e^2\rangle\!\rangle\bigr)\sset
\langle\!\langle\e^2\rangle\!\rangle$.
Hence, the differential $d$ on $A*\C[\e]$
descends to a well-defined differential
on the graded algebra
$(A*\C[\e])\big/\langle\!\langle\e^2\rangle\!\rangle$.
We claim that, under the identification of \S\ref{drin_bar},
we have $\Triv^\hdot{A}=\langle\!\langle\e^2\rangle\!\rangle$,
and therefore
\beq{reduced_bar}
\overline{\B}_\idot{A}=
(A*\C[\e])\big/\langle\!\langle\e^2\rangle\!\rangle=
A*(\C[\e]/\e^2).
\end{equation}

To see this, notice that an element 
of $ A\otimes A^{\otimes n}\otimes A=\B_n{A}$
belongs to $\Triv_n{A}$
if and only if it is a $\C$-linear combination of terms,
each involving
a subexpression like $(\ldots|1_A|\ldots)$.
But such an expresion, when translated into the free product construction,
reads: $(\ldots\e 1_A\e\ldots)=(\ldots\e^2\ldots)$,
and our claim follows.

Observe next that the argument proving acyclicity of the
complex $(A*\C[\e]\,,\,d)$ applies verbatim
to yield acyclicity of
$(A*\C[\e]\big/\langle\!\langle\e^2\rangle\!\rangle\,,\,d)$.
We conclude
\begin{cor}
The reduced bar complex
provides a free $A$-bimodule  resolution of $A$.\qed
\end{cor}
\section{Hochschild homology and cohomology}
\subsection{}
Given an associative $\k$-algebra 
$A$, let $A^\op$ denote
  the opposite algebra, and
write $A\ee=A\otimes A^\op$. There is a canonical
isomorphism $(A\ee)^\op\cong A\ee$. Thus,
an $A$-bimodule is the same thing as a left
$A\ee$-module, and also the same thing as a right $A\ee$-module.
Recall the notation $\bimod A$ for the abelian category of all
$A$-bimodules.  

Given an $A$-bimodule $M$, we write $[M,A]$ for the {\em commutator space},
the $\k$-vector subspace of $M$ spanned by all commutators
$ma-am,\, a\in A,m\in M$. Let $M/[M,A]$ denote the corresponding
 {\em commutator quotient} of $M$.

Equivalently, any $A$-bimodule $M$ may be viewed as either  right
or left
 $A\ee$-module. In particular, we view  $M$ as a right
$A\ee$-module and view $A$ as a left $A\ee$-module, and 
form the tensor product $M\otimes_{A\ee}A$. One has a
 canonical
vector space identification $M\otimes_{A\ee}A=M/[M,A]$.  The
assignment $M\mto  M\otimes_{A\ee}A$ clearly gives a  right exact functor from
 $\bimod A$ to the category of $\C$-vector spaces.
  So, $M\otimes_{A\ee}-$ has left derived functors, the Tor functors
 $\Tor_i^{A\ee}(M,A)$.  For ease of notation, we denote
 $\Tor_i^{A\ee}(M,A)$ by $\HH_i(M)$, the 
\emph{$i^{\text{th}}$ Hochschild homology group} (which is really a $\C$-vector space) of $M$ over $A$.

By definition, for any $A$-bimodule $M$,
we have $$\HH_0(M)=\Tor_0^{A\ee}(M,A)=M/[M,A].$$

In the particular case that $M=A$, we obtain $\HH_0(A)=A/[A,A]$.  Notice that $[A,A]$ is not an ideal in $A$, simply a $\C$-linear subspace of $A$.

Computing
 higher degree Hochschild homology groups requires a  choice of some
projective resolution of $A$ as $A$-bimodules (i.e., left
$A\ee$-modules).  The bar complex provides a canonical choice of such
resolution.
  So, to compute  the groups
 $\HH_i(M)$ we need only apply the functor $M\otimes_{A\ee}-$ to
the bar complex
$$
\B_\idot A\colon\xymatrix{\cdots\ar[r]&A^{\otimes  4}\ar[r]&A^{\otimes  3}\ar[r]&A^{\otimes  2}\ar[r]&0}.
$$
We then tensor this on the left with $M$ over $A\ee$ to yield
$$
M\otimes_{A\ee}\B_\idot\colon\xymatrix{\cdots M\otimes_{A\ee}A^{\otimes  4}\ar[r]&M\otimes_{A\ee}A^{\otimes  3}\ar[r]&M\otimes_{A\ee}A^{\otimes  2}\ar[r]&0}.
$$

To simplify this, pick some $m\otimes(a_0\otimes\cdots\otimes a_n)\in
M\otimes_{A\ee}A^{\otimes   n}$ and write $a_0\otimes
a_1\otimes\cdots 
a_{n-1}\otimes a_n=(a_0\otimes a_n^\op)(1_A\otimes a_1\otimes\cdots
\otimes a_{n-1}\otimes1_A)$. Then
\begin{align*}
m\otimes(a_0\otimes a_1\otimes\cdots\otimes a_{n-1}\otimes a_n)&=m(a_0\otimes a_n^\op)(1_A\otimes a_1\otimes\cdots a_{n-1}\otimes1_A)\\
&=a_nma_0\otimes(1_A\otimes a_1\otimes\cdots\otimes a_{n-1}\otimes1_A).
\end{align*}
By dropping the two $1_A$'s and observing that only scalars on the intermediate $a_j$'s ``commute'' past $\otimes$, we can identify $M\otimes_{A\ee}A^{\otimes n}$ with $M\o A^{\otimes  (n-2)}$.  So, the complex $M\otimes_{A\ee}\B_\idot A$ becomes
$$
\xymatrix{\cdots\ar[r]&M\o A^{\otimes  2}\ar[r]&M\o A\ar[r]&M\ar[r]&0}.
$$
Examining the identification of $M\otimes_{A\ee}A^{\otimes   n}$ with
$M\o A^{\otimes 
(n-2)}$, we find that the differential for this complex is given by (again, using the bar notation)
\begin{align*}
d(m\mid a_1&\mid\cdots\mid a_n)=ma_1\mid a_2\mid\cdots\mid a_n+\\
&\sum_{i=1}^{n-1}(-1)^im\mid a_1\mid\cdots\mid(a_ia_{i+1})\mid\cdots\mid
a_n
+(-1)^na_nm\mid a_1\mid\cdots\mid a_{n-1}.
\end{align*}

We remark that the bar complex, resp., the reduced bar complex,
for $A$ can be recovered from the Hochschild chain complex,
 resp., reduced  Hochschild chain complex,
of the free rank one $A\ee$-module $A\ee=A\otimes A^\op$
(viewed as an $A$-bimodule) as follows
\beq{Hoch_bar}
\B_\idot{A}=C_\idot(A,\, A\otimes A^\op),
\quad\text{resp.,}\quad
\overline{\B}_\idot{A}=\overline{C}_\idot(A,\, A\otimes A^\op).
\end{equation}

\subsection{Hochschild cohomology}\label{S:SquareZero}
As before, we take $A$ to be an associative $\C$-algebra and $M$ is an $A$-bimodule.  
For Hochschild homology, we considered $\Tor_{A\ee}^\hdot(A,M)$.
Now we wish 
to consider $\Ext_{A\ee}^\hdot(A,M)$.  

Recall that $A\ee=A\o A^\op$.  We define 
the \emph{Hochschild cohomology} of $A$ with coefficients in $M$ 
by $\HH^\hdot(M):=\Ext_{A\ee}^\hdot(A,M)$.

\begin{prop}\label{Morita_hoch}
The functors $\HH_\idot$ and $\HH^\hdot$ are both  Morita invariant. In particular,
$$
\HH_\idot(A)=\HH_\idot(\Mat_r(A)),\quad\text{and}\quad
\HH^\hdot(A)=\HH^\hdot(\Mat_r(A))$$
where $\Mat_r(A)$ denotes $r\times r$-matrices over $A$.
\end{prop}
\begin{proof}
This is immediate from the definition:
\[
HH^i(A,A)=\Ext^i_{\Lmod{A^e}}(A,A) \quad\text{and}\quad
HH_i(A,A)=\Tor^i_{\Lmod{A^e}}(A,A),
\]
since, for Morita equivalent algebras,
the corresponding categories of bimodules are equivalent as abelian
categories, hence give rise to isomorphic Ext and Tor groups.

It is instructive, however, to give a direct computational proof
 for the zeroth order homology, that is, to show that
 $\HH_0(A)=\HH_0(\Mat_r(A))$.  
We wish to construct a canonical isomorphism
$$
A/[A,A]\simeq\Mat_r(A)/[\Mat_r(A),\Mat_r(A)].
$$
We will define an isomorphism $\Tr\colon\Mat_r(A)/[\Mat_r(A),\Mat_r(A)]\to A/[A,A]$ by using, as the notation suggests, the standard trace on $\Mat_r(A)$.  In particular, $\Tr$ sends a matrix $(a_{ij})$ to the element
$$
\tr(a_{ij})\bmod[A,A]=\sum_{i=1}^na_{ii}\bmod[A,A].
$$
First, observe that $\Tr(XY)=\Tr(YX)$ for $X,Y\in\Mat_r(A)$, so $\Tr$
factors through a map (also called $\Tr$)
$\Mat_r(A)/[\Mat_r(A),\Mat_r(A)]\to A/[A,A]$.  We wish to show that the
kernel of $\Tr$ is precisely $[\Mat_r(A),\Mat_r(A)]$.  So, choose any
matrix $X=(a_{ij})\in\Ker\Tr$.  We will let $E_{ij}(a)$ denote the
elementary matrix with $a$ in the $ij$-entry and zeroes elsewhere.  Then
for $i\ne j$, we easily find that $E_{ij}(a)=[E_{ij}(a),
E_{jj}(1_A)]$.  So, all matrices with only off diagonal entries lie in $[\Mat_r(A),\Mat_r(A)]$.  Hence we can write
$$
X=(a_{ij})=\diag(a_1,a_2,\ldots,a_r)\bmod[\Mat_r(A),\Mat_r(A)].
$$
Also for $i\ne j$, we can directly compute that
$[E_{ij}(a),E_{ji}(1_A)]=E_{ii}(a)-E_{jj}(a)$.  So, for each
$j=2,\ldots,r$, the matrix $E_{11}(a_j)-E_{jj}(a_j)\in
[\Mat_r(A),\Mat_r(A)]$, hence 
$$
\diag\left(\sum\nolimits_{j=1}^ra_j,0,\ldots,0\right)-
\diag(a_1,\ldots,a_r)\in[\Mat_r(A),\Mat_r(A)].
$$
So, we can write any matrix in $\Mat_r(A)$ as $E_{11}(a)\bmod[\Mat_r(A),
\Mat_r(A)]$ for some $a\in A$.  Since, 
$$
0=\Tr(X)=\Tr(E_{11}(a))=a\bmod[A,A],
$$
we find that  $a\in[A,A]$ and $X\in[\Mat_r(A),\Mat_r(A)]$ to begin with.  So, $\Ker\Tr=[\Mat_r(A),\Mat_r(A)]$, and since $\Tr$ is clearly surjective we obtain that
$$
\Tr\colon\Mat_r(A)/[\Mat_r(A),\Mat_r(A)]\to A/[A,A]
$$
is an isomorphism.
\end{proof}

To calculate this cohomology,
 we will again use the bar complex $\B_\idot A$:
$$
\xymatrix{\cdots\ar[r]^<>(.5){b}&
A^{\otimes4}\ar[r]^<>(.5){b}&A^{\otimes3}\ar[r]^<>(.5){b}&A^{\otimes2}\ar[r]&0}.
$$
Applying the functor $\Hom_{A\ee}(-,M)$ to this complex (and accounting for contravariance), we obtain the sequence
$$
\xymatrix{\Hom_{A\ee}(A^{\otimes4},M)&\Hom_{A\ee}(A^{\otimes3},M)\ar[l]_<>(.5){b}
&\Hom_{A\ee}(A^{\otimes2},M)\ar[l]_<>(.5){b}&\ar[l]0},
$$
Recall that $A^{\otimes2}$ is free of rank one as an $A$-bimodule (which is just a left $A\o A\ee$-module).  So, $\Hom_{A\ee}(A^{\otimes2},M)\simeq M$.  Similarly, if $\f\colon A^{\otimes3}\to M$ is an $A$-bimodule map, then
$$
\f(a_1\otimes a_2\otimes a_3)=a_1\f(1_A\otimes a_2\otimes1_A)a_3.
$$
The association of $\f$ to the $\C$-linear map $A\to M$,
$a\mapsto\f(1_A\otimes a\otimes1_A)$ 
gives an isomorphism between $\Hom_{A\ee}(A^{\otimes3},M)$ and $\Hom_\C(A,M)$.  Indeed, we find that $\Hom_{A\ee}(A^{\otimes n},M)\simeq\Hom_\C(A^{\otimes(n-2)},M)$ for all $n\ge3$.  So, the complex whose cohomology we wish to compute reduces to
$$
\xymatrix{\cdots&\ar[l]_<>(.5){b}\Hom_\C(A^{\otimes2},M)&
\ar[l]_<>(.5){b}\Hom_\C(A,M)&\ar[l]_<>(.5){b}M&\ar[l]0}.
$$

An explicit formula for $b$ in this interpretation is as follows.
  Recall that for the calculation of Hochschild homology, we used the differential
\begin{align*}
&b(a_0\otimes a_1\otimes\cdots\otimes a_n)=a_0a_1\otimes\cdots\otimes
a_n\\
&+\sum_{j=1}^{n-1}a_0\otimes
a_1\otimes\cdots\otimes(a_ja_{j+1})\otimes\cdots\otimes a_n
+(-1)^na_na_0\otimes a_1\otimes\cdots\otimes a_{n-1},
\end{align*}
where we have $a_0\in M$ and $a_j\in A$ for $j\ge1$.  For the first
  differential, $
b\colon M\to\Hom_\C(A,M)$, we find that we must have $b m(a)=am-ma$.
  For $n=1$, we obtain for $f\in\Hom_\C(A,M)$ that $b f(a_1,a_2)=a_1f(a_2)-f(a_1a_2)+f(a_1)a_2$.  Similarly, for all $f\in\Hom_\C(A^{\otimes n},M)$, we find that
\begin{align*}
bf(a_1,a_2,\ldots,a_{n+1})=a_1f(a_2,\ldots,a_{n+1})&+\sum_{i=1}^n(-1)^if(a_1,\ldots,a_ia_{i+1},\ldots,a_{n+1})\\
&+(-1)^{n+1}f(a_1,\ldots,a_n)a_{n+1}.
\end{align*}
We will often write $C^n(A,M):=\Hom_\C(A^{\otimes n},M)$.

With these 
formulae for $b$ in hand, we can now explicitly calculate the first few
cohomology groups.  

\subsection{Interpretation of  $\HH^0$.} For degree zero, we see that
$\HH^0(M)=\Ker b$.  But $m\in\Ker b$ if and only if $bm(a)=am-ma=0$ for
all $a\in A$.  It is natural to call $\{m\in M\st am=ma\}$ the
\emph{center} of the module $A$.  In particular, if $M=A$, 
then we in fact have that $\HH^0(A)=\sZ_A$.  Notice that there is a fair
bit of ``duality'' in the degree zero setting between Hochschild
homology and cohomology.  
For $\HH^0(M)$, we obtain the center of $M$, and for $\HH_0(M)$ we obtain the ``cocenter'' of $M$, namely $M/[A,M]$.

\subsection{Interpretation of  $\HH^1$.} Moving on to $\HH^1(A,M)$,
 we claim that the kernel of $d$ on $C^1(A,M)$ consists precisely of all derivations $A\to M$.  Indeed, if $f\in C^1(A,M)$, then $df(a_1,a_2)=0$ is equivalent to requiring that
$$
a_1f(a_2)-f(a_1a_2)+f(a_1)a_2=0,
$$
that is, that $f(a_1a_2)=a_1f(a_2)+f(a_1)a_2$.  So, $f$ is a derivation.  The image of $d\colon M\to\Hom_\C(A,M)$ is precisely the set of all \emph{inner} derivations, $\Inn(A,M)$, that is, derivations of the form $a\mapsto am-ma$ for some $m\in M$.  So,
$$
\HH^1(M)\simeq\Der(A,M)/\Inn(A,M),
$$
the \emph{outer derivations} from $A$ to $M$.  In particular, 
we can rewrite \eqref{inn_der} in the following way
\beq{derM}
0\to\HH^0(A,M)\map M
\stackrel{\ad}\map \Der(A,M)\map \HH^1(A,M)\to 0.
\end{equation}

Observe that if $M=A$, then $\Der(A)$ is a Lie algebra and $\Inn(A)$ is a Lie 
ideal.  So, $\HH^1(M)$ inherits the structure of a Lie algebra.

\subsection{Interpretation of  $\HH^2$.}
This is given by considering algebra extensions.  
Suppose now that $A$ is an associative algebra and consider an extension of algebras
$$
\xymatrix{0\ar[r]&M\ar[r]&\widetilde{A}\ar[r]&A\ar[r]&0},
$$
where $M$ is an ideal of $\widetilde{A}$ satisfying $M^2=0$.  Then $M$ is an
$A$-bimodule.  Define $a\cdot m$ by choosing a lift $\tilde{a}$ of $a$ in $\widetilde{A}$
and setting $a\cdot m=\tilde{a}m$.  If $\tilde{a}'$ is another such lift, then since
$\tilde{a}-\tilde{a}'$ maps to zero in $A$, it must be an element of $M$.  
Hence $(\tilde{a}-\tilde{a}')m\in M^2=0$, so $\tilde{a}m=\tilde{a}'m$.  

Choose some $\k$-linear splitting (as a vector space)
$c: A\into \widetilde{A}$. Then we get
 a vector space 
direct sum decomposition
$\widetilde{A}\simeq c(A)\oplus
 M$.  Further, for any $a_1,a_2\in A,$ we have 
 $c(a_1)\cdot c(a_2)-c(a_1\cdot a_2)\in M$.
Therefore we can write the product on $\widetilde{A}$ as
$$
(a_1\oplus m_1)(a_2\oplus m_2)= a_1a_2\oplus (a_1m_2+m_1a_2+\beta(a_1,a_2)),
$$
where $\beta$ is an arbitrary bilinear map $A^{\otimes2}\to M$, i.e., an
element of $C^2(A,M)$. 
This formula gives an associative product on
$\widetilde{A}$ if and only if $d\beta=0$.  Indeed,
$$
d\beta(a_1,a_2,a_3)=a_1\beta(a_2,a_3)-\beta(a_1a_2,a_3)+\beta(a_1,a_2a_3)-\beta(a_1,a_2)a_3.
$$
Then by writing out $[(a_1\oplus m_1)(a_2\oplus m_2)](a_3\oplus m_3)$
 and $(a_1\oplus m_1)[(a_2\oplus m_2)(a_3\oplus m_3)]$ and equating
 them, we see that $d\beta$ must equal zero.  We can also check that if
 $\beta_1$ and $\beta_2$ define associative structures on $\widetilde{A}$, then they
 define isomorphic extensions if and only if $\beta_1-\beta_2$
is a Hochschild coboundary.  So, $\HH^2(M)$ classifies extensions of $A$ by $M$.

\begin{rem}
If we consider the
case where $A$ is a unital algebra and we only consider the square-zero
extensions where $\widetilde{A}$ is also a unital algebra, and the homomorphism
$\widetilde{A}\to A$
takes the unit to the unit, then it is easy to check that such extensions
are also classified by $\HH^2(M),$ but where $\HH^2(M)$ is computed using the
{\em reduced} bar complex.
\eer
\subsection{Interpretation of  $\HH^3$.} The group
  $\HH^3(M)$ classifies so-called {\em crossed-bimodules}.  A crossed-bimodule is a map $\f\colon C\to B$ of associative algebras where $B$ is unital, $C$ is nonunital and maps to a $2$-sided ideal of $B$,  the cokernel of $\f$ is $A$, the kernel is $M$, and we require
$$
\f(bcb')=b\f(c)b'\quad\text{and}\quad
\f(c)c'=cc'=c\f(c'),
$$
for all $b,b'\in B$ and $c,c'\in C$.

\subsection{Reduced cochain complex}
One can use the reduced bar complex to compute Hochschild homology
and cohomology.  Specifically, for any $A$-bimodule
$M$, define the reduced Hochschild cochain complex $\overline{C}^\hdot(A,M)$ by
\begin{align*}
\overline{C}^n(A,M):=
\Hom_{\bimod A}(\overline{\B}_n{A},M)&=
\Hom_{\bimod A}(A\otimes\bar A^{\otimes n}\otimes A,M)\\
&\simeq\Hom_\C(\bar A^{\otimes n},M),\quad n=0,1,\ldots\,.
\end{align*}

\section{Poisson  brackets and
 Gerstenhaber algebras}
\subsection{Polyvector fields}
Let $X$ be an affine variety with coordinate ring
$A=\C[X],$ and let $\calE,\scr F,$ be  locally free coherent sheaves on
$X$. We write $E=\Gamma(X,\calE),$ resp., $F=\Gamma(X,\scr F),$
for the corresponding (projective) $A$-modules of global sections.
Then, one has canonical isomorphisms
$$E\otimes_A F\cong \Gamma(X,\,\calE\otimes_{_{\oo_X}}\scr F),
\quad \Lambda_A^p E\cong \Gamma(X,\,\Lambda^p \calE),
\quad \Sym_A^p E\cong \Gamma(X,\,\Sym^p \calE).
$$

In particular, for the tangent sheaf $\calT_X$
we have $\Gamma(X,\,\Sym\calT_X)=\k[T^*X],$
is the algebra
of regular functions on the total space of the cotangent bundle on $X$.

We introduce the notation
 $\Th_p(X):=\Gamma(X,\Lambda^p\calT_X)$
for the  vector space of $p$-polyvector fields on $X$.
The graded-commutative algebra
 $\Th_\idot(X):=\bigoplus_p\,\Th_p(X)=\Gamma(X,\Lambda^\hdot \calT_X)$ 
may be thought of as a {\em odd analogue} of the
commutative algebra $\Gamma(X,\Sym \calT_X)=\C[T^*X]$.

For any polyvector $\pi\in \Th_p$, one defines
a natural {\em contraction operator}
$i_\pi: \Om^k(X)\map \Om^{k-p}(X)\,,\,
\alpha\mapsto i_\pi\alpha,$
where, for $\pi=\xi_1\wedge\ldots\wedge\xi_p$,
the $(k-p)$-form $i_\pi\alpha$,
 is given
by $\eta_1,\ldots,\eta_{k-p}\mto
\alpha(\xi_1,\ldots,\xi_p,\eta_1,\ldots,\eta_{k-p}).$
For $p=1$ this reduces to the standard contraction
of a differential form with respect to a vector field.

Further, there is a natural \emph{Schouten bracket} on $\Th_\idot(X)$:
$$
\{-,-\}\colon\Th_p(X)\times\Th_q(X)\to\Th_{p+q-1}(X).
$$
given by the following formula
\begin{align}\label{Th_bracket}
&\{\xi_1\wedge\cdots\xi_p\,,\,\eta_1\wedge\cdots\eta_q\}=\\
&=\sum_{i,j=1}^{i=p,j=q}(-1)^{i+j}\cdot [\xi_i,\eta_j]\wedge
\xi_1\wedge\cdots\wh{\xi_i}\wedge\cdots\wedge\xi_p\wedge\eta_1\wedge\cdots\wedge\wh{\eta_j}\wedge\cdots\wedge\eta_q.\nonumber
\end{align}
If $p=q=1$, this formula reduces to the usual bracket of vector fields. 

To make the bracket $\{-,-\}$ compatible with the gradings, we note that
$(p+1)$ $+(q+1)-1= (p+q)-1$. Thus,
shifting the natural grading on $\Th_\idot(X)$
by $1$, we obtain a new graded strucre, to be
denoted $\Th_{\bullet-1}(X)$, such that
the bracket \eqref{Th_bracket} is compatible with this
new  graded structure.

Below, we summarize the
various  natural structures on
$\Th_\idot(X)$ and $\Om^\hdot(X)$.
\beq{properties}
\begin{array}{l}
\bullet\quad\pbox{\text{Wedge-product of   polyvector fields  makes $\Th_\idot(X)$  a 
super-commu-}\newline
\text{tative algebra.}}\\
\bullet\quad\pbox{\text{The  bracket \eqref{Th_bracket}
 makes $\Th_{\bullet-1}(X)$  a graded Lie super-algebra.}}\\
\bullet\quad\pbox{\text{The contraction operators $i_\pi,\,\pi\in \Th_\idot(X),$ 
make $\Om^\hdot(X)$ a graded}\newline
\text{module over  $\Th_\idot(X)$,
viewed as a super-commutative algebra.}}\\
\bullet\quad\pbox{\text{The De Rham differential
$d:\Om^\hdot(X)\to\Om^{\bullet+1}(X)$
is an odd derivation}\newline
\text{of the super-commutative algebra
$\Om^\hdot(X)$.}}
\end{array}
\end{equation}

Further, one extends the definition of Lie derivative
from vector fields to polyvector fields.
Specifically, given a polyvector field $\pi\in \Th_p(X),$ 
we define the {\em Lie derivative} operator
$$\L_\pi:=[d, i_\pi] :\ \Om^\hdot(X)\to\Om^{\bullet-p+1}(X).$$
This formula reduces, in the special case $p=1$,
to the classical {\em Cartan homotopy formula}
 for the Lie derivative.

Using Lemma \ref{triv}, one  verifies the following standard
identities
\begin{align}\label{standard_ident}
i_{\pi\wedge \varphi}= i_{\pi}\cdot i_{\varphi},\quad
\L_{\pi\wedge \varphi}= \L_{\pi}\cdot i_{\varphi}
+ (-1)^{\deg \pi} i_{\pi}\cdot \L_{\varphi},\\
[\L_{\pi},\L_{\varphi}]=\L_{\lbrace\pi,\varphi\rbrace},\quad
[i_{\pi},\L_{\varphi}]=i_{\lbrace\pi,\varphi\rbrace},\quad
[i_{\pi},i_{\varphi}]=0.\nonumber
\end{align}

\subsection{Poisson brackets.}
Let $A$ be an associative (not necessarily commutative) algebra.
A  (skew-symmetric)
Lie bracket $\{-,-\}:\ A \otimes A\to A$
is said to be a Poisson bracket if, for any $a\in A$,
the map $\{a,-\}: A\to A$ is a derivation, i.e., the following Leibniz identity
 holds
\beq{Leibniz_poiss}
\{a, b\cdot c\}= b\cdot\{a,  c\} +\{a, b\}\cdot c,\quad \forall
a,b,c\in A.
\end{equation}
 In view of skew-symmetry, the Leibniz identity says
that the bracket
 $\{-,-\}:\ A \otimes A\to A$ is a {\em bi-derivation} on $A$.

In case the  associative product on  $A$ is commutative
we will sometimes call $A$ a {\em commutative Poisson algebra}.
In the non-commutative case, we will say that $A$ has
 a NC-Poisson structure.

Let $A=\C[X]$ be the coordinate ring of  a {\em smooth} affine variety.
Then, it is easy to show that any bi-derivation on
$\C[X]$ is given by a bi-vector, that is,
there is a regular section $\pi\in\Gamma(X, \wedge^2\calT_X)$
such that the bi-derivation has the  form
$(f,g)\mto\langle df\wedge dg,\pi\rangle.$
Thus, giving a Poisson  structure on $\C[X]$ amounts to giving
a  bracket
\beq{bivect}
\{f,g\}=\langle df\wedge dg,\pi\rangle,\quad\text{where}\quad
\pi\in \Gamma(X, \wedge^2\calT_X)\quad\text{is such that}\quad[\pi,\pi]=0
\end{equation}
(the condition  $[\pi,\pi]=0$ on
the Schouten bracket is equivalent to
the Jacobi identity for the Poisson   bracket \eqref{bivect}).
In this case we refer to $\pi$ as a {\em Poisson bivector}.

Fix a smooth variety $X$ with a
A Poisson bivector $\pi$ on a smooth variety $X$ 
gives rise to a canonical Lie algebra structure on
$\Om^1(X)$, the space of 1-forms on $X$.
The corresponding Lie bracket on $\Om^1(X)$
is  given by the formula
\beq{pois_algebroid}
[\alpha, \beta]:= L_{\imath_\pi\alpha}\beta -L_{\imath_\pi\beta}\alpha
-d\langle\alpha\wedge\beta,\pi\rangle.
\end{equation}
This Lie bracket has the following 
properties (which uniquely determine the bracket):

\npb{The De Rham differential
$d: A=\C[X] \to \Om^1(X)$ is  a Lie algebra map, i.e.:}
$$[da,db]=d(\{a,b\}),\quad\forall a,b\in A;
$$
\npb{The following Leibniz identity holds:}
$$[\alpha, f\cdot\beta]=f\cdot[\alpha,\beta] +
(\imath_\pi\beta)f\cdot\alpha, \quad\forall f\in A, \alpha,\beta\in \Om^1(X).$$

\subsection{Gerstenhaber algebras.} The notation of Gerstenhaber
algebra  is an odd analogue of the
notion of Poisson algebra. 
More explicitly, 
\begin{defn}\label{gerst_def}
A Gerstenhaber algebra is a graded
super-commutative algebra $G^\hdot=\bigoplus_i\, G^i$ with a bracket
$$
\{-,-\}\colon G^p\times G^q\to G^{p+q-1}
$$
which makes $G^\hdot$ a Lie super-algebra so that for every 
$a\in G^\hdot$, the map $\{a,-\}$ 
is a super-derivation with respect to the product, i.e., we have
\beq{gerst_axiom}
\{a,\,b\cdot c\}=\{a,b\}\cdot c + (-1)^{(\deg a-1)\deg b}\,
b\cdot \{a,c\}.
\end{equation}
\end{defn}

This definition is motivated by the following basic
\begin{examp} The Schouten bracket makes $\Th_\idot(X)$,
the graded space of
polyvector fields on a smooth manifold $X$,
a Gerstenhaber algebra.
\end{examp}

As another example, let $X$ be a Poisson manifold with
Poisson bivector $\pi$, see \eqref{bivect}. Then, one proves
 the following
\begin{lem}\label{ger_om} The  bracket on 1-forms given
by formula \eqref{pois_algebroid} extends uniquely
to a bracket $[-,-]:\
\Om^i(X)\times\Om^j(X)\to\Om^{i+j}(X),$ such that
the wedge product and the bracket give $\oplus_{i\geq 1}\,\Om^i(X)$
the structure of Gerstenhaber algebra.\qed
\end{lem}

\begin{rem} Suppose $(X,\om)$ is a symplectic manifold.
The isomorphism $\calT_X\iso$ $\calT^*_X,\,\xi\mapsto\imath_\xi\om$
induces a graded algebra isomorphism
$\wedge^\hdot\calT_X\iso\wedge^\hdot\calT^*_X$. Therefore,
one can transport the  Gerstenhaber algebra structure
on $\wedge^\hdot\calT_X$ given by the Schouten bracket
to a  Gerstenhaber algebra structure on $\Om^\hdot(X)$.
The latter one turns out to be the same as
the Gerstenhaber algebra structure of Lemma \ref{ger_om}.
\eer

\begin{rem} Given a manifold $X$ and a bivector 
$\pi\in \Gamma(X, \wedge^2\calT_X),$  such that
$[\pi,\pi]=0$ one can also introduce a {\em different}
Lie bracket $B: \Om^p(X)\times\Om^q(X)\to\Om^{p+q-2}(X),$
which has  degree $-2$. For $\alpha=\alpha_1\wedge\cdots\wedge\alpha_p\in\Om^p(X)$
and 
$\beta=\beta_1\wedge\cdots\wedge\beta_q\in\Om^q(X)$, this new bracket
is defined by
$$
B(\alpha,\beta):=\sum(-1)^{i+j}\langle\alpha_i\wedge\beta_j,
\pi\rangle\alpha_1\wedge\cdots\wedge\wh{\alpha_i}\wedge\cdots\wedge\alpha_p\wedge\eta_1\wedge
\cdots\wedge\wh{\eta_j}\wedge\cdots\wedge\eta_q.
$$
\eer

Let $G^\hdot$ be  a  Gerstenhaber algebra,
and $M^\hdot$ a graded vector space.
We say that $M^\hdot$ is   a  Gerstenhaber module over  $G^\hdot$
if  the square-zero construction, $G^\hdot\sharp M^\hdot$,
is equipped with  a  Gerstenhaber algebra structure
such that $G^\hdot$ is a  Gerstenhaber subalgebra
in  $G^\hdot\sharp M^\hdot$ and, moreover,
we have $M^\hdot\cdot M^\hdot=0=\{M^\hdot,M^\hdot\}$.

\subsection{$\eps$-extension of a  Gerstenhaber algebra.}
Let $\beps$ be the ring of dual numbers.
Given a  Gerstenhaber algebra $G^\hdot=\bigoplus_i\, G^i,$
with operations $(-)\cdot(-)$ and $\{-,-\}$,
one defines, c.f., \cite{TT},
a new  Gerstenhaber algebra $G_\eps^\hdot=\bigoplus_i\, G_\eps^i$,
over the ground ring  $\beps$.
Specifically, put
$$G_\eps^i:=G^i\oplus\eps G^{i-1},$$
and introduce  $\beps$-bilinear operations
 $(-)\cdot_{\eps}(-)$ and $\{-,-\}_{\eps}$,
defined for any homogeneous elements $a,b\in G$ by 
\beq{TT_gerst}
a \cdot_{\eps} b := a\cdot b+ (-1)^{\deg a}\cdot\eps \cdot\{a,b\},\aand
\{a,b\}_{\eps}  := \{a,b\}.
\end{equation}

\begin{rem}
Note that in order to have the new product  $(-)\cdot_\eps(-)$ be
graded commutative, one has the additional term $(-1)^{\deg
a}\cdot\eps \cdot\{a,b\}$
to be ``symmetric'' in the graded sense. Such a construction
would have been impossible for an ordinary Poisson algebra,
where the dot-product is always symmetric and the 
 Poisson bracket is always {\em skew}-symmetric, so that the linear
combination of the two does not have a defininite symmetry.
\eer

Applying the `$\eps$-construction' above
to  $\Th_\idot(X)$, the   Gerstenhaber algebra
of polyvector fields on a manifold $X$,
one obtains a  Gerstenhaber algebra
 $\Th_\idot(X)_\eps.$

The identities in
\eqref{standard_ident}
are conveniently encoded in the  following result.

\begin{prop}[\cite{TT}]\label{TT_prop}
The following formulas
$$(\pi+\eps\varphi)\cdot_{\eps} \alpha:= (-1)^{\deg \pi}\,
i_\pi\alpha,\aand
\{\pi+\eps\varphi,\alpha\}_{\eps}  := \L_\pi\alpha+
\eps\cdot i_\varphi\alpha
$$
make $\Omega^\hdot(X)_\eps$ a Gerstenhaber module
over  $\Th_\idot(X)_\eps.$
\end{prop}

\subsection{Lie algebroids.}\label{algebroid_sec} Let $X$ be an algebraic variety
with structure sheaf $\oo_X$ and tangent sheaf $\calT_X$.
Let $\scr A $ be a coherent sheaf  of $\oo_X$-modules
equipped with a (not necessarily
$\oo_X$-bilinear)  Lie bracket $[-,-]: \scr A \times \scr A  \map \scr A $.
\begin{defn}\label{Liealgebroid} The data of a sheaf $\scr A $ as
above and an $\oo_X$-linear map $\tau: \scr A \to \calT_X,
\,v\mapsto\tau_v,$ (called {\em anchor map})
is said to give a Lie algebroid on $X$ if the following holds:

\pb{The map $\tau$ is a Lie algebra map;}

\pb{We have $[f\cdot v,u]= f\cdot [v,u] +\tau_u(f)\cdot v$
for any $v,u\in \scr A ,\,f\in \oo_X$.}
\end{defn}

\begin{examps} \textsf{(1)}\; The tangent sheaf $\calT_X$ equipped with
the standard Lie bracket on vector fields and with
the identity anchor map $\tau=\id: \calT_X \to\calT_X$
is a Lie algebroid.

\textsf{(2)}\; Any coherent sheaf of $\oo_X$-modules equipped
with an $\oo_X$-{\em bilinear} Lie bracket is a  Lie algebroid
with zero anchor map. In particular, a vector bundle on $X$
whose fibers form an algebraic family of Lie algebras
is a (locally-free)  Lie algebroid on $X$.

\textsf{(3)}\; For any  Lie algebroid $\scr A$, the kernel
of the anchor map $\scr K:=\Ker[\tau:\scr A \to\calT_X]$
is a sub Lie algebroid in $\scr A$. This Lie algebroid
$\scr K$ is of the type decribed in example (2) above.

\textsf{(4)}\; Proposition \ref{der_algebroid} may be conveniently expressed
by saying that,
for any algebra $A$ with center
$\sZ_A$, the space $\Der(A)$ is (the space of global
sections of) a  Lie algebroid on $\Spec \sZ_A$.

\textsf{(5)}\; As a special case of (4) we get:
 the Atiyah algebra of a vector bundle $\calE$ on $X$
 has a natural structure of  Lie algebroid on $X$,
since $\sZ_{_{\End(\calE)}}\simeq \C[X]$. In particular, 
the sheaf ${\scr D}_1$ of (ordinary) first order differential operators
on $X$ is a Lie algebroid on $X$.

\textsf{(6)}\; A Poisson bivector $\pi$ on a smooth
variety $X$ gives the cotangent sheaf
$\calT^*_X$ a Lie algebroid structure. The corresponding
Lie bracket of sections of $\calT^*_X$ is defined by formula
\eqref{pois_algebroid}, and the anchor map
$\calT^*_X\to\calT_X$
is given by  contraction $\alpha\mapsto \imath_\pi\alpha$.
\end{examps}

Fix  a locally free sheaf (a vector bundle)
$\scr A$ on $X$ and an $\oo_X$-linear map $\tau: \scr A\to\calT_X$.
Performing the Symmetric, resp. Exterior,
algebra construction to the sheaf $\scr A$ and to the map
$\tau$
we obtain the following 
  graded algebras and graded algebra morphisms:
$$\Sym^\hdot\scr A\to\Sym^\hdot\calT_X,\quad \Lambda^\hdot\scr A\to
\Lambda^\hdot\calT_X,\quad 
\Lambda^\hdot\scr A^*\to
\Lambda^\hdot\calT_X^*.$$

A Lie algebroid structure on $\scr A$ (with the
 anchor map $\tau$)
gives rise to the following additional structures
on the graded algebras above:

\pb{A Poisson algebra structure on $\Sym^\hdot\scr A$
with Poisson bracket given by}
\beq{kirillov_algebroid}
\{a_1\cdot\ldots\cdot a_p,
b_1\cdot\ldots\cdot b_q\}:=
\sum_{i,j=1}^{i=p,j=q}
[a_i,b_j]\cdot  a_1\cdot\ldots\cdot\wh{a_i}\cdot\ldots\cdot a_p\cdot b_1\cdot
\ldots\cdot\wh{ b_j}\cdot\ldots\cdot b_q.
\end{equation}

\pb{A Gerstenhaber algebra structure on $\Lambda^\hdot\scr A$
 given by formula \eqref{Th_bracket}.}

\pb{A differential graded algebra structure on
$\Lambda^\hdot\scr A^*$ with differential
$d: \Lambda^p\scr A^*\to\Lambda^{p+1}\scr A^*$ defined
as follows.}

For $p=0$:\; the differential $\oo_X=\Lambda^0\scr A^*\to\Lambda^1\scr A^*=\scr A^*$
is  the composite $\oo_X\to \Omega^1_X\to\scr A^*$, where
the first map is the de Rham differential and the second map
is obtained by dualizing the morphism $\tau: \scr A\to \calT_X$;

For $p=1$:\; the differential
$d: \scr A^*\to\Lambda^2\scr A^*$ is given by
$$\langle d\alpha, a\wedge b\rangle
:=\tau_a\langle\alpha,b\rangle-\tau_b\langle\alpha,a\rangle
-\langle\alpha,[a,b]\rangle,\quad
\forall \alpha\in\scr A^*,a,b\in\scr A.
$$

For $p>1$:\;  the differential is given by
$$
d(\alpha_1\wedge\ldots\wedge\alpha_p):=
\sum_{i=1}^p (-1)^{i-1}\cdot
\alpha_1\wedge\ldots\wedge d\alpha_i\wedge\ldots\wedge\alpha_p
$$

Given a Lie algebroid $\scr A$,
we define a {\em right}  $\oo_X$-action
on $\scr A$ by the formula
$u\cdot f:= f\cdot u + \tau_u(f)$.
Lie algebroid axioms insure that
the right action so defined commutes with the
natural left $\oo_X$-action on $\scr A$
and provides $\scr A$ with the
structure of a (not necessarily {\em symmetric})
$\oo_X$-bimodule. This $\oo_X$-bimodule
may be thought of as a coherent sheaf on $X\times X$
set-theoretically concentrated on the
diagonal $X_\Delta\sset X\times X$.

\subsection{Gerstenhaber structure on Hochchild cochains.}
Let $G^\hdot$ be a graded associative, not necessarily commutative,
algebra equipped with  a bracket 
$$
\{-,-\}\colon G^p\times G^q\to G^{p+q-1}.
$$
We say that this bracket makes $G^\hdot$ a {\em noncommutative
Gerstenhaber algebra} provided it  gives $G^\hdot$ the structure
of  Lie super-algebra (in particular, the bracket is
skew-symmetric/symmetric
depending on the parity of its arguments) and
the super-Leibniz identity 
\eqref{gerst_axiom} holds (with the order of factors in
the various dot-products being as indicated in \eqref{gerst_axiom}).

Now, let $A$ be an associative not necessarily commutative algebra
and write
$C^\hdot(A,A)=$\break
$\bigoplus_i\Hom_\C(A^{\otimes i},A)
$ for the Hochschild cochain complex of $A$.

There is a natural associative (non-commutative)
graded algebra structure on $C^\hdot(A,A)$ given by the so-called
{\em cup-product}. It is defined, for any
$ f\in C^p(A,A),\,g\in C^q(A,A),$
 by 
the formula
$$f\cup g:\; a_1,\ldots,a_{p+q}\mto f(a_1,\ldots,a_p)\cdot
g(a_{p+1}\ldots,a_{p+q}).
$$
The  Hochschild differential is a super-derivation
with respect to the cup-product, that is, we have
$$d(f\cup g)= (df)\cup g + (-1)^{\deg f}\cdot f\cup (dg).$$
This formula shows that the cup-product of Hochschild
cocycles is again a
cocycle, and the cup-product of a cocycle and a coboundary is a  
coboundary. Thus the  cup-product
descends to a well-defined associative product on
$\HH^\hdot(A,A)$, the  Hochschild cohomology.
It is not difficult to verify that the resulting graded
algebra structure on $\HH^\hdot(A,A)$
gets identified, under the isomorphism
$\HH^\hdot(A,A)\simeq \Ext_{\bimod A}^\hdot(A,A)$, with the standard
{\em Yoneda product} on the Ext-groups.

In addition to the cup-product, there is a much deeper structure 
on Hochschild
cochains, revealed by the following result

\begin{thm}[Gerstenhaber] \vi
There exists a canonical Lie superalgebra structure 
$\{-,-\}:\ C^p(A,A)\times C^q(A,A)\to C^{p+q-1}(A,A)$, called the
Gerstenhaber
bracket.

\vii The cup-product and the Gerstenhaber
bracket make $C^\hdot(A,A)$ a noncommutative Gerstenhaber algebra.
\end{thm}

\begin{proof}[Proof/Construction]
For $f\in C^p(A,A)$ and $g\in C^q(A,A)$, define
\begin{align}
(f\circ g)(a_1,&\ldots,a_{p+q-1})\\
&=\sum_{i=1}^p(-1)^{(i-1)(q-1)}f(a_1,\ldots,a_{i-1},g(a_i,\ldots,a_{i+q}),\ldots,a_{p+q-1}).
\nonumber
\end{align}
Notice that if we regrade the cochains by $(C')^p:=C^{p-1}$, then the
degrees are additive in this product.  This product is not the
cup-product that we have discussed earlier,
it is not associative. However, one has 
 the following key identity in $C^\hdot(A,A)$
due to  Gerstenhaber
\beq{G_cup} f\cup g - (-1)^{(\deg f)(\deg g)}\cdot
g\cup f = d(f\circ g) - df \circ g -(-1)^{\deg f}\cdot f\circ dg.
\end{equation}

Now, following Gerstenhaber, introduce a bracket on Hochschild cohains by the
formula
$$
\{f,g\}:=f\circ g-(-1)^{(p-1)(q-1)}g\circ f,
$$
It is straightforward to verify that
 this gives us a Lie super-algebra structure claimed in part (i) of
the Theorem. 
\end{proof}

The remarkable fact discovered by 
 Gerstenhaber is that noncommutative Gerstenhaber
algebra structure on Hochschild cochains gives rise
to a {\em commutative}  Gerstenhaber
algebra structure on Hochschild cohomology.
That is, one has

\begin{prop} The cup-product and the Gerstenhaber bracket make\linebreak
$\HH^\hdot(A,A)$ a (super-commutative) Gersenhaber algebra.\qed
\end{prop}
\begin{proof} We need to show that both the cup-product and
 Gerstenhaber bracket descend to  Hochschild cohomology.
The case of cup-product follows from an easy identity
$$d( f\cup g)= df \cup g +(-1)^{\deg f}\cdot f\cup  dg.
$$
Observe  that   formula \eqref{G_cup} insures that
the resulting cup-product  on  Hochschild cohomology
is {\em graded commutative}.

Observe further that the super-commutator on the LHS of \eqref{G_cup} is clearly skew super-symmetric
 with respect to
$f\leftrightarrow g$.
Hence, super-symmetrization of the LHS, hence of the RHS of \eqref{G_cup},
must vanish.
This yields
\beq{gerst_ident}
0=d\{f,g\}-\{df,g\}-(-1)^{\deg f}\cdot \{f, dg\}
\end{equation}
The identity  insures that
the  Gerstenhaber bracket descends 
 to a well-defined bracket on Hochschild cohomology.
\end{proof}

\noindent
{\bf More conceptual approach to the Gerstenhaber bracket.}
We will take $A$ to be finite-dimensional so we can take duals,
otherwise we would need to use coalgebras.  Consider $T(A^*)$.  The
comultiplication $\delta\colon A^*\to A^*\otimes A^*$ extends uniquely
to a superderivation on $T(A^*)$.  Now, he space $\Der(T^*A)$ of all
super-derivations on $T(A^*)$ is a Lie super-algebra.  Now, the Leibniz
rule implies that every superderivation $\theta\in\Der(T(A^*))$ is
determined by where it sends each generator.  
Hence it determines a $\C$-linear map $A^*\to T(A^*)$ (simply see where
each basis element of $A^*$ is sent).  Clearly the correspondence is 
reversible, so
\begin{align*}
\Der(T(A^*))=\Hom_\C(A^*,T(A^*))&=(A^*)^*\otimes T(A^*)\\
&=A\otimes T(A^*)=\Hom_\C(TA,A),
\end{align*}
where we consider the \emph{graded} dual of $T(A^*)$.  Notice that the
$n^{\text{th}}$ degree component of $\Hom_\C(TA,A)$ is given by
$\Hom_\C(A^{\otimes n},A)=C^n(A,A)$.  Since $\Der(T(A^*))$ is a Lie
super-algebra, we obtain a super-bracket on $C^n(A,A)$, 
which is the Gerstenhaber bracket.

\begin{rem} Recall that $\HH^0(A,A)=\sZ_A,$ is the center of
the algebra $A$. For this reason, one may think of
the algebra $\HH^\hdot(A,A)$ as a kind of ``derived center'' of
$A$. The corollary above confirms that this ``derived center''
is indeed a commutative algebra.

Recall  that the center of $A$ may be identified further with
$\Hom(\Id_{\Lmod A},$
$ \Id_{\Lmod A})$, the endomorphism algebra of 
the identity functor on the category of left $A$-modules.
In this spirit, one may think of
$\HH^\hdot(A,A)$ to be $\Ext^\hdot(\Id_{\Lmod A}, \Id_{\Lmod A})$,
the
appropriately defined Ext-algebra  of 
the identity functor.
\eer
\subsection{Noncommutative Poisson algebras.}
In this subsection, given an associative algebra $A$, we always
write $[a,b]:= a\cdot b - b\cdot a$, for the commutator
with respect to the associative product.
To avoid confusion, we will use the notation $\{-,-\}$ for the
Lie bracket on a Lie algebra.

\begin{examp}\label{NCP1}
Let $A$ be any
associative algebra, and $t\in\C$ a fixed number.
For any $a\in A$, the map $[a,-]: A\to A$ is a derivation,
the inner derivation corresponding to $a$.
Hence, the map $t\cdot [a,-]$ is also a derivation.
Therefore, setting $\{a,b\}_t:=t\cdot [a,b]$,
one obtains a Poisson algebra 
structure on $A$. Note that, for $t=0$, the corresponding
Poisson bracket vanishes identically.
\eex

Write $\mathsf {Assoc}$, ${\mathsf {Comm}}$
and $\mathsf {Lie}$  for the operads
of
associative, commutative, and Lie algebras, respectively,
cf. e.g. \cite{GiK} for more information about operads.
Also, let $\mathsf{Poiss}$, and $\mathsf {NC\mbox{-}Poiss}$
denote respectively the  operads
of {\em commutative} and not necessarily commutative Poisson algebras.

Taking the zero-bracket on a commutative associative algebra,
one obtains a functor
$\text{\em commutative}$
$\text{\em algebras}
\map \text{\em commutative Poisson algebras}.$
Further, forgetting the associative product
on a  commutative Poisson algebra, 
gives a functor $\text{\em {Forget}}:\ \text{\em 
Poisson algebras}\map
\text{\em Lie algebras}$.
 These functors give rise
to the following  canonical sequence of morphisms of operads
$$ {\mathsf {Lie}} \too 
{\mathsf{Poiss}}\too {\mathsf{Comm}}.$$
The forgetful functor $\text{\em {Forget}}:\ \text{Poisson algebras}\map
\text{\em Lie algebras}$ has a left adjoint
functor ${}^{\top\!}\text{\em {Forget}}:\ \text{\em Lie algebras}\map
\text{\em Poisson algebras}.$ It is given by
associating to a Lie algebra $\g g$ the symmetric algebra
$\Sym \g g$, equipped with Kirillov-Kostant bracket, cf. section \ref{Kirillov-Kostant}.
Thus, the commutative Poisson algebga $\Sym \g g$
may be thought of as the (commutative) {\em Poisson envelope}
of the Lie algebra $\g g$.

Similarly to the above,
taking the zero-bracket on an associative algebra,
gives a functor $\text{\em associative algebras}
\map \text{\em NC\mbox{-}Poisson algebras}.$
Also,  forgetting the associative product
on an NC-Poisson algebra, 
gives a functor $\text{\em {Forget}}: \text{\em NC\mbox{-}Poisson
algebras}$
$\map
\text{\em Lie algebras}$. These functors give rise
to the following  canonical sequence of morphisms of operads
$$ {\mathsf {Lie}} \too {\mathsf {NC\mbox{-}Poiss}}\too {\mathsf {Assoc}}.$$

\begin{question}\label{NCP2}
Is it true that the operad maps above induce an isomorphism
${\mathsf {NC\mbox{-}Poiss}}/({\mathsf {Lie}})\iso {\mathsf {Assoc}},$
where $({\mathsf {Lie}})$ denotes the operad ideal in
${\mathsf {NC\mbox{-}Poiss}}$ generated by (the image of)
$ {\mathsf {Lie}}$.
\end{question}

It is easy to see that the operad ${\mathsf {NC\mbox{-}Poiss}}$ is {\em
quadratic}. Therefore, since ${\mathsf {Assoc}}^!={\mathsf {Assoc}}$
and ${\mathsf {Lie}}^!= {\mathsf {Comm}},$
the sequence above induces the dual
sequence
$${\mathsf {Assoc}}\too {\mathsf {NC\mbox{-}Poiss}}^!\too {\mathsf {Comm}}.$$

We recall that the operad ${\mathsf{Poiss}}$ is known to be Koszul,
see \cite{GiK} and also \cite{MSS}.

\begin{thm}\label{NCP3} The operad NC-Poiss is Koszul.
\end{thm}
\begin{proof}\footnote{This proof was kindly communicated to me by
Martin Markl.} The theorem would immediately  follow from Theorem 4.5 of
\cite{Ml} once we check that the following {\em distributivity} rule

\beq{dl} [ab,c] = a[b,c] + [a,b]c.
\eeq
that ties up the associative product with the Lie bracket is a
distributive law. This can be done in either of the following two ways.

(1) It is easy to verify directly that \eqref{dl} satisfies the condition of
Definition 2.2 of the above mentioned paper, that is, it is indeed a
distributive law.

(2) A less direct way is the following. It follows from general
theory (see, for example, Theorem 3.2 of \cite{FM})
that \eqref{dl} is a distributive law if, for any vector space $V,$ the
free noncommutative Poisson algebra $NCP(V)$ on $V$ is isomorphic to
$T(L(V)),$ the free associative algebra generated by the free Lie
algebra on $V.$

    It is obvious that $NCP(V) = P_{nc}(L(V))$, the NC-Poisson envelope
of the free Lie algebra $L(V)$ which you defined in your paper. It is
immediate to see that the ideal $I$ defined on page 26 of your paper is
trivial if $g = L(V),$ therefore $P_{nc}(L(V)) = T(L(V))$ and the result
follows.

    Therefore the Koszulity of NC-Poiss follows from the same arguments as
the Koszulity of the usual commutative  ${\mathsf{Poiss}}$.
\end{proof}

The forgetful functor $\text{\em {Forget}}:\ \text{\em NC\mbox{-}Poisson algebras}\map
\text{\em Lie algebras}$ has a left adjoint
functor ${}^{\top\!}\text{\em {Forget}}:\ \text{\em Lie algebras}\map
\text{\em NC\mbox{-}Poisson algebras}.$
In other words, given a Lie algebra $\g g$, there
is a uniquely defined NC-Poisson algebra $\ncP(\g g)$,
called the {\em  NC-Poisson envelope} of $\g g$,
that comes  equipped
with a Lie algebra map $\iota: \g g\to \ncP(\g g)$
(with respect to the Poisson bracket on $\ncP(\g g)$)
and 
such that the following universal property holds:

\npb{For any NC-Poisson algebra $P$ and
a  Lie algebra map $\phi: \g g\to P$ there exists a unique
morphism $\ncP(\phi): \ncP(\g g)\to P,$ of  NC-Poisson algebras,
that makes the following diagram commute}
\beq{univ_poiss}\xymatrix{
\g g \ar[r]^<>(0.5){\iota}\ar[dr]_{\phi}&\ncP(\g g) \ar@{.>}[d]^{\ncP(\phi)}\\
& P
}
\end{equation}

The construction of universal universal NC-Poisson envelope $\ncP(\g g)$
is due to Th. Voronov in \cite{Vo}. It is based on the following

\begin{lem}[\cite{Vo}] Let $P$ be a not necessarily commutative
Poisson  algebra with Poisson bracket $\{-,-\}$.
Then, for any $a,b,c,d,u\in P$, one has
$$
\{a,b\}\cdot u\cdot[c,d]=[a,b]\cdot u\cdot\{c,d\},
$$
where $[x,y]=x\cdot y - y\cdot x$ stands for the commutator for the
associative product.
\end{lem}
\begin{proof} Consider the expression $\{ac, bud\}$.
Compute this in two different ways, first by applying the Leibniz
rule for $\{ac,-\}$, and second, by applying the Leibniz rule
for $\{-,bud\}$.
\end{proof}

Now, given a Lie algebra $\g g$ with Lie bracket $\{-,-\}$, Voronov
considers a two-sided ideal $I\sset T(\g g)$, in the tensor algebra of
the vector space $\g g$, generated by the  elements indicated below
$$ I:=\langle\!\langle\{a,b\}\cdot u\cdot[c,d]-[a,b]\cdot u\cdot\{c,d\}
\rangle\!\rangle_{a,b,c,d\in\g g,\,u\in T(\g g) },
$$
where $[-,-]$ stands for the commutator in the associative algebra
$T(\g g)$ and $\{-,-\}$ stands for the Lie bracket in $\g g$.

\begin{thm}[\cite{Vo}] The assignment
\begin{align*}
(a_1\otimes\ldots \otimes a_k)\,&\times\,(b_1\otimes \ldots\otimes  b_l)\mto
\{a_1\otimes \ldots\otimes  a_k\,,\,b_1\otimes \ldots\otimes  b_l\}:=\\
&\sum_{r=1}^k\sum_{s=1}^l b_1\otimes \ldots\otimes  b_{s-1}\otimes
a_1\otimes \ldots\otimes  a_{s-1}\otimes \{a_r,b_s\}\otimes  a_{r+1}\otimes 
\\
&\hphantom{x}\qquad\qquad\otimes\ldots\otimes  a_k\otimes  b_{s+1}\otimes \ldots b_l
\end{align*}
gives rise to
a well-defined noncommutative Poisson structure on the
associative algebra $T(\g g)/I$. The
universal property \eqref{univ_poiss} holds
for the  NC-Poisson algebra $\ncP(\g g):=T(\g g)/I$
thus defined. \qed
\end{thm}

\begin{examp} Let $\g g$ be a Lie algebra
and  $\U\g g$ its enveloping algebra.
We may consider the associative  algebra $\U\g g$
 as
a  NC-Poisson algebra with Poisson bracket $\{-,-\}=[-,-]$.
Therefore, the canonical Lie algebra
map $\g g \into \U\g g$ gives rise, via the universality property
of  $\ncP(\g g)$, to a natural morphism $\ncP(\g g) \to \U\g g$
of  NC-Poisson algebras. The latter morphism is easily seen to be
surjective.
Thus, $\U\g g$ is a quotient of $\ncP(\g g)$.

A similar argument shows that the commutative
 Poisson algebra, $\Sym \g g$,
 is also a quotient of
 $\ncP(\g g)$.
\eex

\section{Deformation quantization}
\subsection{Star products.}
Let $A$ be an associative $\k$-algebra.  We wish to define a 
twisted ``product'' on $A$.  For $a,b\in A$, define
$$
a\circ_tb=ab+t\beta_1(a,b)+t^2\beta_2(a,b)+\cdots\in A[[t]]
$$
for maps $\beta_j\colon A\times A\to A$.  We wish for this map to be ``associative,'' and ask what conditions this places on the $\beta_j$'s.  Define
$$
a\circ_t(b\circ_tc)=(a\circ_t(bc))+[a\circ_t\beta_1(b,c)]t+[a\circ_t\beta_2(b,c)]t^2+\cdots,
$$
and similarly for $(a\circ_tb)\circ_tc$.  Clearly there is no issue for the constant term.  Consider the coefficients of the $t$ term.  From $(a\circ_tb)\circ_tc$ we obtain $\beta_1(ab,c)+\beta_1(a,b)c$, and from $a\circ_t(b\circ_tc)$ we obtain $a\beta_1(b,c)+\beta_1(a,bc)$.  Equating these two we see that $\beta_1\colon A\otimes A\to A$ must be a Hochschild cocycle in $C^2(A,A)$.
In particular, we  begin with a
commutative algebra $A$ (which we can think of as the affine coordinate
ring of a variety $X$), and we will require the deformed products
$\circ_t$ to be associative only. 
It is interesting to examine how noncommutative $\circ_t$ is, so we define
$$
[a_tb]=a\circ_tb-b\circ_ta=t[\beta_1(a,b)-\beta_1(b,a)]+O(t^2).
$$
Define $\{a,b\}=\beta_1(a,b)-\beta_1(b,a)$, which is clearly a skew-symmetric product on $A$.  So, we can write
$$
[a_tb]=t\{a,b\}+O(t^2).
$$

Now, recall that  commutator in any associative algebra
 satisfies the Jacobi identity.  So, for $(A[[t]],\circ_t)$,
 we have the identities
\begin{align*}
[a_tb\circ_tc]&=[a_tb]\circ_tc+b\circ_t[a_tc]\\
[[a_tb]_tc]&=[[b_tc]_ta]+[[c_ta]_tb].
\end{align*}

Now, we insert the formula $[a_tb]=t\{a,b\}+O(t^2)$ into Leibniz's rule and take only the first order terms.  Since every term of the bracket on $A[[t]]$ introduces a factor of $t$, we need only take the first order terms of each term.  For example, we find that
$$
[a_tb\circ_tc]=[a_tbc]+O(t^2)=t\{a,bc\}+O(t^2).
$$
The other brackets simplify similarly, so we obtain
$$
\{a,bc\}=\{a,b\}c+b\{a,c\}.
$$
Similarly, we examine Jacobi's rule and take the coefficient of the lowest power of $t$, namely $t^2$.  Arguing similarly, we find that
$$
\{\{a,b\},c\}=\{\{b,c\},a\}+\{\{c,a\},b\}.
$$
So, $\{-,-\}$ is a Lie bracket on $A$, and combined with the associative product on $A$ it satisfies the Leibniz rule.  This defines the structure of a \emph{Poisson algebra} on $A$.  So, we can reformulate the problem of classifying deformation quantizations as a problem regarding Poisson algebra structures that can be placed on $A$.

\subsection{}\label{abs_bracket}
The calculations above have one large deficiency: if
$\beta_1$ is identically zero then the bracket $\{\cdot,\cdot\}$ will be
as well.  So, we repeat the above reasoning for the general case of
formal deformations.

Thus, let $A$ be a {\em commutative} associative algebra, and
 $\tilde A$ be a formal flat deformation of $A$,
i.e., a topologically free (not necessarily commutative)
 associative $\C[[t]]$-algebra equiped with an algebra
isomorphism $\tilde A/t\tilde A \cong A.$

Since $A$ is commutative, this implies that for each $\tilde a,\tilde
b\in \tilde A$, 
$\tilde a\tilde b-\tilde b\tilde a\in{t}\tilde A$.  For each pair $\tilde a,\tilde b\in \tilde A$, define a natural number $m(\tilde a,\tilde b)$ to be the maximum integer such that
$$
[\tilde a,\tilde b]\in t^{m(\tilde a,\tilde b)}\tilde A.
$$
If $[\tilde a,\tilde b]$ is contained in every ${t}^i\tilde A$, 
we set $m(\tilde a,\tilde b)=\infty$.  Let
$
N=\min\{m(\tilde a,\tilde b)\;|$
$\;\tilde a,\tilde b\in A\},
$
which is necessarily greater than or equal to $1$.  Now, by Krull's theorem, we know that
$
\bigcap_{i=1}^\infty{t}^i\tilde A=0,
$
so if $N=\infty$ this means that every $[\tilde a,\tilde b]\in{t}^iA$ for each $i\ge1$, hence $[\tilde a,\tilde b]=0$.  This would force $\tilde A$ to be commutative.

Since we are interested in noncommutative deformations, we will only
consider the case $N<\infty$.
  For any $a,b\in A$, choose $\tilde a,\tilde b\in \tilde A$ such that 
\begin{equation}\label{lifts}
\tilde a\bmod{t}\tilde A=a\,,\quad
\tilde b\bmod{t}\tilde A=b.
\end{equation}
Then by the definition of $N$, we know that $[\tilde a,\tilde b]\in
t^N\tilde A$, so 
$t^{-N}(\tilde a\tilde b-\tilde b\tilde a)$ is a well-defined element of
$\tilde A$.  We define 
\begin{equation}\label{bra_def}
\{a,b\}=t^{-N}[\tilde a,\tilde b]\bmod{t}\tilde A.
\end{equation}  
It is then easy to
check that $\{a,b\}$ 
is independent of the choice of $\tilde a$ and $\tilde b$, and that
$\{\cdot,\cdot\}$ gives
 $A$ the structure of a Poisson algebra.

\subsection{A Lie algebra associated to a deformation.} Let now
 $A$ be an  associative, {\em not} necessarily commutative,  algebra
and
write $\sZ_A$ for the {\em center} of~$A$.

Hayashi observed, see [Hay], that given
 $\tilde A$,
  a formal flat deformation of $A$,  the construction of the
previous subsection can be adapted to produce
a Poisson structure on  $\sZ_A$, a commutative algebra.
In more detail, let $a,b \in \sZ_A$,
and choose $\tilde a,\tilde b \in\tilde A$ as in \eqref{lifts}.
Then, put
$$\{a,b\}:=\mbox{$\frac{1}{t}$}[\tilde a,\tilde b]\bmod{t}\tilde A
$$
It is straightforward to verify that
\begin{itemize}
\item $\{a,b\}\in \sZ_A$;
\item $\{a,b\}$ is independent of the choices of $\tilde a,\tilde b
\in\tilde A$;
\item The assignment $a,b\mto \{a,b\},$ makes  $\sZ_A$ a Poisson algebra.
\end{itemize}

\begin{rem}
The reader should be alerted that, in contrast with the case of
\S\ref{abs_bracket}, the Poisson bracket on   $\sZ_A$ thus defined
may turn out to be identically zero.
The reason is that unlike formula \eqref{bra_def},
we now divide by the {\em first} power of $t$ rather than
the $N^{th}$ power (one can check that dividing
by the  $N^{th}$  power does not give rise to a well-defined bracket on $\sZ_A$).
Now, given $a,b \in \sZ_A$, put
$
N(a,b):=\min_{\tilde a,\tilde b\in A}m(\tilde a,\tilde b),
$ where the {\em minimum} is taken over all possible
lifts of $a$ and $b$. Then, it is clear from the  construction
that if $
N(a,b)>1$ then the elements $a,b$ Poisson commute.
Therefore,
the  Poisson bracket on   $\sZ_A$ vanishes whenever
one has $
N(a,b)>1\,,\,\forall a,b\in\sZ_A$.
\eer

The Hayashi construction has been  further  refined in \cite{BD} as follows.
We introduce the following vector subspace
$$
{\tilde A}':=\{\tilde a\enspace\big|\enspace
\tilde a\in\tilde A\enspace\text{such that}\enspace 
\tilde a\bmod{t}\tilde A \in \sZ_A\}.
$$

We put  $\calA:= {\tilde A}'/t\cdot{\tilde A}'$. 
Observe next that, for any $\tilde a, \tilde b\in {\tilde A}',$
we have $[\tilde a, \tilde b]\in t\tilde A,$
We claim that the element
$\frac{1}{t}[\tilde a, \tilde b]$ belongs to ${\tilde A}'$.
Indeed, for any $\tilde c\in\tilde A$, using  Jacobi identity in
$\tilde A$, we find
\begin{align*}\big[\tinv[\tilde a, \tilde b], \tilde c\big]=
\tinv\big[[\tilde a, \tilde c]\,,\tilde b\big]
+\tinv\big[\tilde a,\,[\tilde b,\tilde c]\big]
&\in\tinv[t\tilde A,\tilde b]+\tinv[\tilde a,t\tilde A]\\
&=[\tilde A,\tilde b]+[\tilde a,\tilde A]\sset t\tilde A+t\tilde A=t\tilde A.
\end{align*}
Hence, the expression in the LHS of the top line vanishes modulo
$t$. Thus, we have proved that $\frac{1}{t}[\tilde a, \tilde b]\in
{\tilde A}'$. Further, it is easy to see that 
the class $\frac{1}{t}[\tilde a, \tilde b]\bmod{t}{\tilde A}'$
 depends only on the classes $\tilde a\bmod{t}{\tilde A}'$
and $\tilde b\bmod{t}{\tilde A}'$.

This way, one proves the following
\begin{prop}\label{algebroid}\vi The assignment $\tilde a, \tilde b
\mto\frac{1}{t}[\tilde a, \tilde b]\bmod{t}{\tilde A}'$
induces a Lie algebra structure on $\calA$.

\vii The
projection $\calA\to\sZ_A\,,\,\tilde a\mto\tilde a\bmod{t}\tilde A,$
is a Lie algebra map (with respect to the Hayashi bracket on $\sZ_A$),
that gives rise to a Lie algebra extension:
$$0\too A/\sZ_A \too \calA\too \sZ_A\too 0.
$$

\viii The adjoint action of ${\tilde A}'$ on ${\tilde A}$
descends to a well-defined Lie algebra action of $\calA$ on $A$,
i.e., gives a Lie algebra map $\calA\to\Der(A)\,,\,x\mapsto \partial_x$;
moreover, for any $x\in A/\sZ_A\subset \calA$ and $a\in A$,
we have $\partial_x(a)= [{\tilde x}, \tilde a]\bmod{t}\tilde A,$
where ${\tilde x},\tilde a\in \tilde A$ are any lifts of
$x$ and $a$, respectively.\qed
\end{prop}

\subsection{Example: deformations of the algebra $\End(\calE)$.}
 Let $\calE$ be an algebraic vector bundle
on a affine variety $X$, and $A=\End(\calE)$
the endomorphism algebra of this vector bundle.
The center of this algebra is $\sZ_A=\C[X]\sset \End(\calE)$, the subalgebra
of `scalar' endomorphisms.

Now let $\tilde A$ be a formal deformation of the algebra $A$.
By Hayashi construction, this deformation gives rise
to a Poisson bracket on $\sZ_A=\C[X]$, thus makes
$X$ a Poisson variety. In particular, assigning to
each $z\in \C[X]$ the derivation $\{z,-\}$ yields
a Lie algebra map $\xi: \C[X] \to \Der(\C[X])=\calT(X).$

Recall the Atiyah  algebra $\calA(\calE)=\Der\End(\calE)$
of first order differential operators on $\calE$ with 
scalar principal symbol, introduced in \S\ref{atiyah_alg}.
The following result gives an explicit description of the 
Lie algebra  $\calA$ of Proposition \ref{algebroid} in
the special case at hand.

\begin{prop} There is a natural Lie algebra
map  $\xi_{_{\calA}}: \calA\to \calA(\calE)$
making the second row of the diagram below the
pull-back (via $\xi$) of the standard Lie algebra
extension in the first row of the diagram
$$
\xymatrix{
0\ar[r]&\End(\calE)\ar[r]
&
{\;\calA(\calE)\;}\,\ar[r]
&
\calT(X)\ar[r]
&
0\\
0\ar[r]&\End(\calE)\ar[r]\ar@{=}[u]^{\id}
&\calA\ar[u]^{\xi_{_{\calA}}}\ar[r]&\C[X]\ar[r]\ar[u]^{\xi}
&
0}$$
\end{prop}
\begin{proof} By definition, the Lie algebra $\calA(\calE)$Theorem \ref{der_atiyah},
is the Lie
algebra of derivations of the endomorphism algebra,
which is our algebra $A$. Therefore, producing a map
$\xi_{_{\calA}}: \calA\to \calA(\calE)$
amounts to constructing a Lie algebra map
$\calA\to \Der A$. But the latter map has been
already constructed in part (iii) of Proposition
\ref{algebroid}. It is straightforward to verify, using
part (ii) of that  Proposition
that the map $\xi_{_{\calA}}$ arising in this way indices
the map $\xi$, the vertical arrow on the right of the diagram
above.
\end{proof}
\subsection{}
\begin{defn}
Let $B$ be a $\C[[t]]$-algebra with $t$-linear associative, not necessarily commutative product $\circ_t$ and a $t$-linear Lie bracket $[-,-]_t$ such that
\begin{itemize}
\item{$[b,\cdot]_t$ is a derivation with respect to $\circ_t$; and}
\item{$a\circ_tb-b\circ_ta=t[a,b]_t$.}
\end{itemize}
Then we say that $B$ is a \emph{$t$-algebra}.
\end{defn}

\begin{examp}
\vi If $B$ is a flat deformation $\wh A=A[[t]]$ of an algebra $A$ such that $A$ is commutative, then $[a,b]_t=\frac{a\circ_tb-b\circ_ta}t$ makes $\wh A$ into a $t$-algebra.\\

\vii Let $A$ be a Poisson algebra with Poisson bracket $\{-,-\}$.  Take $\circ_t$ to be the given multiplication and $[-,-]_t=\{-,-\}$, and let $t$ act by zero on $A$.  Then $B$ is a $t$-algebra.
\eex

We let $\scr{MP}_A$ denote the moduli space of flat deformations as $t$-algebra of a given Poisson algebra $A$ (recall that we assume that all Poisson algebras are commutative).  The quantization problem can be summarized as an attempt to understand this space.  Recall the moduli space $\scr M_A$ of deformations of $A$ as an associative algebra.  Giving a Poisson structure on $A$ is equivalent to giving an element of the second Hochschild cohomology of $A$, which we identify with the tangent space of $\scr M_A$ at $A$.  That is, an element of $\scr{MP}_A$ corresponds to an element of $T_A\scr M_A$.  So, we can view $\scr{MP}_A$ is the germ of the blow-up of $\scr M_A$ at the point given by $A$ and tangent direction specified by the given Poisson structure on $A$.

\subsection{Moyal-Weyl quantization.}
Let $(V,\omega)$ be a symplectic manifold (recall that this means that
$\omega\in\Omega^2V$ is nondegenerate and closed, $d\omega=0$). 
 Let $A=\C[V]$.  

The
symplectic form $\omega$ gives rise to a Poisson bracket in the following way.  Since $\omega$ is nondegenerate, it induces an isomorphism of bundles $TV\simeq T^*V$.  Given some section $\alpha\in\Gamma(V,T^*V)$, this isomorphism provides a vector field $\xi_\alpha\in\Gamma(V,TV)$.  So, if $f\in\C[V]$ is a regular function, $df$ is a one-form, that is, $df\in\Gamma(V,T^*V)$.  This gives rise to a vector field $\xi_{df}\in\Gamma(V,TV)$.  So, if we have two functions $f,g\in\C[V]$, we define
$$
\{f,g\}=\xi_{df}g.
$$
This is clearly skew-symmetric.

We consider a special case where $V$ is a symplectic vector space
of dimension $\dim V=2n$,
so $\omega\in\wedge^2V^*$ is a nondegenerate skew-symmetric bilinear
form
$\omega: V\wedge V \to\k$.
Recall that given such a symplectic form $\omega$ on $V$,
 we can find coordinates $p_1,\ldots,p_n,q_1,\ldots,q_n$ on $V$ such that in these coordinates $\omega$ is the standard symplectic form, i.e., 
$$
\omega=\sum_{i=1}^ndp_i\wedge dq_i.
$$
Then if $f$ is a regular function on $V$, we find that the vector field associated to $df$ is given by
$$
\xi_{df}=\sum_{i=1}^n\left[\pder{f}{p_i}\pder{}{q_i}-\pder{f}{q_i}\pder{}{p_i}\right],
$$
so that
$$
\{f,g\}=\sum_{i=1}^n\left[\pder{f}{p_i}\pder{g}{q_i}-\pder{f}{q_i}\pder{g}{p_i}\right].
$$

 Let  $\pi\in\wedge^2V$ be the 
bivector obtained by transporting the 2-form $\om$ via the
isomorphism $\wedge^2 V^*\iso\wedge^2V$ induced by the
symplectic form $\om$. Using $\pi$ we can rewrite the above
 Poisson bracket on $\C[V]$ as $f,g\mapsto \{f,g\}:=\langle df\wedge dg,\pi\rangle$
on $\k[V],$ the polynomial algebra on $V$.
The usual commutative product $m: \k[V]\o \k[V]\to\k[V]$
and the Poisson bracket $\{-,-\}$ make $\k[V]$ a  Poisson algebra.

This Poisson algebra has a distinguished {\em Moyal-Weyl quantization} 
(\cite{M}, see also \cite{CP}).
This is
an associative   star-product 
depending on a  formal quantization parameter $ t$, defined by the formula
\beq{star}
f *_ t g:=m\ccirc e^{\frac{1}{2}  t \pi} (f\o g)\in \k[V][ t],\quad
\forall f,g\in\k[V][ t].
\eeq

To explain the meaning of this formula, view elements of $\Sym V$
as
constant-coefficient differential operators on $V$. Hence,
an element of $\Sym V\o \Sym V$ acts  as a
constant-coefficient differential operator on 
the algebra $\k[V]\o\k[V]=\k[V\times V].$ Now, identify
$\wedge^2V$ with the subspace of skew-symmetric tensors in $V\o V$.
This way, the bivector $\pi\in\wedge^2V\subset V\o V$ becomes
a second order constant-coefficient differential operator
$\pi: \k[V]\o\k[V]\to\k[V]\o\k[V].$ Further, it is clear that
for any element $f\o g\in\k[V]\o\k[V]$ of total degree
$\leq N$, all terms with $d>N$ in the
infinite sum $e^{ t\cdot \pi}(f\o g)=\sum_{d=0}^\infty
\frac{ t^d}{d!} \pi^d(f\o g)$
vanish, so the sum makes sense.
Thus,
the 
symbol $m\ccirc e^{ t\cdot \pi}$ in the right-hand side
of formula \eqref{star} stands for the composition
$$ \k[V]\o\k[V]\stackrel{e^{ t\cdot \pi}}\tooo \k[V]\o\k[V]\o\k[ t]
\stackrel{m\o\Id_{\k[ t]}}\tooo  \k[V]\o\k[ t],
$$
where
$e^{ t\cdot \pi}$
is an infinite-order formal differential operator.

In down-to-earth terms, choose coordinates
$x_1, \ldots, x_n, y_1, \ldots, y_n$ on $V$ such that
the bivector $\pi$, resp., the Poisson bracket $\{-,-\}$, takes 
the canonical form
\beq{pois}
\pi = \sum_i \frac{\partial}{\partial x_i} \o \frac{\partial}{\partial
y_i} - \frac{\partial}{\partial y_i} \o \frac{\partial}{\partial x_i},
\quad\text{resp.,}\quad
\{f,g\}=\sum_i \frac{\partial f}{\partial x_i}\frac{\partial g}{\partial
y_i} - \frac{\partial  f}{\partial y_i} \frac{\partial g}{\partial x_i}.
\eeq
Thus, in canonical coordinates
$x=(x_1, \ldots, x_n), y=(y_1, \ldots, y_n),$ formula \eqref{star} for the Moyal
product
reads
\begin{align}\label{star1}
(f *_ t g)(x,y)&=\sum_{d=0}^\infty
\frac{ t^d}{d!}\left(
\sum_i \frac{\partial}{\partial
x'_i} \frac{\partial}{\partial
y''_i} - \frac{\partial}{\partial y'_i} 
\frac{\partial}{\partial x''_i}
\right)^df(x',y') g(x'',y'')\Big|_{{x'=x''=x}\atop{y'=y''=y}}\nonumber\\
&=\sum_{\mathbf{j},\mathbf{l}\in\Z^n_{\geq 0}}
(-1)^{\mathbf{l}|}\frac{ t^{|\mathbf{j}|+|\mathbf{l}|}}{\mathbf{j}!\,\mathbf{l}!}\cdot
\frac{\partial^{\mathbf{j}+\mathbf{l}}f(x,y)}{\partial x^\mathbf{j}\partial y^\mathbf{l}}
\cdot
\frac{\partial^{\mathbf{j}+\mathbf{l}}g(x,y)}{\partial y^\mathbf{j}\partial x^\mathbf{l}},
\end{align}
where for $\mathbf{j}=(j_1, \ldots,j_n)\in \Z^n_{\geq 0}$
we put $|\mathbf{j}|=\sum_i j_i$ and given
$\mathbf{j},\mathbf{l}\in \Z^n_{\geq 0},$ write
$$\frac{1}{\mathbf{j}!\,\mathbf{l}!}\frac{\partial^{\mathbf{j}+\mathbf{l}}}
{\partial x^{\mathbf{j}}\partial y^{\mathbf{l}}}:=
\frac{1}{j_1!\ldots
j_n!l_1!\ldots l_n!}\cdot\frac{\partial^{|\mathbf{j}|+|\mathbf{l}|}}
{\partial x_1^{j_1}\ldots\partial x_n^{j_n}\partial
y_1^{l_1}\ldots\partial y_n^{l_n}}.
$$

\subsection{Weyl Algebra}
A more conceptual approach to formulas \eqref{star}--\eqref{star1}
is obtained by introducing the {\em Weyl algebra} $A_ t(V)$.
This is a $\k[ t]$-algebra defined by the quotient
$$
A_ t(V):= (TV^*)[ t]/I(u\o u' - u'\o u- t\cdot\langle \pi,
u\o u'\rangle)_{u,u'\in V^*},
$$
where $TV^*$ denotes the tensor algebra of the vector space $V^*$,
and $I(\ldots)$ denotes the two-sided ideal generated by the
indicated set.
 
For instance, if $\dim V=2$ and $p,q$ are canonical coordinates on $V$,
 then
we have
$$
A_t=A_ t(V)=\C\sp<p,q>/(pq-qp=t).
$$

By scaling, there are essentially only two different cases for $t$,
 namely $t=0$ or $t\ne0$.  However, it will be convenient to have a
 continuous parameter.  Also, the algebra $A_1$ will be called simply
 the Weyl algebra.  Notice that even though a monomial does not have a
 well-defined degree (for example, $pq=qp-t$, the left-hand side has
 degree two and the right is not even homogeneous if $t\ne0$), we can
 see that any homogeneous monomial has a highest degree in which it may
 be expressed. 
 We filter $A_t(V)$ by this highest degree for each monomial.

Now, a version of the Poincar\'e-Birkhoff-Witt theorem
says that the natural {\em symmetrization map}
yields a  $\k[ t]$-linear bijection
$\sigma_W: \k[V][ t]\iso A_ t(V)$.
In the special case of a 2-dimensional space $V$,
the linear bijection $\sigma_t\colon\C[p,q]\to A_t$ is given by sending $p^mq^n\in\C[p,q]$ to the average of all possible permutations of the $m$ $p$'s and $n$ $q$'s.  This becomes the identity when we pass to the graded case.  That is, if $\f$ is a homogeneous polynomial, then the principal symbol of $\sigma_t(\f)$ in $\gr A_t$ is precisely $\f$.

Thus, transporting the multiplication map in the Weyl algebra $A_ t(V)$ via
the bijection $\sigma_W,$ one obtains an associative product
$$\k[V][ t]\o_{\k[ t]}\k[V][ t]\to\k[V][ t],
\quad f\o g\mapsto \sigma_W^{-1}(\sigma_W(f)\cdot \sigma_W(g)).
$$
It is known that this  associative product is equal to the one
given by formulas  \eqref{star}--\eqref{star1}.

The easiest way to see the last assertion is to argue heuristically as
follows. We assume $\dim V=2,$ for simplicity.

First of all, one verifies that $f\o g=\sigma_W^{-1}(\sigma_W(f)\cdot
\sigma_W(g))$
admits an expansion with the initial term $f g $ plus terms of higher
order in $t$ whose coefficients are all composed of differential
operators with polynomial coefficients applied to $f$ and $ g $.  
Given this claim, we can then extend the product formula to all smooth
functions,
not just polynomials.  In particular, we choose to take
$$
\f(p,q)=e^{\alpha p+\beta q}\quad\text{and}\quad g (p,q)=e^{\gamma p+\delta q},
$$
where none of $\alpha,\beta,\gamma,\delta$ equal one another and none of them equal zero or one.  Then a differential operator on $\f$ and $ g $ is determined completely by its action on the above $\f$ and $ g $.  

So, we have essentially moved from the problem of computing the product
in $A_t$ to computing its logarithm.  In the case of Lie algebras, 
we can invoke the Campbell-Hausdorff theorem to obtain some partial
information.  In this case, we know that both $e^x$ and $e^y$ are
elements of an associated Lie group, hence so is $e^x\circ_te^y$.  So,
we should be able to express $e^x\circ_te^y=e^{z_t(x,y)}$ for some
function $z_t(x,y)$.  
In particular, we always know the first two terms:
$$
z_t(x,y)=(x+y)+\frac t2[x,y]+O(t^2).
$$
Notice that the first term is $x+y$.  Here is a little exercise: Check that the $\beta_j$'s (the coefficient $t^j$) can be expressed as a differential operator as we claim if and only if the above expansion starts with $x+y$.

We now specialize the general  Campbell-Hausdorff theorem in the case of the
3-dimensional
Heisenberg Lie algebra $\g h$ with basis $\{x,y,c\}$ and 
 commutation
relations
$$ [x,y]=c,\quad [x,c]=[y,c]=0,
$$
in particular,  $c$ is central.

It is clear that the Weyl algebra $A_t$ is a quotient of the enveloping
 algebra
of $\g h$, specifically, we have 
 $A_t=\scr U(\g h)/(c=t)$.  So, we can apply the discussion regarding
 the Campbell-Hausdorff theorem to $\g h$.  Notice that since every
 bracket is central, there can be no nontrivial interated brackets. 
 So, we can check that the Campbell-Hausdorff theorem yields the simple relation
$$
z_t(x,y)=(x+y)+\frac t2[x,y].
$$
We find that
$$
e^x\circ_te^y=e^{x+y+\frac t2[x,y]}.
$$

We can generalize the above to the case $\dim V>2$.  
This will allow us to better understand some of the symmetries of the
situation.  
The Heisenberg Lie algebra $\g h$ for dimension $V$ is given by a central extension of $\C$ by $V$ treated as an abelian Lie group (i.e., $[v,w]=0$ for all $v,w\in V$).  That is,
$$
\xymatrix{0\ar[r]&\C\ar[r]&\g h\ar[r]&V\ar[r]&0}.
$$
The above exact sequence splits as a vector space, so we may write $\g h=\C c+V$, where $c$ is a non-zero element of $\C$.  To determine the Lie bracket on $\g h$, we need only compute $[x,y]$ for $x,y\in V$ (since $c$ is central).  We set $[x,y]=\omega(x,y)c$ for all $x,y\in V$.  Notice that this is invariant under the action of the symplectic group $Sp(V,\omega)$ on $V$.  In our situation, we obtain the formula
$$
e^x\circ_te^y=e^{x+y+\frac t2\omega(x,y)}=e^{\frac t2\omega(x,y)}e^{x+y},
$$
which is a simple scalar correction term.  If we now apply this calculation for $\f(p,q)=e^{\alpha p+\beta q}$ and $\psi(p,q)=e^{\gamma p+\delta q}$, we find
\begin{align*}
\f\circ_t\psi(p,q)&=\exp\left[\frac t2\omega(\alpha p+\beta q,\gamma p+\delta q)\right]\f\psi\\
&=\exp\left\{\frac t2[\alpha\gamma\omega(p,p)+\alpha\delta\omega(p,q)+\beta\gamma\omega(q,p)+\beta\delta\omega(q,q)]\right\}\f\psi\\
&=\exp\left[\frac t2(\alpha\delta-\beta\gamma)\right]\f\psi\\
&=\sum_{n=0}^\infty\frac{t^n}{2^nn!}(\alpha\delta-\beta\gamma)^n\f(p,q)\psi(p,q).
\end{align*}
Of course, $\alpha$ corresponds to differentiating $\f$ with respect to its first argument, $\beta$ is differentiation of $\f$ with respect to the second, etc.  So, for any $f,g\in\C[p,q]$, we obtain
$$
f\circ_tg=\left.\sum_{n=0}^\infty\frac{t^n}{2^nn!}\left(\pder{}{p'}\pder{}{q''}-\pder{}{p''}\pder{}{q'}\right)\f(p',q')\psi(p'',q'')\right|_{\substack{p'=p''=p\\q'=q''=q}}.
$$
The right hand side here is exactly the same expression 
as given by formula \eqref{star1}.

\section{K\"ahler differentials}\label{kahler_diff}
\subsection{}
In order to make some connections to geometry, we are going to discuss a construction for commutative algebras.  So, until further notice, $A$ is a commutative algebra.  K\"ahler differentials for $A$ are $A$-linear combinations of the symbols $db$, where $b\in A$ and $d(ab)=a\,db+b\,da$.  More formally:

\begin{defn}
Set
$\dis
\comO^1(A):=A\o A/(ab\otimes c-a\otimes bc+ac\otimes b).
$
The elements of $\comO^1(A)$ are called the \emph{(commutative) K\"ahler
differentials} of $A$ (intuitively, $a\otimes b$ 
should be thought of as a differential form  $a\,db$).
\end{defn}

Let us examine the connection of $\comO^1(A)$ to Hochschild homology.  Recall that Hochschild homology is the homology of the complex $A^{\otimes3}\to A^{\otimes2}\to A\to0$.  An element $a\otimes b$ of $A\otimes A$ is automatically a $1$-cycle, since $d(a\otimes b)=ab-ba=0$ as $A$ is commutative.  Since 
$$
d(a\otimes b\otimes c)=ab\otimes c-a\otimes bc+ca\otimes b=ab\otimes c-a\otimes bc+ac\otimes b,
$$
we see that the relation for $\comO^1(A)$ is derived precisely from the consideration of $1$-boundaries.  So,
$$
\comO^1(A)\simeq \HH_1(A).
$$

If we think of $\comO^1(A)$ as the free $A$-module generated by symbols $db$ where $d$ acts as a derivation from $A$ to $A$, we find that $d1_A=0$ through the usual argument that a derivation is zero on constants.  Indeed, we find that $d(\lambda1_A)=0$ for all $\lambda\in\C$.

The space of K\"ahler differentials for $A$ plays an important universal role relative to derivations.  Define $\partial\colon A\to\comO^1(A)$ by $\partial a=a\otimes1_A-1_A\otimes a$.  Then $\partial$ is indeed a derivation, and we will often denote it symbolically by $a\mapsto da$.

\begin{thm}\label{T:KaehlerUni}
Let $M$ be an $A$-module and let $\theta\colon A\to M$ be a derivation.
Then the assignment  $\comO^1(\theta): (a\,db)\mapsto a\theta(b)$ gives
an $A$-module map $\comO^1(\theta)\colon\comO^1(A)$
$\to M$,
which is uniquely defined by the requirement  
that the following diagram commutes:
$$
\xymatrix{A\ar[rr]^\partial\ar[dr]_\theta&&
\comO^1(A)\ar[dl]^{\Omega(\theta)}\\
&M}.
$$
\end{thm}

The proof is a routine calculation.

\begin{examp}
Consider the case $A=\Sym V$, the symmetric algebra for the $\C$-vector
space $A$.  
Since a derivation on an algebra is uniquely defined by specifying its
values on generators and applying the Leibniz rule, we find that
$\Der(\Sym V,M)\simeq\Hom_\C(V,M)$.  
Then the theorem asserts that there is an isomorphism between $\Der(\Sym V,M)$ and $\Hom_A(\comO^1(A),M)$, so $\comO^1(A)$ is the unique $\Sym V$-module such that
$$
\Hom_\C(V,M)\simeq\Hom_A(\comO^1(A),M).
$$
Clearly,
 $\comO^1(A)$ is the free $\Sym V$-module with 
$V$ being the space of generators that is, $\comO^1(A)\simeq (\Sym V)\o V$.  
With this definition, we can explicitly calculate that
 $\partial\colon \Sym V=A\to\comO^1(A)=(\Sym V)\o V$ is given by
$$
\partial(v_1\cdots v_n)=\sum_{i=1}^nv_1\cdots \wh v_i\cdots v_n\otimes v_i.
$$
\eex

Suppose now
$A=\C[X]$, the algebra of regular functions on some affine variety $X$.
We  would like to identify $\comO^1(A)$
with   $\calT^*(X)$,
the space of global sections of the cotangent bundle on $X$.
 This can be done in the following way.  

Consider the diagonal embedding $X\subset X\times X$.  Then we can view
$\calT_X$ as the normal bundle to $X$ in $X\times X$, that is,
$TX=T_X(X\times X)$.  Similarly, we view $\calT^*(X)$ as the conormal
bundle $T^*_X(X\times X)$.  
We would like to be able to identify $\comO^1(A)$ with
$\Gamma(X,\,\calT_X^*(X\times X))$, 
the space of regular sections of the corresponding conormal sheaf.

In general, if we have an embedding $X\inj Y$ of affine varieties, then
we obtain an embedding $I(X)\subset\C[Y]$, where $I(X)$ denotes the
ideal of regular functions vanishing on $X$.  We can view $I(X)/I(X)^2$
as a linear form on $T_XY$ which is zero in the 
$X$-direction.  Intuitively, if $f\in I(X)$, then the Taylor expansion
of $f$ around a point $x\in X$ begins with the first derivative.  By
quotienting out $I(X)^2$, we are killing the higher derivative terms,
and only considering the derivative of $f$ applied 
to a tangent vector to $Y$.  So, the definition $\Gamma(X,\,\calT^*_XY)
:=I(X)/I(X)^2$ seems to be a legitimate one.

In our case, $Y=X\times X$ and $X\inj X\times X$ is the diagonal map, 
so $I(X)=\Ker(A\otimes A\to A)$, where the map is multiplication.  
Notice that multiplication is an algebra map if and only if $A$ is 
commutative.  We can now formulate the following proposition.

\begin{prop}
There is a canonical isomorphism of $A$-modules,
$$
\comO^1(A)\simeq I(X)/I(X)^2.
$$
\end{prop}

\begin{proof}
First, we will set $I=I(X)=\Ker(A\otimes A\to A)$ to simplify the
notation.

Given  $A$ a commutative algebra, the multiplication map $m\colon A\otimes A\to A$ is an algebra homomorphism, and we let $I$ denote its kernel, which is therefore an ideal of $A\otimes A$.  Let $M$ be a left $A$-module.  We will view $M$ as an $A$-bimodule by equipping it with the symmetric $A$-bimodule structure, that is, $m\cdot a=am$ for all $m\in M$, $a\in A$.  Then the short exact sequence of $A$-bimodules
$$
\xymatrix{0\ar[r]&I\ar[r]&A\otimes A\ar[r]^m&A\ar[r]&0}
$$
induces a long exact sequence of $\Ext$-groups, the first few 
terms of which we produce below
\begin{align*}
\Hom_{\bimod A}(A,M)&\to\Hom_{\bimod A}(A\otimes A,M)\to
\Hom_{\bimod A}(I,M)\\
&\to\Ext^1_{\bimod A}(A,M)\to\Ext^1_{\bimod A}(A\otimes A,M)\to\cdots.
\end{align*}
Now, since an $A$-bimodule is 
 a left $A\otimes A^\op$-module, and $A$ is commutative, we see that $A\otimes A$ is the free $A$-bimodule of rank one.  Hence $\Ext^1_{\bimod A}(A\otimes A,M)=0$ and $\Hom_{\bimod A}(A\otimes A,M)\simeq M$.  Also, we have by definition that $\Ext^1_{\bimod A}(A,M)\simeq \HH^1(M)$, the first Hochschild cohomology of $A$ with coefficients in $M$, which is precisely the set of all outer derivations from $A$ to $M$.  However, $M$ has a symmetric $A$-bimodule structure, so there are no inner derivations.  Hence, $\Ext^1_{\bimod A}(A,M)\simeq\Der(A,M)$.

Finally, we observe that $\Hom_{\bimod A}(A,M)=\{m\in M\st am=ma\}$,
since this is just the zeroth-degree Hochschild cohomology.  But $M$ is
symmetric, so this is all of $M$.  Therefore, the map $\Hom_{\bimod
A}(A,M)\to\Hom_{\bimod A}(A\otimes A,M)$ is an isomorphism,
 and $\Hom_{\bimod A}(I,M)\to\Der(A,M)$ is a surjection.  A calculation then shows that in fact $\Hom_{\bimod A}(I,M)\simeq\Hom_A(I/I^2,M)$, so we see that $\comO^1(A)\simeq I/I^2$, which is essentially the definition we gave before.

Then given an $A$-bimodule (i.e., an $A\otimes A$-module), the
$A$-bimodule $M/IM$ 
has a symmetric bimodule structure, that is, $am=ma$ for all $a\in A$ and $m\in M$.  Second, we recall that the bar complex is exact, and that it was (initially) given by
$$
\xymatrix{A^{\otimes3}\ar[r]^b&A\otimes A\ar[r]^m&A\ar[r]&0}.
$$
So, $I=\im b$.  So, the following sequence is exact
$$
\xymatrix{0\ar[r]&I\ar[r]&A\otimes A\ar[r]&A\ar[r]&0}.
$$
Now, both $\Tor$ and $\Ext$ are homological functors, so they associate
long exact sequences to short exact sequences such as the one above.
Since $A\otimes A$ is a free $A$-bimodule, 
$\Ext_{\bimod A}(A\otimes A,-)=0$.  By standard results in homological 
algebra, we obtain that for any $A$-bimodule $M$,
$$
\Ext^1_{\bimod A}(A,M)\simeq\Ext^0_{\bimod A}(I,M).
$$
But we know that $\Ext^1_{\bimod A}(A,M)\simeq \HH^1(M)\simeq\Der(A,M)$,
and also that 
$\Ext_{\bimod A}^0(I,M)\simeq \HH^0(I,M)=\Hom_{\bimod A}(I,M)$.  Now recall that $A$ is commutative and $M$ is a symmetric bimodule.  Then
$$
\Hom_{\bimod A}(I,M)=\Hom_{\Lmod A}(I/I^2,M).
$$
Following the line of isomorphisms, we find that indeed $I/I^2\simeq\comO^1(A)$.
\end{proof}

This proof relied heavily on the symmetry of the module actions.  Indeed, this is the first point where we will see noncommutative and commutative geometry diverging.

\begin{rem}
Identify $\comO(A)$ with $I/I^2$ where $I=\Ker[m: A\otimes A\to A]$.
It is easy to check that the map $\comO(\theta)\colon I/I^2\to M$
corresponding to a derivation $\theta\colon A\to M$
is induced by a map $A\otimes A\supset I \to M$ given by a similar 
formula
$$
\comO(\theta):\;\
\sum\nolimits_{i=1}^n\,a_i\otimes b_i
\;\mto \;\sum\nolimits_{i=1}^n\,a_i\theta(b_i).
$$
The Leibniz formula for $\theta$ insures that this map indeed vanishes  on $I^2$.
\end{rem}

\section{The Hochschild-Kostant-Rosenberg Theorem}
\subsection{Smoothness.} Let $X\sset \C^n$ be an algebraic set defined
by a system of polynomial equations
\beq{alg_set}
X=\{x=(x_1,\ldots,x_n)\in\C^n\enspace\big|\enspace f_1(x)=0,\ldots, f_r(x)=0\}.
\end{equation}
Thus, $X$ is the zero variety of the ideal
$I:=\llb f_1,\ldots, f_r\rrb\sset \C[x_1,\ldots,x_n]$.
We call $\C[X]:=\C[x_1,\ldots,x_n]/I$ the (scheme-theoretic) coordinate
ring of $X$. Hilbert's Nullstellensatz says that the ideal
$I$ is {\em radical}, i.e., the algebra $\C[X]$ has no nilpotents,
 if and only if $I$ is equal
to the ideal of all polynomials $f\in \C[x_1,\ldots,x_n]$ that
vanish on the set $X$ pointwise. In this case, the coordinate
ring $\C[X]$ is said to be {\em reduced}, and the algebraic set
$X$ is called an affine algebraic variety.

We would like to discuss the notion of {\em smoothness}
of algebraic varieties. To this end, fix an  algebraic set
 $X$ as in \eqref{alg_set}. For any point
$a\in X$, we introduce the following Jacobian $n\times r$-matrix:
$$J_a(f_1,\ldots,f_r):=\left.\begin{pmatrix}\frac{\partial f_1(x)}{\partial x_1}&\ldots&
\frac{\partial f_1(x)}{\partial x_n}\\
\frac{\partial f_2(x)}{\partial x_1}&\ldots&
\frac{\partial f_2(x)}{\partial x_n}\\
\dots&\dots&\dots\\
\frac{\partial f_r(x)}{\partial x_1}&\ldots&
\frac{\partial f_r(x)}{\partial x_n}
\end{pmatrix}\right|_{x=a}
$$
Further, let $\oo_a$ denote the local ring at $a$ (the localization
of $\C[X]$ with respect to the multiplicative set of all
polynomials $f$ such that $f(a)\neq 0$), and write
$\g m_a\sset \oo_a$ for the maximal ideal of the local ring $\oo_a$.
Thus $\oo_a/\g m_a=\C$, and $\g m_a/\g m_a^{\,2}$ is a finite
dimensional vector space over  $\oo_a/\g m_a=\C$.
We consider the graded algebra
$\Sym^\hdot(\g m_x/\g m_x^{\,2}):=\bigoplus_{i\geq 0}
\Sym^i(\g m_x/\g m_x^{\,2})$, and also 
 the graded algebra $\bigoplus_{i\geq 0}\,\g m^i/\g m_x^{\,i+1}$.
The imbedding $\g m_x/\g m_x^{\,2}\into
\bigoplus_{i\geq 0}\,\g m^i/\g m_x^{\,i+1}$
extends by multiplicativity to a graded algebra
homomorphism
$\Sym^\hdot(\g m_x/\g m_x^{\,2})\map \bigoplus_{i\geq 0}\,\g m^i/\g
m_x^{\,i+1}$.

One has the following basic result.

\begin{thm}\label{smooth_thm} For an irreducible
 algebraic set $X$ the following
properties \vi--\iv are equivalent:

\vi
For any $x\in X,$ the  map
$\Sym^\hdot(\g m_x/\g m_x^{\,2})\map \bigoplus_{i\geq 0}\,\g m^i/\g m_x^{\,i+1}$
is an isomorphism;

\vii The module $\comO^1(\C[X])$ of K\"ahler differentials
is a projective 
$\C[X]$-module;

\viii The coordinate ring  $\C[X]$, viewed as a module over
$\C[X]\ee=\C[X]\otimes\C[X]$, has a finite projective resolution;

\iv For any point $a\in X$ one has $\rk J_a(f_1,\ldots,f_r) =n-\dim X$.
\end{thm}

The algebraic set $X$ satisfying the equivalent conditions
(i)--(iv) of the Theorem is called {\em smooth}.
It is easy to see that condition (i) above implies that
  $\C[X]$ is reduced, i.e., a smooth algebraic set 
is necessarily an affine algebraic variety.

\begin{thm}[HKR]\label{T:HKR} 
Let $A=\C[X]$, where $X$ is a smooth affine variety.  Then
\begin{align*}
\HH_k(A)&=\Gamma(X,\Lambda^k\calT^*(X))=\Lambda_A^k\comO^1(A)\\
\HH^k(A)&=\Gamma(X,\Lambda^k\calT(X))=\Lambda_A^k\Der(A),
\end{align*}
where $\calT(X)$ is the tangent bundle of $X$, $\calT^*(X)$ is the cotangent bundle, and $\Gamma(X,\cdot)$ denotes global sections.
\end{thm}

\begin{rem}
First, observe that the HKR theorem shows that both Hochschild homology and cohomology are commutative algebras (i.e., the total homology $\HH_\idot(A)=\oplus_{n=0}^\infty \HH_n(A)$ is an algebra, as is the cohomology).  In the case of cohomology, this is not too terribly surprising.  Recall that we defined
$$
\HH^\hdot(A)=\Ext_{\bimod A}^\hdot(A,A),
$$
which has a commutative algebra structure induced by the diagonal map $A\mapsto A\o A$.  We can see in a more elementary way that $\HH_\idot(A)=\Tor_\idot^{\bimod A}(A,A)$ has a commutative algebra structure by observing that the multiplication map $A\o A\to A$ is an algebra map if and only if $A$ is commutative.  In particular, we would not expect $\HH_\idot(A)$ to be an algebra if $A$ is not commutative.
\eer

\subsection{From  Hochschild complex to  Chevalley-Eilenberg complex.}
For any {\em commutative} algebra $A$ and
a left $A$-module $M$, viewed as a {\em symmetric} $A$-bimodule, we are going to construct
the following natural maps:
\begin{align}\label{hoch_wedge}
\xymatrix{
{C_\idot(A,M)\;}\ar@<1ex>[rrr]^{\varkappa^M_\idot}&&&
{\;M\otimes_A\,
\Lambda_A^\hdot\Der(A),}\ar@<1ex>[lll]^{\alt^M_\idot}}\\
\xymatrix{
{C^\hdot(A,M)\;}\ar@<1ex>[rrr]^{\varkappa^\hdot_M}&&&
{\;M\otimes_A\,
\Lambda_A^\hdot\comO^1(A).}\ar@<1ex>[lll]^{\alt^\hdot_M}}\nonumber
\end{align}
  The proof involves Lie algebra homology, which we recall here.  If $\g g$ is a Lie algebra and $M$ is a $\g g$-module, then the Lie algebra homology $H^\Lie_p(\g g,M)$ is computed in terms of the following complex: the $p^{\text{th}}$ term $C^\Lie_p(\g g,M)$ is given by $M\o \Lambda^p\g g$, and the differential $d\colon C^\Lie_p(\g g,M)\to C^\Lie_{p-1}(\g, M)$ is given by
\begin{align*}
d_\Lie(m\otimes&(x_1\wedge\cdots\wedge x_p))=\sum_{i=1}^p(-1)^ix_im\otimes(x_1\wedge\cdots\wedge\wh x_i\wedge\cdots\wedge x_p)\\
&+\sum_{j<k}(-1)^{j+k}m\otimes([x_j,x_k]\wedge x_1\wedge\cdots\wedge\wh x_j\wedge\cdots\wedge\wh x_k\wedge\cdots\wedge x_p).
\end{align*}
Now suppose that $A$ is an associative algebra and $M$ is an $A$-bimodule.  Consider $A$ as a Lie algebra under the commutator bracket, and make $M$ a Lie $A$-module via the action $(a,m)\mapsto am-ma$.  Notice that if $A$ is commutative and $M$ is symmetric, then both the Lie algebra structure on $A$ and the Lie $A$-module structure of $M$ are trivial.

Consider the following diagram
$$
\xymatrix{M\o \Lambda_\C^pA\ar[rr]^{\id_M\otimes\alt}\ar[d]_{d_\Lie}&&M\o A^{\otimes p}\ar[d]^d\\
M\o \Lambda_\C^{p-1}A\ar[rr]_{\id_M\otimes\alt}&&M\o A^{\otimes(p-1)}},
$$
where $\alt\colon\Lambda_\C^pA\to A^{\otimes p}$ is the completely alternating map given by
$$
\alt(a_1\wedge\cdots\wedge a_p)=\sum_{\sigma\in S_p}(-1)^{\sgn\sigma}a_{\sigma(1)}\otimes\cdots\otimes a_{\sigma(p)}.
$$
In the case that $A$ is commutative and $M$ symmetric, then the image of the alternating map on $C_p^\Lie(A,M)$ is zero.  In any case, it is tedious but easy to check that $\alt(C_p^\Lie(A,M))\subset Z_p(A,M)$.  By passing to the quotient by $B_p(A,M)$, we obtain a map $\alt\colon M\o \Lambda_\C^pA\to \HH_p(A,M)$.

We will now construct the inverse map and leave it to the 
reader to check the remainder of the proof.  Choose any $m\otimes(a_1\otimes\cdots\otimes a_p)\in C_p(A,M)$.  Define $\pi\colon C_p(A,M)\to M\otimes_A\Lambda_A^p\comO^1(A)$ by
$$
\pi(m\otimes(a_1\otimes\cdots\otimes a_p))=m\otimes da_1\wedge\cdots\wedge da_p.
$$
It is then easy to see that $\pi\circ\alt$ is given by
$$
m\otimes(a_1\wedge\cdots\wedge a_p)\mapsto p!m\otimes(da_1\wedge\cdots\wedge da_p).
$$
We claim that in the case that $M=A$, where $A$ is the coordinate ring of a smooth affine variety, $\pi$ is an isomorphism from $\HH_p(A)$ to $\Lambda_A^p\comO^1(A)$ as desired.  For simplicity's sake, we will let $\comO^p(A)=\Lambda_A^p\comO^1(A)$.

Indeed, we will show that $\pi$ is an isomorphism of $A$-modules.  It is clear that $\comO^p(A)$ is an $A$-module by its definition.  Observe that for any associative algebra, $\HH_p(M)$ has the structure of a $\sZ_A$-module for all $p$ (recall that $\sZ_A$ is the center of $A$).  Indeed, if $z\in \sZ_A$, the action is given by
$$
z\cdot[m\otimes(a_1\otimes\cdots\otimes a_p)]=zm\otimes(a_1\otimes\cdots\otimes a_p).
$$
It is easy to see that the $z$-action commutes with $d$ since $z$ commutes with all elements of $A$.  Since $A$ is commutative, $\sZ_A=A$, hence $\HH_p(M)$ is an $A$-module.  The proof that $\pi$ is an isomorphism then follows a standard argument.  If $\f\colon M\to N$ is a map of $R$-modules for some commutative ring $R$, then $\f$ is an isomorphism if and only if $\f_\g m\colon M_\g m\to N_\g m$ is an isomorphism for every maximal ideal $\g m\subset R$, where $M_\g m$ is $M$ localized at $\g m$.

\subsection{Proof of Theorem \ref{T:HKR}.}
We will only prove the result for the homology, since the result for cohomology is analogous.

Now, recall that $\HH_p(A)=\Tor_p^{A\otimes A}(A,A)$.  
We now localize at some maximal ideal $\g m\subset A$.  We then obtain the localized map
$$
\xymatrix{\Tor_p^{A\otimes A}(A,A)_\g m\ar[r]^<>(0.5){\pi_\g m}
&\left(\comO^p(A)\right)_\g m}.
$$
Ideally, we would hope that $\Tor_p^{A\otimes A}(A,A)_\g m=\Tor_p^{A_\g
m\otimes A_\g m}(A_\g m,A_\g m)$.  Indeed, this is the case.  If $M$ is
any $A$-bimodule, 
then $\Tor_0^{A\otimes A}(A,M)=A\otimes_{A\otimes A}M$.  Then we see that
$$
\Tor_p^{A\otimes A}(A,A)_\g m=(A\otimes_{A\otimes A}M)_\g m=
A_\g m\otimes_{A_\g m\otimes A_\g m}M_\g m=\Tor_p^{A_\g m\otimes\g m}(A_\g m,M_\g m).
$$
Recall that localization is an exact functor.  So, it commutes with
derived functors, in particular, it commutes with the derived functors
of $\otimes$, which are 
precisely the $\Tor$-groups.  So, we are reduced to the case of a map
$$
\xymatrix{\Tor_p^{A_\g m\otimes A_\g m}(A_\g m,A_\g
m)\ar[r]^<>(0.5){\pi_\g m}&\comO^p(A_\g m/\C)}.
$$
So, it suffices to consider a local ring.

\begin{lem}
Let $R$ be a commutative $\C$-algebra, $J\subset R$ an ideal finitely generated by $x_1,\ldots,x_n\in R$, where $x_1,\ldots,x_n$ is a regular sequence (that is, $x_{i+1}$ is not a zero divisor in $R/(x_1,\ldots,x_i)$ for all $i$).  Then
$$
\Tor_1^R(R/J,R/J)\simeq J/J^2.
$$
Indeed we have an isomorphism for all $p$,
$$
\Tor_p^R(R/J,R/J)\simeq\Lambda_{R/J}^p(J/J^2).
$$
\end{lem}

We will omit the proof of this lemma.  With this result in hand, we set $R=A\otimes A$, which we view as $\C[X\times X]$, and we let $J$ be the ideal of regular functions on $X\times X$ vanishing on the diagonal.  If $X$ is smooth, then $J$ is indeed generated by elements forming a regular sequence.  An application of the lemma finishes the proof.

\subsection{Digression: Applications to formality.} In this section we
are going
to use the relation between Hochschild and Chevalley-Eilenberg
complexes to obtain some (non-obvious)  formality results in the
 algebraic geometry.

We fix $X$, a smooth projective variety
with the structure sheaf $\oo_X$. We write
$\Dcoh(X)$ for the bounded derived category of complexes
of $\oo_X$-modules with coherent cohomology sheaves,
cf. e.g. [???]. 

Let $\imath: X \into X\times X$ be the
diagonal imbedding.
Associated to this imbedding,
one has a direct image functor $\imath_*:
\Dcoh(X)\map \Dcoh(X\times X)$,
and an inverse  image functor $\imath^*:
\Dcoh(X\times X)\map \Dcoh(X)$.
We will be interested in the composite functor
$$
\imath^*\imath_*:\
\Dcoh(X)\map \Dcoh(X\times X)\map \Dcoh(X).
$$
It is known, for instance, that for any coherent sheaf $\scr F$ on
$X$,
the cohomology sheaves of the object $\imath^*\imath_*\scr F$
are given by
$$
{\scr H}^i(\imath^*\imath_*\scr F)=\Omega_X^i\otimes_{_{\oo_X}} \scr F.
$$
The Proposition below implies that the object  $\imath^*\imath_*\scr
F\in \Dcoh(X)$
is actually {\em quasi-isomorphic} to a direct sum
of its cohomology sheaves.
\begin{prop}\label{dg_formal}There is a  quasi-isomorphism 
$$\imath^*\imath_*\oo_X\simeq
\bplus_{i=0}^{\dim X}\,\Omega_X^i[i]\quad\text{in}\quad\Dcoh(X).$$
\end{prop}
\begin{proof} Assume first that $X=\Spec A$, is affine.
We have $A\otimes A=\C[X\times X]$.
The kernel  $I:=\Ker(A\otimes A\to A)$ of the multiplication map
may be identified with the
defining ideal of the diagonal $X\sset X\times X$.

We know that the bar complex $\B_\idot A$ provides a resolution of the 
$\C[X\times X]$-module $A=\C[X\times X]/I$ by free
$\C[X\times X]$-modules. For each $i=0,1,\ldots,$
let $\wh\B_i A$ be the
$I$-adic completion of $\B_i A$.
The standard homotopy on the bar complex,
that shows that the complex $\B_\idot A$ is acyclic,
extends to the completions. 
Hence, the completed bar complex provides a resolution
of $\wh{A\ee/I}=\C[X\times X]/I$ of the form
$$ \ldots\to\wh\B_i A\to \wh\B_{i-1} A\to\ldots\to\wh\B_1 A\to
\wh\B_0 A\onto \C[X\times X]/I=\C[X],
$$
where  each term $\wh\B_i A$ is a {\em flat} $\C[X\times X]$-module.

The resolution above globalizes naturally to an arbitrary,
not necessarily affine, variety $X$. Specifically,
let $X$ be any smooth variety. For each $i\geq 0$, we
let $\wh\B_i(X):=\wh{\oo}_{X^{i+2}}$ be  the completion of
the structure sheaf  of the Cartesian power $X^{i+2}$
along the principal diagonal $X\into X^{i+2}$.
The projection $pr_{1,i+2}: X^{i+2}\to X\times X$,
on the first and last factors,
makes $\wh{\oo}_{X^{i+2}}$ a sheaf of $\oo_{X\times X}$-modules.
Thus, we have constructed
 a complex $\wh\B_\idot(X)$
 such that, for
each $i\geq 0$, we have

\noindent
\pb{$\wh\B_i(X)$ is a (not quasi-coherent)
sheaf of  {\em flat} $\oo_{X\times X}$-modules;}

\noindent
\pb{$\wh\B_i(X)$ is  set-theoretically supported on the
diagonal $X_\Delta\sset X\times X$, so may be regarded as a sheaf on $X_\Delta$;}

\noindent
\pb{For any Zariski-open affine subset $U=\Spec A \sset X_\Delta$, we have
$\Gamma(U, \wh\B_\idot(X))=\wh\B_\idot A$.}

We deduce that the object $\imath^*\imath_*\oo_X\in\Dcoh(X)$
is represented by the following
complex of $\oo_{X\times X}$-modules
$$\wh\B_\idot(X)\btimes_{\oo_{X\times X}}\,\oo_{X_\Delta}=
\left[\ldots\to\wh{\oo}_{X^3}\to \wh{\oo}_{X^2}\to \oo_{X_\Delta}\to
0\right].
$$
Now, the assignment sending, for each Zariski-open affine subset $U\sset X_\Delta$,
$a_0\otimes a_1\otimes\ldots\otimes a_i\in 
\Gamma(U, \wh{\oo}_{X^{i+1}})$ to $a_0\wedge da_1\wedge\ldots\wedge
da_i$ yields a quasi-isomorphism of the above complex with
$\bigoplus_{i\geq 0}\,
\Omega^i_X[i].$
\end{proof}

\begin{rem}[Kapranov]
Let $X$ be a smooth algebraic variety.
Consider the complex $RHom$ (over $X\times X$)
from $\oo_\Delta$ to itself, where $\Delta$ is the diagonal.
Its cohomology sheaves, i.e., the Ext's are just
the exterior powers of $\calT_X,$ the tangent bundle. One may
ask if this complex is quasiisomorphic
to the direct sum of its cohomology sheaves.

Further, given a coherent sheaf $F$ on $X,$
let $Quot_h(F)$, be the quot-scheme
that parametrizes quotient sheaves of $F$
with Hilbert polynomial $h$ or, equivalently, subsheaves
$K$ in $F$ with Hilbert polynomial $h_F-h.$ The tangent space to
$Quot_h(F)$
at the point
represented by a sheaf $K$
is the space $Hom_X(K, F/K).$ 
In   \cite{CK1}, the authors have constructed  a  derived quot-scheme,
a smooth
dg-manifold $RQuot$ whose tangent dg-space (i.e.
complex)
at a point $K$ as above is $RHom(K, F/K). $

In the special case where $ h=1$  and $F=\oo_X$, we have
$Quot_h=X$, the quotient sheaves being skyscrapers at
points of $X.$ Yet, the corresponding  derived quot-scheme is different, its
tangent space at $x\in X$ is a complex concentrated in nonnegative
degrees
whose  degree $i\geq 0$ cohomology is
equal to
$\wedge^{i+1} T_x(X)$. According to Proposition \ref{dg_formal},
this dg-manifold is
split, i.e., is quasiisomorphic to  $X$ with the
structure sheaf being the symmetric algebra of
the graded algebra sheaf formed by the direct sum of
$\Omega^i_X$ in degree $(-i+1), i>0.$
\eer

\begin{rem}
For any dg-manifold $ Y$ we have the scheme $\pi_0(Y),$
the spectrum of the $0$th cohomology of $\oo_Y,$
and on $\pi_0(Y)$ we have the sheaf of $\Lie_\infty$ algebras
obtained by restricting the tangent dg-bundle $TY$
onto $\pi_0(Y)$ in the sense of dg-schemes (note that
$\pi_0(Y)$ is a dg-subscheme in $Y$). Denote this sheaf
by $t_Y$ (small $t$ to avoid confusion).
This sheaf of $\Lie_\infty$ algebras defines $Y$ up to
quasiisomorphisms (duality between commutative
and Lie algebras, sheafified along $\pi_0(Y)$). To be more
precise, it determined the formal completion of $Y$
along $\pi_0(Y)$ (which is all that is needed in practice).
 Taking $Y=RQuot$ as before, we get $\pi_0(Y)=X$ and $t_Y =$
the quotient of $R\Hom(\oo_\Delta, \oo_\Delta)$ by the actual
Hom in degree $0,$ then shifted by $1.$ 
\eer
\section{Noncommutative differential forms}
\subsection{}
For an  associative {\em commutative} $\C$-algebra $A$
and  a left $A$-module $M$,
 we may consider the
space $\Der(A,M)$ of all derivations $\theta\colon A\to M$.
We have seen that
the functor $M\mto \Der(A,M)$ is representable
by the $A$-module $\comO(A)$ of K\"ahler differentials.

If $A$ is an associative {\em not necessarily 
commutative} $\C$-algebra,
the space  $\Der(A,M)$ is defined provided $M$ is an $A$-{\em bimodule}.
We are going to show that the functor $M\mto \Der(A,M)$,
defined on the category of $A$-bimodules,
is also representable. 

To this end, we let $m\colon A\otimes A\to A$ denote the multiplication map
viewed as a map of $A$-bimodules.
The kernel of $m$ is a sub-bimodule in $A\otimes A$.

\begin{defn}\label{ncO1}
We denote $\ncO^1(A):=\Ker[m\colon A\otimes A\to A],$ and call it the 
$A$-bimodule of \emph{noncommutative $1$-forms} on $A$. Thus, one has
the fundamental
exact sequence of $A$-bimodules
$$0\too\ncO^1(A)\too A\otimes A \,\stackrel{\tt{mult}}\too\, A\too 0.
$$
\end{defn}
\begin{prop}\label{Om_rep}
For every $M\in\bimod A$, there is a canonical isomorphism
$$
\Der(A,M)\simeq\Hom_{\bimod A}(\ncO^1(A),M).
$$
\end{prop}

{\em Thus, the functor $M\mapsto\Der(A,M)$ is representable
by the $A$-bimodule $\ncO^1(A)$.}

\begin{proof} Recall the bar resolution of $A$:
$$
\B_\idot A:\;\;\xymatrix{\cdots\ar[r]^<>(.5){b}&\B_1
\ar[r]^<>(.5){b}&\B_0\ar[r]^<>(.5){m}&
A\ar[r]&0},
$$
where  $\B_j=\B_jA=
A\otimes A^{\otimes j}\otimes A$, and the final map $\B_0=A\otimes A
\to A$ is the multiplication map.  Recall also that this complex is acyclic.  Now let $\delta\colon A\to M$ be a derivation.  We will define an $A$-bimodule map $\tilde\delta\colon A^{\otimes3}\to M$.  We set
$$
\tilde\delta(a'\otimes a\otimes a'')=a'\delta(a)a''.
$$

We claim that $\delta$ is a derivation if and only if $\tilde\delta$ is
a Hochschild 
2-cocycle.  Indeed,
\begin{align*}
(b\tilde\delta)(a_0\otimes &a_1\otimes a_3\otimes a_4)\\
&=\tilde\delta(a_0a_1\otimes a_2\otimes
a_3)-\tilde\delta(a_0\otimes(a_1a_2)\otimes a_3)
+\tilde\delta(a_0\otimes a_1\otimes(a_2a_3))\\
&=a_0a_1\delta(a_2)a_3-a_0\delta(a_1a_2)a_3+a_0\delta(a_1)a_2a_3.
\end{align*}
Clearly, the latter expression vanishes for all
$a_1,a_2,a_3$ if and only  $\delta$ is a derivation.  

 Since  $\tilde\delta$ is a $2$-cocycle, it is zero on all
 $3$-boundaries, that is, $\tilde\delta|_{d\B_2}=0$.  So, $\tilde\delta$
 descends to a well defined $A$-bimodule homomorphism
 $\tilde\delta\colon \B_1/d\B_2\to M$.  Since the bar complex is exact,
 $\B_1/d\B_2\simeq d\B_1$, and exactness again implies that $d\B_1$ is
 precisely the set of all $2$-cocycles, that is, $d\B_1$ is the kernel
 of the differential $\B_0\to A$, which is precisely the multiplication
 map.  So, 
$$\B_1/d\B_2\simeq d\B_1=\Ker[\B_0\to A]=\Ker[m\colon A\otimes A\to A]=\ncO^1(A).
$$
Hence
 $\tilde\delta\in\Hom_{\bimod A}(\ncO^1(A), M).$
  Clearly we can reverse this argument and produce from every
 $A$-bimodule map $
\ncO^1(A)\to M$ a derivation $A\to M$ in an analogous manner.
\end{proof}

Here is another picture of the $1$-forms on $A$.
Set $I=\Ker[m\colon A\otimes A\to A]$.
Let $\bar A$ denote the (vector space) quotient $A/\C=A/\C\cdot1_A$,
and write $x \mapsto \bar x$ for the projection $A\to\bar A=A/\C.$

\begin{prop}\label{aabar}\vi The map $d\colon
A\to I,\,x\mto dy=y\otimes 1_A-1_A\otimes y$
is a derivation.

\vii The map  $A\o \bar A\to I,
\,x\otimes\bar y\mto x\otimes y-xy\otimes1_A=x\,dy$ is well-defined,
and it is an isomorphisms of left $A$-bimodules.
\end{prop}

\begin{proof} 
 Endow $A\o \bar A$ with an $A$-bimodule structure by 
equipping it with the obvious left action and by setting
$$
(x\otimes\bar y)z=x\otimes\bar{yz}-xy\otimes\bar z.
$$
We let $\psi: x\otimes\bar y\mto x\otimes y-xy\otimes1_A$
denote the map considered in the Proposition.
Observe that $\psi$ is well-defined.  Indeed, 
if $y$ and $y+\lambda1_A$ are two representatives of $\bar y\in\bar A$, then
\begin{align*}
x\otimes(y+\lambda1_A)-x(\lambda1_A+y)\otimes1&=x\otimes y+x\otimes
\lambda1_A-\lambda1_Ax\otimes1_A-xy\otimes1_A\\
&=x\otimes y-xy\otimes1_A+x\otimes\lambda1_A-x\otimes\lambda1_A\\
&=x\otimes y-xy\otimes1_A.
\end{align*}
To show that $\psi$ is surjective, observe that $I$ is generated by
 terms of the form $x\otimes y-xy\otimes1_A$.  Indeed, this follows
 since $I=d\B_1$ (where again $\B_1=A^{\otimes3}$ from the bar complex).
 But then $\psi(x\otimes\bar y)$ is nothing more than $-d(x\otimes
 y\otimes1_A)$.  Since $d$ is an $A$-bimodule map, $d(x\otimes y\otimes
 z)=
d(x\otimes y\otimes1_A)z=-\psi(x\otimes\bar y)z$, so $I$ is 
generated by the image of $\psi$ as an $A$-bimodule.

To see that $\psi$ is an injection, 
consider the multiplication map $m\colon A\o A\to A$.  
Since $\im\psi\subset I=\Ker m$, we can descend to an $A$-bimodule map 
$\bar m\colon A\o A/\im\psi\to A$.  Define $\f\colon A\to
A\otimes A/\im\psi$
 by $\f(a)=a\otimes1_A+\im\psi$.  Then
\begin{align*}
\bar m(\f(a))&=\bar m(a\otimes1_A+\im\psi)=m(a\otimes1_A)=a\\
&\f(\bar m(x\otimes y+\im\psi))=\f(xy)=xy\otimes1_A+\im\psi.
\end{align*}
Since $x\otimes y-xy\otimes1_A\in\im\psi$, we see that
$xy\otimes1_A+\im\psi=
x\otimes y+\im\psi$.  So, $\f$ is an inverse to $\bar m$, hence $\bar m$
is
 an isomorphism.  But then $\im\psi=\Ker m=I$, as claimed.
\end{proof}

Further, combining Propositions \ref{Om_rep} and \ref{aabar},
we obtain the following result, which is completely
analogous to a similar result for K\"ahler differentials
proved in \S\ref{kahler_diff}.
\begin{cor}\label{non_kal} Let $A$ be an associative not necessarily
commutative algebra with unit and $M$ an $A$-bimodule.
For any derivation $\theta: A\to M$,
the assignment $x\,dy\mto x\cdot\theta(y)$
gives a well-defined $A$-bimodule map $\ncO^1(\theta)$
that makes the following diagram commute
$$
\xymatrix{
&&A\ar[dll]_<>(.5){d}\ar[drr]^<>(.5){\theta}&&\\
\ncO^1(A)\ar[rrrr]^<>(.5){\ncO^1(\theta)}&&&&M\qquad\qed
}
$$
\end{cor}

Next we observe that 
\beq{AAee}
H_0(A,A\ee)=A\ee/[A,A\ee]=A\ee\otimes_{A\ee}A=A,
\end{equation}
where the isomorphism $A\ee/[A,A\ee]\iso A$
 is induced, explicitly, by the
assignment $A\ee\ni u\otimes v\mapsto vu$.
One verifies unrevelling the definitions,
that the canonical map $H_0(A,A\ee)\to
H_0(A,A)$ induced by the multiplication morphism
$A\ee\to A$ is nothing but the natural
projection $A\onto A/[A,A]$. 

We  have the following maps
\beq{comp_om}
\xymatrix{
A\o \bar A
\ar[rr]^<>(.5){\text{Prop. \ref{aabar}}}_<>(.5){\sim}&&
\ncO^1(A)=I\ar@{^{(}->}[r]&A\ee
\ar[r]&A\ee/[A,A\ee]=A.
}
\end{equation}
One verifies that the composite map is given by the
formula $x\otimes\bar y\mto yx-xy.$

\subsection{The differential envelope.}
Our goal is to construct an analogue of
de Rham differential on noncommutative differential forms.
To this end, we recall

\begin{defn}
A \emph{differential graded algebra} (DGA) is a graded algebra 
$D=\bigoplus_nD^n$ equipped with a super-derivation 
$d\colon D^n\to D^{n+1}$ such that $d^2=0$.

Given an associative algebra $A$, its \emph{differential envelope}
$D(A)$ is the solution to the following lifting problem.  Consider the
category whose objects consist of all algebra maps from $A$ into the
zeroth degree of some DGA $D^0$.  If $\f\colon A\to D^0$ and $\psi\colon
A\to E^0$ are two such objects, a morphism from $\f$ to $\psi$ is a map
$\theta\colon D\to E$ such that $\theta(D^n)\subset E^n$ and
$d_E\circ\theta=d_D$.  The differential envelope $D(A)$ is then defined
to be the initial object in this category, that is, there is a map of
algebras $i\colon A\to D(A)^0$ such that for any algebra map $\f\colon
A\to D^0$ ($D$ a DGA), there is a unique map of DGA's
$\psi\colon{D}(A)\to D$
 such that $\psi\circ i=\f$.

Explicitly, ${D}(A)$ is generated by the algebra $A$ and all symbols of 
the form $\bar a=da$, where $a\in A$ and $d(ab)=(da)\,b+a\,(db)$.
\end{defn}
\begin{rem}\label{fin_gen_Om} It is clear that if $A$ is generated, as an algebra, by 
$a_1,\ldots,a_n,$ then the elements $a_1,\ldots,a_n, da_1,\ldots,da_n,$
generate $D(A)$ as an algebra. Thus, if $A$ is a finitely generated
algebra
then so is the algebra $D(A)=\oplus_{i\geq 0}\,D^i(A)$. 
Moreover, one verifies similarly that in that case each
homogeneous component $D^i(A)$ is finitely generated as an $A$-bimodule.
\eer

\begin{defn}
Define  the algebra $\ncO^\hdot(A)$ of \emph{noncommutative differential forms} on $A$
to be the tensor algebra (over $A$) of the $A$-bimodule $\ncO^1(A),$
that is
\begin{align*}
\ncO^\hdot(A)&:=T_A\ncO^1(A)=
A\oplus\ncO^1(A)\oplus
T^2_A\ncO^1(A)\oplus T^3_A\ncO^1(A)\oplus\cdots.
\end{align*}
\end{defn}

It turns out that there is a canonical
graded algebra isomorphism
\beq{dommap}
\ncO^\hdot(A)\to{D}(A).
\eeq

To construct this isomorphism, we
first  consider the canonical
algebra map $i\colon A\to{D}(A)$ (this turns out to be an injection).
We also have a derivation $d\colon A\to{D}(A)$ given by $a\mapsto
da=\bar a$.  By our definition of $\ncO^1(A)$, there is then a unique
$A$-bimodule map $\ncO^1(A)\to{D}(A)$.  
This $A$-bimodule morphism extends canonically to an algebra map
\eqref{dommap}.

\begin{prop}\label{DOm}
With the notation as above, the
map $\ncO^\hdot(A)\to{D}(A)$ in  \eqref{dommap} is a
graded algebra isomorphism. 
\end{prop}

Thus, using the isomorphism of the Proposition,
we transport the differential on $D(A)$ to obtain
 a degree one super-differential 
$d\colon\ncO^n(A)\to\ncO^{n+1}(A)$ making $\ncO^\hdot(A)$ a DGA.

\begin{cor} The graded
algebra $\ncO^\hdot(A)$ can be given the structure of DGA with
differential transported from the one on $D(A)$.
\end{cor}

Proposition \ref{DOm} will  be proved later as a special case of 
Theorem \ref{relative_D} below. One also has the following
important result proved
in \cite{CQ1}.

\begin{thm}
There is a canonical isomorphism of left $A$-modules
$$
\ncO^\hdot(A)\cong\bigoplus_{n=0}^\infty A\o \bar A^{\otimes n}
$$
With this isomorphism,
the differential $d\colon\ncO^n(A)\to\ncO^{n+1}(A)$
is given by the formula
$$d(a_0\otimes\bar{a_1}
\otimes\cdots\otimes\bar{a_n})=
1\o\bar{a_0}\otimes\bar{a_1}
\otimes\cdots\otimes\bar{a_n}.$$\end{thm}

 We denote the summand $A\o \bar A^{\otimes n}$
by $\ncO^n(A)$.  
We denote the element $a_0\otimes\bar{a_1}
\otimes\cdots\otimes\bar{a_n}\in\ncO^n(A)$ by $a_0\,da_1\cdots da_n$.
so that the differential reads
$$
d(a_0\,da_1\cdots da_n)=da_0\,da_1\cdots da_n.
$$

\subsection{The universal square-zero extension.}
The algebra imbedding $A\into \ncO^\hdot(A)$ makes,
for each $k$, the graded component $\ncO^k(A)$ an $A$-bimodule.
We define an algebra structure on the vector space $A\oplus \ncO^2(A)$
by the formula
$$
(a,\omega)(a',\omega')=(aa',a\omega'+\omega
a'+c(a,a')),\quad\text{where}\quad c(a,a')=da\,da'.
$$
It is straightforward to check that this way we get
an associative algebra, to be denoted $A{\sharp}_c\ncO^2(A)$. 
Clearly, $\ncO^2(A)$ is a two-sided ideal in  $A{\sharp}_c\ncO^2(A)$.
Moreover,
the natural imbedding $i: \om\mapsto (0,\om)$, and
the projection $(a,m)\mapsto a,$ give rise to a square-zero extension
\begin{equation}\label{E:SqZero}
0\to \ncO^2(A)\stackrel{i}\too A{\sharp}_c\ncO^2(A)\stackrel{p}\too A\to 0.
\end{equation}
\begin{lem} The  square-zero extension \eqref{E:SqZero} is universal,
i.e.,
for any  square-zero extension $I\into \widetilde A \onto A$,
there exists a unique $A$-bimodule map $\varphi$ making the following
diagram commute
$$
\xymatrix{
0\ar[r]&
\ncO^2(A)\ar[r]\ar[d]^<>(.5){\varphi}&A{\sharp}_c\ncO^2(A)\ar[r]\ar[d]^<>(.5){\varphi}&
A\ar[r]\ar@{=}[d]^<>(.5){\id}&0\\
0\ar[r]&
I\ar[r]&\widetilde A\ar[r]&
A\ar[r]&0.}
$$
\end{lem}
\begin{proof}
Choose a $\k$-linear splitting $c: A \into \widetilde A$,
and note that, for any $a_1,a_2\in A$, we have
$c(a_1)\cdot c(a_2)- c(a_1a_2)\in I$. Define the map 
$$\widetilde\varphi\colon 
A\otimes {\bar A}^{\otimes 2}
\map I,
\quad a\otimes a_1\otimes a_2\mto a\cdot[c(a_1)\cdot c(a_2)- c(a_1a_2)].
$$
We use the isomorphism $\colon\ncO^2(A)\cong
A\otimes {\bar A}^{\otimes 2}$ to view the map above
as a map $\widetilde\varphi\colon
\colon\ncO^2(A)\map I,\, a\,da_1\,da_2\mapsto \widetilde\varphi(a\,da_1\,da_2).$
We leave to the reader to verify that this is
an $A$-bimodule map.

Further, write $\widetilde A\cong A \oplus I$ for the
vector space decomposition corresponding to the splitting $c$.
One shows that the following assignment
$$A{\sharp}_c\ncO^2(A)\map 
A \oplus I=\widetilde A,\quad
a'\oplus a\,da_1\,da_2\mto a'\oplus \widetilde\varphi(a\,da_1\,da_2)
$$
gives a map with all the required
properties.
\end{proof}

\subsection{Hochschild differential on non-commutative forms.}
Computing {Hoch}{schild}
 homology of the $A$-bimodule $M=A$ using the bar resolution, we get
$\overline C_\idot(A,A)=\overline \B_\idot(A)\otimes_{A\ee}A$, where 
$\overline \B_\idot(A)$ 
is  the reduced bar complex and $A\ee=A\otimes A^{\op}$.  Thus
$$
{\overline{C}}_n(A,A)=\ncO^n(A)\otimes A\otimes_{A\ee}A\simeq\ncO^n(A).
$$
The Hochschild differential ${\overline{C}}_{n+1}(A,A)\to 
{\overline{C}}_n(A,A)$ is then given by the differential 
$b\colon\ncO^{n+1}(A)\to\ncO^n(A)$ mentioned before.  So, $\HH_\idot(A)=H_\idot(\ncO^\hdot(A),b)$.

\begin{prop} The Hochschild differential $b:
{\overline{C}}_{n+1}(A,A)\map{\overline{C}}_n(A,A)$,
viewed as a map $\colon\ncO^{n+1}(A)\to\ncO^n(A)$,
is given by the following explicit formula
$$
\delta(a_0\,da_1\cdots da_n\otimes a')=
\delta(\alpha\,da_n\otimes a')=(-1)^{\deg\alpha}(\alpha a_n\otimes a'-\alpha\otimes a_na').
$$
\end{prop}
\begin{proof} Direct calculation, see \cite{CQ1}.
\end{proof}

As an application, 
 we are going to construct the following
 exact sequence 
\beq{ex_H_1}
\HH_1(A)\map\ncO^1(A)/[A,\ncO^1(A)]
\stackrel{b}\map [A,A]\map 0.
\end{equation}
where the map $b$ is given by the assignment
\beq{quil_map}
b:\ \ncO^1(A)/[A,\ncO^1(A)]\map A,
\quad u\,dv\mto [u,v]
\end{equation}

We first verify  that the map
$b$ in \eqref{quil_map} is well-defined,
i.e., we have $b([x\,dy,a])=0$ for any
$[x\,dy,a]\in[\ncO^1(A),A]$.    Indeed,
we compute
\begin{align*}
b([x\,dy,a])&=b(x\,dy a)-b(ax\,dy)\\
&=b(x\,d(ya))-b(xy\,da)-[ax,y]\\
&=[x,ya]-[xy,a]-[ax,y]=-([ya,x]+[xy,a]+[ax,y])\\
&=x(ya)-(ya)x-(xy)a+a(xy)-(ax)y+y(ax)=0.
\end{align*}

Now, to construct \eqref{ex_H_1}
is to use a long exact sequence for
 Hochschild homology arising from
the fundamental short 
exact sequence $0\to\ncO^1(A)\to A\ee\to A\to 0$:
$$\ldots\to H_1(A,A\ee)\to
H_1(A,A)\to H_0(A,\,\ncO^1(A))\to H_0(A,A\ee)\to
H_0(A,A)\to 0.
$$
We have $H_0(A,A)= A/[A,A],$
and $H_0(A,\,\ncO^1(A))=\ncO^1(A)/[A,\ncO^1(A)].$
Observe further  that 
$H_1(A,A\ee)=0$, since $A\ee$ is a free $A$-bimodule.
Thus, since $H_0(A,A\ee)=A$ by \eqref{AAee},
the long exact sequence above reduces to a
 short exact sequence as in \eqref{ex_H_1}.

\subsection{Relative differential forms}
 Let $A$ be an associative algebra, $B\subset A$ a subalgebra, and $M$ an
$A$-bimodule.

\begin{defn} A derivation  $f\colon A\to M$  such that $f(b)=0$  for all
$b\in B$
is said to be a derivation
 from $A$ to $M$ \emph{relative to $B$}. We write
$\Der_B(A,M)\sset \Der(A,M) $ for the subspace of all such derivations,
which form a Lie subalgebra.
\end{defn}
Notice that a derivation is a relative one, namely it is relative to
the 
subalgebra $\C1_A\subset A$.  For ease of notation, we will denote this subalgebra simply by $\C$.
\subsection{The Commutative Case}
Let $A$ be a commutative algebra, and led let $M$ be a left $A$-module.
We make $M$ into an $A$-bimodule by equipping it with the symmetric
bimodule structure.  Then $M\mapsto\Der_B(A,M)$ defines a functor
$\Der_B(A,-)\colon\Lmod A\to\Vect$, where $\Vect$ denotes the category of $\C$-vector spaces.  As in the case where $B=\C$, $\Der_B(A,-)$ can be represented by an $A$-module $\comO^1(A/B)$, which we call the \emph{relative K\"ahler differentials}.

If we regard $A$ and $B$ as the coordinate rings of affine varieties and
set $X=\Spec A$ and $Y=\Spec B$, then the embedding $B\inj A$ induces a
surjection $\pi\colon X\to Y$.  
Recall that we view $\Der(A)$ as the global sections of
$\calT_X$, the tangent bundle of $X$, that is, as the
space $\Gamma(X, \calT_X)$ of algebraic
vector fields on $X$. 
 In this geometric picture, 
$\Der_B(A)=
\Der_B(A,A)$ 
corresponds to
 the subspace in $\Gamma(X, \calT_X)$ formed by the
vector fields on $X$ which are parallel to fibers of the projection $\pi
: X\to Y$. More precisely,
the differential of the map $\pi$ may be viewed as
a sheaf map $\pi_*:\ \calT_X\to \pi^*\calT_Y$.
Let $\calT_{X/Y}:=\Ker[\calT_X\to \pi^*\calT_Y],$
be the sheaf whose sections are the vector fields tangent to the fibers.
If $\pi$ has surjective differential, then
 $\pi_*$ is a surjective map of locally free sheaves, hence,
its kernel is again locally free, and we have a short exact
sequence
$$0\too \calT_{X/Y}\too \calT_X\too \pi^*\calT_Y\too 0.
$$

Similarly, we view $\comO^1(A/B)$ as the space of  differential
$1$-forms along the fibers,
that is of sections of 
the {\em relative cotangent bundle}
$\calT^*_{X/Y}$.
Explicitly, dualizing the short  exact
sequence above, one gets the dual  exact
sequence
\beq{rel_cotangent}
0\too \pi^*\calT_Y^*
\too \calT^*_X\too \calT^*_{X/Y}\too 0.
\end{equation}
This presents the space $\comO^1(A/B)=\Gamma(X, \calT^*_{X/Y})$
as a quotient of $\comO^1(A)$.

We can also realize $\comO^1(A/B)$ by relative Hochschild homology.  In particular, it is easy to see (by the same arguments for absolute differential forms) that we can write
$$
\comO^1(A/B)=\Tor_1^{A\otimes_BA}(A,A).
$$
This can be expressed in terms of the kernel $I$of the multiplication
map $A\otimes_BA\to A$.  Indeed, we find that $\comO^1(A/B)\simeq
I/I^2$, as before.  This is defined to be the first relative Hochschild
group, $\HH_1(A;B)$.  It is 
also easy to see (as in the absolute case) that $\comO^1(A/B)=\HH_1(A;B)\simeq A\otimes\bar A$, where now $\bar A=A/B$.  

As in the absolute case, we define $\comO^\hdot(A/B)=T_B\comO^1(A/B)$.  We then obtain the following result, analogous to the absolute case.

\begin{thm}
$\comO^\hdot(A/B)$ is a graded commutative differential envelope of $A$, specifically, it is the universal object in the category whose objects are algebra maps $f\colon A\to D^0\subset D$ such that $d_D\circ f|_B=0$, where $D$ is a graded commutative DGA.
\end{thm}

\subsection{The Noncommutative Case}
We now assume only that $A$ is associative.  Consider the functor 
$$
\Der_B(A,-)\colon\bimod A\to\Vect.
$$  
Again, this functor is representable, and we call its representing object $\ncO^1(A/B)$.  Once again, $\ncO^1(A/B)$ is the kernel of the multiplication map, but this time viewed as a map $A\otimes_BA\to A$.  That is, $\ncO^1(A/B)$ fits in the exact sequence of $A$-bimodules:
$$
\xymatrix{0\ar[r]&\ncO^1(A/B)\ar[r]&A\otimes_BA\ar[r]&A\ar[r]&0}.
$$

It is easy to see that $\ncO^1(A/B)\simeq A\otimes_BA/B$ as a
$B$-bimodule.  
Indeed, the isomorphism is given by sending an element $x\otimes\bar y$
of 
$A\otimes_BA/B$ to $x\otimes y-xy\otimes1_A$ (of course, one needs to
check
 that this is even well-defined).

Recall that in the commutative case, we thought of $\comO^1(A/B)$
as a relative cotangent bundle 
 $\Gamma(X, \T^*_{X/Y})$, where $X=\Spec A$ and $Y=\Spec B$, and 
$\pi\colon X\to Y$ is the surjection induced by the inclusion $B\inj A$. 
Similarly to 
 the exact sequence
\eqref{rel_cotangent}, we have the following exact sequence in the noncommutative case
\begin{align*}
&0\to\Tor_1^B(A,A)\to
A\otimes_B\ncO^1(B)\otimes_BA\to
\ncO^1(A)\to\ncO^1(A/B)\to 0
\end{align*}

There is a connection between $\ncO^1(A/B)$ and differential envelopes of $A$.  Indeed, we have the important following result.

\begin{thm}\label{relative_D}
Let $D^\hdot(A/B)$ denote the relative differential envelope of $A$.  That is, $D^\hdot(A/B)$ is universal in the category whose objects are algebra maps $f\colon A\to D^0\subset D$ where $f|_B=0$ and $D$ is a DGA.  Then
$$
D^\hdot(A/B)\simeq T_A^\hdot\ncO^1(A/B)\simeq A\otimes_BT_B^\hdot(A/B).
$$
\end{thm}

As usual, we set $\ncO^\hdot(A/B)=T_A\ncO^1(A/B)$.

\begin{proof}
We will prove that $D^\hdot(A/B)\simeq\ncO^\hdot(A/B)$.  Of course, this proof will also be valid in the case $B=\C$.  We will complete the proof in two steps.  First, we will observe that $D^1(A/B)$ and $\ncO^1(A/B)$ are isomorphic.  We will then show that $D^\hdot(A/B)\simeq T_AD^1$.

For the first step, notice that $D^1$ is an $A$-bimodule, since we have a map $i\colon A\to D^0$.  If we compose $i$ with the differential $d$ on $D$ (and denote this map by $d$), we obtain an $A$-bimodule derivation $d\colon A\to D^1$.  Since $i$ vanishes on $B$, so does $d$.  We wish to show that $d\colon A\to D^1$ represents $\Der_B(A,-)$, which will show that $D^1$ and $\ncO^1(A/B)$ are isomorphic by the latter's universal property.  Let $\delta\in\Der_B(A,M)$, where $M$ is any $A$-bimodule.  We wish to complete the following diagram
$$
\xymatrix{
A\ar[r]^\delta\ar[d]^d&M\\
D^1\ar@{-->}[ur]}.
$$
We shall perform this by using the square zero construction.  Recall that we set $A{\sharp}M$ to be, as a $\C$-vector space, the direct sum $A\oplus M$.  The product on $A{\sharp}M$ is defined in such a way as to make $A$ a subalgebra and $M^2=0$.  We make $A{\sharp}M$ into a DGA by declaring that $\deg A=0$, $\deg M=1$, and by setting $d(a,m)=(0,\delta m)$.  Define $f\colon A\to A{\sharp}M$ by $f(a)=(a,0)$.  Observe that $d\circ f(b)=0$ for all $b\in B$, hence we obtain a DGA map $f'\colon D(A/B)\to A{\sharp}M$ from the universal property of $D(A/B)$.  Consider the restriction of $f'$ to $D^1(A/B)$.  This must map $D^1(A/B)$ into $M=\{0\}{\sharp}M$, since only $M$ has degree one.  Then using the fact that $f'$ is a DGA map, we find for all $a\in A$, $f'(d(a))=\delta(a)$.  So, $D^1(A/B)$ satisfies the universal property for $\ncO^1(A/B)$, hence these two $A$-bimodules are isomorphic.

Now we wish to show that $D$ is generated as the free algebra over $D^1$.  The inclusion of $D^1\inj D$ as the degree one elements induces an algebra map $\f\colon T_AD^1\to D$.  We need only check that $\f$ is an isomorphism.  Here are two proofs of this fact.  First, we adopt the notation
$$
T_A^nD^1=(D^1)^{\otimes  n}.
$$
Then
$$
T_A^nD^1=T_A^{n-1}D^1\otimes_AD^1\simeq T_A^{n-1}D^1\otimes_A(A\otimes_BA/B),
$$
since we know that $D^1\simeq\ncO^1(A/B)$.  But then
$$
T_A^{n-1}D^1\otimes_A(A\otimes_BA/B)=T^{n-1}_AD^1\otimes_BA/B\simeq D^n.
$$

The second proof is in some way more elementary, as we simply construct an inverse $\psi\colon D\to T_AD^1$ to $\f$.  We recall that $D$ contains $A$ (since $D$ comes with an embedding $A\inj D^0$), and $D$ also contains $\bar A=dA=A/B$.  So, we consider the algebra
$$
D'=T_\C(A+\bar A)/\langle\{a\otimes a'-aa'\otimes1_A,\bar{aa'}=\bar a\otimes a'+a\otimes\bar{a'}\}\rangle.
$$
So $D'$ is simply the free algebra over $A$ and $\bar A$, modulo the relations which give the desired product (as in $\ncO^\hdot(A/B)$) and make $a\mapsto\bar a$ a derivation.  Clearly, $D'$ is precisely $D$.  It is also easy to see that there is a map $D'\to T_AD^1$, and all relations are mapped to zero.  So, $D$ and $T_AD^1$ are isomorphic.
\end{proof}

\section{Noncommutative Calculus}
\subsection{}
Let $A$ be an associative unital $\C$-algebra,
thought of as the coordinate ring of a `noncommutative space'.
Then,  any automorphism $F\colon A\to A$
may be thought of as an  automorphism of that  `noncommutative space'.
Similarly, a derivation $\theta\colon A\to A$
may be viewed as an `infinitesimal
automorphism' of $A$, in the sense that
if the linear map $\theta\colon A\to A$ could have been
exponentiated, i.e., if the
infinite series $\exp(\theta)=\id_A + \theta +
\frac{1}{2}\theta\ccirc\theta
+\frac{1}{3}\theta\ccirc\theta\ccirc\theta+\ldots$
made sense as a map $A\to A$, then a formal
computation shows that 
the map $\exp(\theta)$ would have been an automorphism of $A$.
Geometrically,
one thinks of $\theta$ as a `vector field'
on a noncommutative space; then
$t\mto \exp(t\cdot\theta)$ is the
one-parameter flow of automorphisms of that noncommutative space
generated by our vector field.

Accordingly, an algebra automorphism $F\colon A\to A,$ 
resp.,  a derivation $\theta\colon A\to A$,
induces an automorphism, resp.,  a {\em Lie derivative} endomorphism
$
\L_\theta$,
of 
any of the objects $\Der(A)$, $\ncO^1(A)$ and $\ncO^\hdot(A)$.  

We provide more details. Fix  a derivation $\theta\colon A\to A$.
Since the functor $\Der(A,-)$ is represented by the bimodule
 $\ncO^1(A)$,
we deduce that  there is a unique $A$-bimodule homomorphism
$i_\theta\colon\ncO^1(A)\to A$
corresponding to $\theta$.  Observe further that since
 $\ncO^\hdot(A)=T_A\ncO^1(A),$ is a free algebra of the bimodule
$\ncO^1(A),$ there is a  unique way to extend
the map $i_\theta\colon\ncO^1(A)\to A$ to a degree $(-1)$
{\em super-derivation}
$i_\theta\colon\ncO^\hdot(A)\to\ncO^{\bullet-1}(A)$
of the  algebra $\ncO^\hdot(A)$.
Explicitly, we have
$$
i_\theta(a_0\,da_1\cdots da_n)=\sum_{j=1}^n(-1)^{j-1}a_0\,da_1\cdots\theta(a_j)\cdots da_n.
$$

Recall next that, for any two {\em super-derivations}
$\partial_1,\partial_2$ of {\em odd} degree,
their super-commutator $[\partial_1,\partial_2]:=
\partial_1\ccirc\partial_2+\partial_2\ccirc\partial_1$
is an {\em even} degree  derivation.
In particular, for $\partial_1=d$ and $\partial_2=i_\theta$,
we obtain a degree zero derivation
$d\ccirc i_\theta + i_\theta\ccirc d$.
Similarly, for $\partial_1=\partial_2=d$,
resp.,  for $\partial_1=\partial_2=i_\theta$,
we obtain a degree $(+2)$-derivation $d^2$, resp.,
 a degree $(-2)$-derivation 
$(i_\theta)^2$.

Next we apply Lemma \ref{triv} to $S=\ncO^1(A)\sset \ncO^\hdot(A)=R$.
It is straightforward to get from
definitions that, for any $s\in \ncO^1(A)$ one has $d^2(s)=0=(i_\theta)^2(s).$
It follows that $d^2=0$ and $(i_\theta)^2=0$ identically on
$\ncO^\hdot(A)$.

Similarly,  we define the Lie derivative with respect to
$\theta\in \Der(A)$ as a map
$\L_\theta\colon\ncO^\hdot(A)\to\ncO^\hdot(A)$
given by the formula
$$
\L_\theta(a_0\,da_1\cdots da_n)=\theta(a_0)\,da_1\cdots da_n+
\sum_{j=1}^n a_0\,da_1\cdots d\theta(a_j)\cdots da_n.
$$
A direct calculation shows that
for any $s\in \ncO^1(A)$ one has 
$\L_\theta(s)=[d,i_\theta](s)$. It follows by  Lemma \ref{triv} that
the following Cartan formula holds  on $\ncO^\hdot(A)$
$$
\L_\theta=d\circ i_\theta+i_\theta\circ d.
$$ 
Similarly, one verifies the following identities
\begin{equation}\label{Cart_id}
[\L_\theta,\L_\gamma]=\L_{[\theta,\gamma]},
\quad
[\L_\theta,i_\gamma]=i_{[\theta,\gamma]},
\quad
i_\theta^2=0.
\end{equation}
To prove these identities, observe that in each case,
the identity in question is obvious on the elements of
$\ncO^1(A)$. Hence, using the same argument as above,
we deduce that it holds on the whole of $\ncO^\hdot(A)$.

\subsection{Operations on Hochschild complexes.}
Let $C_\idot(A,A)\,,\,C_k(A,A)=A^{\otimes k},$ be the Hochschild chain complex for $A$.
For any  Hochschild $p$-cochain $c\in C^p(A,A)$ and $k\geq p$,
define a contraction
operator
$i_c: C_k(A,A)\to C_{k-p}(A,A)$ by the formula
$$i_c:\; a_0\otimes \ldots\otimes a_k\mto c(a_1,\ldots,a_p)\cdot
a_{p+1}\otimes \ldots\otimes a_k.
$$

Now suppose we have a derivation $\delta\colon A\to A$.  
Then we can extend $\delta$ to a derivation on each $\B_n$ in the bar complex, namely
$$
\delta(a_1\otimes\cdots a_n)=\sum_{j=1}^na_1\otimes\cdots\otimes\delta(a_j)\otimes\cdots a_n.
$$
Then $\delta\colon\Hom_{\bimod A}(A^{\otimes n},A)\to\Hom_{\bimod
A}(A^{\otimes n},A)$.
  It is not hard to see that $\delta$ commutes with the bar
  differential, hence 
it induces a derivation on the Hochschild cohomology of $A$.

We can generalize the above to all $p\geq 1$,
and define the Lie derivative with respect to
a cochain $c\in C^p(A,A)$ as
an operator $\L_c: C_k(A,A)\to C_{k-p+1}(A,A)$ by the formula
\begin{align*}
\L_c\colon a_0\otimes \ldots\otimes a_k\mto&
\sum_{i=0}^{k-p} (-1)^{(p-1)(i+1)} a_0\otimes \ldots\otimes a_i
\otimes c(a_{i+1},\ldots,a_{i+p})\otimes\ldots\otimes a_k\\
&+\sum_{i=k-p}^k (-1)^{k(j+1)} c(a_{j+1},\ldots,a_0, \ldots)
\otimes a_{p+j-k}\otimes\ldots\otimes a_j.
\end{align*}
\begin{rem} For the 2-cochain $m: a,b \mapsto a\cdot b,$
given by the product in the algebra $A$, the operator
$\L_m$ is nothing but the Hochschild differential on $C_\idot(A,A)$.
\eer

Recall the canonical noncommutative  Gerstenhaber
algebra structure on the Hochschild
cochain complex $C^\hdot(A,A)$.
Recall further that to any (not necessarily commutative)
Gerstenhaber algebra $G^\hdot$ one can 
associate another Gerstenhaber algebra $G^\hdot_\eps,$
called the {\em $\eps$-construction}.
We apply the $\eps$-construction to the Gerstenhaber
algebra $C^\hdot(A,A)$. 

We have the following noncommutative
analogue of Proposition~\ref{TT_prop}.

\begin{thm} \label{TT_thm} For any $b,c\in C^p(A,A)$,
the following formulas
$$
(b+\eps c)\cdot_{\eps} \alpha:= (-1)^{\deg b}\cdot
i_b{\alpha},\aand
\{b+\eps c\,,\,\alpha\}_{\eps}  := \L_b{\alpha}+
\eps\cdot i_c{\alpha}
$$
make $C_\idot(A,A)_\eps$ a Gerstenhaber module
over  $C^\hdot(A,A)_\eps.$
\end{thm}

\subsection{The functor of `functions'}
Let $A$ be a not necessarily commutative associative algebra
thought of as the coordinate ring of a `noncommutative scheme'.
We would like to introduce a vector space $\sr(A)$
playing the role of  `the space of regular functions' on that scheme.
Of course, if $A=\C[X]$ is the coordinate ring of an ordinary
 commutative scheme $X$, then regular  functions
on $X$ are by definition the elements of $A$. So, one might guess that,
 in the noncommutative
case, the equality $\sr(A)=A$ still holds.
This cannot be quite right, however. Indeed, one expects the
space $\sr(A)$ to be a Morita invariant of $A$ since only Morita
invariant notions are `geometrically meaningful'.
Thus, starting from an  ordinary
 commutative scheme $X$, for any $n=1,2,\ldots,$ we may form the
algebra $A=\Mat_n\C\otimes\C[X]$,
which is Morita equivalent to $\C[X]$.
Thus, we would like our definition of the space $\sr(A)$ be such that
the following holds
$$\sr(\Mat_n\C\otimes\C[X])=\sr(\C[X])=\C[X].$$
By Morita invariance of Hochschild homology, see Proposition
\ref{Morita_hoch},
this requirement is satisfied if we introduce the following
\begin{defn}
We define the \emph{space of functions} associated to  an associative
algebra
$A$ to be the
vector space
$$
\sr(A):=A/[A,A].
$$

\subsection{Karoubi-de Rham complex.}
We are going to show that for any associative algebra $A$ with unit,  the
complex $(\ncO^\hdot(A),\,d)$ of noncommutative differential forms
has trivial cohomology:
\beq{om_dr_vanish}
H^i(\ncO^\hdot(A),\,d)=\begin{cases} \k & \text{if}\quad i=0\\
0&\text{if}\quad i>0.
\end{cases}
\end{equation}

To see this, recall the isomorphism
$\ncO^p(A)=A\otimes \overline{A}^{\otimes p},$
where $\oA=A/\k$. The differential $d$
corresponds to the natural projection
$$A\otimes \overline{A}^{\otimes p}\to 
\oA\otimes \overline{A}^{\otimes p}\cong\k\otimes\oA^{\otimes(p+1)}
\sset A\otimes \overline{A}^{\otimes(p+1)}.
$$
The kernel of the latter projection clearly
equals $\k\otimes \overline{A}^{\otimes p}$,
which is exactly the image of the differential
$d: A\otimes \overline{A}^{\otimes p-1}\to
\k\otimes\oA^{\otimes p}$. This proves \eqref{om_dr_vanish}.

Therefore, the differential $d$ on $\ncO^\hdot(A)$ does not give
rise interesting cohomology theory. Things become better with
the following definition.

\begin{defn}
The \emph{noncommutative de Rham complex} of $A$
is a graded vector space defined by
$$
\DR^\hdot(A):=\sr\bigl(\ncO^\hdot(A)\bigr)=
\ncO^\hdot(A)/[\ncO^\hdot(A),\ncO^\hdot(A)],
$$
where $[-,-]$ denotes the super-commutator.  
\end{defn}

The differential $d:\ncO^\hdot(A)\to \ncO^{\bullet+1}(A)$
descends to a well-defined differential  $d: \DR^\hdot(A)
\to \DR^{\bullet+1}(A)$, making the de Rham complex a differential
graded $\C$-vector space.

For example, $\DR^0(A)=A/[A,A]$, since only the 
degree zero terms of $\ncO^\hdot(A)$ contribute.  
Similarly, we have
\beq{HH_DR}
\DR^1(A)=\ncO^1(A)/[A,\ncO^1(A)]=\HH_0(A,\,\ncO^1(A)).
\end{equation}
Hence, from \eqref{ex_H_1}, we obtain a canonical short exact sequence
\beq{HHDR}
0\map \HH_1(A)\map \DR^1(A)\stackrel{b}\map
[A,A]\map 0.
\end{equation}
For $k>1$, the relation between de Rham
complex and Hochschild homology
is more complicated.

One can check that the operations
$d,\L_\theta,$ and $ i_\theta$ on noncommutative
differential forms all    descend to the de Rham complex. 

Let $A$ be a smooth associative and
commutative algebra.  Then, by Hochschild-Kostant-Rosenberg theorem we
have a vector
space (but not necessarily {\em algebra}) isomorphism 
$\comO^\hdot(A)=\HH_\idot(A)$.
Obviously, the Hochschild differential is not the de Rham differential, since the de Rham differential increases degree while the Hochschild homology differential decreases degree.  Indeed, there is another differential (the Connes differential) which yields the de Rham differential (and cyclic homology).

Suppose now that $B\subset A$ is a subalgebra.  The \emph{relative de Rham complex} of $A$ is
$$
\DR(A/B)=\sr(\ncO^\hdot(A/B))=\ncO^\hdot(A/B)/[\ncO^\hdot(A/B),\ncO^\hdot(A/B)],
$$
where now $[-,-]$ is the graded commutator.
\end{defn}

Strictly speaking, to be consistent with this notation, one has to write
$\DR(A/\C)$ rather than $\DR(A)$ in the `absolute' case.

The de Rham complex is Morita
invariant
in the following sense.

\begin{prop}\label{morita_DR}
Let $A$ be an associative algebra.  Then for any $n\in\N$, there is a canonical isomorphism
$$
\DR(\Mat_nA/\Mat_n\C)\simeq \DR(A).
$$
\end{prop}

\begin{proof}
Let $A$ and $B$ be arbitrary associative algebras.  Then it is clear
that there is a canonical isomorphism $\ncO^\hdot\bigl((B\o A)/(1\otimes
B)\bigr)\simeq B\o \ncO^\hdot(A)$.  Also, it is clear that 
$\Mat_nA\simeq\Mat_n\C\o A$.  So, we  find that
$$
\ncO^\hdot(\Mat_nA/\Mat_n\C)\simeq\Mat_n\C\otimes\ncO^\hdot(A)\simeq\Mat_n(\ncO^\hdot(A)).
$$
But our comments regarding $\sr(\Mat_n\C)$ show that
$$
\sr(\Mat_n(\ncO^\hdot(A)))\simeq \sr(\ncO^\hdot(A))=\DR(A).
$$
Thus,
$\sr\bigl(\ncO^\hdot(\Mat_nA/\Mat_n\C)\bigr)=\DR(A),$ and we are done.
\end{proof}

\begin{rem}\label{tr_dr}
One can see by looking through the proof above that the map
$\DR(A) \too$
$ \DR(\Mat_nA/\Mat_n\C)$ that yields the
isomorphism of the
Proposition is induced by the algebra imbedding
$A\into \Mat_nA,\, a\mto \begin{pmatrix}a&0\\0&1_{n-1}\end{pmatrix}.$
The inverse isomorphism $\DR(\Mat_nA/\Mat_n\C)
\iso$
$ \DR(A)$ is induced by the `trace-map' sending
$x=\|x_{ij}\|\in \Mat_nA$ to $\Tr(x)=\sum x_{ii}$.
\end{rem}

Next, let $A=\bigoplus_{i\geq 0}\,A_i$ be a graded
$\k$-algebra. Write $\Gm$ for the multiplicative
group, viewed as an algebraic group over $\k$.

Giving  a $\Z$-grading
on an associative $\k$-algebra $A$  is the same thing as giving an
algebraic
$\Gm$-action on $A$ by algebra automorphisms.
Specifically, given a grading   $A=\bigoplus_{i\in\Z}\,A_i$ 
one defines a $\Gm$-action by 
the formula $\Gm\times A_i\ni t,a \mapsto t^i\cdot a,$
for any $i\in\Z$. Conversely, it is easy to see that
any $\Gm$-action on $A$ by algebra automorphisms
arises in this way from a certain $\Z$-grading on $A$.

Assume now that  $A=\bigoplus_{i\geq 0}\,A_i$
is graded by {\em nonnegative} integers.
Geometrically, this means that the corresponding
$\Gm$-action is a `contraction' of $A$ to the subalgebra $A_0$.

The result below says that de Rham cohomology
is `invariant under contraction'.

\begin{thm}[Poincar\'e lemma]\label{flat1} For a graded algebra
 $A=\bigoplus_{i\geq 0}\,A_i,$ 
the algebra imbedding $A_0\into A$ induces isomorphisms
$$H^j(\DR(A_0))\iso H^j(\DR(A)),
\quad\forall j\geq 0.$$
\end{thm}

\begin{proof} The assignment $A_i \ni a \mapsto i\cdot a,\,
i=0,1,\ldots,$
 gives a derivation of $A$,
called  the \emph{Euler derivation}. This derivation
may be thought of as an infinitesimal generator of
the $\Gm$-action on $A$ corresponding to the grading.

Associated with the Euler derivation,
one has  $\L_\Eu$,
the Lie derivative with respect to $\Eu$,  acting on $\DR(A)$. 
The action of  $\L_\Eu$ on $\DR(A)$ is diagonalizable
with nonnegative integral eigenvalues,
and we write $\DR(A)=\bigoplus_{m\geq 0}\,\DR(A)\langle m\rangle$
for the corresponding eigenspace direct sum decomposition.
It is clear that we have $\DR(A)\langle 0\rangle=\DR(A_0).$

The de Rham differential $d$ commutes with  $\L_\Eu$,
hence preserves the  direct sum decomposition above.
Further, the homotopy formula
$\L_\Eu=d\circ i_\Eu+i_\Eu\circ d$, 
shows that $i_\Eu$ is a chain homotopy between the  map
 $\L_\Eu$ and the zero map.  
Hence, 
 the complex $(\DR(A)\langle m\rangle,\,d)$
is acyclic for all $m$
 except $m=0$. This proves the Theorem.
\end{proof}
\subsection{The Quillen sequence.}
Recall that $
\DR^1(A)=\ncO^1(A)/[\ncO^1(A),A]$, see
\eqref{HH_DR}, and
define a map $b\colon \DR^1(A)\to A$ by $b\colon x\,dy\mapsto[x,y]$.  
We have shown  that this map
 is well-defined, see \eqref{quil_map}.

It is easy to check that $b\circ d=d\circ b=0$.  
Thus $d$ and $b$ yield a $2$-periodic complex
\beq{Quillen_comp}
\xymatrix{A\ar[r]^<>(.5)d&\DR^1(A)\ar[r]^<>(.5)b&A\ar[r]^<>(.5)d&\DR^1(A)\ar[r]&\cdots}.
\end{equation}
This complex gives an approximation to the cyclic homology of
the algebra $A$, see \cite{CQ2} 

Put
$\overline{\DR}^0(A)=A/([A,A]+\C)$, let $\pr: {\overline{A}}
\map  {\overline{A}}/[A,A]\simeq\overline{\DR}^0(A)$ be the natural
projection.

Consider the following sequence of maps, called the {\em Quillen sequence}:
\beq{Quillen_seq}
\xymatrix{0\ar[r]&\overline{\DR}^0(A)\ar[r]^<>(.5)d&
\DR^1(A)\ar[r]^<>(.5)b&{\overline{A}}\ar[r]^<>(.5){\pr}
&\overline{\DR}^0(A)\ar[r]&0},
\end{equation}
where 
 the maps $d$ and $b$ have been introduced above.
Since  $b\circ d=0$ and the image of $b$ is contained
in $[A,A],$ it is clear that the composite of any two consequtive
maps in the sequence is equal to zero, i.e., the sequence
is a {\em complex}.

Consider the case $A=T(V^*)$, where $V$ is a finite-dimensional
$\C$-vector space.  
\begin{lem}
If $A=T(V^*)$, then Quillen's complex \eqref{Quillen_seq}
is an exact sequence.
\end{lem}

\begin{proof} We already know that, for $A=T(V^*)$,
the de Rham complex $\DR^\hdot(A)$ is acyclic in positive degrees,
and $\DR^0(A)=\C$. It follows that the map
$\overline{\DR}^0(A)\to 
\DR^1(A)$ in \eqref{Quillen_seq} is injective.
Further, we have $\im(b)=[A,A]$, hence the complex
is exact at $\overline{A}$ and at $\overline{\DR}^0(A)$.
Thus, it remains only to check that the complex is exact at 
$\DR^1(A)$. 

To this end, it suffices to show that
the class
$$
\big[\overline{\DR}^0(A)\big]-\big[\DR^1(A)\big]+\big[{\overline{A}}\big]-
\big[\overline{\DR}^0(A)\big]
$$
vanishes in the Grothendieck group of graded spaces. 
 The two terms of $\big[\overline{\DR}^0(A)\big]$ cancel to leave
$$
-\big[\DR^1(A)\big]+\big[{\overline{A}}\big].
$$
But $\DR^1(A)\simeq A\otimes V^*=T(V^*)\otimes V^*$, while
${\overline{A}}\simeq T(V^*)\otimes V^*$, as well.  Hence, the
remaining two terms
 cancel as well and the sequence is exact.
\end{proof}

Let  $x_i$ and $x^i$ be dual bases in $V$ and $V^*$.
\begin{prop}
Let $A=T(V^*)$.

\vi
For any $f\in \DR^0(A)$, one has the identity
$
\sum_{i=1}^n\,\left[\pder{f}{x^i},x_i\right]=0,
$
in $[A,A]$.

\vii The space of closed forms in $\DR^2_{closed}\sset\DR^2(A)$
is canonically isomorphic to $[A,A]$.
\end{prop}

\noindent
{\em Proof.\;}
The proof uses the Quillen sequence,
$$
0\to\overline{\DR}^0(A)\to \DR^1(A)\to[A,A]\to0,
$$
which is exact for $A=T(V^*)$.  We also have the following sequence, which is exact by the Poincar\'e lemma,
$$
0\to\overline{\DR}^0(A)\to \DR^1(A)\to \DR^2_{closed}(A)\to0.
$$
Combining these two exact sequences, we have
$$
\xymatrix{0\ar[r]&\overline{\DR}^0(A)\ar@{=}[d]\ar[r]&\DR^1(A)\ar@{=}[d]\ar[r]&[A,A]\ar[r]&0\\
0\ar[r]&\overline{\DR}^0(A)\ar[r]&\DR^1(A)\ar[r]&\DR^2_{exact}\ar[r]&0}.
$$
But then we can use standard diagram chase arguments to construct a map $[A,A]\to \DR^2(A)$, and this must be an isomorphism.  This shows claim (2).  As for (1), observe that
$$
\sum\nolimits_{i=1}^n\,\left[\pder{f}{x^i},x_i\right]=bd(f)=0.\qquad\Box
$$

\subsection{The Karoubi Operator.}
Define the \emph{Karoubi operator} $\kappa: \ncO^\hdot(A)\to\ncO^\hdot(A)$ by 
 $$\kappa(\alpha\,da)=(-1)^{\deg\alpha}da\cdot\alpha,
\quad\text{and}\quad\kappa(\alpha)=\alpha\en\text{if}\en \alpha=a\in \ncO^0(A)=A.
$$

The Karoubi operator provides the following
the  relation between 
$b$ and $d$:
$$
db+bd=\id-\kappa.
$$
If we compare this to differential geometry and think of $b$ as $d^*$,
the adjoint of the de Rham differential with respect
to a euclidean structure on the space of
differential forms,
then $\kappa$ is playing the role of the Laplace operator.  

From a different viewpoint, observe that there is an obvious $n+1$-cycle action
on the bar complex,  namely 
given by the obvious cyclic action on $A^{\otimes(n+1)}$.  Since the
reduced bar complex 
is $A\otimes{\overline{A}}^{\otimes n}\simeq\ncO^n(A)$, 
there is no clear $(n+1)$-cycle action,
and the operation $\kappa$ ``approximates'' an $n$-cycle action.  

\begin{prop}
\pbox{
\vi $\kappa^{n+1}d=d$.\newline
\vii$\kappa^n=\id+b\kappa^nd$.\newline
\viii$\kappa^{n+1}=1-db$.}
\end{prop}

\pf
\vi This is trivial.
To prove \vii we calculate
\begin{align*}
\kappa^n(a_0da_1\ldots da_n)&=(da_1\ldots da_n)a_0\\
&=a_0da_1\ldots da_n +[da_1\ldots da_n\,,\,a_0]\\
&=a_0da_1\ldots da_n +(-1)^n b(da_1\ldots da_n da_0)\\
&=(\id+b\kappa^nd)(a_0da_1\ldots da_n).
\end{align*}
\viii Multiply (ii) on the right by $\kappa$ and observe 
that $\kappa$ commutes with $b$ and $d$.
\begin{align*}
\kappa^{n+1}&=\kappa+b\kappa^nd\kappa=\kappa+b\kappa^{n+1}d\\
&=\kappa+bd\qquad\text{by (1)}\\
&=\kappa+(\id-\kappa-db)=\id-db.\qquad\Box
\end{align*}

\begin{prop}
On $\ncO^n(A)$, we have $
(\kappa^n-1)(\kappa^{n+1}-1)=0.
$
\end{prop}

\begin{proof}
First, we have that $\kappa^n-1=b\kappa^nd$ by (ii).  By (iii),
 $\kappa^{n+1}-1=-db$.  So,
$$
(\kappa^n-1)(\kappa^{n+1}-1)=-b\kappa^nd^2b=0,
$$
since $d^2=0$.
\end{proof}

\subsection{Harmonic decomposition.}
Notice that $(n,n+1)=1$.  So, the polynomial $(t^n-1)\cdot$
$(t^{n+1}-1)$ has
only simple roots except for a double root at one.  The identity
$(\kappa^n-1)({\kappa^{n+1}-~1})=0$
implies that the action of $\kappa$ is 
{\em locally-finite}, and all of its eigenvalues have multiplicity 
one except for  the eigenvalue $1$, which has multiplicity 2.  So,
$$
\ncO^n(A)=\left[\Ker(\kappa-1)^2\right]\;\bigoplus\;
\left[\oplus_{\lambda\in\Spec(\kappa)-\{1\}}\,\Ker(\kappa-\lambda)\right].
$$
The space $\Ker(\kappa-1)^2$ is called the space of \emph{harmonic}
forms, denoted by $\Harm$.  The remaining summand is denoted by
$\Harm^\perp\sset \ncO^\hdot(A)$.  

\begin{prop}
$$
\Harm^\perp=d(\Harm^\perp)\oplus b(\Harm^\perp).
$$
\end{prop}

\begin{proof}
Observe that $db+bd=1-\kappa$ is invertible on $\Harm^\perp$.  Let $G$
be its inverse.  Then $G$ commutes with $b$ and $d$.  So,
$Gdb\colon\Harm^\perp\to d(\Harm^\perp)$
 and $Gbd\colon\Harm^\perp\to b(\Harm^\perp)$ are both projectors.  
This yields the direct sum decomposition.
\end{proof}

\subsection{Noncommutative polyvector fields.}
Recall the notation $\bar{A}=A/\C$; we have the reduced bar complex
$$
\cdots A\otimes\bar A^{\otimes(n+1)}\otimes A\to A\otimes\bar A^{\otimes n}\otimes A\to A\otimes\bar A^{\otimes(n-1)}\otimes A\cdots.
$$

We now consider the reduced cochain complex, that is, for each $n\in\N$ and each $A$-bimodule $M$ we set
$$
\overline C^n(A,M)=\Hom_{\bimod A}(A\otimes\bar A^{\otimes n}\otimes A,M),
$$
where $\overline C$ denotes the reduced cochains.  We let $\overline Z^n(A,M)$ denote the space of reduced $n$-cocycles.

\begin{prop}\label{Z_cocycle}
Suppose $f\in C^n(A,M)$.  Then $f$ is a cocycle if and only if the map $da_1\cdots da_n\mapsto f(1\otimes a_1\cdots a_n\otimes1)$ extends to an $A$-bimodule map $\ncO^n(A)\to M$.
\end{prop}

\begin{proof}
Suppose we are given some $f\in {\overline{C}}^n(A,M)$.  Write $\bar
f\in\Hom_\C(\overline A^{\otimes n},M)$ for the map $\bar
f(\omega)=f(1\otimes\omega\otimes1)$
 for all $\omega\in\overline A^{\otimes n}$.  If $f$ is a cocycle (i.e., if $df=0$), then
$$
f^*(da_1\cdots da_n)=\bar f(a_1\cdots a_n)
$$
extends $\bar f$ uniquely to an $A$-bimodule map
$f^*\colon\ncO^\hdot(A)\to M$. 
 The left $A$-linearity of $f^*$ is trivial, so the cocyclicity
 condition is then becomes equivalent
to right $A$-linearity.
\end{proof}

Write $\overline Z^p(A,M)$ for the $\C$-vector space of cocycles in
the reduced cochain complex $\overline C^p(A,M)$.

\begin{cor} There is a canonical
vector space isomorphism 
$$\overline Z^p(A,M)=\Hom_{\bimod A}(\ncO^p(A),M),$$
i.e., the  
 $A$-bimodule $\ncO^n(A)$ represents the functor 
$M\mto \overline Z^n(A,M)$.
\end{cor}

In particular, if $\theta\in\Der A\subset\Der(T(A^*))$, 
then it extends uniquely to yield a $1$-cocycle $\theta\in Z^1(A,A)$.
  So, every derivation yields a cocycle.
Thus, the above Corollary is a generalization to $p\geq 1$
of the interpretation of the space 
$\overline Z^1(A,M)$ as the space $\Der(A,M)=\Hom_{\bimod
A}(\ncO^1(A),M)$, of derivations from $A$ to $M$
(where the latter equality follows from the universal property of
$\ncO^1(A)$).

\begin{defn} For any $p\geq 1$, we set
$\ncT^p(A):=\Hom_{\bimod A}(\ncO^n(A),A)$,
and call elements of the graded
space $\ncT^\hdot(A):=\bigoplus_{p\geq 1}\,\ncT^p(A)$
\emph{noncommutative polyvector fields} on $A$.
\end{defn}

By definition, there is a natural pairing
$\ncT^p(A) \otimes \ncO^p(A) \too A$.

\begin{prop}\label{grst} For any associative algebra $A$, 
there is  a natural graded Lie super-algebra structure on
$\ncT^{\bullet-1}(A)$, i.e., a super-bracket
such that
$$
[\ncT^p(A)\,,\,\ncT^q(A)]\sset \ncT^{p+q-1}(A).
$$
\end{prop}
\begin{rem}
The degrees given above are precisely those for the usual Schouten bracket.
\eer

\begin{proof}[Proof of Proposition \ref{grst}.]
Inside the Hochschild cochain complex
$C^p(A,A)=$
$\Hom_\C(A^{\otimes
p},A)$ we have a subcomplex of
 reduced cochains $\overline C^p(A,A)=\Hom_\C(\overline A^{\otimes
p},A)$.  Thus, a cochain is a reduced cochain provided it
vanishes whenever  at least one entry is
a scalar.  It follows easily 
 the reduced cochain
complex is preserved by the Gerstenhaber bracket.
  It is also true, although not quite so transparent, that the
  Gerstenhaber 
bracket is compatible with the Hochschild differential.  
In particular, it preserves cocycles.  So, if we let $\overline Z^p
$ denote the reduced $p$-cocycles, we obtain a bracket
$$
\overline Z^p\times\overline Z^q\to\overline Z^{p+q-1}.
$$
But by Proposition \ref{Z_cocycle}
 we can extend a reduced $p$-cocycle $\omega$ to a bimodule map from
$\ncO^p(A)$ to $A$.  An element of $\Hom_{\bimod A}(\ncO^p(A),A)$
is called a \emph{$p$-vector field}.  
The Gerstenhaber bracket then yields a Lie super-algebra structure on $p$-vector fields.
\end{proof}
\begin{question} \vi Does the natural $\Th_\idot(A)$-action on
$\ncO^\hdot(A)$ descend to $\DR^\hdot(A)$~?

\vii Given $c\in\Th_p(A)$, when
does the map $i_c: \ncO^n(A)\to \ncO^{n-p}(A)$
descend to a well-defined map
$\DR^n(A)\to \DR^{n-p}(A)$~?
\end{question} 
\section{The Representation Functor}\label{rep_functor}
\subsection{}

It is believed that 
the noncommutative geometry of
an associative algebra $A$ is `approximated'
(in certain cases) by the
(commutative) geometry of the scheme
${\Rep^A_n}$  of $n$-dimensional representations of $A$.
Moreover, it is expected that this approximation
becomes `better' once the integer $n$ gets larger.

\begin{rem} There exist noncommutative associative algebras
$A$ that do not have {\em any} finite dimensional representations
at all. Thus the idea of looking at finite
 dimensional representations has obvious limitations.
\eer

Suppose $A$ is a finitely generated associative algebra and  $E$ is
a finite dimensional $\C$-vector space.  
Let ${\Rep^A_E}$ denote the affine (not necessarily reduced)
scheme of all $\C$-linear algebra maps $A\to\End_\C E$,
see below for a rigorous   definition of the scheme structure.

Given $a\in A$, for each representation $\rho\in\Rep^A_E$,
the element $\rho(a)$ is a $\C$-linear endomorphism of $E$.
The assignment $\rho \mapsto \rho(a)$ is an
$\End_\C E$-valued regular algebraic function on ${\Rep^A_E}$,
to be denoted $\wh a$. Equivalently, the function $\wh a$
may be viewed as an element of $\End_\C E\otimes\C[\Rep^A_E]$,
a tensor product of the finite dimensional simple algebra 
$\End_\C E$ with $\C[\Rep^A_E]$, the coordinate ring of
the scheme $\Rep^A_E.$

Let $GL(E)$ be the group of invertible
 linear transformations of $E$.  Then we have an action 
of $GL(E)$ on ${\Rep^A_E}$ by conjugation.  That is, if $\f\colon
A\to\End_\C E$ is a representation and $g\in GL(E)$, we define
$(g\cdot\f)(a):=g\f(a)g^{-1}$ for all $a\in A$.

The trivial bundle $E\times \Rep^A_E\to\Rep^A_E$ has a
natural
structure of $GL(E)$-equivariant vector bundle on $\Rep^A_E$
(with respect to the  {\em diagonal}
 $GL(E)$-action on $E\times \Rep^A_E$). We call this vector bundle
the {\em tautological}  vector bundle, to be denoted
$E_\Rep.$ The algebra $\End E_\Rep$, of vector bundle endomorphisms
of $E_\Rep$ is clearly identified with  $\End_\C E\otimes\C[\Rep^A_E]$.

Observe that, for any $a\in A$, the element
$\wh a\in \End_\C E\otimes\C[\Rep^A_E]$
is $GL(E)$-invariant with respect to
the simultaneous $GL(E)$-action on $\End_\C E$ (by conjugation)
and on $\C[\Rep^A_E]$. This way, the assignment $a\mapsto \wh a$
gives a canonical algebra map
\beq{rep_map}
{\mathsf{rep}}: A\map \bigl(\End_\C E\otimes\C[\Rep^A_E]\bigr)^{GL(E)}=
\bigl(\End E_\Rep\bigr)^{GL(E)}.
\end{equation}

To make the scheme structure on ${\Rep^A_E}$ explicit,
we first
 consider  the special case $A=\C\langle x_1,\ldots,x_r\rangle$,
a free algebra on $r$ generators.

An algebra homomorphism $\rho: \C\langle x_1,\ldots,x_r\rangle\map
\End_\C(E)$
is specified by an arbitrary choice of
an $r$-tuple of endomorphisms $X_1:=\rho(x_1),\ldots, X_n:=\rho(x_r)\in
\End_\C(E)$.
Thus, we have an isomorphism of algebraic varieties
\begin{equation}\label{rep_free}
\Rep^{\C\langle x_1,\ldots,x_r\rangle}_E\iso 
\underset{r\text{ factors}}{\underbrace{\End_\C(E)\times\ldots\End_\C(E)}},
\;\rho\mapsto \bigl(\rho(x_1),\ldots,\rho(x_r)\bigr).
\end{equation}

Now, given any finitely generated associative algebra $A$,
choose a finite set $\{x_1,\ldots,x_r\}$ of algebra generators for $A$.
Then, we have
$$
A=\C\langle x_1,\ldots,x_n\rangle/I
$$
for some two-sided ideal $I\sset \C\langle x_1,\ldots,x_n\rangle$.  
It is clear that,
giving an algebra map $A\to \End_\C E$
is the same thing as giving an algebra map $\C\langle
x_1,\ldots,x_n\rangle
\to\End_\C E$ that vanishes on the ideal $I$.
Put another way, the projection
$\C\langle x_1,\ldots,x_n\rangle\onto A$ induces
a closed imbedding $\Rep^A_E\into \Rep^{\C\langle
x_1,\ldots,x_r\rangle}_E$,
and we have
\beq{rep_scheme}
\Rep^A_E= \{f\in  
\Rep^{\C\langle x_1,\ldots,x_r\rangle}_E\cong
\End_\C(E)\times\ldots\End_\C(E)\enspace\big|\enspace f(I)=0\}.
\end{equation}
The RHS of this formula is clearly an algebraic subset 
of a finite dimensional vector space
defined by algebraic equations,
that is, an affine subscheme.
This puts  an affine scheme structure on the LHS of the equation,
that will be shown below to be independent of the choice of the
generators
$x_1,\ldots,x_r$ of the algebra~$A$.

In order to make
the construction of the scheme structure on $\Rep^A_E$
manifestly independent of the choice of generators,
it is coventient to use the {\em functor of points}.

Recall that for any category $\scr C$, an object
$S\in \scr C$ gives rise to a functor
$$S(-):\ \scr C\too \mathsf{Sets}\,,\quad
X \mto S(X):=\Hom_{\scr C}(S,X).
$$
Further, the Yoneda lemma says that, for any $S,S'\in  \scr C$,
every isomorphism of functors $S(-) \iso S'(-)$ is necessarily
induced by an isomorphism of objects $S'\iso S$.
In particular, the functor $S(-)$ determines the object
$S$ uniquely, up to a unique isomorphism.

Recall next that the category  of
affine  schemes of finite type over $\C$ is equivalent, via the `coordinate
ring'
functor, to the category of   finitely generated commutative
$\C$-algebras. Thus, any affine scheme
$S$ is completely determined by the corresponding functor 
\begin{align*}
S(-):\;\;&
\textsf{fin. gen. Commutative Alg.}\too \mathsf{Sets},\\
&B\mto S(B):=\Hom_{\sf{alg}}(\C[S],B)\simeq \Hom_{\textsf{Schemes}}(\Spec B,S).
\end{align*}
The set $S(B)$ is usually referred to as the set of $B$-points
of $S$; for $B=\C[X]$, it is just the set of algebraic maps
$X\to S$.

Now, fix  a finitely generated associative (not necessarily
commutative) algebra $A$,
and an integer $n\geq 1$. Given  a  finitely generated commutative
algebra $B$, write $\Mat_nB$ for the associative algebra of $n\times n$-matrices
with entries in $B$.
We define a functor on the category of
finitely generated {\em commutative} algebras as follows:
\begin{align}\label{rep_points}
\Rep_n^A(-):\;\;&\textsf{fin. gen. Commutative Alg.} \too \mathsf{Sets},\nonumber\\
&B \mto \Rep_n^A(B):=\Hom_{\sf{alg}}(A,\Mat_nB).
\end{align}
For $B=\C[X]$, one may think of the set $\Hom_{\sf{alg}}(A,\Mat_n(\C[X])$
as the set of families $\{\rho_x: A\to \Mat_n\C\}_{x\in X}$
of $n$-dimensional representations of the algebra $A$
parametrized by points of the scheme $X$.
Thus, if  $E=\C^n$ and 
$X$ is a point, we get back to the original definition of
$\Rep^A_n:=\Rep^A_E.$
Observe also that if  $E=\C^n$, then we have
$$\End_\C E\otimes \C[\Rep^A_n]=\Mat_n\C\otimes \C[\Rep^A_n]=
\Mat_n\bigl(\C[\Rep^A_n]\bigr).
$$

The above discussion can be summed up in the following
result that shows at the same time that
the description of the scheme $\Rep^A_E$ given in \eqref{rep_scheme}
is independent
of the presentation $A=\C\langle x_1,\ldots,x_r\rangle/I$.
We restrict ourselves to the case   $E=\C^n$,
in wich case we have
\begin{prop}\label{Prop_rep} \vi For any finitely generated
associative algebra $A$, the corresponding functor
\eqref{rep_points} is {\sf representable} by an affine
scheme, which is called $\Rep^A_n$. 

\vii The coordinate ring $R:=\C[\Rep^A_n]$ is a finitely
generated commutative algebra equipped with a
canonical algebra map ${\sf{rep}}: A \to \Mat_nR,\,
a\mapsto \wh a,$ see \eqref{rep_map},
such that the following universal property holds:

\npb{Given a finitely generated commutative algebra $B$ and
an algebra map $\rho: A \to \Mat_nB$, there exists 
a unique algebra homomorphism $\wh\rho: R\to B$ making the
following diagram commute}
$$
\xymatrix{
A\ar[rr]^<>(0.5){{\sf{rep}}}\ar[drr]_{\rho}&&\Mat_nR\ar[d]^{\Mat_n(\wh\rho)}\\
&&\Mat_nB.
}
$$
\end{prop}
Here and below,
given a linear map $f: V \to U$ of vector spaces,
 we write $\Mat_n(f): \Mat_n(V)\to\Mat_n(U)$
for the map between the spaces of $V$-valued and
$U$-valued matrices
whose entries are all equal to $f$
(observe also that $\Mat_n(V)\simeq\Mat_n\C\otimes V$).

\begin{proof}[Proof of Proposition] To prove the representability of the functor 
one has to choose a set of algebra generators for $A$,
and write $A=\C\langle x_1,\ldots,x_r\rangle/I$.
Then, the set 
$\Hom_{\mathsf{alg}}(A,\Mat_nB)$ may be identified, as in 
\eqref{rep_scheme}, with $B$-points of a  subscheme
in $\Mat_n\C\times\ldots\times\Mat_n\C$.
\end{proof}
\begin{rem} Note that we have defined the scheme
$\Rep^A_n$ through its functor of points, {\em without}
describing the coordinate ring $\C[\Rep^A_n]$. The latter ring
does not have a simple description: below, we will
produce a finite set of generators of the ring
 $\C[\Rep^A_n]$, but determining all the relations among these
generators is a formidable task.
\eer
\begin{cor}[Functoriality] Any algebra map $f: A \to {A'},\,
f(1_A)=1_{A'}$
induces an algebra map $\wh{f}: \k[\Rep^A_n]\to \k[\Rep^{A'}_n],$
hence, a morphism of algebraic varieties
$\wh{f}^*: \Rep^{A'}_n\to\Rep^A_n.$
\end{cor}
\begin{proof}
 Set $B:=\k[\Rep^{A'}_n]$, and apply Proposition \ref{Prop_rep} to 
to the map
$$\rho\colon A\stackrel{f}\map A'\stackrel{{\sf{rep}}}\too 
\Mat_n\bigl(\k[\Rep^{A'}_n]\bigr).
$$
The universal property from the Proposition
implies the existence of an algebra
map $\wh\rho:  \k[\Rep^A_n]\to \k[\Rep^{A'}_n]$ that makes
the following diagram  commute:
$$
\xymatrix{
A\ar[rr]^<>(0.5){{\sf{rep}}_A}\ar[d]^<>(0.5){f}\ar[drr]_{\rho}&&
\Mat_n\bigl(\k[\Rep^A_n]\bigr)
\ar[d]^{\Mat_n(\wh\rho)}\\
A'\ar[rr]^<>(0.5){{\sf{rep}}_{A'}}&&\Mat_n\bigl(\k[\Rep^{A'}_n]\bigr)
}
$$
Thus, we may put $\wh{f}:=\wh\rho$.
\end{proof}

\subsection{Traces.} Let  $A$ be a $\C$-algebra  and
$M$ an $A$-module which is {\em finite dimensional} over $\C$.
Then, the action in $M$ of an element $a\in A$
gives a $\C$-linear map $a: M\to M$. We write $\tr_M(a)$
for the trace of this map.

It is clear that if $N$ is an $A$-submodule in $M$,
then we have
$$\tr_M(a)=\tr_{M/N}(a)+\tr_N(a)\quad\text{\bf (additivity of the trace)}.
$$

Any finite-dimensional $A$-module clearly has finite length,
hence has a finite Jordan-H\"older series
$M=M^0\supset M^1\supset\ldots\supset M^k\supset M^{k+1}=0$.
The isomorphism classes of simple composition factors $M^i/M^{i+1}$ in this
series are defined uniquely, up to permutation.
Hence, the $A$-module $\sss M:=\oplus M^i/M^{i+1}$
is independent of the choice of  Jordan-H\"older series.
By construction, $\sss M$ is a semisimple
$A$-module (i.e., a direct sum of simple $A$-modules),
called the {\em semi-simplification} of $M$.

The additivity property of the trace implies that 
$\tr_M(a)=\tr_N(a)\,,\,\forall a\in A$,
whenever $M$ and $N$ have the same semi-simplification.
Conversely, one has
\begin{thm} Let $A$ be a finitely generated $\C$-algebra
and $M$ and $N$ be finite-dimensional $A$-modules
such that $\tr_M(a)=\tr_N(a)$ for any  $a\in A$.
Then, $\sss M\simeq \sss N$.
\end{thm}

\begin{proof} Both the assumptions and the conclusion 
of the Theorem are unaffected
by replacing $M$ by $\sss M$, and $N$ by $\sss N$.
Therefore, we may reformulate the Theorem as follows:
if $M$ and $N$ are semisimple $A$-modules
such that $\tr_M(a)=\tr_N(a)$ for any  $a\in A$,
then $M\simeq N$.

Thus, from now on we assume that  $M$ and $N$ are  finite-dimensional semisimple
$A$-modules. 

Let ${\mathtt{Ann}}(M\oplus N)\sset A$ be the
annihilator of $M\oplus N$, that is the set of all elements
$a\in A$ that act by zero on $ M\oplus N$.
This is a two-sided ideal in $A$, and
the action of $A$ gives an  algebra
imbedding ${A'}:=A/{\mathtt{Ann}}(M\oplus N)\into
\End_\C(M\oplus N).$
We see that ${A'}$ is a finite dimensional algebra
and that $M\oplus N$ is an ${A'}$-module, which is
moreover semisimple, since it is semisimple as an $A$-module.

Let $\Rad {A'}$ be the radical of ${A'}$,
the intersection of the annihilators of all simple
${A'}$-modules.
Thus,  $\Rad {A'}$ annihilates any
semisimple ${A'}$-module,
in particular,  annihilates  $M\oplus N$.
Furthermore, the structure theory of finite
dimensional algebras over $\C$ says that
$$\overline{{A'}}:={A'}/\Rad {A'}\cong
\oplus_i\, \End_\C(E_i)\quad\text{where}\quad E_i=\C^{r_i}.$$
Thus, each $\End_\C(E_i)\cong \Mat_{r_i}(\C),$ is a simple 
 matrix algebra, and any simple $\overline{{A'}}$-module
is isomorphic to some $E_i$, viewed as an
$\overline{{A'}}$-module
via the projection $\overline{{A'}}\onto \End_\C(E_i).$

We conclude that our semisimple $\overline{{A'}}$-modules
$M$ and $N$ have the form
$M=\oplus E_i$ and $N=\oplus E_j$, respectively.
With this understood, the assumption of the Theorem reads:
$$\tr_M(a)=\tr_N(a)\quad\text{for any}\quad a\in \bplus_i\, \End_\C(E_i).$$
This clearly implies that the direct sums
$M=\oplus E_i$ and  $N=\oplus E_j$ involve the same
summands with the same multiplicities.
Hence $M\simeq N$, and we are done.
\end{proof}

We now fix a finite dimensional vector space $E$ and consider
the scheme $\Rep^A_E$. For any $\rho\in \Rep^A_E$,
let $\OO(\rho)$ denote the $GL(E)$-orbit of the point $\rho$.
The orbit  $\OO(\rho)$ corresponds to the isomorphism class of
the $A$-module $E_\rho$, and we write $\OO(\sss(\rho))$
for the orbit corresponding to the semi-simplification of $E_\rho$.

The following result is a standard application of Geometric Invariant
Theory, see \cite{GIT}.
\begin{thm}\label{git}\vi The orbit $\OO(\rho)$ is closed in
$\Rep^A_E$ if and only if $E_\rho$ is a semisimple $A$-module.

\vii For each $\rho \in \Rep^A_E$, the orbit 
$\OO(\sss(\rho))$
is the  unique closed
$GL(E)$-orbit contained in
$\overline{\OO(\rho)}$, the closure of  $\OO(\rho)$.
\end{thm}

Taking the trace of the $\End_\C E$-valued function $\wh a$,
one  obtains an element 
$\tr(\wh a)\in\C[{\Rep^A_E}]^{GL(E)}$, the $GL(E)$-invariant regular
functions on ${\Rep^A_E}$. 
Equivalently, $\tr(\wh a)$ is the image of $\wh a$ under
the linear map $\tr\otimes\id: \End_\C E\otimes\C[\Rep^A_E]\map
\C\otimes\C[\Rep^A_E]=\C[\Rep^A_E].$

\begin{thm}\label{git2}\vi The functions $\{\tr(\wh a)\,,\,a\in A\}$
generate $\C[\Rep^A_E]^{GL(E)}$ as an algebra.

\vii The  natural map $\Rep^A_E\to \Spec\C[\Rep^A_E]^{GL(E)}$
induces a bijection between the isomorphism classes
of semisimple $A$-modules of the form $E_\rho, \, \rho\in \Rep^A_E$,
and
 maximal ideals of the algebra $\C[\Rep^A_E]^{GL(E)}$.
\end{thm}

\begin{proof} Part (i) can be deduced from
a theorem of LeBruyn-Procesi \cite{LP}.
Part (ii)  follows from the general result,
cf. \cite{GIT}
saying that the elements of the algebra $\C[\Rep^A_E]^{GL(E)}$ separate
closed $GL(E)$-orbits in $\Rep^A_E$.
\end{proof}

\subsection{Noncommutative Rep-scheme (following \cite{LBW}).}
 We have seen that the functor $B\mto
\Hom_{\mathsf{alg}}(A,\Mat_nB)$
is representable, as a functor on the category of
finitely generated {\em commutative} algebras.
L. Le Bruyn and G. Van de Weyer have observed, 
see \cite{LBW}, that this functor is in effect representable
on the much larger category of all finitely generated {\em
non-commutative} associative
algebras. The  non-commutative
algebra, $\sqrt[n]A,$ that represents the functor
$\Hom_{\mathsf{alg}}(A,\Mat_n(-))$
maps surjectively onto $\C[\Rep^A_n]$, the coordinate
ring of the scheme $\Rep^A_n$ constructed in Proposition \ref{Prop_rep}.
Thus,  $\sqrt[n]A$ may be thought of as the
coordinate ring of a non-commutative space,
the canonical  `non-commutative thickening'
of the commutative\linebreak
 scheme $\Rep^A_n$.

The algebra $\sqrt[n]A$ is constructed as
follows, see \cite{LBW} for more details.

First, one forms the free-product algebra $A*\Mat_n\C$,
which contains the matrix algebra $\Mat_n\C$ as a subalgebra.
Let
\begin{align}\label{sqrtA}
\sqrt[n]A&:= \big[A*\Mat_n\C\big]^{\Mat_n\C}\\
&\,=\{r\in A*\Mat_n\C\enspace|\enspace x\cdot r=r\cdot x,
\enspace\forall x\in \Mat_n\C\}.\nonumber
\end{align}
be the centraliser of the subalgebra $\Mat_n\C\sset A*\Mat_n\C.$

\begin{examp}[\cite{LBW}] One can show that
$\sqrt[n]{\C\langle x_1,\ldots,x_r\rangle}\simeq
\C\langle x_{_{11,1}},\ldots,$
$x_{_{rr,n}}\rangle,$
is a free algebra on $n\cdot r^2$ generators.\qed
\eex

Following  \cite{LBW}, we are going to construct, for any (not necessarily
commutative), associative algebra $B$, a
canonical bijection
\beq{lbw}
\Hom_{\mathsf{alg}}(A, \Mat_nB) \stackrel{\sim}\longleftrightarrow 
\Hom_{\mathsf{alg}}(\sqrt[n]A, B)
\end{equation}
The bijection would clearly  yield the representability claim
mentioned above.

The construction of  this bijection is based on the
elementary Lemma below. 
To formulate the Lemma,
fix an associative algebra $R$, 
and  an algebra map
$\imath: \Mat_n\C\to R$, such that $\imath(1)=1$.
Let $R^{\Mat_n}$ denote the centralizer of $\imath\bigl(\Mat_n\C\bigr)$ in $R$.

\begin{lem}\label{free_lem} In the above setting,
the  map
$\Mat_n\C\otimes R^{\Mat_n}\map R\,,\,
m\otimes r\mto \imath(m)\cdot r,$ is an algebra isomorphism.
Moreover, if $R$ is a finitely generated algebra, then
so is the algebra $R^{\Mat_n}$. \qed
\end{lem}

\begin{proof}[Proof (by D. Boyarchenko).]
The map is clearly an algebra homomorphism.
To prove  that
this map is bijective, we form the algebra
$S:=\Mat_n\C\otimes_{_\C}\Mat_n\C^\op$.
This is clearly a simple algebra of dimension $n^4$, which is moreover isomorphic
to $\Mat_{n^2}(\C)$. Thus, by the well-known results
about modules over finite dimensional simple $\C$-algebras,
the algebra $S$ has a unique, up to isomorphism,
simple $S$-module $L$ which has dimension $n^2$.
Furthermore, any faithful $S$-module $M$ (i.e., such that 
$1_S\cdot m=0\;\Longrightarrow\; m=0$, for any $m\in M$)
is isomorphic to a (possibly
infinite) direct sum of copies of $L$,
that is, has the form $L\otimes V$, for some
vector space $V$.

The algebra $\Mat_n\C$ has an obvious
$\Mat_n\C\otimes_{_\C}\Mat_n\C^\op$-module
structure via left and right multiplication. Also, $\dim
\Mat_n\C=n^2$.
We conclude that $L=\Mat_n\C$, as a module over $S=
\Mat_n\C\otimes_{_\C}\Mat_n\C^\op$.

Now, the map  $\iota:\Mat_n\C\to R$ makes $R$ a faithful
(since $\iota(1_{\Mat_n})=1_R$) $\Mat_n\C$-bimodule, hence, a left
$S$-module. Thus, there is an $S$-module isomorphism
$R\simeq S\o V$, where $V$ is some $\C$-vector space.
Under this isomorphism, $R^{\Mat_n}$ clearly corresponds to $(\C\cdot
1_S)\o V$ 
because the center of $S$ is $\C\cdot 1_S$.
 So the statement of the lemma becomes obvious.
\end{proof}

Applying the Lemma to the tautological imbedding $\imath: \Mat_n\C\hookrightarrow
R={A*\Mat_n\C}$, and using the notation \eqref{sqrtA}, 
we obtain a canonical algebra isomorphism
$$\varphi:\ \Mat_n\C\otimes \sqrt[n]A=\Mat_n\C\otimes R^{\Mat_n}\iso
A*\Mat_n\C.
$$

Now, for any algebra $B$, an algebra morphism  $f: \sqrt[n]A\to B$,
induces an algebra morphism $\tilde f: \ A \to \Mat_nB$ given by the
composition:
$$A \into A*\Mat_n\C\stackrel{\varphi^{-1}}\iso
\Mat_n\C\otimes \sqrt[n]A\stackrel{\id\otimes f}\too\Mat_n\C\otimes
B=\Mat_nB.$$
Conversely,
given an algebra map $g: A \to \Mat_nB$, the universal property of
free products yields an algebra map
$g* \jmath :\ A* \Mat_n\C \to \Mat_nB$, where $\jmath :
\Mat_n\C \into \Mat_nB$ is the natural imbedding.
Observe that the centralizer of the image of $\jmath$
is formed by ``scalar matrices'', i.e., we have
$\bigl(\Mat_nB\bigr)^{\Mat_n\C}=B$.
Therefore, the subalgebra $g* \jmath(\sqrt[n]A)$,
the image of the restriction 
of $g* \jmath :\ A* \Mat_n\C \to \Mat_nB$ to the 
centralizer of $\Mat_n\C$ in $A* \Mat_n\C$,
is contained in the subalgebra $\bigl(\Mat_nB\bigr)^{\Mat_n\C}=B$.
This way, one obtains a map $\sqrt[n]g:\ \sqrt[n]A\to B$.

It is straightforward to check that the assignments
$f\mapsto \tilde f$ and $g\mapsto \sqrt[n]g,$
are mutually inverse bijections. This completes the construction of
the bijection in \eqref{lbw}.

\subsection{The Rep-functor on vector fields.}
Below, we are going to relate various `non-commutative' constructions
on $A$
with their  commutative counterparts for $\Rep^A_E$.

First, we claim that any derivation of $A$ gives rise to a $GL(E)$-invariant
vector field
on the scheme $\Rep^A_E$.

\begin{examp}
For instance, let $\theta=\ad a\in \Inn(A)\sset \Der(A)$ be an inner derivation.
The corresponding vector field $\wh\theta$ on $\Rep^A_E$ is then going
to be
tangent to the orbits of the $GL(E)$-action on $\Rep^A_E$,
and it is constructed
as follows. Write $\g g=\Lie^{\,}{GL(E)}=\End_\C E$ for the Lie algebra
of the group  $GL(E)$. To $a\in A$, we have associated 
the element $\wh a\in \End_\C E\otimes \C[\Rep^A_E]$,
which we now view as a $\g g$-valued regular function on
$\Rep^A_E$. For point $\rho\in{\Rep^A_E}$,
let $\wh a(\rho)\in\g g$ be the value of the function
$\wh a$ at $\rho$. Further, the infinitesimal 
$\g g$-action on $\Rep^A_E$ (that is, the differential
of the action map $GL(E)\times \Rep^A_E\to\Rep^A_E$
at the unit element of $GL(E)$)
associates to any $x\in\g g$ and $\rho\in \Rep^A_E$
a vector $x_\rho$ in $T_\rho(\Rep^A_E),$ the tangent space to $\Rep^A_E$ at $\rho$.
In particular, we have the vector
$\wh a(\rho)_\rho\in T_\rho(\Rep^A_E).$
The assignment $\rho\mto\wh a(\rho)_\rho$ gives the required vector field, $\wh\theta$,
on $\Rep^A_E$.
\eex

To study the general case, we first
need to understand the tangent space of ${\Rep^A_E}$ at
some point
$\rho\in{\Rep^A_E}$.  So, $\rho\colon A\to\End_\C E$ is an algebra map,
hence it induces an $A$-module structure on $E$.  We denote this
$A$-module by $E_\rho$.  
Then the tangent space $T_\rho{\Rep^A_E}$ is given by all
$\C$-linear 
 maps $\f\colon A\to\End_\C(E)$ 
such
 that $(\rho+\e\f)\colon A\to\End_{\C[\e]}\bigl(\C[\e]/(\e^2)\o 
E\bigr)$ 
is an algebra map.  Expanding both sides of $(\rho+\e\f)(aa')=(\rho+\e\f)
(a)\cdot(\rho+\e\f)(a')$, we have
$$
\rho(aa')+\e\f(aa')=\rho(a)\rho(a')+\e\f(a)\f(a')+\rho(a)\e\f(a').
$$
Since $\rho$ is an algebra map, $\rho(aa')=\rho(a)\rho(a')$.  Cancelling them and dividing by the common $\e$ factor,  we find that
$$
\f(aa')=\f(a)\rho(a')+\rho(a)\f(a').
$$
If we regard $\End_\C(E_\rho)$ as an $A$-bimodule in the obvious fashion, this equation implies that $\f\in\Der(A,\End_\C(E_\rho))$, i.e., $T_\rho{\Rep^A_E}\simeq\Der(A,\End_\C(E_\rho))$.

Now, suppose we have a derivation $\theta\in\Der(A)$.  We wish to generate a vector field $\wh\theta$ on ${\Rep^A_E}$.  Indeed, for each $\rho\in{\Rep^A_E}$, we set
$$
\wh\theta_\rho(a)=\rho(\theta(a))=\widehat{\theta(a)}(\rho).
$$
It is not hard to check that this is a derivation, hence
$\wh\theta_\rho\in\Der(A,\End_\C(E_\rho))$
$=T_\rho{\Rep^A_E}$. 

\begin{examp}\label{der_free}
By Lemma \ref{der_extend}, a derivation 
$\theta$  of the free algebra $A=\C\langle x_1,\ldots,x_r\rangle$
is determined  by sending each generator $x_i\,,\,i=1,\ldots,r,$ to
an arbitrarily chosen element
$\theta(x_i)=f_i(x_1,\ldots,x_r)\in A,$
which we regard as a non-commutative polynomial in the variables $x_1,\ldots,x_r$.
Given such a  derivation 
$\theta$ we  describe the corresponding
 vector field $\wh\theta$ on $\Rep^A_E$ as follows.

Fix $\rho\in\Rep^A_E$ and
let  $X_1:=\rho(x_1),\ldots, X_n:=\rho(x_r)\in
\End_\C(E)$ be as above.
 We have  $T_\rho\Rep^A_E=\Der(A,\End{E}_\rho)$, where ${E}_\rho$
denotes the $A$-module ${E}$ endowed with the $A$-module structure given
by $\rho$.  
Then composing $\theta$ with $\rho$, we obtain a derivation 
$\rho\circ\theta\colon A\to\End{E}_\rho$. 
Thus, writing  $\wh\theta|_\rho$ for the value of the
vector field $\wh\theta$ at the point $\rho$, 
we have $\wh\theta|_\rho=\rho\circ\theta$.  
Hence, for each $i=1,\ldots,r,$
we find
$\wh\theta|_\rho\colon x_i\mapsto f_i(X_1,\ldots,X_r)$.  
Thus, using the identification \eqref{rep_free} we get
$$
\wh\theta|_{_{(X_1,\ldots,X_r)}}=\bigl(f_1(X_1,\ldots,X_r),\ldots,
f_r(X_1,\ldots,X_r)\bigr),
$$
where the $r$-tuple on the RHS is viewed as a vector
in  $\End_\C(E)\times\ldots\End_\C(E)$.
Now, a vector field $\xi$ on any vector space $R$ is nothing
but a map $R\to R$ that sends each vector $r\in R$ to the 
vector $\xi_r$, the value of $\xi$ at $r$.
With this understood, we see that the vector
field $
\wh\theta$ is  nothing
but the following self-map of  $\End_\C(E)\times\ldots\End_\C(E)$
\begin{equation}\label{selfmap}
\wh{f}:\ \bigl(X_1,\ldots,X_r\bigr) \mto \bigl(f_1(X_1,\ldots,X_r),\ldots,
f_r(X_1,\ldots,X_r)\bigr).
\end{equation}
\eex

\begin{prop}\label{vect_field_map} For any finitely generated associative
algebra $A$, the assignment $\theta\mapsto \wh\theta$ gives a
Lie algebra homomorphism $\Der(A)\to \calT({\Rep^A_E}).$
\end{prop}

\noindent 
{\em Proof (by D. Boyarchenko).\;}
Let $\theta\in \Der(A),$ be a derivation. Set $R:=\C[\Rep^A_n]$,
and recall the notation of Proposition \ref{Prop_rep}.
We have the following characterization of the vector field  $\wh\theta$:

\begin{claim}\label{mitya_claim}
\vi There exists a unique derivation $\widetilde\theta\in \Der(R)$ that makes the
following square commute
$$
\xymatrix{
A\ar[rr]^<>(0.5){\theta}\ar[d]_{{\sf{rep}}}&&A\ar[d]_{{\sf{rep}}}\\
\Mat_nR\ar[rr]^<>(0.5){\Mat_n(\widetilde\theta)}&&\Mat_nR.
}
$$

\vii For any point $\rho\in \Rep^A_n$, the value at $\rho$ of
the vector field on $\Rep^A_n$ induced by the derivation
$\widetilde\theta$
from \vi, is equal to the vector $\wh\theta|_\rho$ constructed above.
Thus, $\wh\theta=\widetilde\theta$ is indeed a regular algebraic
vector field.
\end{claim}

Assuming the Claim, we complete the proof of the Proposition as follows.
Let $\theta,\delta\in\Der(A)$, and $\widetilde\theta,
\widetilde\delta\in \Der(R)$, be the corresponding  derivations of $R$
arising via  Claim
\ref{mitya_claim}.
Then, clearly $[\widetilde\theta,
\widetilde\delta]$ is again a derivation of $R$
and, moreover, the diagrams of part (i) of the Claim for
$\theta$ and $\delta$, respectively, yield
commutativity of the diagram
$$
\xymatrix{
A\ar[rrr]^<>(0.5){[\theta,\delta]}\ar[d]_{{\sf{rep}}}&&&A\ar[d]_{{\sf{rep}}}\\
\Mat_nR\ar[rrr]^<>(0.5){\Mat_n([\widetilde\theta,\widetilde\delta])}&&&\Mat_nR.
}
$$
Now, the uniqueness statement in  Claim
\ref{mitya_claim}(i) combined with (ii) implies
that $\wh{[\theta,\delta]}=[\wh\theta,\wh\delta]$, and the Proposition
is proved.
\smallskip

\begin{proof}[Proof of Claim.] Observe that for $A$ viewed as an
$A$-bimodule,
the square-zero 
construction  produces an algebra $A\sharp A$
which is nothing but $A_\eps:=
A\otimes\beps$. Further, let  $A\into A_\eps$ be
the tautological algebra imbedding.
A linear map $\theta: A\to A$ is a  derivation of $A$ if and only if
the assignment $\theta_\eps: a\mto a+\eps\cdot \theta(a),$
gives an algebra homomorphism $A\to A_\eps$.

Recall the canonical algebra map $\mathsf{rep}: A \to \Mat_nR$.
Tensoring with $\beps$,
one gets a homomorphism $\mathsf{rep}_\eps:A_\eps\to
\bigl(\Mat_nR\bigr)_\eps.$
But we have
$$
\bigl(\Mat_nR\bigr)_\eps=
\Mat_nR\otimes\beps=\Mat_n\bigl(R\otimes\beps\bigr)=
\Mat_n(R_\eps).$$

Now, fix $\theta\in\Der(A)$ and let
$\theta_\eps: A\to A_\eps$ be the corresponding homomorphism.
Composing  $\theta_\eps$ with the 
morphism $\mathsf{rep}_\eps$, we obtain the top row of the following diagram
$$
\xymatrix{
A\ar@{^{(}->}[r]\ar[dr]_{\mathsf{rep}}&A_\eps
\ar[rr]^<>(0.5){\mathsf{rep}_\eps}&&\bigl(\Mat_nR\bigr)_\eps\ar@{=}[d]\\
&\Mat_nR\ar[rr]^<>(0.5){\Mat_n(\phi)}&&\Mat_n(R_\eps).
}
$$
The universal property of $R$ explained in Proposition \ref{Prop_rep}
guarantees the existence and uniqueness of an
algebra map $\phi: R\to R_\eps$
that makes the above diagram commute.
By the discussion at the beginning of the proof, the map
$\phi$ thus defined gives a derivation $\widetilde\theta: R\to R$
such that $\phi=(\widetilde\theta)_\eps.$
This proves the Claim.
\end{proof}

Let $\calD_1^\spad(E_\Rep)$ denote the space of
first order differential operators on $E_\Rep$ with scalar
principal symbol,
cf. \S\ref{atiyah_alg}. This is a Lie algebra, and
the group $GL(E)$ acts naturally on $\calD_1^\spad(E_\Rep)$
by Lie algebra automorphisms.
We consider the Lie subalgebra $\calD_1^\spad(E_\Rep)^{GL(E)}\sset
\calD_1^\spad(E_\Rep),$ of $GL(E)$-invariant elements.

\begin{rem} Since the vector bundle $E_\Rep$ is canonically
and $GL(E)$-equivariantly
trivialized, the symbol map $\sigma: \calD_1^\spad(E_\Rep)\map
\calT(\Rep^A_E)
$ has a natural  $GL(E)$-equivariant splitting $\calT(\Rep^A_E)
\map\calD_1^\spad(E_\Rep).$
This splitting allows to lift the
the Lie algebra map $\theta\mapsto \wh\theta$ of Proposition
\ref{vect_field_map} to a Lie algebra map indicated by the dotted arrow
in the following diagram
$$
\xymatrix{
&&\Der(A) \ar@{.>}[d]^<>(0.5){\Xi}\ar[dr]^<>(0.5){\text{Prop. \ref{vect_field_map}}}&\\
0\ar[r]& \End(E_\Rep)^{GL(E)} \ar[r]&
\calD_1^\spad(E_\Rep)^{GL(E)}\ar[r]& 
\calT(\Rep^A_E)^{GL(E)}\ar[r]&0.
}$$
\end{rem}

\subsection{Rep-functor and the de Rham complex.}
 Given a noncommutative differential
$n$-form, $\omega=$
$a_0\,da_1\cdots da_n\in\ncO^n(A)$, we define
$$
\rho(\omega)=\tr(\wh a_0\,d\wh a_1\cdots d\wh a_n)\in\Omega^n({\Rep^A_E})^{GL(E)},
$$
where the $\Omega^n$ on the right-hand side denotes the usual
differential $n$-forms.  Since the trace is symmetric (that is,
$\tr(ab)=\tr(ba)$), it vanishes on 
$[\ncO^\hdot(A),\ncO^\hdot(A)]$.  Hence the $\rho$ descends to
a map $\rho\colon \DR^\hdot(A)\to\Omega^\hdot({\Rep^A_E})^{GL(E)}$.  

We would like to deal with equivariant cohomology rather than this
invariant part of the (usual) de Rham complex.

\subsection{Equivariant Cohomology}
Let $X$ be an algebraic variety, and let $\g g$ be a Lie algebra.
Suppose we have a Lie algebra map $\g g\to\calT(X)$.  
The $\g g$-equivariant algebraic de Rham complext of $X$ is  the complex
\beq{dr_equiv}
((\Omega^\hdot(X)\otimes\C[\g g])^{\g g},\,d_\g g),
\end{equation}
where $\g g$ acts on $\Omega^\hdot(X)$, by the Lie derivative $\L$,
the
$\g g$-action on $\C[\g g]$ is induced by the  adjoint action of $\g g$
 on, and the differential (called the Koszul differential) $d_\g g$ is defined in the following way.  Choose a basis $e_r$ for $\g g$ and let $e_r^*$ denote the dual basis.  Then we set
$$
d_\g g=d+\sum_{r=1}^ne_r^*i_{e_r},
$$
where $d$ is the usual de Rham 
differential on $\Omega^\hdot(X)$ and $i_{e_r}$ denotes contraction by $e_r$.  

\begin{rem} The cohomology of the complex \eqref{dr_equiv}
are called the equivariant cohomology of $X$. Assume $G$
is a Lie group with Lie algebra $\g g$ such that
the map $\g g\to\calT(X)$ can be exponentiated to
a {\em free} $G$-action on $X$ and, moreover,
 the orbit space
$X/G$ is a well-defined algebraic variety.
Then, the complex  \eqref{dr_equiv} computes the
ordinary de Rham cohomology of $X/G$.
In the general case, the leaves of the vector
fields coming from the image of  $\g g\to\calT(X)$
allow to consider the {\em stack quotient}
$X/\g g$, and the  cohomology of the complex \eqref{dr_equiv}
should be thought of as the de Rham cohomology
of that  stack quotient.
\eer

Choose any element $\omega\otimes f\in(\Omega^\hdot(X)\otimes\C[\g
g])^\g g$. 
 Then $\omega\otimes f(0)\in\Omega^\hdot(X)^\g g$, which yields a surjection
$$
(\Omega^\hdot(X)\otimes\C[\g g])^\g g\to\Omega^\hdot(X)^\g g
$$
that intertwines the differentials $d_\g g$ and $d$.  However, this map
is usually neither surjective nor injective on the level of cohomology.

Now, we let $A$ be an associative algebra and $E$ a finite-dimensional
vector space.  
Let $X={\Rep^A_E}$ be the representation variety of $A$ with 
the natural $GL(E)$-action.  Then we obtain a map
$$
{\mathsf{rep}}_\DR:\ \DR(A)\to \Omega^\hdot(\Rep^A_E)^\g g,
\quad\alpha\mapsto \tr(\wh\alpha),
$$
where $\g g=\g{gl}(E)$, given by $a_0\,da_1\cdots da_n\mapsto\tr(\wh
a_0\,d\wh a_1\cdots d\wh a_n)$, where for each $a\in A$, $\wh a\colon
\Rep^A_E\to\End_\C E$ is given by $\wh a(\f)=\f(a)$.

\begin{prob}
Is there a natural map making the following diagram commute~?
$$
\xymatrix{&\left(\bigl(\Omega^\hdot(\Rep^A_E)\otimes\C[\g g]
\bigr)^\g g ,\,d_\g g\right)\ar[d]\\
(\DR^\hdot(A),d)\ar@{-->}[ur]^{?}\ar[r]&(\Omega^\hdot(\Rep^A_E)^\g g, d)} 
$$
\end{prob}

Here is a simple example of a possible lift (dotted arrow above)
in a very special case.
Let $\alpha=da_1 da_2\in \DR^2(A)$, 
so $\tr(\wh\alpha)=\tr(d\wh a_1 d\wh a_2)$.
Finding the image of
 $\alpha$ under the dotted map amounts to finding
a degree two element $\widetilde\alpha\in
\bigl(\Omega^\hdot(\Rep^A_E)\otimes\C[\g g]\bigr)^\g g$
of the form $\widetilde\alpha= \alpha\otimes 1+\Phi$,
where $\Phi\in \bigl(\C[\Rep^A_E]\otimes\g g^*\bigr)^\g g$.
Since $\alpha$ is closed, the compatibility of the 
 dotted map with the differentials reads
$$i_x\tr(\wh\alpha)= d(\langle\Phi,x\rangle)\quad,\forall x\in\g g.
$$
Using the trace pairing $x,y \mapsto \tr(x\cdot y)$ on
$\g g=\End_C E$, we may identify $\g g$ with $\g g^*$,
and view $\Phi$ as an element of $\C[\Rep^A_E]\otimes\g g$.
With these identifications, we have $\langle\Phi,x\rangle=
\tr(\Phi\cdot x).$ It is easy to check that the equation
above is satisfied if
one puts $\Phi:= \wh a_1\cdot \wh a_2$ (product in 
the associative algebra $\k[\Rep^A_E]\otimes \End_\C E$).

\section{Double-derivations and the double-tangent bundle.}
\subsection{}
Given an associative algebra $A$, we define the $A$-bimodule
of {\em double-derivations} of $A$ by
\beq{double_der}
\dder(A):=\Der(A,A\ee) = \Hom_{A\ee}(\ncO^1(A), A\ee),
\end{equation}
where the $A$-bimodule structure comes from
the $A\ee$-action on the entry $A\ee$ by {\it right} multiplication.

The fundamental exact sequence \eqref{ncO1}
gives rise to 
an exact sequence
\beq{double_fund}
0\to\End_{A\ee}\bigl(\ncO^1(A)\bigr)\to\dder(A)\to\Der(A)\to
\Ext^1_{A\ee}\bigl(\ncO^1(A),\,\ncO^1(A)\bigr).
\end{equation}

 \begin{rem} Note that
the space $\End_{A\ee}\bigl(\ncO^1(A)\bigr)$ on the left is an associative
algebra, in particular, a Lie algebra with respect to the commutator
bracket. The space $\Der(A)$ is also a Lie algebra.
Furthermore, the sequence above resembles the exact sequence
for the Atiyah algebra, cf. \S\ref{atiyah_alg},
 of the `vector bundle' $\ncO^1(A).$
\eer

 We observe that if $\ncO^1(A)$ is a projective
$A\ee$-module (in which case the algebra $A$ is said to
be {\em formally smooth}, see \S22), the $\Ext^1$-term on the right of
\eqref{double_fund} vanishes, hence the map
$\dder(A)\to\Der(A)$ becomes surjective.

It is immediate to check that the map
\beq{delta}
\Delta: A\map A\ee=A\otimes A,\quad a\mto a\otimes 1 -
1\otimes a
\end{equation}
 gives a derivation, i.e.,
we have $\Delta\in \dder(A)$. The derivation $\Delta$
maps to zero under the projection
 $\dder(A)\to\Der(A)$, cf. \eqref{double_fund},
hence,  $\Delta$ may be identified with an element of
$\End_{A\ee}\bigl(\ncO^1(A)\bigr)$.
It is easy to verify that the latter element
is nothing but the identity map
$\Id_{\ncO^1(A)}$.

\begin{examp}\label{dere_curve} Let $A=\C[X]$ be the coordinate ring of an
affine algebraic variety. Then, $A\ee=\C[X\times X]$,
and the bimodule $\ncO^1(A)$ is the ideal 
of the diagonal $X_\Delta\sset X\times X.$
We see that if $\dim X >1$, then
$X_\Delta$ has codimension $\geq 2$ in $X\times X$,
It is easy to deduce,
by Hartogs type theorem, that
in this case  the map $\ad: A\ee\to \dder(A)$
is an isomorphism.

Assume now that $\dim X =1$, i.e., $X$ is a curve.
Then, we have $\dder(A)=\Gamma(X\times X,\,\oo_{X\times X}(X_\Delta)),$
is the space of rational functions on $X\times X$ with simple
pole along $X_\Delta.$  On the other hand, if $X$ is a curve ($\dim X=1$),
then the exact sequence in \eqref{double_fund} reduces to 
\beq{extx}
0\to\C[X\times X]\to \dder(A)\to \calT(X)\to 0.\qquad\lozenge
\eeq
\end{examp}

The ``meaning'' of the double-derivation bimodule
may be seen
in terms of the $\Rep$-functor, via the following construction  due to
Kontsevich-Rosenberg \cite{KR}.

Fix a finite dimensional vector space $E$, and
let $\rho: A \to \End_\C E$ be a representation of $A$ in $E$
(to be denoted $E_\rho$). 
The action map $A\o E_\rho\to
E_\rho$ 
is surjective since $A$ contains the unit. This gives a surjective map of $A$-modules
${\mathtt{act}}: A \o E \onto E_\rho$,
where $ A \o E $ is regarded as a free left $A$-module 
generated by the vector space $E$.
We set $K_\rho:=$ ${\Ker(A \o E \onto E_\rho),}$
a left $A$-module. 

For any left $A$-module $M$, we have
$\Hom_A(A \o E,\,M)=\Hom_\C(E_\rho,M)$ and
$\Ext_A^1(A \o E,\,M)=0$. Hence,
the long exact sequence 
of Ext-groups arising from the short exact sequence
$ K_\rho\into A \o E \onto E_\rho$ reads:
\begin{align}\label{ext_rep_kr}
0\to \Hom_A(E_\rho,M)\stackrel{{\mathsf{ev}}}\map\Hom_\C(E_\rho,M)
\map & \Hom_A(K_\rho,M)\\
&\map\Ext^1_A(E_\rho, M)\to 0.\nonumber
\end{align}
It is easy to see that if $M$ is finite dimensional over $\C$,
then the space $\Ext^1_A(E_\rho, M)$ is also  finite dimensional over $\C$.
It follows that 
\beq{fin_dim_K}\dim_\C M<\infty\quad\Longrightarrow\quad
\dim_\C \Hom_A(K_\rho,M)<\infty.
\end{equation}

\subsection{}
We consider
the scheme $\Rep^A_E$. 
To simplify notation, write
 $\End:=\End_\C E$, a simple associative
algebra isomorphic
to $\Mat_n\C$ where $n=\dim E$.
We set $R:=\C[\Rep^A_E]$,
the coordinate ring of the scheme $\Rep^A_E$,
and write $\C_\rho$  for the 1-dimensional $R$-module
such that $f\in \C[\Rep^A_E]$ acts in $\C_\rho$ 
via  multiplication by the number $f(\rho)\in\C$.
We identify the algebra $R\ee=R\otimes R$
with $\C[\Rep^A_E\times\Rep^A_E].$
We have the tautological algebra map
$\rep: A\to R\otimes\End,\,a\mapsto \wh a$. 

Further, we consider the algebra $R\otimes\End\otimes R$.
There are two algebra maps
$\rep_l,\rep_r: A\to R\otimes\End\otimes R,$
given by the formulas
$a\mapsto \wh a\otimes 1,$ and $a\mapsto 1\otimes\wh a,$
respectively. These two maps make
$R\otimes\End\otimes R$ an $A$-bimodule.

We  define the double-tangent $R$-bimodule
by
\beq{d_tang}
\calT\ee_E(A):=\Der(A,R\otimes_\k\End_\k E\otimes_\k R).
\end{equation}
This is clearly a (not necessarily free)
 $\C[\Rep^A_E\times\Rep^A_E]$-module
whose geometric  fiber at a point
$(\rho,\varphi)\in \Rep^A_E\times\Rep^A_E$
is defined as
$$\calT\ee_E(A)|_{(\rho,\varphi)}:=(\C_\rho\otimes\C_\varphi)
\bigotimes\nolimits_{R\otimes
R}\calT\ee_E(A).
$$
The  corresponding
coherent sheaf 
on $\Rep^A_E\times\Rep^A_E$ will be referred to as the {\em
double-tangent bundle} on $\Rep^A_E$.

\begin{lem}\label{double_T} For any
$(\rho,\varphi)\in \Rep^A_E\times\Rep^A_E$, one has a canonical
isomorphism
$$\calT\ee_E(A)|_{(\rho,\varphi)}\;\simeq\;\Hom_A(K_\rho,E_\varphi).
$$
\end{lem}

\begin{proof}
We observe that this  $A$-bimodule structure on $R\otimes\End\otimes R$
commutes with the obvious $R$-bimodule structure.
It follows that the latter  induces
an $R$-bimodule structure
on each Hochschild cohomology
group $\HH^p(A, R\otimes\End\otimes R)\,,\,p\geq 0.$
In particular, the sequence \eqref{derM} becomes, in our
present situation, the following exact sequence of 
$R$-bimodules (to be compared with \eqref{ext_rep_kr}):
\begin{align}\label{der_double_HH}
0\to \HH^0(A,R\otimes &\End\otimes R)\to R\otimes\End\otimes R\\
&\to \Der(A,R\otimes\End\otimes R)
\to \HH^1(A,R\otimes\End\otimes R)\to 0.\nonumber
\end{align}

We leave to the reader to verify that,
for any two representations $\rho,\varphi\in \Rep^A_E$,
 the geometric fibers at $(\rho,\varphi)$ of the 
$R$-bimodules occurring in \eqref{der_double_HH}  are given by
\begin{align*}
&\HH^0(A,R\otimes\End\otimes
R)\big|_{(\rho,\varphi)}=\Hom_A(E_\rho,E_\varphi),\\
&\;(R\otimes\End\otimes
R)\big|_{(\rho,\varphi)}\enspace=\enspace\Hom_\C(E_\rho,E_\varphi),\\
&\Der(A,R\otimes\End\otimes R)\big|_{(\rho,\varphi)}=\Hom_A(K_\rho,E_\varphi),\\
&\HH^1(A,R\otimes\End\otimes
R)\big|_{(\rho,\varphi)}=\Ext_A^1(E_\rho,E_\varphi).
\end{align*}

The statement  of the Lemma is now clear.
\end{proof}

Next, we consider the following maps
$$
\xymatrix{
A\otimes A
\ar[rr]^<>(0.5){\rep_l\otimes\rep_r}&&
(R\otimes\End)\otimes(\End\otimes R)
\ar[rr]^<>(0.5){{\id_R\otimes m\otimes\id_R}}
&&
R\otimes\End\otimes R,
}
$$
where $m: \End\otimes\End\to\End$ denotes the multiplication in the
algebra $ \End$.
We observe that the two maps above are morphisms of $A$-bimodules,
i.e., of $A\ee$-modules. Let $\tau$ denote the composite
map $A\ee=A\otimes A\too R\otimes\End\otimes R,$
which is an $A\ee$-module map again.

Thus, from formula \eqref{d_tang} we deduce 
that the $A\ee$-module map $\tau$ gives rise to a canonical
map
$$\rep\ee:\;\dder(A):=\Der(A,A\ee)\stackrel\tau\too
\Der(A,R\otimes\End\otimes R)=:\calT\ee_E(A).
$$

\subsection{Double-derivations for a free algebra.}
Throughout this subsection we let $A=\C\langle
x_1,\ldots,x_r\rangle$ be
a free associative algebra
on $r$ generators. For each $i=1,\ldots,r$,
we introduce an element $D_i\in \dder(A)$
uniquely defined by the following conditions
$$D_i(x_i)=1\otimes 1\quad\text{and}\quad
D_i(x_j)=0 \quad\forall j\neq i.
$$

Now, let $F: A\to A$ be an algebra homomorphism,
such that
$x_1\mapsto F_1,\ldots,x_r\mapsto F_r,$
where $F_1,\ldots,F_r\in \C\langle x_1,\ldots,x_r\rangle$.
The $r$-tuple
$F_1,\ldots,F_r$ determines $F$ uniquely,
and we will think of the elements
$F_i\in \C\langle x_1,\ldots,x_r\rangle$ as  functions $F_i=F_i(x_1,\ldots,x_r)$
in $r$ non-commutative variables $x_1,\ldots,x_r$.

We define the {\em Jacobi matrix} for the map $F$ to be
the following $A\otimes A$-valued $r\times r$-matrix
\beq{jac}
DF=\|D_i(F_j)\|_{i,j=1,\ldots,r}\in \Mat_r(A\otimes A),
\end{equation}
where $D_i(F_j)$ denotes the image of $F_j\in A$ under the derivation
$D_i: A\to A\otimes A$ introduced above.

The following result is proved by a straightforward computation,
see \cite[Lemma 6.2.1]{HT}.
\begin{prop}[Chain rule] For any two algebra homomorphisms
$F,G: A\to A,$ in $\Mat_r(A\otimes A)$ one has
$$ D(G\ccirc F)= DG \ccirc DF,\quad\text{where}\quad
|DG \ccirc DF|_{kl}:=\sum_i D_k(G(F(x_i)))\cdot D_i(x_l).\quad\Box
$$
\end{prop}

Next, we fix a finite dimensional vector space $E$
and consider the scheme $\Rep^A_E$, which is canonically
isomorphic to  $\End_\C(E)\times\ldots\End_\C(E)$ via
the identification, cf. \eqref{rep_free} 
$$\Rep^A_E\ni \rho\mto \bigl(X_1:=\rho(x_1),\ldots,
X_n:=\rho(x_r)\bigr).
$$

Analogously to the case of derivations considered in Example
\ref{der_free},
any automorphism $F: A \to A$ induces an automorphism 
$\wh{F}$ of the vector space $\Rep^A_E$. The map $\wh{F}$
is given, in terms of the $r$-tuple
$F_1,\ldots,F_r\in \C\langle x_1,\ldots,x_r\rangle,$
 by the formula
\beq{hat_F}\wh{F}:\,\bigl(X_1,\ldots,
X_n\bigr)\mto \bigl(F_1(X_1,\ldots,
X_n),\ldots, F_r(X_1,\ldots,X_n)\bigr).
\end{equation}

For each $i=1,\ldots,r$, and any $a\in A$, we have an
element $D_i(a)=\sum D'_i(a)\otimes D''_i(a)\in A\otimes A$.
Thus, given a point $\rho\in\Rep^A_E$,
 we have a well-defined element
$\rho(D_i(a))=\sum \rho(D'_i(a))\otimes \rho(D''_i(a))\in \End_\C(E)\otimes\End_\C(E)$.
In particular, for each $i,j=1,\ldots,r$, there is an
element $\rho(D_i(F_j))\in \End_\C(E)\otimes\End_\C(E)$.
These elements
give rise to a linear map
$$(DF)_\rho:\ \Rep^A_E\too \Rep^A_E,$$
given by the formula
\begin{align}\label{jac_map}
(Z_1,\ldots,&Z_r)\mto \\
&\Bigl(\sum \rho(D'_1(F_j))\cdot Z_j\cdot \rho(D''_1(F_j))\,,\,\ldots\,,\,
\sum \rho(D'_r(F_j))\cdot Z_j\cdot \rho(D''_r(F_j))\Bigr).\nonumber
\end{align}

The differential of the map \eqref{hat_F} can be read off from the Jacobi
matrix \eqref{jac} by means of the following result, proved in 
\cite[Lemma 6.2.2]{HT}.

\begin{prop} The differential of the map $\wh{F}$
at a point $\rho=(X_1,\ldots,
X_n)\in \Rep^A_E$ is a linear map $d_\rho\wh{F}:
T_\rho\Rep^A_E\to T_\rho\Rep^A_E$
given by formula \eqref{jac_map}, that is, we have
$d_\rho\wh{F}=(DF)_\rho\,,\,\forall \rho \in\Rep^A_E.$\qed
\end{prop}

\subsection{The Crawley-Boevey construction.}\label{CrB}
Recall next that we have an isomorphism $(A\ee)^{\op}\cong A\ee$.
Therefore,
 {\em right} multiplication in the
algebra $A\ee$ makes $A\ee$, hence $\dder(A)$ a 
{\em left} $A\ee$-module,
that is, an $A$-bimodule.

We consider
 $T_A\dder(A)$, the tensor algebra
of the  $A$-bimodule $\dder(A)$,
and view the derivation $\Delta$, see \eqref{delta},
as an element of $T^1_A\dder(A)=\dder(A).$
Given an element
$c\in A$, viewed as an element of  $T^0_A\dder(A)=A$,
we introduce, following  \cite{CrB}, an associative algebra
$\Pi^c(A):=T_A\dder(A)/\llb \Delta-c\rrb,$
where $\llb \Delta-c\rrb$ stands for the two-sided
ideal generated by $\Delta-c$.

We first consider the case: $c=0$.
To this end, we use
formula \eqref{derM} in the special case $M=A\ee$
and get a canonical exact
sequence of $A$-{\em bimodules}:
\beq{double_der_seq}
0\map\HH^0(A,A\ee)\map A\ee
\stackrel{\ad}\map\dder(A)\map \HH^1(A,A\ee)\map 0.
\end{equation}

This way, one proves
\begin{prop}\label{DD_prop} \vi The above map $A\ee
\map\dder(A)$ sends the element $1\in A\ee$ to the
derivation $\Delta$, see \eqref{delta}. Thus,
$\HH^1(A,A\ee)$ is isomorphic to a quotient
of 
$\dder(A)$ by the sub $A$-bimodule generated by $\Delta$. 

\vii The isomorphism in 
\vi induces  a
graded algebra isomorphism \hfill\break
$\Pi^0(A)=T_A\dder(A)/\llb \Delta\rrb\cong
T_A\HH^1(A,A\ee)$.\qed
\end{prop}

In the general case of an arbitrary $c\in A$, the grading on the
tensor algebra induces a natural increasing filtration
on the algebra $T_A\dder(A)/\llb \Delta-c\rrb$.
Furthermore, part (ii) of Proposition \ref{DD_prop}
yields a canonical surjective graded algebra
map
$$T_A\HH^1(A,A\ee)\onto \gr T_A\dder(A)/\llb \Delta-c\rrb.
$$
It would be interesting to find some sufficient conditions
for the map above to be an isomorphism (a 
Poincar\'e-Birkhoff-Witt type isomorphism).
 
The case $c=1$ is especially important. 
To explain the geometric meaning of the algebra
$\Pi^1(A)$, we return to the setup of Example \ref{dere_curve} and
let $A=\C[X]$ be the coordinate ring of a smooth
affine algebraic curve.

\begin{thm}\label{BBthm} The algebra $\Pi^1(\k[X])$
is canonically isomorphic, as a filtered algebra, to
$\scr D(X,\Om^{1/2}_X)$,
the filtered algebra of {\sf{twisted differential operators}}
acting on {\sf{half-forms}} on $X$.
\end{thm}

\begin{rem} It has been shown in \cite{CrB}, that
there is a natural graded algebra isomorphism
 $\Pi^0(\k[X])\cong \k[T^*X],$ where
$T^*X$ stands for the total space
of the cotangent bundle on the curve $X$.
\end{rem}

We refer to \cite{BB} for the basics of the theory
of  twisted differential operators
on an algebraic variety. According to this
theory, sheaves of twisted differential operators
on a smooth affine algebraic variety $X$
are classified by the group
$H^1_\DR(X, -)$. This group vanishes if
$X$ is a curve. Thus, in the case at hand,
one can replace the algebra
$\scr D(X,\Om^{1/2}_X)$ in Theorem \ref{BBthm}
by the algebra $\scr D(X,\oo_X)$
of usual (not twisted)  differential operators
on $X$. In this form, 
an algebra isomorphism
$\Pi^1(A)\cong\scr D(X,\oo_X)$
 has been already established
in \cite{CrB}.

The advantage of our present version of
 Theorem \ref{BBthm}, involving
{\sf{twisted}} differential operators,
is that the isomorphism of the Theorem
becomes canonical. In particular,
both the statement and proof of the
Theorem generalize easily to
the case of {\em sheaves} of algebras
of twisted differential operators
on a {\em not necessarily affine}
 smooth curve $X$. To explain this,
observe that in the case of an affine curve
we have $A\ee=\k[X\times X]$,
hence
 the bimodule $\ncO^1 A\sset A\otimes A$ is the ideal 
of the diagonal divisor $X_\Delta\sset X\times X.$ Therefore, we have 
\beq{dderD}
\dder A=\Hom_{A^e}(\ncO^1 A, A\otimes A)=
\Gamma(X\times X,\,\oo_{X\times X}( X_\Delta )),
\eeq
is the space of regular functions on $(X\times X)\sminus X_\Delta$ with 
at most simple
pole along the diagonal divisor $ X_\Delta .$

Now, in the case of an arbitrary smooth,
 not necessarily affine, curve $X$ it is natural to define the sheaf
$\ncO^1(\oo_X)$ to be the ideal sheaf of the diagonal divisor
 $X_\Delta\sset X\times X,$
and to put
$$\dder(\oo_X)=\Hom_{\oo_{X\times X}}(\ncO^1(\oo_X),\,\oo_{X\times X})
\cong\oo_{X\times X}( X_\Delta ),$$ 
see \eqref{dderD}. Repeating the definitions
above, on constructs
a sheaf  $\Pi_X$ of filtered algebras on $X$ that
  corresponds,
locally in the Zariski topology, to the algebra $\Pi(A)$.
The sheaf-theoretic version of
Theorem \ref{BBthm} says that there is a canonical
isomorphism 
between  $\Pi_X$, viewed as a
 sheaf  of filtered algebras on $X$
in the  Zariski topology,
and 
the sheaf of twisted differential operators
on $X$
acting on half-forms.
We note that,
for any complete
curve of genus $\neq 1$,
the sheaf of twisted differential operators
acting on half-forms is
{\em not} isomorphic
to the sheaf of  usual, not twisted,
  differential operators.

\subsection{Sketch of proof of Theorem \ref{BBthm}.} First of all,
for any associative not necessarily commutative algebra $A$
we have the tautological $A$-bimodule imbedding $\ncO^1A\into A\o A$
and an  $A$-bimodule map $\ad: A\o A \dder A$. Composing these two maps,
yields a canonical $A$-bimodule morphism
\beq{self}
\phi:\
\ncO^1A\too \dder A=\Hom_{A\ee}(\ncO^1A, A\ee).
\eeq
We observe that the morphism above is {\em self-dual}, i.e.,
applying the functor $\Hom_{A\ee}(-,A\ee)$ to $\phi$ one gets
the same morphism $\phi$ again.

\begin{lem}\label{phi} Assume that the bimodule
$A\o A$ has trivial center, that is for $x\in A\o A$ we have
$ax=xa,\,\forall a\in A\en\Longrightarrow\en x=0.$ Then, the map
$\phi$ in \eqref{phi} is injective and one has a canonical $A$-bimodule
isomorphism
$$\Pi_{\leq 1}(A)\cong \Coker(\phi)=\dder A/\ncO^1A.$$
\end{lem}
\begin{proof} For any three objects $U\sset V\sset W$ of an abelian 
category, one has a short exact sequence
\beq{x}0\map
V\stackrel{f}\map W\oplus (V/U) \stackrel{g}\map W/U\map 0,
\eeq
where the map $f$ is the direct sum of the imbedding
$V\into W$ with the projection $V\onto V/U$
and the map $g$ is given by $g(w\oplus (v\mres U))=-w\mres U.$

Now, let $A$ be an associative algebra satisfying the assumptions of the
 lemma.
We let the triple
 $U\sset V\sset W$ to be the following triple of $A$-bimodules
$U=\ncO^1A\sset V=A\o A\sset\dder A,$ where the rightmost imbedding
 induced
by the map $\ad$, which is injective by our assumptions.
From the fundamental short exact sequence
we get $A=A\o A/\ncO^1A=V/U$. Further, the bimodule  $\Pi^1_{\leq 1}(A)$
is by
definition the
quotient $\dder A\oplus A= W\oplus V/U$ by the image
of $A\o A=U$. Hence, the short exact sequence \eqref{x} yields
the desired isomorphism
$$\Pi^1_{\leq 1}(A)\cong(\dder A\oplus A)/A\o A
\cong (W\oplus V/U)/V\cong W/U=\dder A/\ncO^1A.$$
\end{proof}

We can now proof  Theorem \ref{BBthm}.
By standard arguments, see e.g., \cite{BB},
to prove the theorem it suffices to construct
a canonical isomorphism of {\em Atiyah algebras}:

\beq{atiyah}
\xymatrix{
0\ar[r]&A=\Pi_{\leq 0}(A)\ar[r]\ar@{=}[d]^<>(0.5){\Id}&
\Pi_{\leq 1}(A)\ar[r]\ar@{=}[d]^<>(0.5){\Phi}&
\Der A\ar[r]\ar@{=}[d]^<>(0.5){\Id}&0
\\
0\ar[r]&\k[X]\ar[r]&
{\scr D}_{\leq 1}(X,\Om^{1/2}_X)\ar[r]&
\calT(X)\ar[r]&0
}
\eeq

The extension in the bottom row of this diagram can
be computed explicitly.
Specifically, one shows
that this extension is canonically isomorphic,
to (the spaces of global sections of) the
following extension of sheaves 
$$
0\to \oo_{X\times X}/\oo_{X\times X}(-X_\Delta)
\to\oo_{X\times X}(X_\Delta)/\oo_{X\times X}(-X_\Delta)
\to\oo_{X\times X}(X_\Delta)/\oo_{X\times X}
\to0.
$$
The quotient sheaves on both sides are nothing but
the cotangent and tangent sheaf on $X_\Delta=X$, respectively.
So the above extension reads
\beq{result}
0\to \oo_X\to\oo_{X\times X}(X_\Delta)/\oo_{X\times X}(-X_\Delta)
\to\calT_X\to 0.
\eeq

Further, 
we have the following diagram of
 short exact sequences, cf.  \eqref{extx}.
\beq{vg2}
\xymatrix{
0\ar[r]&A\otimes A\ar[r]^<>(0.5){j}\ar@{=}[d]&
\Der(A,A\o A)\ar[r]\ar@{=}[d]^<>(0.5){\Psi}&
\Der A\ar[r]\ar@{=}[d]^<>(0.5){\Phi}&0
\\
0\ar[r]&\Gamma(\oo_{X\times X})\ar[r]&
\Gamma(\oo_{X\times X}(X_\Delta))\ar[r]&
\Gamma(\oo_{X\times X}(X_\Delta)/\oo_{X\times X})
\ar[r]&0
}
\eeq

In the  bottom row of the diagram above we have used shorthand notation
$\Gamma(-)$ for $\Gamma(X\times X,-)$; this row is obtained by applying
 the global sections functor to
the natural extension of sheaves on $X\times X$. 
The vertical isomorphism $\Phi$, in the diagram, follows from the
identification
$\Der A=\calT(X)$, with the space of regular
vector fields on $X$. 
The vertical isomorphism $\Psi$ comes from
\eqref{dderD}.

Observe  that the function $1\in A\otimes A$
corresponds under the  above identifications
to the element
$\Id_\Om\in \Hom_{A^e}(\ncO^1A,\ncO^1A)$. Therefore,  in the diagram we have
$j(1)=\Delta$, and the map  $j$ is nothing but
the imbedding 
$\ad :  A\otimes A\into  \Der(A,A\otimes A),$
of inner derivations.

Thus, we use the isomorphism
$\calT(X)=\Der A=\dder A/A\ee,$ see \eqref{vg2},
to identify the  short exact sequence in \eqref{result}
with the canonical short exact sequence
$$
0\to A\to \dder(A)/\ncO^1A\to \Der A\to 0.
$$
This  short exact sequence can be identified with
the extension in the top row of  diagram \eqref{atiyah}
using Lemma \ref{phi}.
\qed

\section{Noncommutative Symplectic Geometry}
\subsection{}
Let $A$ be an associative algebra and $\omega\in \DR^2(A)$ a
noncommutative $2$-form. Contraction with $\omega$ gives a linear
map $i_\omega: \Der(A)\to \DR^1(A)\,,\,\theta\mapsto i_\theta\omega.$
\begin{defn}
The pair $(A,\omega)$ is called a \emph{noncommutative symplectic
manifold} 
if  $d\omega=0$ in $\DR^3(A)$, i.e., the $2$-form $\omega\in \DR^2(A)$
is closed, and furthermore, $\omega$ is \emph{nondegenerate},
 i.e., the map $i_\omega: \Der(A)\to \DR^1(A)$  is a bijection. 
\end{defn}

Fix a noncommutative symplectic
manifold $(A,\omega)$.

\begin{defn}
A derivation  $\theta\in\Der(A)$ is called \emph{symplectic} if
$\L_\theta\omega=0$.
  We denote by $\Der_\omega(A)$ the Lie algebra of all symplectic derivations, i.e., 
$\Der_\omega(A)=\{\theta\in\Der(A)\st\L_\theta\omega=0\}$.
\end{defn}
 The space $\Der_\omega(A)$ inherits the Lie algebra structure from
$\Der(A)$ 
given by commutators.

\begin{lem}\label{der_om}
A derivation $\theta$ is symplectic if and only if $i_\theta\omega$ is
closed in $\DR^1(A)$.
\end{lem}

\begin{proof}
This is simply an application of Cartan's formula, namely
$$
\L_\theta\omega=i_\theta\circ d\omega+d\circ i_\theta\omega=d(i_\theta\omega),
$$
since $d\omega=0$ by assumption.
\end{proof}

Recall that, viewing $\sr(A)$ as $\DR^0(A)$, we have a map $d\colon \sr(A)\to \DR^1(A)$.  For every $f\in \sr(A)$, $df$ is exact in $\DR^1(A)$, hence closed.  The previous lemma shows that the closed forms in $\DR^1(A)$ are identified with the symplectic vector fields.  We let $\theta_f$ denote the symplectic vector field associated to $df$ under this identification.  As in the classical theory, $\theta_f$ is called the Hamiltonian vector field associated to $f$.  We then define a Poisson bracket on $\sr(A)$ by
$$
\{f,g\}=i_{\theta_f}(dg).
$$
Notice that since $dg\in \DR^1(A)$, $\{f,g\}$ is indeed contained in $\DR^0(A)=\sr(A)$.  It is clear that we have the following several equalities
$$
\{f,g\}:=i_{\theta_f}(dg)=i_{\theta_f}i_{\theta_g}\omega=-i_{\theta g}i_{\theta f}\omega=-i_{\theta_g}(df)=-\L_{\theta_f}g=-\L_{\theta_g}f.
$$

\begin{thm}
(i) $\{-,-\}$ makes $\sr(A)$ a Lie algebra.\\
(ii) The map $f\mapsto\theta_f$ is a Lie algebra homomorphism from $\sr(A)$ to $\Der_\omega(A)$.
\end{thm}

\begin{proof}
The skew symmetry of $\{-,-\}$ is immediate.  We will first establish (ii), which will then prove (i).  Let $\theta,\gamma$ be arbitrary derivations from $A$ to $A$.  Then by our standard identities we have
$$
i_{[\theta,\gamma]}=\L_\theta i_\gamma-i_\gamma\L_\theta=di_\theta i_\gamma+i_\theta di_\gamma-i_\gamma di_\theta-i_\gamma i_\theta d.
$$
Now we specialize to the case $\theta=\theta_f$ and $\gamma=\theta_g$ and consider $i_{[\theta_f,\theta_g]}\omega$.  By the definition of the Hamiltonian vector field, we know that $i_{\theta_f}\omega=df$ and $i_{\theta_g}\omega=dg$.  Also, $d\omega=0$ by definition, so we are left with
$$
i_{[\theta_f,\theta_g]}\omega=di_{\theta_f}(dg)+i_{\theta_f}d(dg)-i_{\theta_g}d(df).
$$
The latter two terms are zero, since $d^2=0$.  Therefore,
$$
i_{[\theta_f,\theta_g]}\omega=di_{\theta_f}(dg)=d\{f,g\}.
$$
But by definition, $\theta_{\{f,g\}}$ is the unique symplectic
derivation such 
that $i_{\theta_{\{f,g\}}}\omega$ $=d\{f,g\}$.  Hence $\theta_{\{f,g\}}=[\theta_f,\theta_g]$, which establishes (ii).  Since $f\mapsto\theta_f$ is an isomorphism (\emph{a priori} only of vector spaces), this shows that $\{f,g\}$ must satisfy the Jacobi identity (finishing (i)).  Indeed, if we choose $h\in \sr(A)$, then the equality $\theta_f\theta_gh-\theta_g\theta_fh=\theta_{\{f,g\}}h$ becomes $\{f,\{g,h\}\}-\{g,\{f,h\}\}=\{\{f,g\},h\}$, which is precisely Jacobi's identity after some rearranging.  
\end{proof}

\begin{examp}\label{symp_free}
Let  $(V,\omega_V)$ be a finite dimensional symplectic vector space.
We claim that the symplectic structure on $V$ induces  a
noncommutative  symplectic structure on
$A=T^\hdot (V^*)$, the tensor algebra. 
Explicitly, let $x_1,\ldots,x_n,y_1,\ldots,y_n\in V^*$
be  canonical coordinates in $V$, so that $\omega_V=\sum_{i=1}^n\,
dx_i\wedge dy_i$.
We will see in the next section that one has
$$\Der(A)=A\otimes V,\quad\DR^1(A)=A\otimes V^*,\quad
\ncO^i(A)=A\otimes T^i(V^*)\otimes A,\,\forall i\geq 0.
$$
We put
 $\om_A:=\sum_{i=1}^n\,1 \otimes (x_i\otimes y_i)\otimes 1\in A\otimes T^2(V^*)\otimes A=
\ncO^2(A)$.
Further, 
the nondegeneracy of $\omega_V$ implies that the assignment
$v\mapsto\omega_V(v,-)$ yields a vector space
isomorphism $V\iso V^*$. The latter induces
an isomorphism 
$$\id_A\otimes\omega_V:\ 
\Der(A)=A\otimes V\iso A\otimes V^*=\DR^1(A).
$$
It is easy to verify that the last isomorphism
is nothing but the map  $\theta\mapsto i_\theta\omega_A$
arising from the noncommutative 2-form $\omega_A\in \DR^2(A)$.
Thus, $(A=T(V^*), \om_A)$ is  a
noncommutative  symplectic structure.
\eex

\begin{question} \vi Given a noncommutative symplectic structure on an associative 
algebra $A$, can one define a Lie super-algebra structure
on $\oplus_{i\geq 1}\,\ncO^i(A)$ which is a noncommutative analogue
of the Lie super-algebra structure of Lemma \ref{ger_om} ?

\vii In case of a positive answer to part \vi, does the Lie
 super-algebra structure
on $\oplus_{i\geq 1}\,\ncO^i(A)$ combined with the
associative product on $\ncO^\hdot(A)$ give rise to the structure of
noncommutative Gerstenhaber algebra ?
\end{question}

\subsection{Noncommutative `flat' space}
In  Noncommutative Geometry, the free associative
algebra $\C\langle x_1,\ldots,$
$x_n\rangle$ plays the role
of  coordinate ring of an $n$-dimensional affine space.
It will be convenient to introduce
 an $n$-dimensional $\C$-vector space $V$
with coordinates $x_1,\ldots,x_n$ (thus, the elements
 $x_1,\ldots,x_n$ form a base in $V^*$, the dual space).
This allows to
adopt a
`coordinate free' point of view and to identify
the algebra $\C\langle x_1,\ldots,x_n\rangle$
 with $A=T^\hdot (V^*)$,
the tensor algebra of $V^*$.

  The
derivations from $A$ to an $A$-bimodule $M$
are specified precisely by a linear map from 
$V^*$ to $M$--it is then extended to a derivation uniquely by the
Leibniz rule. 
 So, $\Der(A,M)\simeq V\o M$.  Therefore, the functor
$\Der(A,-)$ is represented by the free $A$-bimodule
generated by the space $V^*$. Thus, we find
$$
\ncO^1(A)\simeq A\o V^*\o A\simeq A\bigotimes
(\oplus_{i>0}\,T^i(V^*))\simeq A\o \overline A,
$$
Hence, for any $p\geq 1,$ we obtain
$$\ncO^p(A)=T^p_A\ncO^1(A)= A\o V^*\o A\o \ldots \o 
V^*\o A
\quad\text{($p$ factors $V$)}.
$$
Observe  that the assignment
$a_1\otimes v_1\otimes\ldots \otimes v_p\otimes a_{p+1}
\mto 
a_1* v_1*\ldots * v_p* a_{p+1}$
gives an imbedding
$T^p_A\ncO^1(A)\into  A* A$ (= free product of two copies of $A$).
Further, 
an element of the form
$\ldots v_1\otimes 1_A\otimes v_2\otimes 1_A\otimes\ldots 1_A\otimes
v_m\otimes\ldots$ goes under this imbedding to
the element $\ldots*(v_1\otimes v_2\otimes\ldots\otimes
v_m)*\ldots$. Thus we deduce that the imbedding above yields
a graded algebra isomorphism
$$
\ncO^\hdot(A)\iso A* A= T(V^*)\,*\, T(V^*),$$
where the grading on the left-hand side corresponds
to the total grading with respect to the second
factor $A=T(V^*)$ on the right-hand side.

To describe $\DR^0(A)$, it is 
instructive to
identify $A=T(V^*)$ with $\C\langle x_1,\ldots,$
$x_n\rangle$.
Then, $A$ may be viewed
as  $\C$-vector space whose basis is 
formed by all possible  words in the alphabet formed by the $x_i$'s.
The algebra structure is given by concatenation of words.
Further, let $A_{\mathsf{cyclic}}\sset A$ be the
$\C$-linear span of all  cyclic words. It is clear that
the composite map $A_{\mathsf{cyclic}}\into A\onto A/[A,A]$
is a bijection. Thus, we may identify
$$\DR^0(A)=\sr(A)=A/[A,A]=A_{\mathsf{cyclic}}$$
via this bijection.

  It is more difficult to analyze the $k^{\text{th}}$ degree of $\DR(A)$
where $k>0$.  
In general we can only identify it as some quotient of $A\otimes TV$.  
However, in the particular case $k=1$, we find that
$$
\DR^1(A)=\ncO^1(A)/[\ncO^1(A),A]\simeq A\o V^*,
$$
as a complex vector space,
since no other combinations of $\ncO^k(A)$'s in the ``denominator'' will
yield degree one.

So, for all $x\in V^*$, we have an element $dx=1\otimes x\in \DR^1(A)$.  Similarly, for every $v\in V$, we have $\partial_v\colon \DR^0(A)\to A$ given by
$$
\partial_v(f)=df(v)
$$
for all $f\in \DR^0(A)=\sr(A)=A/[A,A]$.  If $x_i$ is a basis of $V^*$ and $x^i$ is the corresponding dual basis of $V$, consider the map $d$ which sends every $f\in \DR^0(A)$ to 
$$
\sum_{j=1}^n\pder{f}{x^j}\otimes dx_j,
$$
where we write $\pder{}{x^j}=\partial_{x^j}$.

\begin{examp} Let  $\dim_\C V=2$, and equip $V$ with 
 symplectic basis $x$ and $y$ and the standard symplectic 2-form  $dx\wedge dy$.
Let us  calculate the Poisson bracket $\{f,g\}$ of two elements of
$\sr(A)$, 
where $A=TV$.  
View an element $f\in \sr(A)=TV/[TV,TV]$ as a cyclic 
word in $x$ and $y$. We have already seen that $\DR^1(A)=A\otimes V$, so we can write $df=f_x\otimes x+f_y\otimes y$, where $f_x,f_y\in A$.  This defines two maps $\pder{}{x},\pder{}{y}\colon \sr(A)\to \sr(A)$, the \emph{partial derivative} maps, given by
$$
\pder{}{x}\colon f\mapsto f_x\bmod[A,A]\qquad\text{and}\qquad\pder{}{y}\colon f\mapsto f_y\bmod[A,A].
$$

Now, the correspondence $A\otimes V=\DR^1(A)\to\Der(A)=A\otimes V^*$ given by the symplectic structure $\omega\in \DR^2(A)$ is nothing more than the canonical map $A\otimes V\to A\otimes V^*$ that is the identity on $A$ and given by the identification $V\simeq V^*$ given by the symplectic structure $\omega_V$ on $V$.  Hence, 
$$
\theta_f=\pder{f}{x}x^*+\pder{f}{y}y^*,
$$
where $x^*,y^*$ is the basis dual to $x,y$ under the correspondence induced by $\omega_V$ (that is, $x^*(v)=\omega(x,v)$ for all $v\in V$, and similarly for $y^*$).  Then
\begin{align*}
&\{f,g\}=\theta_f(dg)\\
&=\left[\pder{f}{x}x^*+\pder{f}{y}y^*\right]\left[\pder{g}{x}x+\pder{g}{y}y\right]
=\pder{f}{x}\pder{g}{y}\omega_V(x,y)+\pder{g}{x}\pder{f}{y}\omega_V(y,x).
\end{align*}
Since $\omega_V=dx\wedge dy$, $\omega_V(x,y)=1$ and $\omega_V(y,x)=-1$.  So we get the familiar formula
$$
\{f,g\}=\pder fx\pder gy-\pder fy\pder gx.
$$
\eex 
The next Proposition gives a non-commutative analogue 
of the  classical Lie algebra exact sequence:

$$
0\too \left[\begin{array}{c}\mbox{\it constant}\\
\mbox{\it functions}\end{array}\right]  \too
 \left[\begin{array}{c}\mbox{\it regular}\\
\mbox{\it functions}\end{array}\right]   \too
 \left[\begin{array}{c}\mbox{\it 
symplectic}\\
\mbox{\it  vector fields}\end{array}\right] \too 0\,,
$$
 associated
with a connected and simply-connected symplectic manifold.

Let $\omega=\sum_i dx_i\otimes dy_i.$ This is a symplectic
 structure on  $A=T(E^*).$
\begin{prop}\label{Lie}
There is a natural Lie algebra
central extension:
$$0\too \C\too  A /[ A , A ]\too \Der_\omega(A)\too 0\,.$$
\end{prop}

\begin{proof} It is immediate from Lemma \ref{der_om}
and Theorem \ref{flat1}
that for the map: $f\mapsto \theta_f$ we have:
$\,\Ker\lbrace{A /[ A , A ]\too \Der_\omega(A)\rbrace}\,$
$=\,\Ker\,d$. Further, by Theorem \ref{flat1}
we get: $\Ker\,d=\C$, and
every closed element in $\DR^1A$ is exact.
This yields surjectivity of the map:
$A /[ A , A ]\to \Der_\omega(A).$
\end{proof}

\section{Kirillov-Kostant Bracket}\label{Kirillov-Kostant}

In this section, we fix a  finite-dimensional Lie algebra
 $\g g$ over $\C$. We are going to define a Poisson bracket
on the polynomial algebra $\C[\g g^*]\cong\Sym(\g g)$.

\subsection{Coordinate formula.}
Fix a basis $e_1,\ldots,e_n$ of $\g g$, and let $c_{ij}^k$ denote the structure constants for this basis.  That is, for all $i,j=1,\ldots,n$, we have that
$\,
[e_i,e_j]=\sum_{k=1}^nc_{ij}^ke_k.
\,$

Since $\g g$ is finite-dimensional, it is isomorphic to its second dual.  So, each $e_i$ gives rise to a linear functional on $\g g^*$, which we denote by $x_i$ (that is, for all $\f\in\g g^*$, $x_i(\f)=\f(e_i)$).  Now, we can identify $\Sym(\g g)$ with the polynomial algebra $\C[\g g^*]$.  Then if $\f,\psi\in\C[\g g^*]$, we have that
$$
\{\f,\psi\}=\sum_{i,j,k=1}^nc_{ij}^k\pder{\f}{x_i}\pder{\psi}{x_j}x_k.
$$
Notice that in this case the Poisson bracket reduces the degree by $1$
(that is, $\deg\{\f,\psi\}=\deg\f+\deg\psi-1$) since two derivatives are
taken and a factor of $x_k$ is multiplied in.  In the Poisson bracket
associated to the Weyl quantization, the degree is reduced by two
because no factor of $x_k$ is introduced.

There is also an explicit 
coordinate free way of writing the  Poisson bracket.
To this end, it is convenient to use the 
$\Sym(\g g)$-realization of our algebra. Specifically,
given two monomials $a=a_1\cdots a_n\in \Sym^n\g g$ and $b=b_1\cdots b_m\in\Sym^m\g g$, we have
\begin{equation}\label{KK_bracket}
\{a,b\}=\sum_{i,j}\,[a_i,b_j]\cdot a_1\cdot\ldots\cdot
\wh a_i\cdot\ldots\cdot a_n\cdot b_1\cdot\ldots\cdot\wh b_j\cdot\ldots\cdot b_m.
\end{equation}

\subsection{Geometric approach.}
Choose two functions $\f,\psi\in\C[\g g^*]$ (or even two smooth function
on $\g g^*$).  Fix some $\lambda\in\g g^*$.  Then $d_\lambda\f$ and
$d_\lambda\psi$ are linear functionals from $T_\lambda\g g^*\simeq\g
g^*$ to $\C$.  Since $\g g$ is finite-dimensional, we identify $\g g$
and $\g g^{**}$.  By abuse of notation, we let $d_\lambda\f$ denote the
element of $\g g$ corresponding to the linear function $d_\lambda\f$ on
$\g g^*$ under this identification.  Then we set
\beq{geom_kir}
\{\f,\psi\}(\lambda)=\langle\lambda,[d_\lambda\f,d_\lambda\psi]\rangle.
\end{equation}
That is, we take the elements $d_\lambda\f$ and $d_\lambda\psi$ of $\g g$ and compute their Lie bracket.  We then evaluate the linear functional on $\g g^*$ associated to this element of $\g g$ on $\lambda$.  

\subsection{Symplectic structure on coadjoint orbits.}
Let $G$ denote any connected Lie group such that $\g g=\Lie(G)$ (in the
future, we will call this a Lie group associated to $\g g$).  Consider
the adjoint action of $G$ on $\g g$.  By transposing, this gives rise to
the coadjoint action on $\g g^*$.  We can then decompose $\g g^*$ into
the disjoint union $\g g^*=\sqcup {\mathbb{O}}$ of $G$-orbits.  

 According to a theorem of Kirillov and
Kostant, every coadjoint orbit ${\mathbb{O}}$
admits a canonical symplectic structure.
Explicitly, for any $\la\in{\mathbb{O}}_k$ we have a natural isomorphism
$T_\la{\mathbb{O}}=\g g/\g g(\la)$,
where $\g g(\la)$ denotes the Lie algebra of the isotropy group
of the point $\la$. Define the pairing
\beq{bra_kir}
\g g/\g g(\la)\times\g g/\g g(\la)\to\C
\quad\text{by}\quad
(x,y)\mapsto\langle\la, [x,y]\rangle,\quad\forall x,y\in\g g.
\end{equation}
\begin{prop}[Kirillov-Kostant]\label{kir}
The pairing above descends to a 
well-defined skew-symmetric 2-form $\om_{\mathbb{O}}$ on
the coadjoint orbit ${\mathbb{O}}$.\qed
\end{prop}

The symplectic form $\omega_{\mathbb{O}}$ gives rise to a Poisson
bracket $\{-,-\}_{\mathbb{O}}$ on the space of functions on the orbit ${\mathbb{O}}$.  
Since $\g g^*$ is the disjoint union of the ${\mathbb{O}}$'s,  the formula
$$
\{\f,\psi\}|_{\mathbb{O}}=\{\f|_{\mathbb{O}},\psi|_{\mathbb{O}}\}\quad
\text{for any coadjoint orbit}\enspace {\mathbb{O}}\sset\g g^*
$$
defines
a  Poisson
structure on the whole of $\g g^*$. It is clear from
\eqref{bra_kir} that this  Poisson 
structure is nothing but the one given by formula
\eqref{bra_kir}.

\begin{examp}
Let $E$ be a finite dimensional vector space
and  $\g g=\End{E}$ the Lie algebra of endomorphisms of $E$
with the commutator bracket.  The trace provides an invariant
bilinear form on $\g g$, and this allows us to identify $\g g$ with its
dual $\g g^*$.  Each conjugacy class ${\mathbb{O}}\sset\g g$ becomes a
coadjoint orbit in $\g g^*$.  For 
$p\in{\mathbb{O}}\sset\g g$, the Lie algebra $\g g(p)$
is nothing but the {\em centralizer} of $p$ in $\g g=\End{E}$.
Thus, Kirillov-Kostant 2-form on the tangent space at $p$ is given by
$$
\om_p:\ \g g/\g g(p)\times\g g/\g g(p)\to\C,\quad
(x,y)\mapsto\Tr(p\cdot[x,y]),\quad\forall x,y\in\g g.
$$
\eex

We now give an interesting example of a {\em noncommutative}
 Kirillov-Kostant structure.
Let $A=\C[e]/(e^2-e)=\C1\oplus\C e$.  
It is then possible to calculate $\ncO^\hdot(A)$ and $\DR(A)$
concretely.  
In particular, we find that
$$
H^j(\DR(A))=\begin{cases}
\C,&\qquad\text{if $j$ is even;}\\
0,&\qquad\text{if $j$ is odd.}
\end{cases}
$$

Let $
\omega=e\,de\,de\in\ncO^2(A).
$ be a
 noncommutative  2-form on $A$.
\begin{lem} The pair $(A,\omega)$ 
 is a noncommutative  symplectic structure.\qed
\end{lem}
This   symplectic structure may be thought of as a
`universal' noncommutative
Kirillov-Kostant  structure.
Indeed, fix a vector space $E$ of dimension $\dim E=n$,
 and let $\Rep^A_E$ denote the variety of all
algebra homomorphisms $A\to \End_\C{E}.$
of $A$.  Each representation is uniquely determined by the image of $e$,
which maps to 
some idempotent, whose image has some rank.  So, we can write
\begin{align*}
\Rep^A_E=\bigsqcup_{k\le n}\{p=p^2\in\End{E}\st\dim(\im p)=k\}
=\bigsqcup_{k\le n}{\mathbb{O}}_k\subset\End{E},
\end{align*}
where each ${\mathbb{O}}_k$ denotes the conjugacy class (under
$GL(E)$) of idempotents 
with rank $k$.
 It is not 
difficult to see that the canonical map 
$\DR^2(A)\too\comO^2(\Rep^A_n)$ sends $e\,de\,de$  to the
ordinary
Kirillov-Kostant 
form (Proposition \ref{kir})  on ${\mathbb{O}}_k\,,\,k=0,1,2,\ldots$.

\subsection{The algebra $\cal O_A$}
Let $A$ be an associative algebra, and as usual let $\sr(A)=A/[A,A]$,
which is a vector space.  It is sometimes  useful to form
the commutative algebra $\cal O_A=\Sym \sr(A)$.  
Consider ${\Rep^A_E}$, the variety of all representations of $A$ on the vector space $V$.  Then we have defined a map $\sr(A)\to\C[{\Rep^A_E}]$, which we denote by $a\mapsto\Tr\wh a$.  Extend this to a map $\cal O_A\to\C[{\Rep^A_E}]$. 

Suppose now that $A$ has a noncommutative symplectic structure $\omega$.
This makes $\sr(A)$ into a Lie algebra, as we have seen before.  
Then $\cal O_A=\Sym \sr(A)$ becomes a Poisson algebra with respect to
the Kirillov-Kostant bracket
\eqref{KK_bracket}.

\subsection{Drinfeld's bracket.}
This section is taken from Drinfeld's paper \cite{Dr}.

Let $\g a$ be a Lie algebra.  Define
$$
\sr(\g a)=\g a\otimes\g a/\C\text{-span of}\langle\{x\otimes y-y\otimes x,[x,y]\otimes z-x\otimes[y,z]\st x,y,z\in\g a\}\rangle.
$$

As the notation suggests, $\sr(\g a)$ is meant to be a Lie analogue of $\sr(A)$ for an associate algebra $A$.  The set we are quotienting out by is meant to reproduce the key properties of the span of the commutators $[A,A]$ in the definition of $\sr(A)$ for $A$ associative.  Indeed, there is a striking similarity embodied in the following observation.  Suppose $\f$ is a linear functional $\sr(\g a)\to\C$.  Then this induces a symmetric invariant bilinear form.  Simply define $\tau_\f\colon\g a\times\g a\to\C$ by $\tau_\f(x,y)=\f(x\otimes y)$.  Then $\tau_\f(x,y)=\tau_\f(y,x)$ and $\tau_\f([x,y],z)=\tau_\f(x,[y,z])$ both follow directly from the definition of $\sr(\g a)$.

This is analogous to the associative case.  Suppose we are given a linear functional $\Tr\colon \sr(A)\to\C$ (here $A$ is associative).  Then this induces a symmetric bilinear form $(a,a')=\Tr(aa')$.  Indeed, since $\sr(A)$ consists of all cyclic words in $A$, $(a,b)=\Tr(ab)=\Tr(ba)=(b,a)$.  It is also clear that $(ab,c)=\Tr((ab)c)=\Tr(a(bc))=(a,bc)$.  So, we see that the second condition in the definition of $\sr(\g a)$ replaces associativity.

Let us apply these considerations to the case where $\g a$ is the free Lie algebra on $m$ generators $x_1,\ldots,x_m$.  Then $\sr(\g a)$ can be realized as pairs of words $(w_1,w_2)$, where each $w_j$ is composed of ``Lie expressions'' in the $x_i$'s.  That is, each $w_j$ is a nested collection of Lie brackets of various generators $x_i$.  

\medskip
\hskip-\parindent\textbf{Claim.} $\sr(\g a)$ has a natural Lie algebra structure.
\medskip

Indeed, the Lie algebra structure is a noncommutative version of the Kirillov-Kostant Poisson bracket on $\Sym(\g g)$, where $\g g$ is a Lie algebra with a nondegenerate, invariant inner product $(-,-)$.  In this context, $\Sym(\g g)$ plays the role of $\sr(\g a)$.  

As mentioned, we let $\g g$ be a Lie algebra with nondegenerate, invariant
inner product $(-,-)$.  This inner product induces an isomorphism
$\g g\times\g g\simeq\g g^*\times\g g^*$, which is the Lie algebra
associated to $(\g g\times\g g)^*$.  A pair of Lie words in $\g g$ then
yield a formula for a Lie word in the Lie algebra $\mathrm{Lie}((\g
g\times\g g)^*)$ via this isomorphism.  This rule determines a formula for
a Poisson bracket $\{-,-\}$ on $\sr(\g a)$.

We now define maps $\pder{}{x_i}\colon \sr(\g a)\to \sr(\g a)$ in the
following way.  Let $f=(w_1,w_2)$ be a pair of Lie words in the basis
$x_n$ of $\g a$--this is the form of an element of $\sr(\g a)$.  Consider
the substitution $x_j\mapsto x_j+z$.  This yields some expression in $z$
and the $x_j$'s, and take the $z$-linear part of it.  Using the properties
of $\sr(\g a)$, we can rewrite this linear part as
$$
(z,w_3),
$$
where $w_3$ is a word in the $x_j$'s.  This word $w_3$ is defined to be
$\pder{f}{x_j}$.

\begin{lem}[Poincar\'e Lemma]
Let $f\in \sr(\g a)$, and $P_1,\ldots,P_n\in \sr(\g a)$.

\noindent
\vi One has $\dis
\sum\nolimits_{j=1}^n\;\left\{x_j\,,\,\pder{f}{x_j}\right\}=0.
$
\vskip 3pt

\noindent
\vii If $
\sum_{j=1}^n\;\{x_j\,,\,P_j\}=0,
$
then there exists $f\in \sr(\g a)$ such that $P_j=\pder{f}{x_j}$.
\end{lem}

\section{Review of (commutative) Chern-Weil Theory}
\subsection{}
Let $\g g $ be a finite-dimensional Lie algebra.  Then the exterior algebra $\Lambda\g g^*$ is equipped with a differential $d$ of degree $+1$.  This differential is called the Koszul-Chevalley-Eilenberg differential, but we will usually shorten the name to the Koszul differential.  The Koszul differential is defined on $\Lambda^1(\g g^*)$ by $d(\lambda)(x,y)=\lambda([x,y])$.  It is then extended to all of $\Lambda(\g g^*)$ by the Leibniz rule.  Consider the Weyl algebra of $\g g$, namely
$$
W(\g g)=\Sym(\g g)^*\otimes\Lambda\g g^*.
$$
Thus  $
W(\g g)$ is a super-commutative
algebra, a tensor product of a commutative algebra ($\Sym(\g g)^*$) with a
super-commutative algebra ($\Lambda\g g^*$).  We wish to equip $W(\g g)$ with a differential.  To accomplish this, consider the graded Lie super-algebra $\tilde{\g g}=\g g_0\oplus\g g_{-1}$, where $\g g_0=\g g_{-1}=\g g$ as vector spaces, and the subscript indicates the degree of each component.  The Lie super-algebra structure is given in the following way.  We set $[x,y]_{\tilde{\g g}}=[x,y]_{\g g_0}$ for all $x,y\in\g g_0$.  We declare $\g g_{-1}$ to be an abelian subalgebra, that is, $[z,w]_{\tilde{\g g}}=0$ for all $z,w\in\g g_{-1}$.  So, we need only define $[x,w]_{\tilde{\g g}}$, where $x\in\g g_0$ and $w\in\g g_{-1}$.  Define a map $\partial\colon\tilde{\g g}\to\tilde{\g g}$ to be zero on $\g g_0$ and the identity isomorphism $\g g_{-1}\to\g g_0$.  Then define $[x,w]_{\tilde{\g g}}=[x,\partial w]_{\g g_0}$.

With $\tilde{\g g}$ defined as above, we find that
$$
\Lambda(\tilde{\g g}^*)=\Sym(\g g)^*\otimes\Lambda\g g^*=W(\g g).
$$
Indeed, this is an isomorphism of graded algebras, if we give $W(\g g)$ the grading
$$
W_n(\g g)=\bigoplus_{n=2p+q}\Sym^p(\g g^*)\otimes\Lambda^q(\g g^*).
$$
So, if $\lambda\colon\g g\to\C$ is a linear functional, its image in $\Sym(\g g)^*$ has degree two, while its image in $\Lambda\g g^*$ has degree one.  We can use the identification $W(\g g)=\Lambda\tilde{\g g}^*$ to define a differential on $W(\g g)$.  Let $\partial\colon\Lambda\tilde{\g g}^*\to\Lambda\tilde{\g g}^*$ denote the extension of the transpose of the map $\partial\colon\tilde{\g g}\to\tilde{\g g}$ defined above, and let $d\colon\Lambda\tilde{\g g}^*\to\Lambda\tilde{\g g}^*$ denote the Koszul differential on $\Lambda\tilde{\g g}^*$.  Then we define $d_W=d+\partial$.

A useful alternative picture of $\tilde{\g g}$ is provided by the
following.  Define $\g g_\varepsilon=\g g\otimes(\C[\e]/(\e^2))=\g g\oplus\e\g g$.
Clearly, $\g g$ corresponds to $\g g_0$ and $\e\g g$ corresponds to $\g
g_{-1}$ in $\tilde{\g g}$.  Define $\partial_\varepsilon\colon\g g_\varepsilon\to\g g_\varepsilon$ by
$\partial_\varepsilon(\e)=1$ and $\partial_\varepsilon(x)=0$ for all $x\in\g g\oplus\{0\}$.  
\begin{rem}
Notice that formally, $\partial_\varepsilon=\frac{\partial}{\partial \varepsilon}$.
\eer
  Clearly, $(\g g_\varepsilon,\partial_\varepsilon)$ is isomorphic to
$(\tilde{\g g},\partial)$, and we will from here on
 identify these two constructions.

\begin{prop}
For the cohomology of the complex $(W(\g g),d_W)$, we have
$$
H^j(W(\g g))=\begin{cases}
\C,&\qquad\text{if $j=0$;}\\
0,&\qquad\text{if $j>0$.}
\end{cases}
$$
\end{prop}
\begin{proof} Since $d_W=d+\partial$ is a sum
of two anti-commuting differentials,
we may view $(W(\g g)\,,\,d+\partial)$ as a bicomplex.
Therefore, we can compute the cohomology
of the total differential via the standard spectral sequence
for a bicomplex. Now, the differential
$\partial$ is induced, by definition,
by the  differential on ${\widetilde{\g g}}$, which is the identity map
$\id: \g g\to\g g,$ by definition. The two term complex
$\g g\to\g g,$ given by this latter map is clearly acyclic.
It follows that the  induced differential $\partial:
W(\g g)\to W(\g g)$ is  acyclic as well.
Hence, the  spectral sequence implies that
the total differential $d_W=d+\partial$
has trivial cohomology in all positive degrees.
\end{proof}

We now introduce a useful piece of notation.
Let $\lambda\in\g g^*$ be a linear functional.  Then we can view
$\lambda$ as both an element of $\Sym(\g g^*)$ and an element of
$\Lambda\g g^*$.  We denote the image of $\lambda$ in $\Sym(\g g^*)$ by
$\lambda_+$, and its image in $\Lambda\g g^*$ is denoted by $\lambda_-$.

Clearly, the elements $\lambda_+$ and $\lambda_-$ generate $\Sym(\g g^*)$ and $\Lambda\g g^*$, respectively, so if we can calculate the action of $d_W$ on them we have a complete description of $d_W$ on $W(\g g^*)$.  Indeed, we find that 
$$
d_W\lambda_-=\lambda_++\lambda_-([-,-])\in\Sym^1\g g^*\oplus\Lambda^2\g g^*=W_2(\g g),
$$
and
$$
d_W\lambda_+=\sum_{j=1}^n\ad^*x_j(\lambda)\otimes x_j^*,
$$
where $\{x_1,\ldots,x_n\}$ is a basis of $\g g$, and $x_j^*$ is the dual basis of $\g g^*$.

Using this characterization of $W(\g g)$ and $d_W$, we are able to deduce the following.

\begin{prop}[Universal property]
If $D$ is a super-commutative DGA, and $\f\colon\g g^*\to D$ is a $\C$-linear map, then there exists a unique map of super-commutative DGA's $\f_W\colon W(\g g)\to D$ such that $\f=\f_W|_{\Lambda^1\g g^*}$.
\end{prop}

\begin{proof}
Extend $\f$ to $\tilde\f\colon\Lambda\g g^*\to D$ as an algebra map.  Notice that $\tilde\f$ will not, in general, commute with the differential.  But it is then possible to extend $\tilde\f$ to $\Sym(\g g^*)$ so that it kills this difference.  This extension $\f_W$ then commutes with differentials, and is the desired extension.  Uniqueness is clear.
\end{proof}

For all $x\in\g g$, we define the \emph{Lie derivative} with respect to $x$, $\L_x\colon W(\g g)\to W(\g g)$ by $L_x=\ad^*x$.  This is a super-derivation of degree zero.  We also define \emph{contraction} $i_x\colon W(\g g)\to W(\g g)$ by $i_x(f)=0$ for all $f\in\Sym(\g g^*)$ and $i_x\alpha$ is simply the contraction of $\alpha$ by $x$ for all $\alpha\in\Lambda(\g g^*)$.  Notice that this is a degree $-1$ map (hence it should be zero on $\Sym(\g g^*)$ since these have only even degrees in $W(\g g)$).  Indeed, $i_x$ is a super-derivation.

\begin{prop}
The Cartan formula holds on $W(\g g)$, that is,
$$
\L_x=d_W\circ i_x+i_x\circ d_W,\quad\text{for all $x\in\g g$.}
$$
\end{prop}

\begin{proof}
This is immediate to verify on the generators of $W(\g g)$.
The result then follows from Lemma \ref{triv}.
\end{proof}

\begin{defn}
An element $u\in W(\g g)$ is \emph{basic} if $\L_xu=i_xu=0$ for all $x\in\g g$.  We let $W(\g g)_\basic$ denote the set of all basic elements of $W(\g g)$
\end{defn}

\begin{lem}
$W(\g g)_\basic=(\Sym(\g g^*))^\g g$, where $\g g$ acts on $\Sym(\g g^*)$ by $\L_x$.
\end{lem}

\begin{proof}
Both $\L_x$ and $i_x$ vanish on $(\Sym(\g g^*))^\g g$ by definition for all $x\in\g g$.  Conversely, $i_x\alpha=0$ for all $x\in\g g$ and some $\alpha\in\Lambda\g g^*$ forces $\alpha=0$, hence $\bigcap_{x\in\g g}\Ker i_x\subset\Sym(\g g^*)$.
\end{proof}

\subsection{Connections on $G$-bundles.}
We will now place the Weyl algebra in a geometric context.  
Suppose we have a principal $G$-bundle 
$$
\xymatrix{
P\ar[r]^G&M},
$$
where $G$, $P$, and $M$ are all connected, $G$ is a Lie group, and the Lie algebra of $G$ is $\g g$.  A \emph{connection} on $P$ is a $\g g$-valued $1$-form $\nabla$ on $P$, i.e., it is an element of $\Om^1(P)\otimes\g g$, satisfying
\begin{itemize}
\item{$\nabla$ is $\g g$-equivariant with respect to the diagonal action on $\Om^1(P)$ and $\g g$, i.e., $\L_x\nabla=0$;}
\item{For each $x\in\g g$, let $\xi_x$ be the vector field on $P$ associated to $x$ by the $G$-action.  Then $\nabla(\xi_x)=x$.}
\end{itemize}

We call an element of $\Om^\hdot(P)$ \emph{basic} if $\L_x\omega=i_x\omega=0$.  We then have the following lemma.

\begin{lem} Assume the group $G$ is connected.
Then, the
 pullback along the bundle map yields a canonical isomorphism $\Om^\hdot(P)_\basic\simeq\Om^\hdot(M)$.
\end{lem}
\begin{proof} It is clear that a differential form on the total space
$P$ descends to a well-defined  differential form on $M$ if and only if
it is $G$-invariant, and annihilates all vectors tangent to the fibers
of the
bundle projection. But if $G$ is connected,  the $G$-invariance of
$\alpha$ 
is equivalent to $\L_x\alpha=0,$ for any  $x\in\g g$.
\end{proof}

We observe further 
that a connection gives rise to a linear map 
$\Phi^\nabla\colon\g g^*\to\Om^1(P)$.  Namely, $\Phi^\nabla(\lambda)=\lambda\circ\nabla$.

By the universal property of $W(\g g)$, $\Phi^\nabla$ extends to a DGA map $\Phi_W^\nabla\colon W(\g g)\to\Om^\hdot(P)$.  By the connection conditions placed on $\nabla$, we see that $\Phi_W^\nabla$ commutes with $\L_x$ and $i_x$.  So,
\begin{align*}
\Phi_W^\nabla(W(\g g)_\basic)\subset\Om^\hdot(P)_\basic\simeq\Om^\hdot(M).
\end{align*}
But $\Phi_W^\nabla(W(\g g)_\basic)=\Phi_W^\nabla((\Sym(\g g^*))^\g g)$,
and 
$(\Sym(\g g^*))^{\g g}\simeq\C[\g g]^{\g g}$.  The map $\Phi_W^\nabla$ is called the \emph{Chern character map}.  We will also denote this map $\Phi_W^\nabla$ by $\ch$.

\begin{defn}
Let $\nabla$ be a connection.  Then the \emph{curvature} of $\nabla$, $K(\nabla)$ is defined by
$$
K(\nabla)=d\nabla+\frac12[\nabla,\nabla]\in\Om^2(P)\otimes\g g.
$$
If $\nabla=\sum_{j=1}^n\nu_j\otimes x_j$, then 
$$
[\nabla,\nabla]=\sum_{j=1}^n\sum_{k=1}^n\nu_j\wedge\nu_k\otimes[x_j,x_k]
\quad\text{and}\quad\frac12[\nabla,\nabla]=\sum_{1\leq j <k\leq n}
\nu_j\wedge\nu_k\otimes[x_j,x_k].$$
\end{defn}

It follows from the definition of curvature and the Chern character that
for all $\lambda\in\g g^*$, $\ch(\lambda_+)=\lambda\circ K(\nabla)$.

Next, we consider
 the ideal $(\Sym^1(\g g^*)\otimes1)W(\g g)$ of $W(\g g)$.  For simplicity, we let $\g g^*_+=\Sym^1\g g^*\otimes1\subset W(\g g)$.  Since $d_W(\g g^*_+)\subset g^*_+W(\g g)$, we see that $\g g^*_+W(\g g)$ is a differential ideal of $W(\g g)$.

\begin{defn}
The \emph{Hodge filtration} of $W(\g g)$ is a {\sl decreasing}
filtration $W(\g g)\supset g^*_+W(\g g)\supset \ldots
\supset F^{p-1}W(\g g)\supset 
F^pW(\g g)$ given by
$$
F^pW(\g g)=(\g g^*_+)^pW(\g g).
$$
\end{defn}

Notice that
  $W(\g g)/F^1W(\g g)=\Lambda\g g^*\subset W(\g g)$.

\begin{lem}
For each $p$,
$$
H^{2p}(F^pW(\g g))=(\Sym^p\g g^*)^{\g g}.
$$
\end{lem}

\subsection{Transgression map.}
The short exact sequence
$$
\xymatrix{0\ar[r]&F^1W(\g g)\ar[r]&W(\g g)\ar[r]&W(\g g)/F^1W(\g g)\ar[r]&0}
$$
gives rise to a long exact sequence in cohomology.  Recall that the cohomology of $W(\g g)$ is the cohomology of a point (i.e., $\C$ in degree zero, and zero elsewhere).  The inclusion of $F^pW(\g g)\subset F^1W(\g g)$ yields a homomorphism $H^{2p}(F^pW(\g g))\to H^{2p}(F^1W(\g g))$.  Recalling that $H^{2p}(F^pW(\g g))=(\Sym^p\g g^*)^\g g$, we obtain the exact sequence
$$
\xymatrix{(\Sym^p\g g^*)^{\g g}\to H^{2p}(F^1W(\g g))\to H^{2p-1}(W(\g g)/F^1W(\g g))=H^{2p-1}(\Lambda\g g^*,d)},
$$
where $d$ denotes the Koszul differential.

But then $H^{2p-1}(\Lambda\g g^*,d)$ is isomorphic to
 the $2p-1$-Lie cohomology of $\g g$, $H^{2p-1}_{\mathrm{Lie}}(\g g)$, 
which is isomorphic to $H^{2p-1}(G)$.  It is well know that the cohomology of $G$ can be calculated using only invariant differential forms on $G$.

A geometric meaning behind this can be found by considering the universal $G$-bundle $EG\to BG$.  Then $EG$ is contractible, and $W(\g g)\simeq\Om^\hdot(EG)$ and $W_\basic(\g g)\simeq\Om^\hdot(BG)$.  Then the $2p$-cohomology of $F^1W(\g g)$ is calculating $H^{2p}(BG)$, so we obtain a homomorphism $H^{2p}(BG)\to H^{2p-1}(G)$.

\subsection{Chern-Simons formalism.}
Let $D$ be a DGA, with differential $d\colon D^n\to D^{n+1}$.  We set $\sr(D)=D/[D,D]$, where $[-,-]$ here denotes the super-commutator, i.e.,
$$
[x,y]=xy-(-1)^{(\deg x)(\deg y)}yx.
$$
Then $d$ descends to a super-differential $d$ on $\sr(D)$.  Fix any $a\in D^1$.  We define its \emph{curvature} $F:=da+a^2$.

\begin{prop}[Bianchi Identity]
With $D$, $a$, and $F$ as above, the following identity holds in $D$:
$$
(d+\ad a)F=0.
$$
\end{prop}

\noindent
{\em Proof.\;}
Observe that $dF=d(da+a^2)=d^2-a\,da+da\cdot a=a\,da+da\cdot a$.  Also,
\begin{align*}
(\ad a)F=[a,F]&=aF-(-1)^{(\deg a)(\deg F)}Fa\\
&=aF-Fa=a(da+a^2)-(da+a^2)a\\
&=a\,da-da\cdot a.\qquad\Box
\end{align*}

Using the Bianchi identity, we see that in $\sr(D)$ we have
$$
d(F^n)=-[a,F^n]=0,
$$
since all super-commutators are zero in $\sr(D)$.  Hence the elements $\frac{F^n}{n!}$ are cocycles in $\sr(D)$.

Consider the algebra $D[t]=D\otimes\C[t]$.  Take an element $a_t\in D^1[t]=D^1\otimes\C[t]$.  Define $F_t=da_t+a_t^2$, the curvature of $a_t$.  Then a simple computation shows that
$$
\pder{}{t}\left(\frac{F_t^n}{n!}\right)=d\left[\frac1{(n-1)!}\pder{a_t}{t}F_t^{n-1}\right].
$$

In particular, take $a_t=ta$.  Then $F_t=t\,da+t^2a^2$.  Then we can integrate the above  equation in $t$.  Indeed, this yields that 
$$
\frac{F^n}{n!}=d\cs_{2n-1}(a),
$$
where
$$
\cs_{2n-1}(a)=\int_0^1a\frac{F_t^{n-1}}{(n-1)!}\,dt.
$$
The element $\cs_{2n-1}(a)$ is called the $2n-1$-\emph{Chern-Simons class} of $a$.  The map $\frac{F^n}{n!}\mapsto\cs_{2n-1}(a)$ is a trangression map.

\subsection{Special case: $\g g=\g{gl}_n$.}
Let $G=\oper{GL}_n$. Then, giving a principal $G$-bundle
$P\to M$ is the same thing as giving
an ordinary vector bundle on $M$.
Let  $\nabla$ be a  connection on $P$.  
Then we define $\ch_k=\frac1{k!}\Tr(K(\nabla)^k)$, which is
$$
\ch_k=\frac1{k!}\Tr((d\nabla+\frac12[\nabla,\nabla])^k).
$$

Now, suppose that $\nabla_0,\ldots,\nabla_N$ are $N+1$ connections on the same bundle $P$.  We wish to show that $\ch_k$ is independent of the connection used.  This will follow essentially from the fact that the space of connections is convex.  Let $\Delta$ be the standard $N$-simplex in $\R^{N+1}$, that is,
$$
\Delta=\{(t_0,\ldots,t_N)\in\R^{N+1}\st t_j\ge0, \sum_{j=0}^nt_j=1\}.
$$
For each $t\in\Delta$, define $\nabla(t)=\sum_{j=0}^Nt_j\nabla_j$.  Now, it is a well-known fact that a connection $\nabla_j$ may be written in a local trivialization of $P$ as a sum of the usual differential and a matrix of $1$-forms, $\nabla_j=d+A_j$.  Define
$$
\ch_m^N(\nabla_0,\ldots,\nabla_N)=\int_\Delta\Tr(K(\nabla(t))^k)\,dt.
$$
This is a complicated expression in $A_0,\ldots,A_N$ and
$dA_0,\ldots,dA_n$.

\subsection{Quantized Weil algebra.}
Let $(\g g,B)$ be a Lie algebra 
equipped with an invariant, nondegenerate bilinear form $B
\colon\g g\times\g g\to\C$.  Then this form induces a canonical isomorphism $\g g\simeq\g g^*$.  In the standard Chern-Weil theory, we set
$$
W(\g g)=\Sym\g g\otimes\Lambda\g g.
$$
(The above formula is correct, since $\g g\simeq\g g^*$.)

Following  Alekseev-Meinrenken \cite{AM}, we would like to
`quantize' the algebra $W(\g g)$ by
replacing
 $\Sym\g g$ by the universal enveloping algebra and $\Lambda\g g$ 
by the Clifford algebra ${\clif}\g g$.  Recall that 
$$
{\clif}\g g=T\g g/(xy+yx-2 B(x,y)).
$$

We set
$$
\scr W(\g g)=U\g g\otimes{\clif}\g g.
$$
It turns out that $\scr W(\g g)$ has a differential.  It is this algebra that acts as the quantized  version of $W(\g g)$.

This quantized Weyl algebra is connected to the work of  Alekseev-Meinrenken
in the following fashion.  Consider the family of Lie algebras $\g g_t$,
where $\g g_t=\g g$ as a vector space and $[x,y]_t=t[x,y]$ for all $x,y\in\g g_t$.  Define
$$
\scr W_t(\g g)=\scr W(\g g_t)=U\g g_t\otimes{\clif}(\g g_t).
$$
Then $\scr W_t(\g g)$ is a flat family of DGA's, and $\scr W_0(\g
g)\simeq W(\g g)$.  Based on our previous work with deformations of
associative algebras, 
we conclude that the family $\scr W_t(\g g)$ induces a Poisson bracket 
on $W(\g g)$ making it a Poisson DGA.  

Let $\Delta\in\Sym^2(\g g)\subset W(\g g)$ be the
 canonical element corresponding to the inverse of the nondegenerate
bilinear form $B$ on $\g g$.

Using the deformation argument above, one proves

\begin{thm}\label{W_bracket}
An invariant, nondegenerate bilinear form on  $\g g$
gives rise to a Poisson (super) algebra structure $\{-,-\}$
on $W(\g g)$ such that
 $$d_W(u)=\frac12\{\Delta,u\},\quad\forall u\in W(\g g).$$
\end{thm}

\section{Noncommutative Chern-Weil theory}
\subsection{}
In the previous section, we discussed quantized Chern-Weil theory, which
could be considered a part of noncommutative geometry ``in the small.''
That is, it is simply a deformation of the usual Chern-Weil theory.  
We now want to begin with a noncommutative algebra which will replace the Lie algebra $\g g$.

Let $A$ be a (possibly noncommutative) associative algebra.  Let $A^*$ denote the dual of $A$.  This is a coalgebra, with comultiplication map $\Delta\colon A^*\to A^*\otimes A^*$.  For the sake of simplicity, we will use the Sweedler notation, that is, we write
$$
\Delta(\lambda)=\sum\lambda'\otimes\lambda''
$$
for all $\lambda\in A^*$.  

Following the paper \cite{C},
we define $\ncW(A)=T(A^*_+\oplus A^*_-)$, where $A_+^*=A^*_-=A^*$.  We
make $\ncW(A)$ graded  algebra (under the usual multiplication for the tensor algebra)
by taking 
$$\ncW(A)_p:=\bigoplus\nolimits_{2n+m=p}\,(A^*_+)^{\otimes
n}\otimes(A^*_-)^{\otimes m}.$$

Define a differential $d_W$ on $\ncW(A)$ by
\begin{align*}
d_W\lambda_-=\lambda_++\sum\lambda_-'\otimes\lambda_-''\quad\text{and}\quad
d_W\lambda_+=\sum(\lambda_+'\otimes\lambda_-''-\lambda_-'\otimes\lambda_+'').
\end{align*}
As usual, $\lambda$ is a linear functional in $A^*$, and $\lambda_+$ (respectively, $\lambda_-$) represents its image in $A^*_+$ (respectively $A^*_-$).  This differential makes $\ncW(A)$ a DGA.

Similarly to the construction
of Bar-complex as a free product, it is sometimes useful to
have the following alternative definition
\beq{W_eps}
\ncW(A)=T\bigl((A_\eps)^*\bigr),\quad\text{where}\quad A_\eps:=A\otimes \C[\e]/(\e^2).
\end{equation}
Let $\partial\colon\ncW(A)\to\ncW(A)$ be the $\C$-linear map 
sending $\e\mapsto1$ and $1_A\mapsto0$.  Define
$$
d_W=\partial+\delta,
$$
where $\delta$ is essentially the differential dual to the
multiplication map, as before.
The two differentials $\partial$ and $\delta$ anti-commute,
hence $d_W\ccirc d_W=0$. 
Observe that, in the presentation \eqref{W_eps},
the differential $\partial$ can be suggestively
written as $\partial=\partial/\partial\e$.

As in the commutative and quantized cases, there is a Poincar\'e lemma for $\ncW(A)$. 

\begin{lem}
The cohomology of the complex $(\ncW(A),d_W)$ are given by
$$
H^j(\ncW(A))=\begin{cases}
\C,&\qquad\text{if $j=0$;}\\
0,&\qquad\text{if $j>0$.}
\end{cases}
$$
\end{lem}

\begin{proof}
Suppose that the multiplication on $A$ is trivial.  Then it is easy to see that $H^j(\ncW(A))$ is zero for $j>1$, and $\C$ for $j=0$.  This follows, since the differential $d_W$ involved the transpose $\delta$ of the multiplication map, which would also be zero.

We construct an isomorphism of DGA's from $(\ncW(A),d_W)$ to $(\ncW(A),\tilde d)$, where $\tilde d$ is the differential given when the multiplication is trivial.  Indeed, observe that $\ncW(A)$ is independent of the multiplication on $A$.  The map is given by
$$(\ncW(A)\,,\,d_W)\too(\ncW(A)\,,\,\tilde d)\,,\quad
\lambda_-\mapsto\lambda_-\,,\,\lambda_+\mapsto\lambda_++\sum\lambda'_-\otimes\lambda_-''.
$$
It is easy to check this is a DGA isomorphism.

An alternative proof of the Lemma can be obtained as follows.
It is clear that the  two term complex
 $\partial/\partial_\e: \e A^* \to A^*$
has trivial homology. Thus, the spectral sequence
associated to the bicomplex given by the
differentials  $\partial$ and $\delta$ implies
the result.\end{proof}

\begin{lem}[Universal Property]
Given a DGA $D$ and a $\C$-linear map $\f\colon A^*\to D$, 
there is a unique extension $\f_W\colon\ncW(A)\to D$ such that $\f_W|_{A^*}=\f$.
\end{lem}

As usual, we set $R\ncW(A)=\ncW(A)/[\ncW(A),\ncW(A)]$, where the
commutators are graded.  This still has trivial cohomology (i.e., the
cohomology of a point).  
The same calculation as before suffices--we simply take $A$ to have the trivial product.

\begin{thm}[\cite{C}]
There exists a canonical transgression map
$$
\xymatrix{H^j(\sr(F^p\ncW(A)))\ar[r]^{\sim\hphantom{stuff}}&H^{j-1}(\sr(\ncW(A)/F^p\ncW(A)))}.
$$
\end{thm}

\begin{defn}
Define $I_+=\ncW(A)A_+^*\ncW(A)\subset\ncW(A)$.  Then this is a differential ideal of $\ncW(A)$, i.e., $dI_+\subset I_+$.  The \emph{Hodge filtration} of $\ncW(A)$ is given by $F^p\ncW(A)=I_+^p$.
\end{defn}

\subsection{Example:  case $A=\C$.}
In this case, $\ncW(A)$ is the free algebra with generators $\alpha$ in degree one and $\partial\alpha$ in degree two.  Extend $\partial$ by the Leibniz rule to be a super-derivation of $T(\alpha,\partial\alpha)$ such that $\partial^2=0$.  Let $F=\partial\alpha+\alpha^2$.  Then one finds that $\partial F=-[\alpha,F]$.  On $\ncW(A)$, we have the Hodge filtration $F^p\ncW(A)$, the two-sided ideal generated by terms involving $\partial\alpha$ at least $p$ times.  As before, one sees that $\partial(F^p\ncW(A))\subset\ncW(A)$.  Finally, we consider the algebra $\sr(\ncW(A)/F^p\ncW(A))$.

\begin{thm}
Let $H^\hdot(\sr(\ncW(A)/F^p\ncW(A)))$ denote the cohomology of the 
complex\break
 $\sr(\ncW(A)/F^p\ncW(A))$ with differential
$d_W=\partial+\delta$ 
(where $\delta$ is the differential dual to the multiplication).  Then
$H^\hdot(\sr(\ncW(A)/F^p\ncW(A)))$ has a $\C$-basis formed by 
the Chern-Simons classes $\cs_{2n-1}$ for $n\ge p$.
\end{thm}

\begin{proof}
A $\C$-basis of $\sr(\ncW(A)/F^p\ncW(A))$ is given by $\alpha^{2k-1}$, $k\ge1$, and $(\partial\alpha)^\ell$, $1\le\ell<p$.  Calculating the $\delta$-cohomology, we see that $a^{2k-1}$ transgresses to $(\partial\alpha)^k$ for $1\le k\le p-1$.  The remaining $\alpha^{2k-1}$'s are sent to the Chern-Simons classes, $\cs_{2k-1}$ ($k\ge p$).
\end{proof}

\subsection{Gelfand-Smirnov bracket.} We fix
 a finite-dimensional
associative algebra $A$ equipped with an invariant trace $\Tr\colon A\to
\C$, that is, $\Tr(a_1a_2)=\Tr(a_2a_1)$ for all $a_1,a_2\in A$.
Write $\Delta\in \Sym^2(A^*)$ for the canonical
element corresponding to the bilinear form,
i.e., such that $B(x,y)=\langle\Delta, x\otimes y\rangle,$
for any $x,y\in A$.

The following result is a noncommutative analogue of
Theorem \ref{W_bracket}

\begin{thm}\label{ncW_bracket}
 If the trace pairing $a_1\times a_2\mto  \Tr(a_1a_2)$ is
non-degenerate,
then the graded vector space
$\sr(\ncW(A))$ has a canonical Lie super-algebra structure
 such that
 $$d_W(u)=\frac12\{\Delta,u\},\quad\forall u\in \sr(\ncW(A)).$$
\end{thm} 

\begin{proof}[Hint of Proof:] The trace pairing on $A$
induces an isomorphism of vector spaces
$A\iso A^*,\,a\mapsto \Tr(a\cdot(-)).$ Hence,
we get an isomorphism $\kappa:\  A\oplus A \iso A^*\oplus A^*$.
Further, the trace pairing also gives rise to a
non-degenerate skew-symmetric $\C$-bilinear 2-form $\om$ on the vector space $A\oplus A$
defined by the formula
$$(a\oplus a')\,\times\, (b\oplus b')\mto \Tr(ab')-\Tr(a'b).$$
Transporting this 2-form from $A\oplus A$ to $A^*\oplus A^*$
via the isomorphism $\kappa$ makes $A^*\oplus A^*$ a symplectic
vector space. Therefore, the tensor algebra
$T(A^*\oplus A^*)$ acquires a natural structure of noncommutative
symplectic manifold, see Example \ref{symp_free}.
Thus, we get a Lie bracket on $\sr\bigl(T(A^*\oplus A^*)\bigr)$.

Now, the algebra $\ncW(A))$ is just $T(A^*\oplus A^*)$,
as an associative algebra. So, the above construction
can be adapted, by inserting suitable signs (due to the
fact that the first copy  $A^*\sset \ncW(A)$ is placed in degree 1
and the second  copy is placed in degree 2),
to produce the required Lie super-algebra structure on 
$\sr\bigl(\ncW(A)\bigr)$.
\end{proof}

Gelfand and Smirnov considered in \cite{GeSm}
a very special case of this situation
where  $A=\C\oplus\cdots\oplus\C$ ($n$ copies)
is a semisimple algebra equiped with the natural
trace $\Tr: (x_1\oplus\cdots\oplus x_n) \mto \sum_i x_i$.

It is clear that for 
$A=\C\oplus\cdots\oplus\C$, one has
$$\ncW(A)=\C\langle a_1,\ldots,a_n,b_1,\ldots,b_n\rangle,
$$
is  the free graded algebra 
 on $2n$ generators, $a_1,\ldots,a_n$ and $b_1,\ldots,b_n$ where $\deg
a_j=1$ and $\deg b_j=2$.

The differential $d_W$ is  
 a super-derivation such that on generators we have
$d_W(a_j)=b_j$, and $d_W(b)=0.$

Set $\sr^\hdot(\ncW(A))=\ncW(A)/[\ncW(A),\ncW(A)]$, where
$[\ncW(A),\ncW(A)]$ is the linear span of all super-commutators of
elements of $\ncW(A)$, and the grading on
the RHS is induced from that on $ \ncW(A)$.
 We can view $\sr(\ncW(A))$ as all cyclic words in
$a_j$'s and $b_j$'s (where $\deg a_j=1$ and $\deg b_j=2$).

Define maps $\pder{}{a_j}\colon \sr(\ncW(A))\to \sr(\ncW(A))$ in the following way.  Choose any word $x_1\cdots x_k$ in $\ncW(A)$.  Set $\delta_{a_j}(x_1\cdots x_k)$ to be $x_2\cdots x_k$ if $x_1=a_j$ and zero otherwise.  Then $\pder{}{a_j}(x_1\cdots x_k\bmod[\ncW(A),\ncW(A)])$ is defined to be the sum of $\delta_{a_j}$ applied to all cyclic permutations of $x_1\cdots x_n$.  Maps $\pder{}{b_j}\colon \sr(\ncW(A))\to \sr(\ncW(A))$ are defined analogously.  Then set
$$
\{P,Q\}=\sum_{j=1}^n\pder{P}{a_j}\pder{Q}{b_j}+(-1)^{(\deg P)(\deg Q)}\pder{Q}{a_j}\pder{P}{b_j}.
$$

\begin{prop} 
\vi The above map 
$$\{-,-\}\colon \sr^p(\ncW(A))\times \sr^q(\ncW(A))\too \sr^{p+q-3}(\ncW(A))$$ 
is
well-defined, and gives a Lie super-algebra structure on $\sr(\ncW(A))$.

\vii  The above bracket is equal to the one arising from Theorem \ref{ncW_bracket} 
in the special case $A=\C\oplus\cdots\oplus\C$. in particular,

\viii
The differential $d$ on $\ncW(A)$ descends to a well-defined Lie
(super)algebra super-derivation of $\sr(\ncW(A))$ such that,
for all $P\in \sr(\ncW(A))$, one has
$$
dP=\{b_1^2+\cdots+b_n^2,P\}.
$$
\end{prop}

Assume, for simplicity, that  $n=1$, 
so $\ncW(A)=\C\langle a,b\rangle$, where $b=da$.  
If $a$ were giving a connection on a bundle, then 
we would consider the quantity $da+a^2=b+a^2$, which is the curvature of the connection.  This motivates the following definition

\begin{defn}
For $\ncW(A)=\C\langle a,b\rangle$ as above, define for each $k=0,1\ldots$
$$
\ch_k=(a^2+b)^k\in \sr(\ncW(A)).
$$
\end{defn}

Then we have the following lemma.

\begin{lem}
The following identities hold in $\sr(\ncW(A))$.
\begin{enumerate}
\item[(a)]{$d(\ch_k)=0$.}
\item[(b)]{$\{\ch_k,\ch_l\}=0$ for all $k,l\in\N\cup\{0\}$.}
\end{enumerate}
\end{lem}

By the Poincar\'e lemma for $\sr(\ncW(A))$, 
 every closed element of $\sr(\ncW(A))$  of degree $>0$ is exact.  Hence, there is some
element $\ch^1_k\in \sr(\ncW(A))$ such that $d(\ch^1_k)=\ch_k$.  Indeed, an
explicit formula for $\ch_k^1$ can be found.  To simplify its
expression, we introduce the following notation.  For any $x,y\in \ncW(A)$ and
$k,l\in\N\cup\{0\}$, we set $\sigma_{k,l}(x,y)$ to be (the image under
the quotient of) the sum of all noncommutative 
words in elements of $\ncW(A)$ that contains exactly $k$ symbols $x$ and $l$
symbols 
$y$.  Then we have the following result.

Following \cite{GeSm}, for each $k$, set
\begin{align*}
\ch_k^1:=\frac{a}{(k-1)!}\cdot\Bigg[\frac{b^{k-1}}{k}
+\frac{\sigma_{1,k-2}(a^2,b)}{k-1}+\frac{\sigma_{2,k-3}(a^2,b)}{k+2}+\ldots
+\frac{a^{2(k-1)}}{2k-1}\Bigg].
\end{align*}
\begin{prop}
For any $k$,  in $\sr(\ncW(A))$ one has $d(\ch_k^1)=\ch_k$.
\end{prop}

\section{Chern Character on $K$-theory}
\subsection{Infinite matrices.} Fix an  associative algebra $A$.
For each integer $n\geq 1$ we have the algebra $\Mat_nA$ of
$n\times n$-matrices with entries in $A$. The
assignment
$$\Mat_nA\too\Mat_{n+1}A,\quad
x\mto\begin{pmatrix}x&0\\0&0\end{pmatrix}
$$
gives an algebra imbedding (note that the unit $1_n\in\Mat_nA$ does {\em not} 
go to $1_{n+1}\in \Mat_{n+1}A$).
We let $\Mat_\infty(A):=\underset{n\to\infty}\lim \Mat_nA$ denote the
corresponding direct limit under the ``upper left hand corner''
inclusions.
 Thus $\Mat_\infty(A)$ is an associative
algebra
without unit that can be identified with the algebra of infinite
matrices with finitely many nonzero entries.

Further, let $\GL_n(A)\sset \Mat_nA$ be the group of {\em invertible} 
$n\times n$-matrices with entries in $A$.
The map
$$g\mto \begin{pmatrix}g&0\\0&1\end{pmatrix}
$$
gives a group imbedding $\GL_n(A)\into\GL_{n+1}(A)$
(note that this time the map does take the unit into unit).
We let $\GL_\infty(A):=\underset{n\to\infty}\lim \GL_n(A)$ denote the
corresponding direct limit.  Thus $\GL_\infty(A)$ is a group that can be
identified with the 
group  of infinite
matrices $g=\|g_{ij}\|$ such that the matrix
$g-\Id=\|g_{ij}-\delta_{ij}\|$
has only  finitely many nonzero entries (here $\delta_{ij}$ denotes the
Kronecker
delta).

For $i,j\in\{1,\ldots,n\}$ and $a\in A$, we let $E_{ij}(a)$ 
denote the elementary $n\times n$ matrix with $a$ in the $ij$-position and zero elsewhere.

Also, for any group $G$, let $[G,G]\sset G$ denote the (normal)
subgroup 
 generated by the elements
$ghg^{-1}h^{-1}\,,\,g,h\in G$.  We will denote $ghg^{-1}h^{-1}$ by $[[g,h]]$.  
\begin{lem} The group
$[GL_\infty(A),GL_\infty(A)]$ is generated by matrices of the form $E_{ij}(a)$ where $i\ne j$.
\end{lem}

\begin{proof}
First, observe that
$$
[[E_{ij}(a),E_{k\ell}(b)]]=\begin{cases}
1,&\qquad\text{if $j\ne k$, $i\ne\ell$;}\\
E_{i\ell}(ab),&\qquad\text{if $j=k$, $i\ne\ell$;}\\
E_{kj}(-ba),&\qquad\text{if $j\ne k$, $i=\ell$.}
\end{cases}
$$
It is an easy computation to see that
$$
\begin{pmatrix}X&0\\0&X^{-1}\end{pmatrix}=\begin{pmatrix}1&X\\0&1\end{pmatrix}\begin{pmatrix}1&0\\-X&1\end{pmatrix}\begin{pmatrix}0&-1\\1&0\end{pmatrix}.
$$
Hence $\begin{pmatrix}X&0\\0&X^{-1}\end{pmatrix}$ is a product of elementary matrices
$$
\begin{pmatrix}[[Y,Z]]&0\\0&1\end{pmatrix}=\begin{pmatrix}Y&0\\0&Y^{-1}\end{pmatrix}\begin{pmatrix}Z&0\\0&Z^{-1}\end{pmatrix}\begin{pmatrix}(ZY)^{-1}&0\\0&ZY\end{pmatrix}.
$$
Thus, $[[Y,Z]]$ is a product of elementary matrices, hence any element
in $[GL_\infty(A),$
$GL_\infty(A)]$ may be written as a product of elementary matrices.
\end{proof}
\subsection{The group $K^0(A)$.} Fix an associative algebra $A$.
Recall that  $K^0(A)$ is defined to be an abelian group
which is a quotient of the free  abelian group
$\Z$-generated by the isomorphism classes $[P]$ of all
finite rank projective (left) $A$-modules $P$ modulo
the subgroup generated by the following relations:
$$[P]+[Q]-[P\oplus Q].$$
In other words,  $K^0(A)$ is the Grothendieck group of (the
exact category of) finite rank projective $A$-modules,
equipped with a  semigroup structure by direct sum.

Each finite rank projective $A$-module $P$ is a direct summand of a free
$A$-module $A^n$, that is, one has
an $A$-module direct sum decomposition
$A^n=P\oplus Q$. The projection to $P$ along $Q$ gives
a map
$A^n=P\oplus Q\map P\into A^n,$ 
which is given by an $n\times n$-matrix, i.e., by 
 an element $e\in \Mat_nA$. It is clear that $e$ is an
idempotent,
i.e., $e^2=e$. 
 Direct sum of modules corresponds to direct sum of idempotents, where
$$
e\oplus e':=\begin{pmatrix}e&0\\0&e'\end{pmatrix}.
$$

Let $\scr P(A)$ denote the set of idempotents in $\Mat_\infty(A)$.
We view any $n\times n$-matrix as being imbedded into $\Mat_\infty(A)$,
thus any idempotent $e\in \Mat_nA$ becomes an idempotent in
 $\Mat_\infty(A)$.
\begin{lem}\label{conj_invar} Two idempotents $e, e'\in\Mat_\infty(A)$  give rise to isomorphic
 projective $A$-modules if and only if $e'=geg^{-1}$ for some $g\in
\GL_\infty(A).$
\end{lem}
\begin{proof}
Let $e\in\Mat_nA$ and $e'\in\Mat_mA$ be two idempotents
such that one has an $A$-module isomorphism
$A^ne\simeq A^me'$. Put $P_1=A^ne,\,Q_1=A^n(1-e),$
and $P_2=A^me',\,Q_2=A^m(1-e').$ Thus, 
$A^n=P_1\oplus Q_1,\, A^m=P_2\oplus Q_2,$ and
we are given an isomorphism
$\varphi: P_1\iso P_2.$ We consider the following chain
of isomorphisms
\begin{align*}
A^{n+m}= &A^n\oplus A^m\iso P_1\oplus Q_1\oplus P_2\oplus Q_2
\xymatrix{{}\ar[rrr]_<>(0.5){\sim}^<>(0.5){\varphi\oplus
\Id_{Q_1}\oplus\varphi\inv\oplus\Id_{Q_2}}&&&{}}\\
&P_2\oplus Q_1\oplus P_1\oplus Q_2\xymatrix{\ar[r]_<>(0.5){\sim}^<>(0.5){P_1
\leftrightarrow P_2}&{}}P_1\oplus Q_1\oplus
P_2\oplus Q_2\iso
A^n\oplus A^m=A^{n+m}.
\end{align*}
The composite map $A^{n+m}\iso A^{n+m}$ is an isomorphism,
hence, it is given by an invertible matrix $g\in\Mat_{n+m}A$.
It is clear that $g\inv(e\oplus 0_m)g=0_n\oplus e',$
and we are done.
\end{proof}

We introduce an equivalence relation  $e\sim e'$ on $\scr P(A)$ by
$e'=geg^{-1}$ for some $g\in
\GL_\infty(A)$, 
and denote by $[e]$ the equivalence class of $e$.  We define a 
semigroup structure on the equivalence classes of idempotents 
by $[e]+[e']=[e\oplus e']$ and $[0]=0$.  
This way, one  can rephrase the definition of
 $K^0(A)$ in terms of idempotents as follows:
$$ K^0(A) \cong \bigl(\scr P(A)/\!_\sim\,,\,\oplus\bigr).$$

\subsection{Chern class on $K^0$ and $K^1$.}
Recall that for any associative algebra $A$
 there is a trace map $\tr\colon\Mat_nA\to A/[A,A]$ given by
$$
\tr(a_{ij})=\sum_{i=1}^na_{ii}\bmod[A,A].
$$
Since $\tr(xy)=\tr(yx)$, we see that if $e$ and $e'$ 
are equivalent idempotents (suppose $e'=geg^{-1}$), then
$$
\tr(e')=\tr(geg^{-1})=\tr(e).
$$
So, $\tr$ descends to a well defined map $\tr\colon\scr P(A)/\!_\sim\,\to
A/[A,A]$.  Notice that $A/[A,A]=\DR^0(A)$.  This map is additive.

\begin{prop} The assignment $[e]\mto  \tr(e)$
extends to a group homomorphism:
$$\cch_0\colon K^0(A)\too \Ker\left[\DR^0(A)\to \DR^1(A)\right],$$ 
called the \emph{Chern
character}.
\end{prop}
\begin{proof}
Additivity of the map is clear.
So, we must show that for all $[e]\in K^0(A)$, $\cch_0[e]$
is a closed form.  Indeed, choose some representative $e\in\Mat_nA$.
Then $e^2=e$.  Applying $d$ to both sides yields $e\,de+(de)e=de$. This
way, one proves
$$
e\,de=de(1-e)\quad\text{and}\quad (de)e=(1-e)\,de.
$$
Therefore, we calculate
\begin{align*}
\tr(e\,de)&=\tr(e^2\,de)=\tr(e(de)(1-e))=\tr((1-e)e(de))=0,
\end{align*}
since $(1-e)e=e-e^2=e-e=0$.  Similarly, we see that $\tr((1-e)\,de)=0$.  So,
$$
\tr(de)=\tr(e\,de)+\tr((1-e)\,de)=0.
$$
Clearly, $\tr$ and $d$ commute, so $\cch_0([e])$ is closed.
\end{proof}

Next, we define 
$$
K^1(A):=GL_\infty(A)/[GL_\infty(A),GL_\infty(A)].
$$ 

We are going to construct a Chern character for $K^1(A)$.
\begin{prop} There is a natural group homomorphism
$$
\cch_1\colon K^1(A)\to\Ker\left[\DR^1(A)\stackrel{b}\longrightarrow A\right].
$$
\end{prop}

\begin{proof}
Choose any $[g]\in K^1(A)$.  Choose some representative $g$ of $[g]$, and define
$$
\cch_1[g]=\tr(g^{-1}\,dg).
$$
First, let us check that this is a group homomorphism.  Indeed
\begin{align*}
\tr((g_1g_2)^{-1}d(g_1g_2))&=\tr[(g_1g_2)^{-1}dg_1\cdot g_2]+\tr[g_2^{-1}g_1^{-1}g_1\,dg_2]\\
&=\tr(g_2^{-1}g_1^{-1}\,dg_1\cdot g_2)+\tr[g_2^{-1}\,dg_2]\\
&=\tr[g_1^{-1}\,dg_1\cdot g_2g_2^{-1}]+\tr(g_2^{-1}\,dg_2)\bmod[GL_\infty(A),GL_\infty(A)]\\
&=\tr[g_1^{-1}\,dg_1]+\tr[g_2^{-1}\,dg_2],
\end{align*}
as desired.  Again, $\cch_1[g]$ is a $b$-cycle, since
$$
b\left[\tr(g^{-1}\,dg)\right]=\tr(g^{-1}g)-\tr(gg^{-1})=0.
$$
\end{proof}
\subsection{Chern classes via connections.} Given a finite rank
projective
(left) module $M$ over an associative algebra $A$, one can associate
to $M$ its de Rham characteristic classes $\ch_k(M)\in\DR^{2k}(A),
\,k=1,2,\ldots,$ as follows.

Choose a direct sum decomposition $M\oplus N=A^r$ and let
$e\in\Mat_rA$ be the corresponding projector $A^r\to M$.
Then, one has a well-defined non-commutative 1-form
$de\in \ncO^1(\Mat_rA)$. 
For each $k=1,2,\ldots,$ we consider the differential form
\beq{ede}
e(de)^{2k}:=e\,\underbrace{de\,\,de\,\ldots\,de}_{2k\enspace\text{factors}}\;\in\;
\ncO^{2k}(\Mat_rA),
\end{equation}
and the corresponding class in $\DR^{2k}(\Mat_rA/\Mat_r\C)$. 
Let $\Tr\bigl(e(de)^{2k}\bigr)$ be the image of that class under
the canonical `trace'-isomorphism
$\DR^{2k}(\Mat_rA/\Mat_r\C)\iso$
$\DR^{2k}(A)$, cf. Remark \ref{tr_dr}.

\begin{prop}\label{ede_prop} \vi The class $\Tr\bigl(e(de)^{2k}\bigr)$
is independent of the choice of presentation of $M$ as a direct
summand in a free $A$-module, hence is intrisically attached to $M$. The assignment 
$$[M]\mto \ch_k([M):=\mbox{$\frac{1}{k!}$}\Tr\bigl(e(de)^{2k}\bigr)\in \DR^{2k}(A)$$
gives a group homomorphism $K^0(A)\map \DR^{2k}(A)$.

\vii The class $\Tr\bigl(e(de)^{2k}\bigr)\in \DR^{2k}(A)$ is
closed, i.e., $d\bigl[\Tr\bigl(e(de)^{2k}\bigr)\bigr]=0.$
\end{prop}
\begin{proof} It is clear that if $e'=geg\inv, \, g\in\GL_\infty(A),$ is another projector
then $\Tr\bigl(e(de)^{2k}\bigr)=\Tr\bigl(e'(de')^{2k}\bigr)$,
due to the invariance of the trace. Therefore, Lemma
\ref{conj_invar} implies independence of  presentation of $M$ as a direct
summand in a free $A$-module. 

Further, given two idempotents $e_1,e_2\in\Mat_\infty(A),$ we clearly have
$d(e_1\oplus e_2)=(de_1)\oplus(de_2)$,
hence, $(e_1\oplus e_2)\cdot (d(e_1\oplus e_2))^{2k}=
e_1(de_1)^{2k}\oplus e_2(de_2)^{2k}$.
Therefore, additivity of the trace implies that,
for any finite rank projective $A$-modules $P$ and $Q$,
one has $\ch_k([P]\oplus[Q])=\ch_k([P])+\ch_k([Q])$.
This completes the proof of (i). 

Part (ii)
may be verified by a direct computation. Instead of doing so,
below
we will give an alternative, more conceptual, construction of
the characteristic classes in terms of {\em connections}, and then
prove an analogue of Proposition \ref{ede_prop}(ii) in that
more general framework.
\end{proof}

Following A. Connes
\cite{Co},
one introduces
\begin{defn}\label{connection}
A connection on a left $A$-module $M$ is a linear
map $\nabla: M\to \ncO^1(A)\otimes_A M$ such that
$$\nabla(am)= a\cdot\nabla(m) + da\otimes m,\;\forall a\in A,m\in M.$$
\end{defn}

\begin{lem}\label{conn_exist} A left $A$-module $M$ admits a connection if and only if
it is projective.
\end{lem}
\begin{proof} Assume $M$ is projective and write it as a direct
summand of a free $A$-module $A\otimes E$,
for some $\C$-vector space $E$. Thus there is an
$A$-module imbedding $i: M\into A\otimes E$ and a projection
$A\otimes E\onto M$ such that $p\ccirc i=\Id_M$.
We define $\nabla$ to be the following composite map
\beq{levi_cevi}
M\stackrel{i}\into A\otimes E\stackrel{d\otimes\Id_E}\map
\ncO^1(A)\otimes E\iso\ncO^1(A)\otimes_A(A\otimes E)
\stackrel{\Id_\Om\otimes p}\map\ncO^1(A)\otimes_A M.
\end{equation}
It is easy to see that this map gives a connection,
sometimes called the {\em Grassmannian connection}
induced from $A\otimes E$.

Conversely, let $M$ be any left  $A$-module.
Observe that since $\ncO^1(A)\cong \bar{A}\otimes A$ is a free
right  $A$-module,
the functor $(-)\otimes_AM$ takes the fundamental exact
sequence $0\to\ncO^1(A)\to A\ee\to A\to 0$
to an {\em exact} sequence that looks as follows:
\beq{OM}
0\map \ncO^1(A)\otimes_AM\stackrel{j}{\map} A\otimes M
\stackrel{\ac}\map M\map 0.
\end{equation}
Here the map $j$ takes $da\otimes m$ to $a\otimes m-1\otimes (am)$,
and the map $\ac: A\otimes M\to M$ is
the action map. 

Now, given a connection $\nabla: M\to \ncO^1(A)\otimes_AM$,
we define a map 
$$s: M\to A\otimes M,\; m\mto 1\otimes m - j\ccirc\nabla(m).$$
For any $a\in A$, using the definition of $j$ and
Definition \ref{connection}, we compute
\begin{align*}
s(am)-a\cdot s(m)&=1\otimes (am)-j\ccirc\nabla(am)-a\cdot[1\otimes m - j\ccirc\nabla(m)] \\
&=1\otimes (am)-a\otimes m -j[a\cdot\nabla(m) + da\otimes m]
+j[a\cdot\nabla(m)]
\\
&=j[da\otimes m]-j[a\cdot\nabla(m) + da\otimes m]
+j[a\cdot\nabla(m)]=0.
\end{align*}
Hence, $s$ is an $A$-module map, moreover,
one finds 
$$\ac\ccirc s(m)=\ac[1\otimes m - j\ccirc\nabla(m)]=
m- (\ac\ccirc j)\ccirc\nabla(m)= m-0=m.$$
We see that the map $s$ provides a splitting of the projection $\ac$ in
\eqref{OM}, therefore $M$ is a direct summand in
$A\otimes M$, hence it is projective.
\end{proof}

It is easy to see that any connection $\nabla: M\to \ncO^1(A)\otimes_A
M$
has a unique extension to a map $\nabla: \ncO^\hdot(A)\otimes_A M
\to\ncO^{\hdot+1}(A)\otimes_A M$ such that
\beq{con_OM}\nabla(\alpha\cdot \mu)=d\alpha\cdot
\mu+(-1)^{\deg\alpha}\cdot\alpha\cdot\nabla(\mu),\;\forall \alpha\in
\ncO^\hdot(A),
\mu\in \ncO^\hdot(A)\otimes_A M.
\end{equation}
In particular, one defines the curvature of the connection 
$\nabla$ as the composite map 
\beq{curvature}
R_\nabla:\
\ncO^\hdot(A)\otimes_A M\stackrel{\nabla}\map
\ncO^{\hdot+1}(A)\otimes_A M\stackrel{\nabla}\map
\ncO^{\hdot+2}(A)\otimes_A M.
\end{equation}
From formula \eqref{con_OM}, we compute
\begin{align*}
R_\nabla(\alpha\cdot \mu)&=\nabla\ccirc\nabla(\alpha\cdot \mu)=
\nabla[d\alpha\cdot
\mu+(-1)^{\deg\alpha}\cdot\alpha\cdot\nabla(\mu)]\\
&=d^2(\alpha)\cdot
\mu+(-1)^{\deg d\alpha}\cdot d\alpha\cdot\nabla(\mu)
+(-1)^{\deg\alpha}\cdot d\alpha\cdot\nabla(\mu)\\
&\enspace+(-1)^{\deg\alpha+\deg d\alpha}\cdot
\alpha\cdot\nabla\ccirc\nabla(\mu)
=\alpha\cdot R_\nabla(\mu).
\end{align*}
Thus, the curvature is an $\ncO^\hdot(A)$-linear map.

We can now proceed to the construction of characteristic classes.
Fix a   finite rank projective left
$A$-module $M$ and choose an imbedding of $M$ into
a free module $A^r$ as a direct summand. Let
$e\in \Mat_rA$ be the idempotent that projects $A^r$ to $M$,
and let $\nabla_e:=e\ccirc d\ccirc e$ be the corresponding
Grassmannian connection on $M$. 
We may view the curvature
$R_{\nabla_e}= e\ccirc d\ccirc e\ccirc d\ccirc e$ as an element of $\Mat_rA\otimes
\ncO^{2}(A)$. 
For each $k=1,2,\ldots,$ we consider the class of the
element  $(R_{\nabla_e})^k\in \Mat_rA\otimes
\ncO^{2k}(A)$ in $\DR^{2k}(\Mat_nA/\Mat_n\C)$.
Thus, applying the  canonical isomorphism $\Tr: \DR^\hdot(\Mat_nA/\Mat_n\C)\iso \DR^\hdot(A)$
of Proposition \ref{morita_DR} to $\frac{1}{n!}(R_{\nabla_e})^k$
we obtain an element
\beq{chern_max}
\ch_k(M,\nabla_e):=\mbox{$\frac{1}{n!}$}\Tr(R_{\nabla_e})^k\in \DR^{2k}(A),
\end{equation}
called the $k$-th de Rham Chern character class.

\section{Formally Smooth Algebras}
\subsection{} We are going to study the concept of `smoothness' in
noncommutative geometry. Throughout this section $A$ denotes a finitely
generated associative algebra.

Recall that  a two-sided ideal $I$ of an  associative algebra $B$ 
is said to be nilpotent if there exists $n>0$ such that
$b_1\cdot\ldots\cdot b_n=0,$ for any $b_1,\ldots,b_n\in I$.
\begin{defn}
A finitely generated associative algebra $A$ is called \emph{formally
smooth} 
if the following
lifting property holds.
  For every algebra $B$ and a  nilpotent two-sided ideal $I\subset B$,
 given a map $A\to B/I$ there 
is a lift $A\to B$ such that the following diagram commutes:
\beq{lif}
\xymatrix{&B\ar[d]\\
A\ar@{.>}[ur]\ar[r]&B/I,}
\end{equation}
where $B\onto B/I$ is the quotient map.
\end{defn}

To build some intuition for formally smooth algebras
we consider commutative case first.

\begin{thm} For  a finitely-generated commutative $\C$-algebra
the following conditions are equivalent.
\begin{enumerate}
\item[(a)]{Let $m\colon A\otimes A\to A$ be the multiplication map.  Then $\Ker m$ has the locally complete intersection property (this is basically used in the proof of the Hochschild-Kostant-Rosenberg theorem).}
\item[(b)]{$\comO^1(A)$ is a projective $A$-module.}
\item[(c)]{$A$ satisfies the lifting property \eqref{lif} for any
commutative algebra $B$.} \qed
\end{enumerate}
\end{thm}

Recall that ${\Rep^A_E}$ denotes the 
algebraic variety of all representations of $A$ on~$E$.

\begin{prop}
If $A$ is formally smooth and finitely generated, then for every
finite-dimensional $\C$-vector space $E$, the scheme ${\Rep^A_E}$ is smooth.
\end{prop}

\begin{proof}
Take $E=\C^n$, 
 and  let $\Rep_n^A:=\Rep^A_{\C^n}={\Rep^A_E}$.  
For any scheme $X$ and any finitely generated commutative algebra $B$, we set
$$
X(B)=\Hom_{\sf{alg}}(\C[X],B).
$$
Giving an element of $X(B)$ is equivalent to giving an algebraic map $\Spec B\to X$.  The elements of $X(B)$ are called the $B$-points of $X$.  In the case of $\Rep_n^A$, we see that for any such $B$,
$$
\Rep_n^A(B)=\Hom_{\sf{alg}}(A,\Mat_nB).
$$
Observe that if $I\subset B$ is a nilpotent ideal, then so is $\Mat_n(I)\subset\Mat_nB$.  Let $R=\C[\Rep_n^A]$.  We will check that $R$ is formally smooth.  In other words, we wish to see if the obvious map
$$
\Hom_{\sf{alg}}(R,B)\to\Hom_{\sf{alg}}(R,B/I)
$$
is a surjection.  By definition, $\Hom_{\sf{alg}}(R,B)\simeq\Hom_{\sf{alg}}(A,\Mat_nB)$.  By the formal smoothness of $A$,
$$
\Hom_{\sf{alg}}(A,\Mat_nB)\to\Hom_{\sf{alg}}(A,\Mat_nB/\Mat_n(I))
$$
is surjective since $\Mat_n(I)\subset\Mat_nB$.  The proof then follows.
\end{proof}

\begin{prop}\label{smooth2} The algebra
$A$ is formally smooth if and only if $\ncO^1(A)$ is a projective
$A$-bimodule.
Equivalently, $A$ is formally smooth if and only if
the functor $\Der(A,-)$, on the category
of $A$-bimodules,  is  exact.
\end{prop}

\begin{lem}\label{dim1} If $\ncO^1(A)$ is a projective $A\ee$-module, then
the categories $\Lmod A$ and $\Lmod {A\ee}$ both have
 homological dimension less than or equal to
1, i.e., we have
$$
\Ext_{\Lmod A}^j(M,N)=0\quad\text{resp.}\quad
\Ext_{\bimod A}^j(M,N)=0,\quad\text{for all}\quad j>1.
$$
\end{lem}

\begin{proof} We prove the statement for left $A$-modules;
the proof for $A$-bimodules is similar.

  Recall the fundumental exact sequence
$$0\map \ncO^1(A)\map A\ee\map A\map 0.
$$
Now, if $\ncO^1(A)$ is a projective $A\ee$-module, then
 $\ncO^1(A)$ is projective as a left and right module.
So tensoring the fundamental sequence by any left $A$-module  $M$
preserves exactness, hence yields an exact
sequence of left $A$-modules 
$$
0\to\ncO^1(A)\otimes_AM\to A\o M\to M\to0.
$$
Here, the left $A$-module $A\o M$ is projective since
it is free.  
Further, the projectivity of $\ncO^1(A)$ implies that
 $\ncO^1(A)$  is projective as a left module.  Indeed, if $\ncO^1(A)$ is a direct summand of the free bimodule $A\otimes E\otimes A$, then $\ncO^1(A)\otimes_AM$ is a direct summand of the free left module $A\otimes E\otimes M$.
Thus, we have constructed a length two resolution of
$M$ by projective $A$-modules.
Computing the Ext-groups via this resolution we conclude that
$
\Ext_{\Lmod A}^j(M,N)=0$ for any $j>1$.
\end{proof}

We observe that if $I^2=0$, then $B=A{\sharp}I$, the square zero construction, in the lifting problem if $A$ is formally smooth.

\begin{lem}\label{proj_lif}
The following are equivalent.
\begin{enumerate}
\item{$\ncO^1(A)$ is projective.}
\item{Lifting property holds for any square zero extension.}
\item{$\HH^2(A,M)=0$ for any $A$-bimodule $M$.}
\end{enumerate}
\end{lem}

\begin{proof}[Proof of Lemma \ref{proj_lif}.]
We will prove that each of the claims (1) and (2) is equivalent to (3).

Recall that the square zero extensions the algebra $A=B/I$ by $I$ are
classified by $\HH^2(B/I,I)$. 

Suppose that $\HH^2(A,M)=0$ for all $M$.  Suppose we wish to lift a map
$\alpha\colon A\to B/I$.  
The  pull-back  of extension $I\to B\to B/I$ via $\alpha$ gives a
commutative diagram:
$$
\xymatrix{0\ar[r]&I\ar@{=}[d]_<>(0.5){\id}\ar[r]&E\ar[d]\ar[r]&A\ar[r]\ar[d]^\alpha&0\\
0\ar[r]&I\ar[r]&B\ar[r]^\beta&B/I\ar[r]&0,}
$$
where the left-hand vertical map is the identity.  The algebra $E$ is given by the fiber product of $A$ and $B$, that is, $E=\{(b,a)\in B\oplus A\st\alpha(a)=\beta(b)\}$.  Since $\HH^2(A,M)$, the top row is split by some $\sigma\colon A\to E$.  Then, letting $E\to B$ be denoted by $\tau$, the composition $\tau\circ\sigma$ is the desired lift.

Now, suppose that the lifting property for square zero extensions holds.  We wish to show that $\HH^2(A,M)=0$.  Now, an element of $\HH^2(A,M)$ gives a square zero extension $0\to M\to E\to A\to0$.  Now, since lifting holds for the square zero case, we can lift the identity map $A\to A$ to a map $A\to E$, which splits the extension.  So, (2) and $\HH^2(A,M)=0$ are equivalent.
 
Finally,  an $A$-bimodule $P$ is projecive if and only if
$\Ext_{\bimod A}^1(P,M)=0$ for all $A$-bimodules $M$. Hence, $\ncO^1(A)$
is projecive if and only if
$\Ext_{\bimod A}^1(\ncO^1(A),M)=0$ for all $A$-bimodules $M$.
But the long exact sequence for  $\Ext$
arising from the fundamental short exact 
sequence $\ncO^1(A)\to A\ee\to A$ yields
$$
\Ext_{\bimod A}^1(\ncO^1(A),M)=\Ext_{\bimod A}^{2}(A,M)\simeq
\HH^2(A,M).
$$
Thus, $\ncO^1(A)$ is projecive if and only if
$\HH^2(A,M)=0$.
 This completes the proof.
\end{proof}

\begin{proof}[Proof of Proposition \ref{smooth2}]
Let $I\subset B$ be any nilpotent ideal, $I^n=0$, of an algebra $B$.  We will proceed by induction on $n$.  Consider the exact sequence
$$
0\to I^{n-1}/I^n\to B/I^n\to B/I^{n-1}\to0.
$$
This is then a square zero extension of $B/I^{n-1}$ if $n\ge2$.  Take a map $A\to B/I^{n-1}$.  Then let $E$ be the fiber product of $A$ and $B/I^n$ to obtain the commutative diagram (with exact rows)
$$
\xymatrix{0\ar[r]&I^{n-1}/I^n\ar@{=}[d]_<>(0.5){\id}\ar[r]&E\ar[d]\ar[r]&A\ar[d]\ar[r]&0\\
0\ar[r]&I^{n-1}/I^n\ar[r]&B/I^n\ar[r]&B/I^{n-1}\ar[r]&0}.
$$
Then by assumption (and the lemma), there is a splitting $A\to E\to B$.  So, by inducting on $n$ until $I^n=0$, we obtain lifting.

The implication that there exists a lifting implies $\ncO^1(A)$ is
projective follows
 from Lemma \ref{proj_lif}.
\end{proof}

\subsection{Examples of formally smooth algebras.} Here are a few examples:
\begin{enumerate}
\item{The free associative algebra $\C\langle x_1,\ldots,x_n\rangle$.}
\item{$\Mat_n\C$.}
\item{$\C[X]$ where $X$ is a smooth affine curve.}
\item{The path algebra of a quiver.}
\item{The upper triangular matrices. This is a special case of (4),
since the  algebra of upper triangular matrices  is nothing but
the path algebra of the
quiver 
$\bullet\to\bullet\to\cdots\to\bullet$.}
\item{If $A$ and $B$ are formally smooth, then so are  $A\oplus B$ and $A*B$.}
\end{enumerate}

The reader should be warned that the (commutative) polynomial algebra
$\C[x_1,$
$\ldots,x_n]$ is {\em not} formally smooth, for any $n>1$.

\subsection{Coherent modules and algebras.} It is perhaps clear from 
discussion in the
previous sections that a formally smooth finitely generated associative
(not necessarily commutative) algebra $A$ should be viewed
as a `noncommutative analogue' of the coordinate ring
of a smooth affine algebraic variety $X$. Accordingly,
the category of finitely generated $A$-bimodules  should be viewed
as a `noncommutative analogue' of the abelian category
$\Coh(X)$. 

An immediate problem that one encounters with 
such an  analogy is that a
finitely generated formally smooth algebra $A$ is 
typically {\em  not}
Noetherian, hence, neither the category
of finitely generated left $A$-modules nor the
 category
of finitely generated $A$-bimodules, are abelian categories,
in general.
This is so, for instance, in the `flat' case  where
$A=\C\langle x_1,\ldots,x_n\rangle,$
is a free associative algebra on $n$ generators.
This difficulty can be dealt with by replacing the notion
of a finitely generated module by a more restrictive notion of
{\em coherent module}.

In general, let $A$ be an associative algebra. We introduce, see
\cite{Po}.
\begin{defn} A (left) $A$-module $M$ is called {\em coherent}
if $M$ is finitely generated and, moreover,
the kernel of any $A$-module map $A^r\to M$ is also
a finitely generated $A$-module.
\end{defn}

It is straightforward to verify that if $f: M\to N$ is a morphism
between two coherent modules, then both the kernel and cokernel
of $f$ are again coherent modules. Thus, coherent $A$-modules
form an abelian category. 

In order for the concept of coherent module to be useful
one has to know that, for a given algebra $A$, there are
``sufficiently many'' coherent modules. We will see below that
this is indeed the case for  formally smooth algebras.

First, we recall that an algebra $A$ is called \emph{hereditary} if $\Ext^j_{\Lmod A}(M,N)=0$
for all $A$-modules $M$ and $N$ and $j\ge2$,
in other words,  $A$ is hereditary if
the category $\Lmod A$ has homological dimension less than or equal to
1,
c.f. Lemma \ref{dim1}.

\begin{lem}
An algebra $A$ is hereditary if and only if every submodule of a 
projective $A$-module is again projective.
\end{lem}

\begin{proof}
Let $P$ be a projective $A$-module, and let $P'\subset P$ be a submodule.  Then we have the short exact sequence
$$
0\to P'\to P\to P/P'\to0,
$$
where the middle term is projective.  If $M$ is any $A$-module, 
we obtain a long exact sequence in $\Ext$ groups:
$$
\ldots\to\Ext^1(P,M)\to\Ext^1(P',M)\to\Ext^2(P/P',M)\to\Ext^2(P,M)\to\ldots.
$$
Since $P$ is projective, $\Ext^1(P,M)=\Ext^2(P,M)=0$, hence we find that
$$
\Ext^1(P',M)\simeq\Ext^2(P/P',M)
$$
for all $A$-modules $M$.  But $A$ is hereditary, hence $\Ext^2(P/P',M)=0$.  Therefore, $\Ext^1(P',M)=0$ for all $M\in\Lmod A$, hence $P'$ is projective.

Conversely, let $M$ be any $A$-module.  Then we have the resolution
$$
0\to P\to A^{\oplus s}\to M\to0,
$$
where $P=\Ker(A^{\oplus s}\to M)$.  Since $A^{\oplus s}$ is free, it is projective.  Therefore $P$ is also projective by assumption.  This shows that every $A$-module $M$ has a projective resolution of length at most two, hence $\Ext^j(M,N)=0$ for all $j\ge2$.
\end{proof}
 It is known that $A$ is hereditary if and only if
 all left ideals of $A$ are projective.

We observe next that if $A$ is formally smooth, 
then by Lemma \ref{dim1} the category
$\Lmod A$ has homological dimension less than or equal to 1,
so $A$ is a \emph{hereditary} algebra. 
Furthermore, the proposition below insures that the category
of coherent $A$-modules, reps.,  coherent $A$-bimodules,
is sufficiently ``large''.
\begin{lem} Let $A$ be a finitely generated hereditary
algebra. Then, any finite rank free 
$A$-module is  coherent.
\end{lem}
  
\begin{proof}[Proof (by D. Boyarchenko)]  Let $M= A^n$ be a free left $A$-module of
finite rank $n$. Clearly, $M$ is finitely generated, so we only have to prove
that if $f : A^r \to M$ is a homomorphism, then $\Ker(f)$ is also finitely
generated. Let $K = \Ker(f),$ and $ Q = \im(f).$ We have a short exact sequence
$0 \to K \to A^r \to Q \to 0.$
By construction,
$Q$ is a submodule of the free module $M.$ Hence $Q$ is projective,
since  $A$ is hereditary.
Hence
the above exact sequence splits. In particular, this yields a surjection 
$A^r \onto K,$ which implies that $K$ is finitely generated.
It follows that any  finite rank free left $A$-module is  coherent.
\end{proof}

\begin{cor} If $A$ is a formally smooth algebra, then
the cokernel of any $A$-module map
$A^m\to A^n$ is a coherent $A$-module.
Also, any finite rank free $A\ee$-module, is  coherent.
\end{cor}

\begin{proof} The first claim follows from the previous Lemma, since a
formally smooth algebra is hereditary.

Further, one proves easily that if $M$ is a coherent left $A$-module
and $N$  is a coherent left $B$-module, then
$M\otimes N$  is a coherent left $A\otimes B$-module.
The claim on $A\ee$-modules follows.
\end{proof}
\subsection{Smoothness via torsion-free connection}\label{tors}
For any associative algebra $A$, we define
$$
D(A)=T(A+{\overline{A}})/(\overline{ab}=a\bar b+\bar ab, a\otimes a'=aa'\otimes1),
$$
where as usual ${\overline{A}}=A$ as a vector space.  Then $a\mapsto\bar a$ is a differential.  Earlier, we showed that $D(A)\simeq\ncO^\hdot(A)\simeq T_A\ncO^1(A)$.

Let us consider the commutative situation.  Let $X$ be a smooth affine algebraic variety.  Then we have the two following relations:
\begin{align*}
\Omega^\hdot(X)=\Lambda^\hdot\Omega^1(X)\quad\text{and}\quad
\C[TX]=\Sym\Omega^1(X),
\end{align*}
where we use the convention $\Omega^1(X)=\comO^1(\C[X])$.
  In the noncommutative case, there is no difference between the
exterior 
and symmetric powers.  Hence we see that the single differential graded algebra $D(A)=\ncO^\hdot(A)$ can be simultaneously thought of as noncommutative differential forms and as functions on the ``noncommutative tangent bundle.''  When we wish to stress the latter interpretation, we will write $\check D(A)$ instead of $D(A)$.

Fix a finite dimensional vector space $E$
 and consider the variety $\Rep^{\check D(A)}_E$ of 
algebra maps $\check D(A)\to \End_\k{E}.$

\begin{prop}
The varieties $\Rep^{\check D(A)}_E$ and $T\Rep_E^A$ are isomorphic.
\end{prop}

\begin{proof}
An element of $\Rep^{\check D(A)}_E$ is a homomorphism $\check D(A)\to\End{E}$.  But a homomorphism from $\check D(A)$ can be specified by giving the image of $a$ and $\bar a$ for each $a\in A$.  So, let $\rho(a)$ be the image of $a$, and $\f(a)$ the image of $\bar a$.  Then we can easily check that $\rho\colon A\to\End{E}$ must be a homomorphism, while $\f\colon A\to\End{E}$ is a derivation.  This is precisely a point of $T\Rep_E^A$.  Given a point $T\Rep_E^A$, we can clearly reverse the arguments made above to construct a homomorphism $\check D(A)\to\End{E}$.
\end{proof}

\begin{thm}
An associative algebra $A$ is formally smooth if and only if the natural map $A\to{\overline{A}}$ can be extended to a derivation of $\check D(A)$ of degree~$+1$.
\end{thm}

Suppose we are in the commutative case, that is, $A=\C[X]$ for some 
affine variety $X$.  Then we view $\check D(A)$ as the coordinate ring
of the total space of the tangent bundle $TX$ of $X$.  A derivation
$\nabla$ of $\C[TX]$ is a 
vector field on $TX$.  This gives a connection on $TX$.  Indeed, since
$\nabla$ is a derivation, this connection has no torsion (that is,
$\nabla_{[\xi,\eta]}
=\nabla_\xi\eta-\nabla_\eta\xi$).  

\begin{proof}[Proof of Theorem]
We already know that $A$ is formally smooth if and only if every square zero extension
$$
0\to M\to E\to A\to0
$$
splits. To check that  every square zero extension splits,
it suffices to check this for the universal 
square zero extension, $A{\sharp}_c\ncO^2(A)$, see \ref{E:SqZero}.
A splitting of the latter extension is provided
by an algebra map $\psi\colon A\to E$ 
such that $\psi\circ j=\id$, where $j\colon A\to E$ is the inclusion map.  So, write $\psi(a)=(a,-\phi(a))$ for some function $\phi\colon A\to\ncO^2(A)$ (we are not asserting any properties for $\phi$ except linearity, which is obvious).  Then we see that
\begin{align*}
\psi(a_1)\psi(a_2)=(a_1a_2,-a_1\phi(a_2)-\phi(a_1)a_2-da_1\,da_2)&
=(a_1a_2,-\phi(a_1a_2))\\
&=\psi(a_1a_2).
\end{align*}
So, we find that
$
\phi(a_1a_2)-a_1\phi(a_1)-\phi(a_1)a_2=da_1\,da_2.
$
\end{proof}

Now we can define a derivation of $\check D(A)$ by $a\mapsto\bar a$ 
and $\bar a\mapsto\phi(\bar a)$.  The proof of the opposite implications follows from the universality of equation~\eqref{E:SqZero}.

But the left hand side of the above equation is precisely $\delta\phi$, where $\delta$ is the Hochschild differential.  We write the right hand side as $d\otimes d$ to obtain the equation
$$
\delta\phi=d\otimes d.
$$

Now, we claim that extending the map $a\mapsto\bar a$ to a derivation $+1$ on $\check D(A)$ is equivalent to giving a map ${\overline{A}}\to\ncO^2(A)$.  

\begin{lem}
Giving a map $\phi$ satisfying $\delta\phi=d\otimes d$ is equivalent to giving an $A$-bimodule splitting of the sequence
$$
0\too \ncO^1(A)\stackrel{j}\too\ncO^1(A)\otimes A
\stackrel{m}\too\ncO^1(A)\too 0
$$
where $m$ is right multiplication, $m(\omega\otimes a)=\omega\cdot a$, and
$$
j(\alpha\,da)=\alpha a\otimes1-\alpha\otimes a,\quad\text{for all}
\quad\alpha\in\ncO^1(A)\,,\,a\in A.
$$
\end{lem}

(By splitting, we mean a 
map $p\colon\ncO^1(A)$ $\otimes A\to\ncO^2(A)$).
\begin{proof}
Giving such a $\phi$ is equivalent to giving a map ${\overline{A}}\to\ncO^2(A)$, and we use this along with bilinearity to see that this is equivalent to giving a map
$$
p\colon\ncO^1(A)\otimes A=A\otimes{\overline{A}}\otimes A\to\ncO^2(A).
$$
Since $\phi$ satisfies $\delta\phi=d\otimes d$, we  can check that $p$ is a splitting, i.e., $pj=\id$.  Conversely, if $pj=\id$, we see that
\begin{align*}
pj(da_1\,da_2)&=p(da_1\cdot a_2\otimes1-da_1\otimes a_2)\\
&=p(d(a_1a_2)\otimes1-a_2\,da_2\otimes1-da_1\otimes a_2)=da_1\,da_2,
\end{align*}
forces $\delta p=d\otimes d$.
\end{proof}

\section{Serre functors and Duality}
\subsection{}
We write $\Vect$ for the category of finite dimensional
vector spaces, and $V\mapsto V^*=\Hom_\C(V,\C),$  for the
obvious duality functor on $\Vect$.

In this section, we will freely use the language of derived categories.
We write $[n]$ for the shift by $n$ in a triangulated category.
A functor $F: D_1 \to D_2$ between triangulated categories
is said to be a {\em triangulated functor}
if it takes distinguished triangles into  distinguished  triangles,
and commutes with the shift functors.

A $\C$-linear category $D$ is said to be {\em Hom-finite}
if,  for any two objects $M,N\in D$,
the space $\Hom_D(M,N)$ has finite dimension over $\C$.

\begin{defn}[Bondal-Kapranov]\label{serre_def}
Let $D$ be a $\Hom$-finite triangulated category.
An exact functor $\se: D \to D$ is called a {\em Serre functor}
if there are functorial vector space isomorphisms
$$\Hom_D(M,N)^*\simeq \Hom_D(N, \se(M)),\quad\text{for any}\quad
M,N\in D.
$$
\end{defn}

Let  $D$ be a  $\Hom$-finite category.
For   any  object $M\in D$, we consider the composite functor
$\Hom_D(M,-)^*: D \to \Vect\,,\, N\mapsto 
\Hom_\C\bigl(\Hom_D(M,N),\,\C\big)=\Hom_D(M,N)^*$.
Observe that, if the category $D$ has a Serre functor
$\se$, then the functor  $\Hom_D(M,-)^*$ is, by definition,
{\em represented} by the object $\se(M)$.
Yoneda Lemma \ref{yoneda} insures that
 the object representing this functor is unique, if exists,
up to (essentially unique) isomorphism. Using this, it is not difficult to deduce
that any two Serre functors on a $\Hom$-finite  category must
be isomorphic to each other.

The  definition of Serre functor is motivated by the following geometric
example.

\subsection{Serre duality.}\label{Serre_duality} Let $X$ be a smooth projective
algebraic variety of dimension $d$. Write 
$\Dcoh(X)$ for the bounded derived category
of sheaves of $\oo_X$-modules on $X$ with coherent cohomology sheaves. 
As a consequence
of {\em completeness} of $X$, 
the category $\Dcoh(X)$  is $\Hom$-finite.
This follows from  the well-known result,
saying that $\dim \left(\bigoplus_i\, H^i(X, \cal F)\right)<\infty$,
for any 
 coherent sheaf $\cal F$ on $X$.

Let $K_X:=\Omega^d_X$ be
the {\em canonical line bundle} on $X$, the line bundle
of top-degree differential forms on $X$, viewed as an invertible
sheaf on $X$.
\begin{prop}\label{Serre_duality_prop}
The functor $M\mto 
M\otimes K_X[d]$ is a Serre functor on $\Dcoh(X)$.
\end{prop}

\begin{rem} \vi
The condition that $X$ is smooth is not very essential here.
In the non-smooth
case, one has replace $K_X[d]$ by the {\em dualising complex}
$\kappa_{_X}$, and to restrict oneself to the
category  $\Perf(X)\sset\Dcoh(X)$ of {\em perfect complexes},
that is, the full triangulated subcategory in $\Dcoh(X)$ formed
by complexes which are  quasi-isomorphic to
bounded complexes of locally free sheaves.
For example, if $X$ is a Cohen-Macaulay projective scheme, then
the functor $M\mto 
M\otimes \kappa_{_X}$ gives a Serre functor on $\Perf(X)$.
In general, the Serre functor is an equivalence
between the category  $\Perf(X)$ and the category of
bounded coherent complexes with finite injective dimension.

\vii 
According to Kontsevich, any  $\Hom$-finite triangulated category with
a Serre functor should be thought of as  the category $\Dcoh(X)$
for some  {\em complete}
`noncommutative space' $X$,  possibly singular.
\eer

Below, we will interchangeably
use the words ``locally free sheaf'' and ``vector bundle''
and, given such a vector bundle $E$, write $E^*$ for the dual vector bundle.

The Proof of  Proposition \ref{Serre_duality_prop} uses the following important result
\begin{thm}[Grothendieck]\label{groth}
An algebraic variety $X$ is smooth if and
only
if every coherent sheaf on $X$ has a finite resolution by
locally free sheaves, equivalently,
if the (shifts of)  vector bundles on $X$ generate the category
 $\Dcoh(X)$.\qed
\end{thm}

\begin{proof}[Proof of  Proposition \ref{Serre_duality_prop}]
 The result is essentially a reformulation of the
standard Serre duality. The latter says that, for any 
vector bundle $E$ on $X$, one has
$$
H^i(X,E)^*\cong H^{d-i}(X,K_X\otimes E^*),\quad
\forall i\in\Z
\quad(\text{\bf Serre duality}).
$$

Now, for any two vector bundles $F_1,F_2,$ on $X$, we 
have $\Ext^\hdot(F_1,F_2)=H^\hdot\bigl(X,$
$\scr H\!om(F_1,F_2)\bigr)=
H^\hdot\bigl(X,F_1^*\otimes F_2),$
where $\Ext^\hdot(-,-)$ stands for the Ext-group in the
abelian category $\Coh(X)$, and $\scr H\!om(F_1,F_2)\in \Coh(X)$ stands for
the internal Hom-sheaf which, for vector bundles,
is isomorphic to $F_1^*\otimes F_2$. Using this, we
compute
\begin{align*}
\Hom_{\Dcoh(X)}(F_1,F_2[i])^*&=\Ext^i(F_1,F_2)^*=H^i(X,\,F_1^*\otimes F_2)^*\\
\text{(by Serre duality)}\quad&=H^{d-i}\bigl(X\,,\,
K_X\otimes (F_1^*\otimes F_2)^*\bigr)\\
&=H^{d-i}(X,\,K_X\otimes F_1\otimes F_2^*)\\
&=\Ext^{d-i}(F_2\,,\,K_X\otimes F_1)\\
&=\Hom_{\Dcoh(X)}(F_2[i]\,,\,F_1\otimes K_X[d])\\
&=\Hom_{\Dcoh(X)}\bigl(F_2[i]\,,\,\se(F_1)\bigr).
\end{align*}
Thus, we have checked the defining property of Serre functor
in the special case of vector bundles
(more precisely, the chain of isomorphisms above may be refined 
to yield  a morphism
of functors $\Hom(-, N)^* \to \Hom(N, \se(-))$,
which we have shown to be an isomorphism 
for locally free sheaves).
The general case now follows from Proposition \ref{groth}.
\end{proof}

We keep the above setup, and for any integer $n\geq 1$,
set $\se^n=\se\ccirc\se\ccirc\ldots\ccirc\se$ ($n$ times).
Clearly, we have $\se^n: M\mapsto M\otimes K_X^{\otimes n}[dn]$.
We see that any global section $s\in \Gamma(X,K_X^{\otimes n})$
gives, for each $M\in \Dcoh(X)$,
a morphism $\Phi_s: M\to M\otimes K_X^{\otimes n}=\se^n(M)[-dn]\,,\, m\mapsto m\otimes s$. 
Thus, we have a morphism of functors $\Phi_s: \Id_{\Dcoh(X)}
\to \se^n[-dn].$ 
This way, we get a linear map of vector spaces
\beq{orlov}\Gamma(X,K_X^{\otimes n})\too \Hom(\Id_{\Dcoh(X)}\,,\,\se^n[-dn]),
\quad s\mapsto \Phi_s,\end{equation}
which is easily seen to be an isomorphism.

We apply this to prove the following interesting result,
first due to Bondal-Orlov.
\begin{thm} Let $X$ and $Y$ be smooth projective varieties
such that the canonical  bundles $K_X$ and $K_Y$ are
both ample  line bundles on $X$ and $Y$, respectively.
Then, any trianulated equivalence $\Dcoh(X)\iso \Dcoh(Y)$
implies an isomorphism $ X\simeq Y$,
of algebraic varieties.
\end{thm} 
\begin{rem} 
The Theorem says that
a smooth projective variety with ample canonical
class is completely determined by the corresponding
triangulated category  $\Dcoh(X)$.
In particular, with the assumptions above,
one has
$\,\Dcoh(X)\cong  \Dcoh(Y)\quad\Longrightarrow\quad
X\cong Y.$
Such an implication is definitely false for varieties
with non-ample, e.g. with trivial, canonical bundles.

On the other hand, it is not difficult to show
that, for {\it any} algebraic varieties
$X$ and $Y$, an equivalence
$\Coh(X)\cong \Coh(Y),$ of {\it abelian} categories,
does imply an isomorphism $X\simeq Y$
(Hint: the assignment sending a point
$x\in X$ to the sky-scrapper sheaf at $x$
sets up a bijection between the set $X$ and
the set of (isomorphism classes of) simple objects
of the category $\Coh(X)$. This way, one recovers
$X$ from  $\Coh(X)$, as a set. A bit more efforts
allow to recover $X$ as an algebraic variety, as well.)
\eer

\begin{proof}[Proof of the Theorem (after Kontsevich)]
Put $D:=\Dcoh(X)$. We are going to 
 give a canonical procedure of
reconstructing the variety $X$ from the triangualted category
$D$.

To this end, 
observe  that, for each $n,m\geq 0$, the obvious sheaf morphism
 $K_X^{\otimes n} \otimes K_X^{\otimes
m}\to K_X^{\otimes (n+m)}$ induces
a linear map
$$
\Gamma(X,K_X^{\otimes n}) \otimes \Gamma(X,K_X^{\otimes m})
\too \Gamma(X,K_X^{\otimes (n+m)}).
$$

Similarly, for each $n,m\geq 0$,
there is a composition of morphisms of functors defined as follows
\begin{align*}
&\Hom(\Id_{D}\,,\,\se^n[-dn]) \, \otimes\,
\Hom(\Id_{D}\,,\,\se^m[-dm])\\ 
&{\xymatrix{
\ar[rr]^<>(0.5){\Id\otimes\se^n[-dn]}&&{\Hom(\Id_{D}\,,\,\se^n[-dn]) \, \otimes\,
\Hom(\se^n[-dn]\,,\,\se^{n+m}[-d(n+m)])}}}\\
&\hphantom{x}\hskip 60mm\too
\Hom(\Id_{D}\,,\,\se^{n+m}[-d(n+m)]).
\end{align*}
(here, we put $\se^0:=\Id_{D}$, by definition).
This way, from \eqref{orlov} we deduce a
graded algebra isomorphism
$$\bigoplus\nolimits_{n\geq 0}\,\Gamma(X,K_X^{\otimes n})\simeq
\bigoplus\nolimits_{n\geq 0}\,\Hom(\Id_{D}\,,\,\se^n[-dn]).
$$

Now, if $K_X$ is very ample, then the variety $X$ may be
obtained from the graded algebra on the LHS above via the
standard $\Proj$-construction, that is, we have
\beq{proj}
X=\Proj\left(\bigoplus\nolimits_{n\geq 0}\,\Gamma(X,K_X^{\otimes
n})\right)\simeq
\Proj\left(\bigoplus\nolimits_{n\geq
0}\,\Hom(\Id_{D}\,,\,\se^n[-dn])\right).
\end{equation}
If $K_X$ is ample but not very ample, we replace $K_X$
in this formula by
a sufficiently large power of $K_X$.
Thus,  formula  \eqref{proj} gives way to reconstruct the variety out of
the corresponding derived
category $\Dcoh(X)$.

We observe next that the integer $d=\dim X$ can be characterized
as follows:

\npb{$d$ is the unique integer with the property that there exists
an object $M\in \Dcoh(X)$ such that $\se(M)\simeq M[d]$.}

Now, let $X,Y$ be two smooth projective varieties such that
$\Dcoh(X)\cong  \Dcoh(Y)$.
Then, by the characterization above,
one must have $\dim X=\dim Y$.
Hence, the integer $d$ in the RHS of \eqref{proj} is equal
to the one in a similar formula for $Y$.
Further,  the uniqueness of the Serre functor
mentioned after Definition \ref{serre_def}
implies that the Serre functor on
$\Dcoh(X)$ goes under the equivalence
$\Dcoh(X)\cong  \Dcoh(Y)$
to the  Serre functor on
$\Dcoh(Y)$.
Hence, the  equivalence
yields a  graded algebra isomorphism
\begin{align*}
\bigoplus\nolimits_{n\geq
0}\,\Hom\bigl(\Id_{\Dcoh(X)}\,,\,&\se^{n+m}[-d(n+m)]\bigr)\\
&\simeq
\bigoplus\nolimits_{n\geq
0}\,\Hom\bigl(\Id_{\Dcoh(Y)}\,,\,\se^{n+m}[-d(n+m)]\bigr).
\end{align*}
Therefore, the corresponding Proj-schemes are isomorphic,
and we are done.
\end{proof}
\begin{rem} A similar result (with similar proof)
holds in the case where the varieties have ample
anti-canonical classes $(K_X)^{-1}$ and $(K_Y)^{-1}$.
\eer

\subsection{Calabi-Yau categories.}
Recall that a smooth variaty $X$ is called a {\em Calabi-Yau}
manifold if it has trivial canonical bundle, $K_X\cong\oo_X$.
Assuming in addition that $X$ is {\em projective} and has dimension $d$,
 the Calabi-Yau
property
can be reformulated as an isomorphism of functors
$\se(-)\simeq (-)[d]$. Motivated by this, one introduces the following
\begin{defn}\label{CY} A $\Hom$-finite triangulated category
$D$, with Serre functor $\se$, is said to be
a {\em Calabi-Yau category} of dimension $d$,
if there is an isomorphism of functors $\se(-)\simeq (-)[d]$.
In such a case, we write $d=\dim D$.
\end{defn}

We observe further that it makes sense to consider
Calabi-Yau categories of {\em fractional dimension}.
Specifically, we say that $\dim D=m/n,$ provided
there is an isomorphism of functors $\se^n(-)\simeq (-)[m]$.

\begin{examp}[Kontsevich] Let $A$ be the associative algebra
of upper-triangular $n\times n$-matrices (with zero diagonal entries).
Then, one can show that the category $D^b(\Lmod A)$ has 
dimension $\frac{n-1}{n+1}$.
This category is, in effect, related to the category
of coherent sheaves on the orbifold with
$\mathbf{A_n}$-type isolated singularity,
i.e., with
singularity of the form $\C^2/(\Z/n\Z)$.
\eex

The following important result is due to
\cite[Lemma 2.7]{BK}, cf. also \cite{BKR}.
\begin{thm} Let $D$ and $D'$ be two triangulated categories with Serre
functors, and let $F: D\to D'$ be an exact functor that
intertwines  the Serre functors on $D$ and $D'$. Assume in addition that

\npb{$F$ has a left adjoint $F^\top: D'\to D$ and
the adjunction morphism $\Id_D\to F\ccirc F^\top$ is an isomorphism.}

\npb{The category $D$ is indecomposable, i.e., there is no nontrivial
decomposition $D=C_1\oplus C_2$, into triangulated subcategories.}

\npb{$D$ is a Calabi-Yau category.}

Then, the functor $F$ is an equivalence of triangulated categories.\qed
\end{thm}
\subsection{Homological duality.}
Given an associative algebra $A$, let
$D(\Lmod A)$ be the derived category of
all (unbounded) complexes of left $A$-modules.

\begin{defn}\label{perf} Let $\Perf(A)$ be the full
subcategory in $D(\Lmod A)$ formed by the
complexes $C$, such that $C$ is quasi-isomorphic to a bounded
complex of finite rank projective $A$-modules.
\end{defn}

The following Lemma provides an interesting,
 purely category-theoretic interpretation of the
category $\Perf(A)$.

\begin{lem} An object $M\in D(\Lmod A)$ belongs to
$\Perf(A)$ if and only if $M$ is {\sl{compact}},
i.e.,  if the functor $\Hom_{D(\Lmod A)}(M,-)$
commutes with arbitrary direct sums, cf. Definition \ref{compact}.\qed
\end{lem}

Observe that, for any left $A$-module $M$, the space $\Hom_A(M,A)$
has a natural {\em right} $A$-module structure induced
by  right multiplication of $A$ on itself.
This gives a functor $\Hom_A(-,A): \Lmod A \to \Rmod A=\Lmod {A^\op}$.
Similarly, one has a functor  $\Hom_{A^\op}(-,A): \Lmod {A^\op} \to \Lmod A.$
Observe further that if $M$ is a projective
left  $A$-module, then $\Hom_A(M,A)$ is a projective
right  $A$-module, and vice versa.
Therefore, the functor  $\Hom_A(-,A)$
gives rise to  well-defined derived functors
$R\Hom_A(-,A): \Perf(A) \rightleftarrows
\Perf(A^\op)$.

Below, we will frequently use the following
canonical isomorphisms

\begin{align}\label{can_iso_perf}
&\!\!\!\!\bullet\enspace
M\iso R\Hom_{A^\op}\bigl(R\Hom_A(M,A)\,,\,A\bigr),\quad\forall M\in \Perf(A),\nonumber\\
&\!\!\!\!\bullet\enspace
R\Hom_A(M,A)\lotimes_A\,N \iso R\Hom_A(M,N),\quad\forall M,N\in \Perf(A).\\
&\!\!\!\!\bullet\enspace
A\lotimes_{A\ee}(L \o R)\simeq R\lotimes_A L\,,\quad\text{for any}\quad
L\in\Perf(A)\,,\,R\in\Perf(A^\op).
\nonumber
\end{align}

\noindent
Here, each of the isomorphisms is clear for finite rank free
modules, hence, holds for  finite rank projective
modules. This yields the result for arbitrary objects of $\Perf(A).$

Recall next that an abelian category $\scr C$ is said to
have {\em homological dimension} $\leq d$ if,
for any objects $M,N\in \scr C$, we have
$\Ext^i_{\scr C}(M,N)=0,$ for all $i>d$. For a smooth
variety of dimension $d$, the category $\Coh(X)$ is
known to have dimension $\leq d$ (and the inequality is in effect 
an equality).

\begin{rem}  We recall that if  $A$ and $B$ are both
formally smooth algebras
then  $A\otimes B$ is not necessarily  formally smooth.
Similarly, if $\Lmod A$ and $\Lmod B$  both have finite
 homological dimension then this is  not necessarily 
so for $\Lmod {(A\otimes B)}$, e.g., take
$A=B=\mathbb{K}$, a field of infinite transcendence degree over $\C$
(this example is due to Van den Bergh).
\eer

Further, assume that the algebra $A$ is left Noetherian, and let 
$D^b(\Lmod A)$  be the full 
subcategory in $D(\Lmod A)$ formed by
the complexes $C\in D(\Lmod A)$ such that

\pb{Each cohomology group $H^i(C)$ is a finitely generated $A$-module;}

\pb{$H^i(C)=0$ for all but finitely many $i\in \Z$.}

\noindent
It is easy to show that $D^b(\Lmod A)$ is a triangulated
subcategory that contains $\Perf(A)$. Furthermore, one proves

\begin{lem}\label{perf_equiv} For a  left Noetherian
algebra $A$, the following conditions are equivalent:

\vi The  category $\Lmod A$ has finite homological dimension;

\vii The inclusion $\Perf(A)\into D^b(\Lmod A)$ is an equivalence;

\viii Any finitely-generated left $A$-module has
a finite resolution by finitely-generated projective $A$-modules.
\qed
\end{lem}

\subsection{Auslander-Reiten functor.}
Below, it will be helpful for us to observe that 
a left $A\ee$-module is the same thing
as an $A$-bimodule, and also
 is the same thing
as a right  $A\ee$-module. This
may be alternatively explained by the existence of
the canonical algebra isomorphism $(A\ee)^\op\simeq A\ee,$
given by the flip.

Now, the object $A\ee=A\otimes A$
is clearly both  a left  and  right $A\ee$-module.
The left  $A\ee$-action on $A\ee$
corresponds to the `outer' $A$-bimodule
structure on $A\otimes A$, explicitly given by
$(a',a'') : x\otimes y \mto (a'x)\otimes (ya'').$
We will indicate this `outer action' by writing
$A\ee= 
{}_{_A}A\otimes A_{_A}.$ More generally,
given an $A$-bimodule $M$, we will use
the notation ${}_{_A}M$, resp., $M_{_A}$,
whenever we want to emphasize that $M$ is viewed
as a left, resp. right, $A$-module.
With these notations,
the right $A\ee$-action on $A\ee$
corresponds to the `inner' $A$-bimodule
structure: $A\ee =A_{_A}\,\otimes\, {}_{_A}A$, explicitly given by
$(a',a'') : x\otimes y \mto (xa')\otimes (a''y).$

View $A$ and $ A\ee$ as a left $A\ee$-modules,
 and put 
$$\su:= R\Hom_{D(\Lmod {A\ee})}(A,A\ee).$$
The right $A\ee$-module structure on $A\ee$ induces
one on $R\Hom_{D(\Lmod {A\ee})}(A,A\ee)$. This makes 
$\su$ a complex of right  ${A\ee}$-modules.
As we have explained, any  right $A\ee$-module
may be as well viewed as a left  $A\ee$-module.
Thus, we may (and will)
regard $\su$ as an object of $D(\Lmod {A\ee})$.

\begin{examp}\label{free_sample} Let $V$ be a finite dimensional
vector space.
The tensor algebra $A=TV$ is homologically smooth,
and has a standard $A$-bimodule resolution:
$$
0\map TV \otimes V\otimes TV \stackrel\varkappa\too
TV \otimes TV \stackrel{\tt{mult}}\too TV\map 0,
$$
where the map $\varkappa$ is given by
$\varkappa:\ a\otimes v\otimes b\mapsto (a\cdot v)\otimes b- a\otimes (v\cdot b)$
(this is a special
case of Koszul bimodule-resolution for a general
 Koszul algebra, cf. e.g. \cite{BG},\cite{VdB1}).
Therefore, we find that, for $A=TV,$
the object $\su$ is represented by the  following two-term complex
$$ TV \otimes TV \too
TV \otimes V^*\otimes TV,
\quad a\otimes b\mto \sum_{i=1}^r \big[
(a\cdot v_i)\otimes {\check v}_i\otimes b -
a\otimes {\check v}_i\otimes (v_i\cdot b)\big],
$$
where $\{v_i\}$ and $\{{\check v}_i\}$ are dual bases of $V$ and $V^*$,
respectively.
\eex

\begin{defn} The functor
$$\Perf(A)\to D(\Lmod A),\quad\text{resp.,}\quad
\Perf(A\ee)\to D(\Lmod {A\ee}),
$$
given by $M\mapsto \su\lotimes_A M,$
will be called the {\em Auslander-Reiten functor}.
\end{defn}
Auslander and  Reiten considered a similar functor (on certain abelian
categories) in their study of representation theory 
of finite dimensional algebras.

Next, we extend the notion of Hochschild homology and cohomology
to objects of $D(\Lmod {A\ee})$ by the formulas
$$\HH^i(A,M):=\Hom_{\Perf(A\ee)}(A,M[i])\aand
\HH_i(A,M):=H^{-i}(A\lotimes_{A\ee}\,M)
$$
where the negative sign is chosen in order to make the
present definition compatible with the standard
definition of Hochschild homology of a bimodule,
as given in \S5.

With these definitions we have the following result, see
\cite{VdB2}.

\begin{prop}[Duality]\label{vdb_duality} For any $M\in \Perf(A\ee)$, and $i\in\Z$,
there is a natural isomorphism $\HH^i(A,M)\simeq \HH_{-i}(A,
\su\lotimes_{A}\,M).$
\end{prop}
\begin{proof}
For any $M \in \Perf(A\ee),$ we have
$\HH^i(A,M)=\Hom_{\Perf(A\ee)}(A,M[i])$.
Therefore, by the second formula in
\eqref{can_iso_perf}, we find
$$\HH^i(A,M)=H^i\bigl(\Hom_{\Perf(A\ee)}(A,A\ee)\,\lotimes_{A\ee}\,\,
 M\bigr)= H^i(\su\lotimes_{A\ee}\, M).
$$
Observe that, for any right $A\ee$-modules $R$
and left $A\ee$-module $L$, 
the object $R\lotimes_{A} L$ carries
an $A$-bimodule structure, equivalently, a left
 $A\ee$-module structure; furthermore, an analogue of the
third isomorphism in \eqref{can_iso_perf} says
$R\lotimes_{A\ee} L \simeq
A\lotimes_{A\ee} \bigl(R\lotimes_A L\bigr).$
Using this formula, and the previous calculation, we find
\begin{align*}
\HH^i(A,M)=H^i(\su\lotimes_{A\ee}\, M)=
H^i\bigl(A\lotimes_{A\ee} & (\su\lotimes_A M)\bigr)=
\HH_{\!-i}(A,\,\su\lotimes_A\,M).
\end{align*}
\end{proof}

We say that an
associative (not necessarily commutative) algebra $A$ is {\em Gorenstein}
of dimension $d$ if one has
$$\Ext_{\bimod A}^i(A,A\ee)\cong
\begin{cases}
A&\text{if}\enspace i=d\\
0&\text{otherwise}.
\end{cases}
$$
Clearly, for a Gorenstein algebra $A$, in $\Perf(A\ee)$
one has $\su\cong A[d]$. Hence, from Proposition \ref{vdb_duality} we deduce
\begin{cor} Given a  Gorenstein algebra $A$ of dimension $d$, 
for any $M\in \Perf(A\ee)$ there are canonical isomorphisms
$\HH^\hdot(A,M)\simeq \HH_{d-\idot}(A,M).
$\qed
\end{cor}

Assume next that $A$ is a {\em finite-dimensional} algebra
of finite homological dimension. 
Then any  finite rank projective $A$-module is finite dimensional,
hence any object of $\Perf(A)$, resp., $\Perf(A^\op)$,
is quasi-isomorphic to a complex
of $\C$-finite dimensional $A$-modules.
Thus, taking the $\C$-linear dual of (a complex of) finite dimensional
vector spaces induces a functor
$\Perf(A^\op)\to D^b(\Lmod A)\iso \Perf(A)\,,\,M$
$\mapsto M^*:=\Hom_\C(M,\C)$.

\begin{prop}\label{fin_dim_serre} Let $A$ be a finite-dimensional algebra of
finite homological dimension. Then

\vi The composite functor
$$
\xymatrix{
\Perf(A)\ar[rr]^<>(0.5){R\Hom_A(-,A)}&& \Perf(A^\op)
\ar[rr]^<>(0.5){\Hom_\C(-,\C)}&&\Perf(A),}
$$
$M\mapsto R\Hom_A(M,A)^*,$  is a Serre functor on $\Perf(A)$.

\vii The functor $M\mto \su\lotimes_A M$ is an inverse
to the Serre functor $\se$; furthermore, one
has a functorial isomorphism $\su\lotimes_A M\simeq
R\Hom_A(A^*,M).$
\end{prop}

Our proof of the Proposition given below exploits the following
useful isomorphism that holds for any  finite dimensional
algebra $A$:
\beq{useful}
R\Hom_{A\ee}(A, E^*)\simeq
(A\lotimes_{A\ee} E)^*\,,\quad\text{for any}\quad E\in\Perf(A\ee).\quad\Box
\end{equation}

\begin{proof}[Proof of Proposition \ref{fin_dim_serre}.]
 Recall that, for any left $A$-modules
$M,N$, the vector space $\Hom_\C(M,N)$ has a natural
$A$-bimodule structure,
and we have a natural isomorphism
\beq{useful2}\Hom_A(M,N)=\Hom_{\bimod A}\bigl(A,\,\Hom_\C(M,N)\bigr).
\end{equation}
Using this formula, for $\C$-finite dimensional
  left $A$-modules $M,N$, we  compute
\begin{align*}
\Hom_A\bigl(N,\,\Hom_A(M,A)^*\bigr)&=
\Hom_{\bimod A}\bigl(A,\, N^*\o \Hom_A(M,A)^*\bigr)\\
&=\Hom_{\bimod A}\bigl(A,\, (\Hom_A(M,A)\o N)^*\bigr).
\end{align*}

Therefore,
 for any finite rank projective $A$-modules
$M,N$, using the previous calculation, we find
\begin{align*}
\Hom_A\bigl(N,\,\Hom_A(M,A)^*\bigr)&=
\Hom_{\bimod A}\bigl(A,\, (N\o \Hom_A(M,A))^*\bigr)\\
&=\bigl(A\lotimes_{A\ee} \bigl(N\o \Hom_A(M,A)\bigr)\bigr)^*\\
\text{by \eqref{useful}}\quad&=\bigl(\Hom_A(M,A)\otimes_A N\bigr)^*\\
\text{by \eqref{can_iso_perf}}\quad&=\bigl(\Hom_A(M,N)\bigr)^*.
\end{align*}

Thus, we have established a natural isomorphism
$$\Hom_A(M,N)^*\simeq \Hom_A\bigl(N,\,\Hom_A(M,A)^*\bigr),$$
for any finite rank projective $A$-modules.
It follows that a similar isomorphism holds for any
objects $M,N\in \Perf(A)$. Thus, we have
verified the defining property of Serre functor,
and part (i) is proved.

We claim next that the functor  $M\mto \su\lotimes_A M$ is
a left adjoint of the Serre functor, that is,
one has a functorial isomorphism
\beq{adj1}
\Hom_{\Perf(A)}(N, \se(M)) \cong \Hom_{\Perf(A)}(\su\lotimes_A N,\,M).
\end{equation}

To establish the isomorphism above, we first use the
canonical adjunction isomorphism for tensor products.
This says 
$$R\Hom_A(\su\lotimes_A N,\,M)\simeq R\Hom_{A\ee}(\su,\, M \otimes_{_\C}N^*)=
R\Hom_{A\ee}\bigl(\su,\, \Hom_\C(N,M)\bigr).
$$

Now, using the second formula in \eqref{can_iso_perf} we compute
\begin{align*}
R\Hom_{A\ee}\bigl(\su,& \Hom_\C(N,M)\bigr)=
R\Hom_{A\ee}(\su,A\ee)\lotimes_{A\ee}\Hom_\C(N,M)\\
&=R\Hom_{A\ee}(R\Hom_{A\ee}(A,A\ee)\,,\,A\ee)
\lotimes_{A\ee}\Hom_\C(N,M)\\
\text{by \eqref{can_iso_perf}}\quad&=A\lotimes_{A\ee}\Hom_\C(N,M)\\
\text{by \eqref{useful}}\quad&=\bigl(R\Hom_{A\ee}\bigl(A,
\Hom_\C(N,M)\bigr)\bigr)^*\\
&=\bigl(R\Hom_A(M,N)\bigr)^*=R\Hom_A(N,\se(M)),
\end{align*}
where the last equality holds by the definition of  Serre functor.

Thus we have proved our claim that,
for the left adjoint functor ${}^{\top\!}\se$, we have
${}^{\top\!}\se(-)
= \su\lotimes_A (-).$
But the explicit form
of the Serre functor provided by part (i) clearly shows that this
functor is an equivalence. Hence, its left adjoint functor
must be an inverse of $\se$, and the first statement
of part (ii) follows. 

To prove the last statement
 we compute, cf. \cite[\S2]{CrB}:
$$R\Hom_A(A^*,M)\cong R\Hom_A(A^*,A)\lotimes_A M \cong R\Hom_{A\ee}(A,
A\otimes A)\lotimes_A M,
$$
where the first isomorphism is due to the second formula in \eqref{can_iso_perf}
and the second  isomorphism is due to \eqref{useful2}.
\end{proof}

\subsection{Homologically smooth algebras.} 

Recall that the category of $A$-bimodules has a natural
monoidal structure $M,N\mto M\otimes_A N$,
where $M\otimes_A N$ is again an $A$-bimodule.
This gives, at the level of derived categories,
a monoidal structure $\Perf(A\ee)\bigotimes
\Perf(A\ee)\map \Perf(A\ee)\,,\,M,N\mto M\lotimes_A N$,
where we identify objects of $\Perf(A\ee)$
with complexes of  $A$-bimodules. 
In a similar way, the category
$\Perf(A)$ is a {\em module category}
over the monoidal category  $\Perf(A\ee)$,
with the module structure
$\Perf(A\ee)\bigotimes\Perf(A)\map \Perf(A)$
given by the derived tensor product over $A$.

Fix an
algebra $A$.
We consider $A$ as an $A$-bimodule, that is,  as an object of
$D(\Lmod {A\ee})$. 

Following Kontsevich, we introduce
\begin{defn} The algebra $A$ is called {\em homologically smooth}
if $A$ is a \emph{compact} object of $D(\Lmod {A\ee})$,
equivalently,
if  $A\in \Perf(A\ee)\sset D(\Lmod {A\ee})$, that is, if 
$A$ has a finite resolution
by finitely-generated
projective (left) $A\ee$-modules.
\end{defn}

This definition is motivated by the following result
\begin{lem} An affine algebraic variety $X$ is smooth if and
only if its coordinate ring, $\C[X]$, is a  homologically smooth
algebra (more generally, a scheme $X$ is smooth  if and
only if $\oo_{X_\Delta}$, the structure sheaf of the diagonal
$X_\Delta\sset X\times X,$ is a compact object 
in $D(\Lmod {\oo_{X\times X}})$).
\end{lem}
\begin{proof} If $X$ is smooth, then so is $X\times X$.
Hence, by Grothendieck's theorem \ref{groth},
the structure sheaf $\oo_{X_\Delta}$
has a finite locally-free resolution.
Conversely, let $\scr F_d\to\ldots\to\scr F_1\to\scr F_0\onto \oo_{X_\Delta}$
be a resolution of $\oo_{X_\Delta}$ by locally free  sheaves on
 $X\times X$. Each term of the resolution, as well as the sheaf
 $\oo_{X_\Delta}$ itself, is {\em flat} relative to 
the first projection $X\times X\to X$. Hence,
for any point $x\in X$, the resolution restricts
to an exact sequence
$$
\scr F_d|_{\{x\}\times X}\to\ldots\to\scr F_1|_{\{x\}\times X}\to\scr
F_0|_{\{x\}\times X}
\onto \C_x,
$$
where $\C_x$ denotes the sky-scrapper sheaf at the point $x$.
This exact sequence provides a bounded 
resolution of  $\C_x$ by
locally free sheaves on $X$. Hence, by the standard
regularity criterion, cf. e.g. \cite{Eis},
$X$ is smooth at the point $x$. Thus, $X$ is smooth at every 
point, and we are done.
\end{proof}

Thus, homologically smooth algebras should be thought of as
coordinate rings of 
{\em smooth}
`noncommutative spaces'.
We remark that the property of being "homologically smooth"
is  weaker 
than that of being "formally smooth":
many algebras, such as polynomial algebras, universal enveloping
algebras, etc., are homologically smooth, but
not formally smooth. On the other hand, we have
\begin{lem} \vi 
Any formally smooth algebra is homologically smooth.

\vii If $A$ and $B$ are homologically smooth, then so is $A\otimes B$,
and $A^\op$.
\end{lem} 
\begin{proof} If $A$ is formally smooth, then
we have a length two projective resolution
$\ncO(A)\to A\ee\onto A$. Hence, $A\in \Perf(A\ee)$,
and  $A$ is  homologically smooth.
Part (ii) is straightforward, and is left to the reader.
\end{proof}

\begin{lem} Let $A$ be an associative algebra
such that $A\ee$ is Noetherian.

\vi The following conditions on $A$ are equivalent

\noindent
\pb{The algebra  $A$ is  homologically smooth;}

\noindent
\pb{The category $\Lmod {A\ee}$ has finite
homological dimension;}

\noindent
\pb{There exists $d\gg 0$ such that, for any $A$-bimodule $M$,
one has\linebreak
 $\HH^i(A,M)=0,$ for all $i>d$.}

\noindent
\vii If $A$ is  homologically smooth, then
the equivalent conditions {\sf{(i)-(iii)}} of Lemma \ref{perf_equiv}
hold for $A$.
\end{lem}
\begin{proof} See \cite{VdB2}.
\end{proof}

Assume now that  $A$ is a homologically smooth algebra, so $A\in
\Perf(A\ee).$
Then
$A$ clearly plays the role of unit for the monoidal structure
$(-)\lotimes_A(-),$ on $\Perf(A\ee).$
An object $R\in \Perf(A\ee)$ is said to be
{\em invertible} if there exists another object
$R'\in \Perf(A\ee)$, such that in $\Perf(A\ee)$ one has
$$R\lotimes_A R'\simeq A\simeq R'\lotimes_A R.
$$
In this case, one calls $R'$ an inverse of $R$,
which is uniquely defined up to isomorphism.

It
 is  straightforward to see from the definitions that the object
$\su\in D(\Lmod {A\ee})$
belongs to $\Perf(A\ee)$ and, moreover,
 the derived tensor product with $\su$ preserves the category
$\Perf(A)$, that is,
gives a functor
\beq{R_functor}
\su\lotimes_A(-):\ \Perf(A) \too\Perf(A),\quad M\mto
\su\lotimes_A\, M.
\end{equation}
\begin{question}[Kontsevich]\label{invertible}
Is it true that  $\su\in \Perf(A\ee)$ is an invertible object,
for any  homologically smooth
algebra $A$ ?
\end{question}

It is not known (to the author) whether the two-term complex representing
the object $\su$ for the free algebra 
$A=TV, \dim V>1$, see Example \ref{free_sample}, is an
invertible object in
$\Perf\bigl((TV)\ee\bigr),$ i.e., whether
Question \ref{invertible} has a positive answer for  free
associative algebras with more than one generator.

Assume that $A$ is a  homologically smooth
algebra such that $\su$ is invertible, and write
$\sd\in \Perf(A\ee)$ for the inverse of $\su$ (which is
well-defined up to isomorphism). It is clear that the functors
 $$\Perf(A)\to\Perf(A),\quad M\mapsto \su\lotimes_A\, M,
\quad\text{resp.,}\quad M\mapsto \sd\lotimes_A\, M,
$$
are mutually inverse auto-equivalences of $\Perf(A)$.
Furthermore,  one can prove the following

\noindent
\pb{There are algebra quasi-isomorphisms:}
$$A^\op\stackrel{_\sim}\to R\Hom_{\Perf(A)}({}_{_A}\sd,\,{}_{_A}\sd),
\aand A\stackrel{_\sim}\to R\Hom_{\Perf(A^\op)}(\sd_{\!_A},\,\sd_{\!_A}).
$$
\pb{In $\Perf(A\ee)$, there is an isomorphism:}
$$ \sd\simeq R\Hom_{\Perf(A\ee)}(A,\,{}_{_A}\sd\otimes \sd_{\!_A}).
$$
\begin{rem}
The above properties show that,
if  $\su$ is an invertible object,
then its inverse, $\sd\in \Perf(A\ee)$,
is the {\sl{rigid dualizing complex}} for $A$,
as defined by Van den Bergh \cite{VdB1}. This is known
to be the case, for instance, if  $A=\C[X]$ is the coordinate ring of
a smooth variety $X$. Then, $\Perf(\C[X])=\Dcoh(X)$,
and the functor
$\sd\lotimes_A\,(-)$ reduces to 
$M\mapsto M\lotimes_A\, \comO(A)$
(this is not a Serre functor because the category
$\Perf(\C[X])$ is not $\Hom$-finite, since $X$ is
not compact.)
\eer

\section{Geometry over an Operad}\subsection{}
 We recall that, for any
 $\Z/(2)$-graded vector space $M=M_\ev\oplus M_\odd$,
one defines 
the parity 
reversal operator,  $\Pi$, such that $\Pi{M}=M_\odd \oplus M_\ev$.
Thus  $\Pi{M}$ is
a  $\Z/(2)$-graded  vector space again.

 Let 
$\PP=\{\PP(n),\,n=1,2,\ldots\}$ be a $\C$-linear quadratic operad
with $\PP(1)=\C$, see
[GiK]. Let $\s_n$ denote the Symmetric group on $n$ letters.
Given $\mu\in \PP(n)$ and a $\PP$-algebra $A$,
we will write: $\mu_A(a_1,\ldots,a_n)\,$
for the image of $\mu\otimes a_1\otimes\ldots\otimes a_n$
under
the structure map: $\PP(n)\otimes_{_{\s_n}} A^{\otimes n} \too A\,.$
Following [GiK, \S1.6.4], we introduce
an {\it enveloping algebra} $\upa$, the associative unital
$\C$-algebra  such that
the abelian category of (left) $A$-modules is equivalent to the
category of left modules over  $\upa$, see [GiK, Thm. 1.6.6]. 
The algebra $\upa$ is generated by the symbols:
$u(\mu, a)\,,\, \mu\in \PP(2), a\in A,$ subject to certain
relations, see [Ba, \S1.7].

A $\PP$-algebra in the monoidal category of
$\Z/(2)$-graded, (resp. $\Z$-graded) super-vector spaces, 
see [GiK, \S1.3.17-1.3.18],
will be referred to as
a  $\PP$-{\it super-algebra}, (resp. graded
 super-algebra). Any $\PP$-algebra may be regarded as a
$\PP$-superalgebra concentrated in degree zero.
Following \cite[\S5]{Gi}, we define a  free graded $\PP$-algebra 
(resp. super-algebra) generated by $V$ by
$$
\T^{^{_\idot}}\pp\!{V} :=\bigoplus\nolimits_{i\geq 1}\;\;
\PP(i) \otimes_{_{\s_i}}V^{\otimes i}\;\aand\;
\Tc^{^{_\idot}}\pp{_{\!}}{V} :=\bigoplus\nolimits_{i\geq 1}\;\;
\PP(i) \otimes_{_{\s_i}}(\Pi{V})^{\otimes i}
.$$

Fix a $\PP$-algebra $A$. Following \S\ref{tensor_algebra},
we consider the category
$A\textsf{-algebras},$ whose objects are pairs $(B,f)$, where
$B$ is  a $\PP$-algebra and $f: A\to B$ is
a $\PP$-algebra morphism.
Arguing as in \S\ref{tensor_algebra},
 we get a functor $A\textsf{-algebras}\too A\textsf{-modules}$.
By \cite[Lemma 5.2]{Gi}, this functor
has a right adjoint, i.e., we have

\begin{lem}\label{TM} Given a $\PP$-algebra $A$, there
is a functor: $M \mapsto T_A^{^{_\idot}}M\,,$
(resp. $M \mapsto \tc^{^{_\idot}}_AM$) assigning to a left
$A$-module $M$ a graded 
 $\PP$-algebra 
$\,T^{^{_\idot}}_AM= \bigoplus_{i\geq 0}\;T^i_AM$ (resp. graded 
 $\PP$-superalgebra $\,\tc^{^{_\idot}}_AM
= \bigoplus_{i\geq 0}\;T^i_A(\Pi{M})\,$)
equipped with
a canonical $\PP$-algebra isomorphism 
$\imath: A \iso T^0_AM$.
Moreover, for any $\PP$-algebra map: $A\to B$, one has
 a natural adjunction isomorphism:
$$\Hom_{_{A{\sf{\tdash}modules}}}(M,B)\iso 
\Hom_{_{\PP{\sf{\tdash}algebras}}}(T^{^{_\idot}}_AM,B).\qquad\Box$$
\end{lem}

An ideal $I$ in a $\PP$-algebra $A$ will be called $N$-nilpotent if,
for any $n\geq N\,,\,\mu\in\PP(n)$, and $a_1,\ldots,a_n\in A\,,$ one
has: $\mu_A(a_1,\ldots,a_n)=0,$ whenever at least $N$ among the elements
$a_1,\ldots,a_n$ belong to $I$.

Given a left $A$-module $M$, one defines the
{\em square zero} extension $A\sharp M$,
cf.  \cite[Definition 3.2.6]{Ba}, or \cite[Lemma 5.1]{Gi}.
The following useful reformulation of the notion
of a left $A$-module is due to [Ba,~1.2], [Qu]:

\begin{lem}\label{La-mod}
Giving a left $A$-module structure on a vector space $M$ is equivalent
to giving a $\PP$-algebra structure on $\,A\sharp M:= A\oplus M\,$ 
such that the  following
conditions hold:

\vi The imbedding: $a\mapsto a\oplus 0$
makes $A$ a $\PP$-subalgebra in $A\sharp M$.

\vii $\;M$ is a 2-nilpotent ideal in $A\sharp M$. \quad\qed
\end{lem}

\begin{defn}\label{der}
 A $\C$-linear map $\theta:
A\to M$ is called a {\it derivation} if
the map: $\,a\bigoplus m \mapsto a\,\bigoplus\, \theta(a)\!+\!m\,,$
is an automorphism of the $\PP$-algebra $A\sharp M$.
\end{defn}

Let
$\derp(A,M)$ denote the $\C$-vector space of all derivations
from $A$ to $M$. It is straightforward to see that 
 the ordinary commutator makes $\derp(A,A)$ 
a Lie algebra.\medskip

Next we define, following [Ba, Definition 4.5.2], 
an $A$-module of {\it K\"ahler differentials}
 as the left $\upa$-module, $\Om^1\pp{A}$, generated by the symbols
$da$, for $a\in A$, subject to the relations:

\noindent
$\bullet\quad${$d(\lambda_1 a_1 +\lambda_2 a_2) = \lambda_1 da_1+\lambda_2
da_2,\enspace\forall \lambda_1,\lambda_2 \in \C;$}

\noindent
$\bullet\quad${$
 d(\mu(a_1,a_2)) = u(\mu,a_1)\otimes da_2 + u(\mu^{(12)}, a_2)\otimes
da_1,\enspace
\forall \mu\in \PP(2),a_1,a_2\in A,$}

\noindent
where $u(\mu,a)$
denote the standard generators of $\upa$, see [Ba].
\medskip

By construction, $\Om\pp^1\!A$ is a left $A$-module,
and the assignment $a\mapsto da$ gives a derivation
$d\in \derp(A\,,\,\Om\pp^1\!A)$.
Moreover, this derivation is universal in the following
sense.
Given any left $A$-module $M$ and a derivation
$\theta: A\to M$, there exists an 
$A$-module morphism  $\,\Om^1\theta: \Om\pp^1\!A \to M$,
uniquely determined by the condition that 
$\,(\Om^1\theta)(da)=\theta(a)\,.$ It follows 
that the 
$A$-module of K\"ahler differentials
represents the functor
$\derp(A,-)$, i.e., we have (see [Ba, Remark 4.5.4]):

\begin{lem}\label{om-der} 
For any left $A$-module $M$ there is a natural isomorphism:
$$\derp(A,M) \;\simeq\; \Hom_{_{A{\sf{\tdash}mod}}}(\Om^1\pp{A}, M)\,.\quad\square$$
\end{lem}

\begin{lem}\label{B}
There is a natural $A$-module morphism
$\delta: \Om^1(A\sharp\Om^1\pp{A})\too$
$ \Om^1\pp{A}\bigotimes_{_{\upa}}\,\Om^1\pp{A}$.
\end{lem}

\begin{lem}\label{A&M-module} Let $A$ be a $\PP$-algebra, and $M$ a
left $A$-module. Giving the structure of a left $A\sharp M$-module
on $N$ is equivalent to giving a $\PP$-algebra structure
on $A\oplus M\oplus N$ such that the  following
conditions hold:

\vi The imbedding: $a\mapsto a\oplus 0\oplus 0$
makes $A$ a $\PP$-subalgebra in $A\oplus M\oplus N$.

\vii $\;M$ is a 2-nilpotent ideal in $A\sharp M$. \quad\qed
\end{lem}

Next, we define $\Om\bulp A$, the {\it differential envelope}
of a $\PP$-algebra $A$,
 as the  
graded $\PP$-{\it super}-algebra:
$\,\Om\bulp A = \tc^{^{_\idot}}_A(\Om^1\pp{A})\,.$
The canonical derivation $d: A \to \Om^1\pp{A}$
extends, via the Leibniz rule,
 to a $\PP$-superalgebra derivation 
$d\,:\Om\pp^{^{_\idot}}\!A \to \Om\pp^{\bullet+1}\!A$, such that $d^2=0$.
Thus, $\Om\bulp A$ is a differential
graded $\PP$-{\it super}-algebra.
Further, there is a natural $\PP$-superalgebra imbedding
$j: A=\Om\pp^0A \into \Om\bulp A$. Lemmas \ref{TM} and \ref{om-der} 
imply that this imbedding has the following universal
property, cf. [Ko2] and [KR]: given a differential
graded $\PP$-superalgebra $B$ and a graded $\PP$-algebra
morphism $f: A\to B$, there exists
a unique morphism $\,\Om(f): \,\Om\bulp A \too B\,$
of differential
graded $\PP$-superalgebras, such that: $\,f=\Om(f)\ccirc j\,.$

\begin{prop}[\cite{Gi}, Proposition 5.6]\label{d^2=0} \vi There is a natural
super-differential $\,d: \Om\bulp A \too \Om^{\bullet+1}\pp{A}$, $
d^2=0$,
such that its restriction $A=\Om^0\pp{A}\to \Om^1\pp{A}$
coincides with the canonical $A$-module
derivation $d: A \to \Om^1\pp{A}$. 

\vii The differential graded $\PP$-superalgebra $\,(\Om\bulp A\,,\,d)\,$
is the differential envelope of $A$.\qed
\end{prop}

We set $\Th\pp\bul{A} = 
\Hom_{_{A{\sf{\tdash}mod}}}(\Om^{^{_\idot}}\pp{A}\,,A)$,
and call the $\C$-vector space $\Th_pA := $\linebreak
$\Hom_{_{A{\sf{\tdash}mod}}}(\Om^p\pp{A}\,,
A)$ the space of $\;p${\it-polyvectors} on $A$.\medskip

\begin{defn}\cite{Ba}\label{right}
A $\C$-vector space  is called a {\it right}
$A$-module if it is a right $\,\up A$-module.
\end{defn}\medskip

Assume from now on that $\PP$ is a Koszul operad, 
 and write $\PP^!$ for the quadratic dual operad, see [GiK].
Following
[GiK, \S4.2], [KV] and [Ba, \S4],
one defines Hochschild homology, $H_\idot^\PP(A,N)$,
of  $A$ with coefficients\footnote{in 
[GiK, \S4.2] only the case of trivial coefficients has been
considered} 
 in a {\it right} $A$-module $N$.
To this end,  write 
$\sgn$ for the  1-dimensional sign-representation,
of $\s_n$, and  given an $\s_n$-module $E$, set
$\,E^\vee := \Hom_\C(E,\C)\otimes \sgn$.
We define $H_\idot^\PP(A,N)$ as the homology groups
 of the differential graded vector space:
\begin{equation}\label{GiK1}
C^\PP_\idot(A,N)= \bigoplus_{n\geq 1}\;C^\PP_n(A,N)\quad,\quad
C^\PP_n(A,N)\;= \;N\,
\otimes_{_{\Bbbk}}\,\PP^!(n)^\vee\,
\otimes_{_{\s_n}}\,
A^{\otimes n}\,,
\end{equation}

As was pointed out in [Ba, Prop. 4.5.3], for any
right $A$-module $N$, there is a natural
isomorphism: $\,H_1^\PP(A,N)\simeq N\otimes_\upa \Om^1_{_{\!\PP\!}}A\,$.
In particular, for $N=\upa$,
the enveloping algebra of $A$ regarded as a right $A$-module,
one obtains an isomorphism of {\it left} $A$-modules:
\begin{equation}\label{om1=H}
\Om^1_{_{\!\PP\!}}A\;\simeq \;H_1^\PP(A,\,\upa)\,,
\end{equation}
where the left $A$-module structure on the RHS is induced by the
left $A$-module structure on $\upa$.

The right $A$-module structure on the enveloping algebra $\upa$
enables us to  form the chain complex, see (\ref{GiK1}) and [Ba, \S5]:
\begin{equation}\label{W}
\bpa = C^\PP_\idot(A\,,\,\upa)= \bigoplus_{n\geq 1}\enspace
\left(\upa\,
\otimes_{_{\Bbbk}}\,\PP^!(n)^\vee\,
\otimes_{_{\s_n}}\,
A^{\otimes n}\right)\,.
\end{equation}
Further, the left $A$-module structure on $\upa$ makes
$\bpa$ into the following augmented complex of {\it left} $A$-modules:
\begin{equation}\label{bar_complex}
\ldots \B^\PP_iA \too \B^\PP_{i-1}A\too\ldots\too
\B^\PP_2A \too \B^\PP_1A\onto \Om\pp^1A\;,
\end{equation}
\noindent
Here the augmentation: $\B^\PP_1A\onto \Om\pp^1A$ is induced
by the tautological map: $\, \B^\PP_1A = \upa\,
\otimes_{_{\Bbbk}}\,\PP^!(1)^\vee\,
\otimes_{_{\Bbbk}}\,A\,\iso \upa\otimes A\,$, using the observation
that the image of the morphism: $\,\B^\PP_2A \too \B^\PP_1A\,$
is the submodule of $\,\B^\PP_1A=\upa\otimes A\,$
generated by all the elements of the form:
$\,d(\mu(a_1,a_2)) - u(\mu,a_1)\otimes da_2 - u(\mu^{(12)}, a_2)\otimes
da_1\,$, which is exactly the kernel
of the projection: $\,\upa\otimes A \onto \Om^1\pp{A}\,$, see relations
of condition (ii) 
in the definition of $\Om\pp^1A$.

In the associative case the complex $\bpa$ is essentially
the standard Bar-resolution of the algebra $A$ viewed as an
$A$-bimodule. In particular, for an associative
algebra $A$, the bimodule
$\,{\mathtt{Coker}}(\B^\PP_2A \too \B^\PP_1A)\,$ coincides with:
$\,{\mathtt{Coker}}(A^{\otimes 4} \too A^{\otimes 3})\,,$ which is 
equal, due to the exactness of the bar-resolution, to
$\,{\mathtt{Ker}}(A\otimes A \too A)\,.$ 
For this reason, given an algebra $A$ over a general operad $\PP$,
 we will refer to 
$\bpa$, or to the complex (\ref{bar_complex}),
as the {\it Bar-complex} of $A$. \medskip

Similarly, there is a notion of {\it co}-homology, $H^{^{_\idot}}_\PP(A,M),$
with coefficients in a {\it left} $A$-module $M$, see [KV],
and one has:
$\der\pp(A,A) \simeq H^1_\PP(A,A)$.

\begin{prop}\label{bal_result}
For any $\PP$-algebra $A$ and an $A$-module $M$ there is a canonical
isomorphism: 
$$\;H^i\pp(A,M) \;\simeq\;H^i\bigl(\Hom_{_{A{\sf{\tdash}mod}}}(\bpa, M)\bigr)
\quad,\quad\forall i\geq 1\,.$$
\end{prop}

\begin{proof} See [Ba, Proposition 5.2].\qed\medskip
\end{proof}
\begin{conjecture}\label{bpa_conjecture}
If $\PP$ is a Koszul operad, then
$\bpa$ is a  resolution of the left $A$-module $\Om\pp^1A$;
equivalently, and one has:
 $\,H_i(\bpa)= 0\,,
\,\forall i>1\,.$ 
\end{conjecture}

The Conjecture would imply the following result.

\begin{cor}\label{tor_conjecture} We have:
$\,{\mathtt  {Tor}}^{^\upa}_i(\upa,\,A)=0\,,$ for all $i>0$, and also:
$$H^\PP_{\bullet+1}(A, N)\;
\simeq\; {\mathtt  {Tor}}^{^\upa}_{\bullet}(N, \Om^1\pp{A})\quad,\quad
H^{\bullet+1}\pp(A,M)\;
\simeq\; {\mathtt
{Ext}}^{^{_\idot}}_{_{A{\mathsf{{\tdash}mod}}}}(\Om^1\pp{A}, M)\,,$$
for any right $A$-module $N$ and left $A$-module $M$.\qed
\end{cor}

\begin{conjecture}\label{forms}
There is a natural differential graded space morphism:
$\Om^{^{_\idot}}_{_{\!\PP\!}}A \to C^\PP_\idot(A, A).$
\end{conjecture}

Following Grothendieck and Quillen, see [CQ], we introduce

\begin{defn}\label{smooth}
A $\PP$-algebra $A$ is said to be {\textsl  formally-smooth} if it satisfies the
lifting property with respect to all nilpotent ideals, that is any
$\PP$-algebra homomorphism: $A\to R/I$ can be lifted to a 
$\PP$-algebra homomorphism: $A\to R$, provided $I$ is a nilpotent ideal
in $B$.
\end{defn}

\begin{prop}\label{smothness_properties}  An
 algebra $A$ over a Koszul operad $\PP$ is
formally-smooth if and only if the following two conditions hold, cf. [Gr]:

\vi $\;\;\Om^1\pp{A}$ is a projective left $A$-module;

\vii For any presentation $A=R/I$, where $R$ is a free $\PP$-algebra,
the morphism $\jmath$ below, induced by the
canonical  map: $I \into R \stackrel{d}{\too} \Om\pp^1R$,
is
 injective, i.e., the
following canonical sequence is exact:
$$0\too I/I^2 \stackrel{\jmath}{\too} 
\upa\otimes_{_{\U^{\!^{_\PP}\!}\!R}} \Om\pp^1R
\too \Om\pp^1A \too 0\,.$$
\end{prop}

\begin{proof} See [Gr].\end{proof}

Note that it has been shown in [CQ] that condition \vii above automatically holds
in the associative case (of course, it is not automatic in the
associative
commutative case, see [Gr]).

\begin{prop}\label{quasi-smooth}
The following 3 properties of a $\PP$-algebra $A$ are equivalent:

\vi The left $A$-module $\Om^1_{_{\!\PP\!}}A$ is projective.

\vii Any $\PP$-algebra extension: $I \into R \onto A\,,$
where $I$ is a 2-nilpotent ideal in\par
 $R$ has a $\PP$-algebra
splitting $\,A\into R$, that is:
$R\simeq A\sharp I$.

\viii For any left $A$-module $M$ we have: $\,H^2\pp(A, M)=0$.
\end{prop}

\begin{proof} It has been shown in [Ba, Theorem 3.4.2] that 
2-nilpotent $\PP$-algebra extensions: $I \into R \onto A\,$
are classified by the 2-d cohomology group: $H^2\pp(A, I)$.
This proves that: (ii)$\;\Longleftrightarrow\;$(iii).
Further, by Corollary \ref{tor_conjecture} we have:
$H^2\pp(A, I) = {\mathtt{Ext}}^1_{_{{\sf A\tdash{mod}}}}(\Om\pp^1A, M)\,.$
Thus, (i)$\;\Longleftrightarrow\;$(iii).
\end{proof}

\begin{prop}\label{tensor_smooth}
If $A$ is a formally-smooth $\PP$-algebra, and $M$ a projective
$A$-module,
then the  $\PP$-algebra $T_AM$ is formally-smooth.
\end{prop}
\begin{proof}
 Copy the proof in [CQ, Proposition 5.3(3)].
\end{proof}

\parindent=0pt

\footnotesize{

}
\end{document}